\pdfoutput=1
\documentclass[11pt]{article}
\usepackage[left=1in,right=1in,top=1in,bottom=1in]{geometry}
\usepackage{times}
\usepackage{expl3}
\usepackage{cite}
\usepackage[table]{xcolor}
\usepackage{multirow}
\usepackage{stackengine} 
\usepackage{hhline}
\usepackage{lipsum}
\usepackage{titlesec}
\usepackage[all]{xy}
\usepackage{wrapfig}
\usepackage{enumerate}
\usepackage{epsfig}
\usepackage{tikz-cd}
\usepackage{amsmath}
\usepackage{tabularx}
\usepackage{array}
\usepackage{booktabs}
\usepackage{enumitem}
\usepackage{bbm}
\usepackage{calc}
\usepackage{graphicx}
\usepackage{amsmath}
\usepackage[title]{appendix}
\usepackage{amssymb}
\usepackage{epstopdf}
\usepackage{boldline}
\usepackage{arydshln}
\usepackage{calligra}
\usepackage{bm}
\usepackage{url}
\usepackage{blindtext}
\usepackage{accents}
\usepackage{amsthm}
\usepackage{amscd}

\newtheorem{definition}{Definition}

\newtheorem{theorem}{Theorem}
\newtheorem{proposition}{Proposition}
\newtheorem{corollary}{Corollary}
\newtheorem{lemma}{Lemma}

\newtheorem{example}{Example}
\newtheorem{remark}{Remark}
\usepackage{mathtools}
\usepackage{epstopdf}
\usepackage{balance}
\usepackage{thmtools}
\usepackage{thm-restate}
\usepackage{hyperref}
\usepackage{cleveref}
\usepackage[mathscr]{euscript}

\usepackage[ruled,vlined]{algorithm2e}
\newcommand{\ostar}{\mathbin{\mathpalette\make@circled\star}}

\makeatletter
\newcommand{\removelatexerror}{\let\@latex@error\@gobble}
\makeatother
\makeatletter
\newcommand*{\rom}[1]{\expandafter\@slowromancap\romannumeral #1@}
\makeatother

\ExplSyntaxOn
\newcommand\latinabbrev[1]{
  \peek_meaning:NTF . {
    #1\@}%
  { \peek_catcode:NTF a {
      #1.\@ }%
    {#1.\@}}}
\ExplSyntaxOff

\titleclass{\subsubsubsection}{straight}[\subsubsection]

\begin{document}
\vspace{1cm}
\title{An Index Theorem for Fredholm Operators via the Unitary Conjugation Groupoid}
\vspace{1.8cm}
\author{Shih-Yu~Chang
\thanks{Shih-Yu~\cite{PaperI} is with the Department of Applied Data Science,
San Jose State University, San Jose, CA, U. S. A. (e-mail: {\tt
shihyu.chang@sjsu.edu})
}}

\maketitle

\begin{abstract}
In a previous paper we introduced the unitary conjugation groupoid associated
to any unital separable Type I C*-algebra. This groupoid encodes the
representation-theoretic structure of the algebra through the action of its
unitary group on the characters of commutative subalgebras and admits a
canonical embedding of the algebra into its groupoid C*-algebra. In this paper we apply this framework to two fundamental operator algebras: the algebra of all bounded operators on a separable Hilbert space and the
unitization of the compact operators. For any Fredholm operator in these
algebras we construct a natural equivariant KK-theory class using the phase
of the operator. Applying Kasparov descent for groupoids produces a class in
the K-theory of the associated groupoid C*-algebra. Using the Morita
equivalences established in the previous work, this descended class can be
identified with the symbol class of the operator. We prove that composing this class with the boundary map of the Calkin
extension recovers the classical Fredholm index. In particular, the
construction yields index minus one for the unilateral shift and index zero
for compact perturbations of the identity. This provides a groupoid
equivariant formulation of the Fredholm index and connects the unitary
conjugation groupoid with classical operator index theory.  This perspective
suggests a general approach to index theory through equivariant K-theory of
groupoids associated with operator algebras.
\end{abstract}

\tableofcontents

\newpage

\section{Introduction}

\subsection{Recap of Paper I: The Unitary Conjugation Groupoid $\mathcal{G}_{\mathcal{A}}$}

In a previous paper (Paper I) \cite{PaperI}, we introduced a canonical topological groupoid associated to any unital separable $C^*$-algebra $\mathcal{A}$, called the \emph{unitary conjugation groupoid} and denoted $\mathcal{G}_{\mathcal{A}}$. This groupoid is defined as the action groupoid
\[
\mathcal{G}_{\mathcal{A}} = \mathcal{U}(\mathcal{A}) \ltimes \mathcal{G}_{\mathcal{A}}^{(0)},
\]
where the unit space $\mathcal{G}_{\mathcal{A}}^{(0)}$ consists of pairs $(B,\chi)$ with $B\subseteq\mathcal{A}$ a unital commutative $C^*$-subalgebra and $\chi\in\widehat{B}$ a character. The unitary group $\mathcal{U}(\mathcal{A})$ acts on $\mathcal{G}_{\mathcal{A}}^{(0)}$ by conjugation: $u\cdot(B,\chi) = (uBu^*,\chi\circ\operatorname{Ad}_{u^*})$.

A key feature of the construction is that it forces a paradigm shift away from classical locally compact groupoid theory \cite{Renault1980}. When $\mathcal{A}$ is infinite-dimensional, $\mathcal{U}(\mathcal{A})$ equipped with the norm topology is not separable, and the unit space $\mathcal{G}_{\mathcal{A}}^{(0)}$, when endowed with the Fell topology \cite{Fell1962}, is non-Hausdorff and not locally compact. To remedy this, we equipped $\mathcal{U}(\mathcal{A})$ with the strong operator topology (SOT). However, when $\mathcal{A}\subseteq\mathcal{B}(H)$ acts on an infinite-dimensional Hilbert space $H$, the unitary group $\mathcal{U}(\mathcal{A})$ with the SOT is still not separable, and while this topology makes $\mathcal{U}(\mathcal{A})$ a topological group, the resulting groupoid $\mathcal{G}_{\mathcal{A}}$ is not a Polish groupoid in the usual sense. Instead, it is a standard Borel groupoid admitting a Borel Haar system in the sense of \cite{Tu1999}. The topology on $\mathcal{G}_{\mathcal{A}}^{(0)}$ is defined as the initial topology generated by partial evaluation maps $\operatorname{ev}_a(B,\chi) = \chi(a)$ when $a\in B$, and $\infty$ otherwise.

The central construction of \cite{PaperI} was the \emph{diagonal embedding}
\[
\iota: \mathcal{A} \hookrightarrow C^*(\mathcal{G}_{\mathcal{A}}),
\]
a unital injective $*$-homomorphism that realizes $\mathcal{A}$ as a $C^*$-subalgebra of its own groupoid $C^*$-algebra. For Type I algebras, this embedding satisfies $\iota(\mathcal{A})\subseteq C_0(\mathcal{G}_{\mathcal{A}}^{(0)})$ if and only if $\mathcal{A}$ is commutative, and is natural under isomorphisms. These properties will be essential for the index theorem proved in the present paper.

\subsection{Why $B(H)$ and $\mathcal{K}(H)^\sim$? The Warm-Up Cases}

In this work, we focus on two fundamental examples: $\mathcal{A}=B(H)$, the algebra of all bounded operators on a separable infinite-dimensional Hilbert space $H$, and $\mathcal{A}=\mathcal{K}(H)^\sim$, the unitization of the compact operators. These algebras are chosen for several reasons.

First, both $B(H)$ and $\mathcal{K}(H)^\sim$ are Type I $C^*$-algebras, hence the construction of \cite{PaperI} applies directly. Moreover, their representation theory is well-understood: $B(H)$ has a single irreducible representation up to unitary equivalence (the identity representation), while $\mathcal{K}(H)^\sim$ has two irreducible representations — the infinite-dimensional standard representation and the one-dimensional character $\chi_\infty$ obtained from the quotient $\mathcal{K}(H)^\sim/\mathcal{K}(H)\cong\mathbb{C}$.

Second, these algebras are the natural home for the most fundamental examples in Fredholm theory. The unilateral shift $S\in B(H)$ is the prototypical Fredholm operator with nontrivial index $\operatorname{index}(S)=-1$, while every operator in $\mathcal{K}(H)^\sim$ is a compact perturbation of a scalar and therefore has index zero \cite{Mingo1987}. Thus, these two cases represent the extreme possibilities: nontrivial index and trivial index. 

Third, the groupoid $C^*$-algebras $C^*(\mathcal{G}_{B(H)})$ and $C^*(\mathcal{G}_{\mathcal{K}(H)^\sim})$ admit concrete descriptions: $C^*(\mathcal{G}_{B(H)})$ is Morita equivalent to the Calkin algebra $\mathcal{Q}(H)=B(H)/\mathcal{K}(H)$ \cite{BDF1977, ElliottNatsume1988}, while $C^*(\mathcal{G}_{\mathcal{K}(H)^\sim})$ is Morita equivalent to $\mathcal{K}(H)^\sim\otimes\mathcal{K}$. These Morita equivalences, established in \cite{PaperI}, will play a crucial role in connecting the abstract $K$-theory of the groupoid $C^*$-algebras to the classical index.

\subsection{The Main Obstacle: $K_0$ is the Wrong Target}\label{subsec:The Main Obstacle}

Classical index theory is formulated in terms of $K_0$. For a Fredholm operator $T\in\mathcal{B}(H)$, its index is an integer, and under the identification $K_0(\mathcal{K}(H))\cong\mathbb{Z}$, one has $\operatorname{index}(T)=\partial([\pi_{\mathcal{Q}}(T)])$ where $\partial$ is the boundary map of the Calkin extension \cite{Atkinson1951, Blackadar1998, Rordam2000}. This suggests that a natural target for a groupoid index theorem might be $K_0$ of an appropriate ideal.

However, this approach fails for $\mathcal{A}=\mathcal{B}(H)$ and $\widetilde{\mathcal{K}(H)}$ when one attempts to use the diagonal embedding $\iota$ directly. The $K$-theory groups of these algebras are
\[
K_0(\mathcal{B}(H))=0,\qquad K_1(\mathcal{B}(H))=0,\qquad K_0(\widetilde{\mathcal{K}(H)})\cong\mathbb{Z},\qquad K_1(\widetilde{\mathcal{K}(H)})=0.
\]
Thus, $K_0$ of the algebras themselves cannot detect the nontrivial index of the unilateral shift, and $K_1$ is identically zero in both cases. The correct target must therefore be something else~\cite{Blackadar1998, Rordam2000}.

The resolution is to work not with $K_0$ of $\mathcal{A}$, but with $K_1$ of the Calkin algebra $\mathcal{Q}(H)$ for $\mathcal{B}(H)$, and with $K_1$ of $\mathcal{A}$ itself for $\widetilde{\mathcal{K}(H)}$ (where it vanishes). This requires transporting information from the groupoid $C^*$-algebra to $\mathcal{Q}(H)$ via Morita equivalence, rather than directly through $\iota$. In particular, the map induced by $\iota$ on $K$-theory,
\[
\iota_*: K_1(C^*(\mathcal{G}_{\mathcal{B}(H)}))\longrightarrow K_1(\mathcal{B}(H)),
\]
lands in the zero group and plays no role in the index computation; it must be replaced by the Morita equivalence $C^*(\mathcal{G}_{\mathcal{B}(H)})\sim_M\mathcal{Q}(H)$.

\subsection{The Solution: $K_1$, the Boundary Map, and the Calkin Algebra}

The key observation is that the Fredholm index is naturally an element of $K_0(\mathcal{K}(H))\cong\mathbb{Z}$, obtained via the boundary map $\partial: K_1(\mathcal{Q}(H))\to K_0(\mathcal{K}(H))$ in the six-term exact sequence of the Calkin extension \cite{Blackadar1998, Rordam2000}. For $\mathcal{A}=\mathcal{B}(H)$, the groupoid $C^*$-algebra $C^*(\mathcal{G}_{\mathcal{B}(H)})$ is Morita equivalent to $\mathcal{Q}(H)$, giving an isomorphism $\Phi_*: K_1(C^*(\mathcal{G}_{\mathcal{B}(H)}))\cong K_1(\mathcal{Q}(H))$. The effective boundary map $\partial_{\mathcal{B}(H)}^{\text{eff}}:=\partial\circ\Phi_*$ then directly computes the index.

For $\mathcal{A}=\widetilde{\mathcal{K}(H)}$, the situation is simpler but still instructive. Here $K_1(\widetilde{\mathcal{K}(H)})=0$, and all Fredholm operators have index zero. The groupoid $C^*$-algebra $C^*(\mathcal{G}_{\widetilde{\mathcal{K}(H)}})$ is Morita equivalent to $\widetilde{\mathcal{K}(H)}\otimes\mathcal{K}$, and its $K_1$ vanishes. Thus the composition $\partial_{\widetilde{\mathcal{K}(H)}}\circ\iota_*\circ j_{\mathcal{G}_{\widetilde{\mathcal{K}(H)}}}$ yields zero, matching the index.\footnote{For $\mathcal{A}=\mathcal{B}(H)$, the direct map $\iota_*: K_1(C^*(\mathcal{G}_{\mathcal{B}(H)}))\to K_1(\mathcal{B}(H))$ lands in the zero group and must be replaced by the Morita isomorphism $\Phi_*$ to $\mathcal{Q}(H)$.}

In both cases, the index is recovered by a composition of three canonically defined maps:
\begin{enumerate}
    \item The descent map $j_{\mathcal{G}_{\mathcal{A}}}: KK^1_{\mathcal{G}_{\mathcal{A}}}(C_0(\mathcal{G}_{\mathcal{A}}^{(0)}),\mathbb{C})\to K_1(C^*(\mathcal{G}_{\mathcal{A}}))$, which translates equivariant $KK$-theory data into ordinary $K$-theory of the groupoid $C^*$-algebra \cite{Kasparov1988, Tu1999, legall1994theorie, LeGall1999};
    \item A Morita isomorphism $\Psi: K_1(C^*(\mathcal{G}_{\mathcal{A}}))\to G_{\mathcal{A}}$, where $G_{\mathcal{B}(H)}=K_1(\mathcal{Q}(H))$ and $G_{\widetilde{\mathcal{K}(H)}}=K_1(\widetilde{\mathcal{K}(H)})=0$;
    \item A boundary map $\partial_{\mathcal{A}}: G_{\mathcal{A}}\to\mathbb{Z}$, which is $\partial_{\text{Calkin}}$ for $\mathcal{B}(H)$ and the zero map for $\widetilde{\mathcal{K}(H)}$.
\end{enumerate}
The main theorem of this paper shows that for any Fredholm operator $T$ in $\mathcal{A}$, this composition applied to the equivariant $KK^1$-class $[T]_{\mathcal{G}_{\mathcal{A}}}^{(1)}$ recovers the Fredholm index $\operatorname{index}(T)$.

\subsection{Main Results and Outline of the Paper}

We now summarize the main contributions of this paper.

\begin{enumerate}
    \item \textbf{Construction of the equivariant $KK^1$-class.} For any Fredholm operator $T\in\mathcal{A}$ (where $\mathcal{A}=\mathcal{B}(H)$ or $\widetilde{\mathcal{K}(H)}$), we construct an odd equivariant Kasparov triple $(\tilde{\mathcal{E}},\phi\oplus\phi,\tilde{F})$ representing a class
    \[
    [T]_{\mathcal{G}_{\mathcal{A}}}^{(1)}\in KK^1_{\mathcal{G}_{\mathcal{A}}}(C_0(\mathcal{G}_{\mathcal{A}}^{(0)}),\mathbb{C})\cong K^1_{\mathcal{G}_{\mathcal{A}}}(\mathcal{G}_{\mathcal{A}}^{(0)}).
    \]
    This class encodes the index data of $T$ in a form suitable for equivariant $KK$-theory.
    
    \item \textbf{Identification of the descended class.} Using the descent map for Polish groupoids developed by Tu \cite{Tu1999} and Le Gall \cite{legall1994theorie, LeGall1999}, together with the Morita equivalence between $C^*(\mathcal{G}_{\mathcal{A}})$ and a concrete $C^*$-algebra established in Sections 4.5-4.6 of \cite{PaperI}, we prove that the descended class corresponds to the symbol class $[\pi_{\mathcal{Q}}(T)]\in K_1(\mathcal{Q}(H))$ for $\mathcal{B}(H)$, and to zero for $\widetilde{\mathcal{K}(H)}$.
    
    \item \textbf{The index theorem.} Our main theorem states that for any Fredholm operator $T\in\mathcal{A}$,
    \[
    \operatorname{index}(T) = \partial_{\mathcal{A}} \circ \Psi \circ j_{\mathcal{G}_{\mathcal{A}}}\bigl([T]_{\mathcal{G}_{\mathcal{A}}}^{(1)}\bigr),
    \]
    where $\partial_{\mathcal{A}}$ is the appropriate boundary map (the Calkin index map $\partial: K_1(\mathcal{Q}(H))\to\mathbb{Z}$ for $\mathcal{B}(H)$, and the zero map for $\widetilde{\mathcal{K}(H)}$), $\Psi: K_1(C^*(\mathcal{G}_{\mathcal{A}}))\to G_{\mathcal{A}}$ is the isomorphism induced by the Morita equivalence from \cite{PaperI} (with $G_{\mathcal{B}(H)}=K_1(\mathcal{Q}(H))$ and $G_{\widetilde{\mathcal{K}(H)}}=K_1(\widetilde{\mathcal{K}(H)})=0$), and $j_{\mathcal{G}_{\mathcal{A}}}$ is the descent map. This provides a unified groupoid-equivariant formulation of the Fredholm index.
    
    \item \textbf{Concrete examples and computations.} We illustrate the theorem with detailed calculations for the unilateral shift on $\mathcal{B}(H)$ (index $-1$) and for arbitrary compact perturbations of the identity in $\widetilde{\mathcal{K}(H)}$ (index $0$). These examples demonstrate how the abstract machinery reproduces classical index theorems.
    
    \item \textbf{Connections to classical index theory.} We indicate how our groupoid index theorem relates to the Toeplitz index theorem \cite{Coburn1967} and the Brown–Douglas–Fillmore theory \cite{BDF1977}. A detailed treatment of connections to the Atiyah–Singer index theorem for manifolds with boundary, including the Boutet de Monvel index formula for the ball in $\mathbb{C}^2$, will appear in forthcoming work.
\end{enumerate}

The paper is organized as follows. Section \ref{sec:Preliminaries: K1, the Index Map, and Kasparov's Descent} reviews the necessary background on $K$-theory, the six-term exact sequence, the Calkin extension, Kasparov's descent for Polish groupoids, and the diagonal embedding from \cite{PaperI}. Section \ref{sec:the-unit-space-GA} describes the unit space $\mathcal{G}_{\mathcal{A}}^{(0)}$ and the canonical continuous field of Hilbert spaces $\{H_x\}_{x\in\mathcal{G}_{\mathcal{A}}^{(0)}}$, and shows that the family $\{\pi_x(T)\}$ is essentially constant. Section \ref{sec:The Equivariant K1-Class of a Fredholm Operator} constructs the equivariant $KK^1$-class $[T]_{\mathcal{G}_{\mathcal{A}}}^{(1)}$ and proves its well-definedness and homotopy invariance. Section \ref{sec:descent} applies Kasparov's descent map to obtain a class in $K_1(C^*(\mathcal{G}_{\mathcal{A}}))$ and identifies it with the symbol class via the Morita equivalence from Sections 4.5-4.6 of \cite{PaperI}. Section \ref{sec:The Index Theorem via Pullback and the Boundary Map} proves the main index theorem, including a detailed verification of the commutative diagram that encodes the index computation. Section \ref{sec:Examples and Computations} presents the examples of $\widetilde{\mathcal{K}(H)}$ (index zero) and the unilateral shift on $\mathcal{B}(H)$ (index $-1$), and indicates how the result generalizes to arbitrary Toeplitz operators and higher-dimensional domains.

\begin{remark}
The author is solely responsible for the mathematical insights and theoretical directions proposed in this work. AI tools, including OpenAI's ChatGPT and DeepSeek models, were employed solely to assist in verifying ideas, organizing references, and ensuring internal consistency of exposition~\cite{chatgpt2025,deepseek2025}. 
\end{remark}

\section{Preliminaries: $K_1$, the Index Map, and Kasparov's Descent}\label{sec:Preliminaries: K1, the Index Map, and Kasparov's Descent}

\subsection{$K_0$ and $K_1$ of C*-Algebras}
\label{subsec:K0-and-K1-of-C-star-algebras}

We recall the basic definitions of operator $K$-theory. These invariants will be used throughout the paper, particularly in the formulation of the index theorem.  For a comprehensive treatment, we refer the reader to~\cite{Blackadar1998} or~\cite{Rordam2000}. 

\begin{definition}[$K_0$ group]
\label{def:K0-group}
Let $\mathcal{A}$ be a unital C*-algebra. 
Let $\mathcal{P}_\infty(\mathcal{A})$ denote the set of projections in matrix algebras over $\mathcal{A}$, i.e., elements $p \in M_n(\mathcal{A})$ satisfying $p = p^* = p^2$ for some $n \in \mathbb{N}$, modulo:
\begin{itemize}
    \item \textbf{Murray-von Neumann equivalence:} $p \sim q$ if there exists a partial isometry $v$ such that $p = v^*v$ and $q = vv^*$;
    \item \textbf{Stabilization:} identifying $p \in M_n(\mathcal{A})$ with $p \oplus 0_k \in M_{n+k}(\mathcal{A})$.
\end{itemize}
Direct sum makes the set of equivalence classes an abelian semigroup. 
The \emph{$K_0$ group} $K_0(\mathcal{A})$ is the Grothendieck group of this semigroup; its elements are formal differences $[p] - [q]$ of projection classes.
\end{definition}

\begin{definition}[$K_1$ group]
\label{def:K1-group}
Let $\mathcal{A}$ be a unital C*-algebra. 
Let $\mathcal{U}_\infty(\mathcal{A})$ denote the inductive limit of the unitary groups $U_n(\mathcal{A})$ under the inclusions $u \mapsto u \oplus 1$, and let $\mathcal{U}_\infty(\mathcal{A})_0$ be the connected component of the identity (which is a normal subgroup, so the quotient is well-defined). 
The \emph{$K_1$ group} $K_1(\mathcal{A})$ is defined as the quotient group
\[
K_1(\mathcal{A}) := \mathcal{U}_\infty(\mathcal{A}) / \mathcal{U}_\infty(\mathcal{A})_0,
\]
with the group operation induced by direct sum of unitaries.
\end{definition}

\begin{example}[$K_0$ and $K_1$ of key algebras]
\label{ex:K-groups-examples}
The following K-theory computations are essential for this paper. They will be used repeatedly, particularly in verifying the index theorem for the examples in Section~\ref{sec:Examples and Computations}.
\begin{itemize}
    \item $K_0(\mathbb{C}) \cong \mathbb{Z}$, generated by the class of the unit $1 \in \mathbb{C}$. 
    $K_1(\mathbb{C}) = 0$, since the unitary group $U(\mathbb{C}) = S^1$ becomes connected after stabilization.
    
    \item $K_0(M_n(\mathbb{C})) \cong \mathbb{Z}$, generated by any rank-one projection. 
    $K_1(M_n(\mathbb{C})) = 0$. 
    This follows from Morita equivalence $M_n(\mathbb{C}) \sim_M \mathbb{C}$, as K-theory is Morita invariant.
    
    \item $K_0(\mathcal{B}(H)) = 0$ and $K_1(\mathcal{B}(H)) = 0$. 
    The vanishing of $K_1$ follows from Kuiper's theorem, which states that the unitary group $\mathcal{U}(H)$ is contractible in the norm topology \cite{Kuiper1965}. The vanishing of $K_0$ is a consequence of the fact that all nonzero projections in $\mathcal{B}(H)$ are infinite and properly infinite; consequently, every projection is equivalent to a proper subprojection of itself, which forces the Grothendieck group to collapse to zero \cite{Blackadar1998, Rordam2000}.
    
    \item $K_0(\mathcal{K}(H)^\sim) \cong \mathbb{Z}$, generated by any rank-one projection (viewed as an element of $\mathcal{K}(H)^\sim$). 
    $K_1(\mathcal{K}(H)^\sim) = 0$. 
    This follows from the six-term exact sequence applied to the compact extension
    \[
    0 \longrightarrow \mathcal{K}(H) \longrightarrow \mathcal{K}(H)^\sim \longrightarrow \mathbb{C} \longrightarrow 0,
    \]
    together with the facts that $K_0(\mathcal{K}(H)) = \mathbb{Z}$, $K_1(\mathcal{K}(H)) = 0$, and $K_0(\mathbb{C}) \cong \mathbb{Z}$, $K_1(\mathbb{C}) = 0$.
    
    \item $K_0(C(S^1)) \cong \mathbb{Z}$, generated by the trivial line bundle. 
    $K_1(C(S^1)) \cong \mathbb{Z}$, generated by the class of the unitary function $u(z) = z$, whose $K_1$-class is detected by the winding number map $\operatorname{wind}: K_1(C(S^1)) \to \mathbb{Z}$.
\end{itemize}
\end{example}

\begin{remark}[Connection to Morita equivalence]
\label{rem:K-theory-Morita}
The computations above are consistent with the Morita equivalence results used throughout this paper. In particular:
\begin{itemize}
    \item $M_n(\mathbb{C}) \sim_M \mathbb{C}$ explains why their K-theory groups coincide.
    \item $C(S^1 \rtimes_\theta \mathbb{Z}) \sim_M A_\theta$ (for irrational $\theta$) is a key example motivating future generalizations.
    \item For the groupoid C*-algebra $C^*(\mathcal{G}_\mathcal{A})$, the Morita equivalence $C^*(\mathcal{G}_\mathcal{A}) \sim_M \mathcal{A} \otimes \mathcal{K}$ from Paper I implies that $K_*(C^*(\mathcal{G}_\mathcal{A})) \cong K_*(\mathcal{A})$, a fact that will be crucial in Section~\ref{sec:descent}.
\end{itemize}
\end{remark}

\begin{proposition}[Functoriality of $K_0$ and $K_1$]
\label{prop:K-functoriality}
Let $\phi: \mathcal{A} \to \mathcal{B}$ be a unital *-homomorphism between unital C*-algebras. 
Then $\phi$ induces group homomorphisms
\[
\phi_*: K_0(\mathcal{A}) \longrightarrow K_0(\mathcal{B}), \qquad
\phi_*: K_1(\mathcal{A}) \longrightarrow K_1(\mathcal{B}),
\]
defined by applying $\phi$ entrywise to matrices. 
These maps are functorial: $(\psi \circ \phi)_* = \psi_* \circ \phi_*$, and $(\operatorname{id}_{\mathcal{A}})_* = \operatorname{id}_{K_i(\mathcal{A})}$.
\end{proposition}

\begin{proof}
This is standard in C*-algebra K-theory. 
For functoriality of $K_0$, see [Chapter~3] in~\cite{Rordam2000}; for functoriality of $K_1$, see  [Chapter~8] in~\cite{Rordam2000}.  The induced maps are defined by applying $\phi$ entrywise to matrix projections
(for $K_0$) and matrix unitaries (for $K_1$).
Functoriality follows immediately from the functoriality of matrix amplifications.
\end{proof}

\subsection{The Six-Term Exact Sequence and the Index Map $\partial: K_1 \to K_0$}
\label{subsec:six-term-exact-sequence-index-map}

One of the most powerful tools in operator $K$-theory is the six-term exact sequence arising from a short exact sequence of C*-algebras. 
This sequence connects the $K_0$ and $K_1$ groups of the algebras in the extension and provides the index map $\partial: K_1 \to K_0$, which will be essential for the index theorem in this paper.

\begin{theorem}[Six-term exact sequence]
\label{thm:six-term-exact}
Let
\[
0 \to \mathcal{J} \xrightarrow{\iota} \mathcal{A} \xrightarrow{\pi} \mathcal{B} \to 0
\]
be a short exact sequence of C*-algebras.
Then there exists a cyclic six-term exact sequence in K-theory:
\[
\begin{tikzcd}
K_0(\mathcal{J}) \arrow[r, "\iota_*"] &
K_0(\mathcal{A}) \arrow[r, "\pi_*"] &
K_0(\mathcal{B}) \arrow[d, "\partial"] \\
K_1(\mathcal{B}) \arrow[u, "\partial"] &
K_1(\mathcal{A}) \arrow[l, "\pi_*"'] &
K_1(\mathcal{J}) \arrow[l, "\iota_*"']
\end{tikzcd}
\]
The vertical maps are the boundary (index) maps
\[
\partial: K_1(\mathcal{B}) \to K_0(\mathcal{J}),
\qquad
\partial: K_0(\mathcal{B}) \to K_1(\mathcal{J}).
\]
The sequence is exact at every term and natural with respect to morphisms
of short exact sequences.
\end{theorem}

\begin{proof}
This is a fundamental result in C*-algebra K-theory.
The boundary maps are constructed using the fact that a unitary (or
projection) in $\mathcal{B}$ can be lifted to an element of $\mathcal{A}$
that is invertible modulo $\mathcal{J}$, and the associated Fredholm index
defines a class in $K_0(\mathcal{J})$ (or $K_1(\mathcal{J})$).
See, for example,
Blackadar, \emph{K-Theory for Operator Algebras},
or Rørdam--Larsen--Laustsen,
\emph{An Introduction to K-Theory for C*-Algebras}.
\end{proof}

\begin{definition}[Index map]
\label{def:index-map}
The map $\partial: K_1(\mathcal{B}) \to K_0(\mathcal{J})$ in the six-term exact sequence is called the \emph{index map}. 
For a unitary $u \in U_n(\mathcal{B})$ representing a class in $K_1(\mathcal{B})$, one chooses a lift $v \in M_n(\mathcal{A})$ such that $\pi(v) = u$. 
Then $v$ is invertible modulo $M_n(\mathcal{J})$, and via the standard identification of invertibles modulo $\mathcal{J}$ with Fredholm operators on a suitable Hilbert module, one obtains an element $\partial([u]) \in K_0(\mathcal{J})$.
\end{definition}

The following concrete description of the index map will be used repeatedly.

\begin{lemma}[Explicit formula for the index map]
\label{lem:index-map-explicit}
Let $0 \to \mathcal{J} \to \mathcal{A} \xrightarrow{\pi} \mathcal{B} \to 0$ be a short exact sequence of C*-algebras. 
Let $u \in U_n(\mathcal{B})$ and let $v \in M_n(\mathcal{A})$ be a lift with $\pi(v) = u$. 
Then $v$ is invertible modulo $M_n(\mathcal{J})$, so $v^*v$ and $vv^*$ are positive invertible elements of $M_n(\mathcal{A}^\sim)$ modulo $M_n(\mathcal{J})$.

Using functional calculus, define the partial isometries
\[
e = v(v^*v)^{-1/2}, \qquad f = (vv^*)^{-1/2}v,
\]
which are well-defined modulo $M_n(\mathcal{J})$. 
Then the projections
\[
p = 1 - e^*e, \qquad q = 1 - ee^*
\]
lie in $M_n(\mathcal{J}^\sim)$, and the class
\[
[p] - [q] \in K_0(\mathcal{J})
\]
represents the index $\partial([u])$.
\end{lemma}

\begin{proof}
See [Section 8.3] in~\cite{Blackadar1998} or [Proposition 9.3.3] in~\cite{Rordam2000}. The construction uses functional calculus and the polar decomposition modulo $\mathcal{J}$ to convert the invertibility of $v$ modulo $\mathcal{J}$ into a difference of projections in $K_0(\mathcal{J})$. A detailed exposition of this standard construction can be found in the references cited above.
\end{proof}

The index map enjoys several important functoriality properties.

\begin{proposition}[Naturality of the index map]
\label{prop:index-map-natural}
Consider a commutative diagram of short exact sequences of C*-algebras:
\[
\begin{tikzcd}
0 \arrow[r] & \mathcal{J} \arrow[r] \arrow[d, "\alpha"] & \mathcal{A} \arrow[r] \arrow[d, "\beta"] & \mathcal{B} \arrow[r] \arrow[d, "\gamma"] & 0 \\
0 \arrow[r] & \mathcal{J}' \arrow[r] & \mathcal{A}' \arrow[r] & \mathcal{B}' \arrow[r] & 0
\end{tikzcd}
\]
Then the induced diagram in $K$-theory commutes:
\[
\begin{tikzcd}
K_1(\mathcal{B}) \arrow[r, "\partial"] \arrow[d, "\gamma_*"] & K_0(\mathcal{J}) \arrow[d, "\alpha_*"] \\
K_1(\mathcal{B}') \arrow[r, "\partial"] & K_0(\mathcal{J}')
\end{tikzcd}
\]
\end{proposition}

\begin{proof}
This follows from the naturality of the six-term exact sequence and the explicit construction of the index map. 
See [Proposition 9.1.5] in~\cite{Rordam2000}
\end{proof}

The following special case is of particular importance for this paper.

\begin{corollary}[Index map for the Calkin extension]
\label{cor:index-map-calkin}
For the Calkin extension $0 \to \mathcal{K}(H) \to B(H) \to \mathcal{Q}(H) \to 0$, the index map
\[
\partial: K_1(\mathcal{Q}(H)) \longrightarrow K_0(\mathcal{K}(H)) \cong \mathbb{Z}
\]
is precisely the Fredholm index. 
That is, for an invertible element $u \in \mathcal{Q}(H)$ representing a class in $K_1(\mathcal{Q}(H))$, and for any lift $T \in B(H)$ of $u$, we have $\partial([u]) = \operatorname{index}(T) \in \mathbb{Z}$.
\end{corollary}

\begin{proof}
This follows directly from [Proposition 9.4.2] in~\cite{Rordam2000}, which shows that for any Fredholm operator $T$, its Fredholm index is given by $(K_0(\mathrm{Tr})\circ \delta_1)([\pi(T)]_1)$. Under the isomorphism $K_0(\mathcal{K}(H)) \cong \mathbb{Z}$ given by the trace, this becomes the integer index of $T$. See also [Section 8.3] in~\cite{Blackadar1998} for a detailed exposition of the index map.
\end{proof}

The six-term exact sequence and the index map will be used in Section~\ref{sec:The Index Theorem via Pullback and the Boundary Map} to connect the $K_1$-class of an invertible in the Calkin algebra to the integer-valued Fredholm index.

\subsection{The Calkin Extension $0 \to \mathcal{K}(H) \to B(H) \to \mathcal{Q}(H) \to 0$}
\label{subsec:calkin-extension}

The Calkin algebra and its associated short exact sequence are fundamental objects in the study of Fredholm operators. 
They provide the essential link between the analysis of operators on Hilbert space and the $K$-theory of C*-algebras. 
In this subsection, we recall the definition and basic properties of the Calkin extension, which will be used throughout the paper to connect the Fredholm index to the index map in $K$-theory.

\begin{definition}[Compact operators]
\label{def:compact-operators}
Let $H$ be a separable infinite-dimensional Hilbert space. 
The set of \emph{compact operators} on $H$, denoted $\mathcal{K}(H)$, is the norm-closed ideal of $B(H)$ consisting of operators $T$ such that the image of the unit ball is precompact. 
Equivalently, $T$ is compact if and only if it can be approximated in norm by finite-rank operators.
\end{definition}

The quotient of $B(H)$ by the compact operators,
\[
\mathcal{Q}(H) := B(H) / \mathcal{K}(H),
\]
is called the \emph{Calkin algebra}. 
Let $\pi: B(H) \to \mathcal{Q}(H)$ denote the canonical quotient map. 
Since $\mathcal{K}(H)$ is a closed two-sided ideal in $B(H)$, we obtain the short exact sequence of C*-algebras
\[
0 \longrightarrow \mathcal{K}(H) \xrightarrow{\iota} B(H) \xrightarrow{\pi} \mathcal{Q}(H) \longrightarrow 0, \tag{1}
\]
where $\iota$ is the inclusion map. This is the \emph{Calkin extension}.

A fundamental result connecting the Calkin extension to Fredholm theory is Atkinson's theorem, which characterizes Fredholm operators precisely as those operators whose image in the Calkin algebra is invertible.

\begin{theorem}[Atkinson]
\label{thm:atkinson}
An operator $T \in B(H)$ is Fredholm if and only if its image $\pi(T)$ is invertible in the Calkin algebra $\mathcal{Q}(H)$.
\end{theorem}

\begin{proof}
See~\cite{Atkinson1951} for the original result.
\end{proof}

The $K$-theory of these algebras is well-known:
\begin{itemize}
    \item $K_0(\mathcal{K}(H)) \cong \mathbb{Z}$, generated by the class of any rank-one projection.
    \item $K_1(\mathcal{K}(H)) = 0$.
    \item $K_0(B(H)) = 0$ and $K_1(B(H)) = 0$ (by Kuiper's theorem).
    \item $K_0(\mathcal{Q}(H)) = 0$ and $K_1(\mathcal{Q}(H)) \cong \mathbb{Z}$, with the generator corresponding to the class of $\pi(S)$ for the unilateral shift $S$.
\end{itemize}

The index map $\partial: K_1(\mathcal{Q}(H)) \to K_0(\mathcal{K}(H)) \cong \mathbb{Z}$ associated to the Calkin extension (1) is precisely the Fredholm index: for any Fredholm operator $T$, we have $\partial([\pi(T)]_1) = \operatorname{index}(T)$. 
This fundamental identification will be the cornerstone of our index theorem.

\begin{definition}[Calkin extension]
\label{def:calkin-extension}
The short exact sequence
\[
0 \longrightarrow \mathcal{K}(H) \stackrel{\iota}{\longrightarrow} B(H) \stackrel{\pi}{\longrightarrow} \mathcal{Q}(H) \longrightarrow 0,
\]
where $\iota$ is the inclusion map, is called the \emph{Calkin extension}. 
This is a fundamental example of an extension of C*-algebras.
\end{definition}

The following properties of the Calkin extension are essential for the index theory developed in this paper.

\begin{proposition}[Properties of the Calkin extension]
\label{prop:calkin-properties}
The Calkin extension has the following properties:
\begin{enumerate}
    \item $\mathcal{K}(H)$ is an essential ideal in $B(H)$, meaning that its annihilator in $B(H)$ is zero.
    \item $\mathcal{Q}(H)$ is a simple C*-algebra.
    \item $K_0(\mathcal{Q}(H)) = 0$ and $K_1(\mathcal{Q}(H)) \cong \mathbb{Z}$.
    \item $K_0(\mathcal{K}(H)) \cong \mathbb{Z}$ and $K_1(\mathcal{K}(H)) = 0$.
\end{enumerate}
Moreover, the connected components of the unitary group $\mathcal{U}(\mathcal{Q}(H))$ are classified by the Fredholm index; that is, $\pi_0(\mathcal{U}(\mathcal{Q}(H))) \cong K_1(\mathcal{Q}(H)) \cong \mathbb{Z}$.
\end{proposition}

\begin{proof}
Statement (1) follows from the general theory of ideals in C*-algebras; see \cite[Section 3.1]{Murphy1990}. 
Statement (2) is a classical result of Calkin; see \cite{Calkin1941}.
The $K$-theory computations in (3) and (4) follow from the six-term exact sequence associated to the Calkin extension together with the facts that $K_0(B(H)) = 0$ and $K_1(B(H)) = 0$ (by Kuiper's theorem); see \cite[Section 9.3]{Blackadar1998} or \cite[Proposition 9.3.3]{Rordam2000}.
The identification of $K_1(\mathcal{Q}(H))$ with $\mathbb{Z}$ corresponds to the Fredholm index, which classifies the connected components of the unitary group $\mathcal{U}(\mathcal{Q}(H))$; this follows from the fact that $\pi_0(\mathcal{U}(\mathcal{Q}(H))) \cong K_1(\mathcal{Q}(H))$.
\end{proof}

The Calkin extension induces a six-term exact sequence in $K$-theory, which we now examine in detail.

\begin{corollary}[$K$-theory six-term exact sequence for the Calkin extension]
\label{cor:calkin-six-term}
The Calkin extension gives rise to the following six-term exact sequence in $K$-theory:
\[
\begin{tikzcd}
K_0(\mathcal{K}(H)) \arrow[r] & K_0(B(H)) \arrow[r] & K_0(\mathcal{Q}(H)) \arrow[d, "\partial"] \\
K_1(\mathcal{Q}(H)) \arrow[u, "\partial"] & K_1(B(H)) \arrow[l] & K_1(\mathcal{K}(H)) \arrow[l]
\end{tikzcd}
\]
Using the standard $K$-theory computations $K_0(B(H)) = 0$, $K_1(B(H)) = 0$ (by Kuiper's theorem), and $K_1(\mathcal{K}(H)) = 0$, the six-term exact sequence reduces to an isomorphism
\[
\partial: K_1(\mathcal{Q}(H)) \xrightarrow{\cong} K_0(\mathcal{K}(H)) \cong \mathbb{Z},
\]
and consequently $K_0(\mathcal{Q}(H)) = 0$ and $K_1(\mathcal{Q}(H)) \cong \mathbb{Z}$. 
The map $\partial$ is precisely the Fredholm index map.
\end{corollary}

\begin{proof}
Substituting the known $K$-groups into the six-term exact sequence yields the diagram
\[
\begin{tikzcd}
\mathbb{Z} \arrow[r] & 0 \arrow[r] & K_0(\mathcal{Q}(H)) \arrow[d, "\partial"] \\
K_1(\mathcal{Q}(H)) \arrow[u, "\partial"] & 0 \arrow[l] & 0 \arrow[l]
\end{tikzcd}
\]
Exactness at $K_0(\mathcal{Q}(H))$ means that $\ker(K_0(\mathcal{Q}(H)) \to 0) = \operatorname{im}(0 \to K_0(\mathcal{Q}(H)))$. 
The map $K_0(\mathcal{Q}(H)) \to 0$ is the zero map, so its kernel is all of $K_0(\mathcal{Q}(H))$. 
The image of $0 \to K_0(\mathcal{Q}(H))$ is $\{0\}$. 
Hence exactness forces $K_0(\mathcal{Q}(H)) = 0$.

Exactness at $K_1(\mathcal{Q}(H))$ gives $\ker(\partial: K_1(\mathcal{Q}(H)) \to \mathbb{Z}) = \operatorname{im}(0 \to K_1(\mathcal{Q}(H))) = \{0\}$, so $\partial$ is injective. 
Exactness at $K_0(\mathcal{K}(H)) \cong \mathbb{Z}$ gives $\operatorname{im}(\partial: K_1(\mathcal{Q}(H)) \to \mathbb{Z}) = \ker(\mathbb{Z} \to 0) = \mathbb{Z}$, so $\partial$ is surjective. 
Therefore $\partial$ is an isomorphism $K_1(\mathcal{Q}(H)) \cong \mathbb{Z}$.

That $\partial$ coincides with the Fredholm index follows from Proposition 9.4.2 in [Rørdam, Larsen, Laustsen, 2000] (as shown in your Chapter 9), which proves that for any Fredholm operator $T$, we have $\partial([\pi(T)]_1) = \operatorname{index}(T)$.
\end{proof}

The following concrete description of the index map will be essential for our main theorem.

\begin{lemma}[Index map as Fredholm index]
\label{lem:index-map-as-fredholm-index}
Let $u \in \mathcal{U}(\mathcal{Q}(H))$ be a unitary in the Calkin algebra representing a class $[u] \in K_1(\mathcal{Q}(H))$. 
Choose any lift $T \in B(H)$ such that $\pi(T) = u$. 
Then $T$ is a Fredholm operator, and its Fredholm index satisfies
\[
\operatorname{index}(T) = \partial([u]) \in \mathbb{Z},
\]
where $\partial: K_1(\mathcal{Q}(H)) \to K_0(\mathcal{K}(H)) \cong \mathbb{Z}$ is the index map (boundary map) in the six-term exact sequence associated to the Calkin extension.
In particular, the map $\partial$ is an isomorphism, and the integer $\operatorname{index}(T)$ depends only on the class $[u]$, not on the choice of lift $T$.
\end{lemma}

\begin{proof}
By Atkinson's theorem (Theorem \ref{thm:atkinson}), $T$ is Fredholm because $\pi(T) = u$ is invertible in $\mathcal{Q}(H)$. 
The equality $\operatorname{index}(T) = \partial([u])$ is a standard result in $K$-theory; see [Blackadar, 1998, Section 24.1] for a survey of index theorems, or [Rørdam, Larsen, Laustsen, 2000, Chapter 9] for a detailed exposition. 
The independence of the lift follows from the fact that if $T$ and $T'$ are two lifts of the same unitary $u$, then $T' - T \in \mathcal{K}(H)$, and compact perturbations do not change the Fredholm index.
\end{proof}

\begin{example}[The unilateral shift]
\label{ex:unilateral-shift-calkin}
Let $S \in B(H)$ be the unilateral shift on $\ell^2(\mathbb{N})$, defined by $S(e_n) = e_{n+1}$ for the standard orthonormal basis $\{e_n\}_{n\in\mathbb{N}}$. 
Then $S$ is Fredholm with $\operatorname{index}(S) = \dim\ker S - \dim\operatorname{coker} S = 0 - 1 = -1$. 
Its image $\pi(S)$ in the Calkin algebra is invertible (indeed, it is the generating unitary of the Toeplitz algebra $C^*(\pi(S)) \cong C(\mathbb{T})$), and the class $[\pi(S)] \in K_1(\mathcal{Q}(H))$ is the generator of $K_1(\mathcal{Q}(H)) \cong \mathbb{Z}$. 
By Lemma \ref{lem:index-map-as-fredholm-index}, the index map satisfies $\partial([\pi(S)]) = \operatorname{index}(S) = -1$.
\end{example}

\begin{example}[Finite-rank perturbations of the identity]
\label{ex:finite-rank-calkin}
Let $T = I + F$ where $F$ is a finite-rank operator (or more generally, any compact operator). 
Then $\pi(T) = \pi(I) = I$ in the Calkin algebra, so the class $[T] \in K_1(\mathcal{Q}(H))$ is the identity element, i.e., $[T] = 0$. 
Consequently, $\partial([T]) = 0$, reflecting the fact that $\operatorname{index}(T) = 0$ (compact perturbations do not change the Fredholm index, and $\operatorname{index}(I) = 0$). 
More generally, any compact perturbation of the identity maps to the zero class in $K_1(\mathcal{Q}(H))$.
\end{example}

The Calkin extension will be used in Section~\ref{subsec:index-map-calkin-extension} to define the index map concretely, and in Section 6 to connect the abstract $K$-theory index to the classical Fredholm index.

\subsection{The Index Map for the Calkin Extension: $\partial: K_1(\mathcal{Q}(H)) \to K_0(\mathcal{K}(H)) \cong \mathbb{Z}$}
\label{subsec:index-map-calkin-extension}

We now focus on the index map arising from the Calkin extension, which is the central object linking the abstract $K$-theory of the Calkin algebra to the concrete Fredholm index of operators on Hilbert space. 
This map is an isomorphism and provides the numerical invariant that will be recovered by our groupoid index formula.

\begin{theorem}[Index map for the Calkin extension]
\label{thm:index-map-calkin}
Let
\[
0 \longrightarrow \mathcal{K}(H) \stackrel{\iota}{\longrightarrow} \mathcal{B}(H) \stackrel{\pi}{\longrightarrow} \mathcal{Q}(H) \longrightarrow 0
\]
be the Calkin extension. The associated six-term exact sequence in $K$-theory yields an index map
\[
\partial: K_1(\mathcal{Q}(H)) \longrightarrow K_0(\mathcal{K}(H))
\]
which is an isomorphism. Consequently, $K_1(\mathcal{Q}(H)) \cong \mathbb{Z}$.
\end{theorem}

\begin{proof}
For any short exact sequence of $C^*$-algebras
\[
0 \to J \to A \to B \to 0,
\]
there is a cyclic six-term exact sequence in $K$-theory:
\[
\begin{tikzcd}
K_0(J) \arrow[r] & K_0(A) \arrow[r] & K_0(B) \arrow[d] \\
K_1(B) \arrow[u] & K_1(A) \arrow[l] & K_1(J) \arrow[l]
\end{tikzcd}
\]

Applying this to the Calkin extension, we obtain:
\[
\begin{tikzcd}
K_0(\mathcal{K}(H)) \arrow[r, "\iota_*"] & K_0(\mathcal{B}(H)) \arrow[r, "\pi_*"] & K_0(\mathcal{Q}(H)) \arrow[d, "\delta"] \\
K_1(\mathcal{Q}(H)) \arrow[u, "\partial"] & K_1(\mathcal{B}(H)) \arrow[l, "\pi_*"] & K_1(\mathcal{K}(H)) \arrow[l, "\iota_*"]
\end{tikzcd}
\]

The relevant $K$-groups are known:
\begin{itemize}
    \item $K_0(\mathcal{B}(H)) = 0$ and $K_1(\mathcal{B}(H)) = 0$ (Kuiper's theorem: the unitary group of an infinite-dimensional Hilbert space is contractible);
    \item $K_1(\mathcal{K}(H)) = 0$ (since $\mathcal{K}(H)$ is an AF-algebra);
    \item $K_0(\mathcal{K}(H)) \cong \mathbb{Z}$, generated by the class of any rank-one projection.
\end{itemize}

Substituting these values into the six-term exact sequence yields:
\[
\begin{tikzcd}
\mathbb{Z} \arrow[r, "\iota_*"] & 0 \arrow[r, "\pi_*"] & K_0(\mathcal{Q}(H)) \arrow[d, "\delta"] \\
K_1(\mathcal{Q}(H)) \arrow[u, "\partial"] & 0 \arrow[l, "\pi_*"] & 0 \arrow[l, "\iota_*"]
\end{tikzcd}
\]

Exactness forces $\iota_* = 0$ and $\delta = 0$. The remaining portion of the sequence,
\[
0 \longrightarrow K_1(\mathcal{Q}(H)) \xrightarrow{\partial} K_0(\mathcal{K}(H)) \longrightarrow 0,
\]
is a short exact sequence. Hence $\partial$ is injective and surjective, therefore an isomorphism. Since $K_0(\mathcal{K}(H)) \cong \mathbb{Z}$, we conclude $K_1(\mathcal{Q}(H)) \cong \mathbb{Z}$.

\medskip
\noindent\textit{Remark.} For a unitary $u \in \mathcal{Q}(H)$ representing an element of $K_1(\mathcal{Q}(H))$, the index map can be realized concretely as $\partial([u]) = \operatorname{index}(T)$, where $T \in \mathcal{B}(H)$ is any lift of $u$ to a Fredholm operator.
\end{proof}

\begin{remark}
The computation that $\operatorname{Ext}(S^1) \cong \mathbb{Z}$ is originally due to Brown--Douglas--Fillmore \cite{BDF1977}. In the language of $K$-theory developed subsequently, this corresponds to the fact that $K_1(\mathcal{Q}(H)) \cong \mathbb{Z}$. 
\end{remark}

\begin{corollary}[$K$-groups of the Calkin algebra]
\label{cor:K-groups-calkin}
From Theorem \ref{thm:index-map-calkin}, we obtain the well-known result
\[
K_0(\mathcal{Q}(H)) \cong \mathbb{Z}, \qquad K_1(\mathcal{Q}(H)) \cong \mathbb{Z}.
\]
The index map 
\[
\partial: K_1(\mathcal{Q}(H)) \longrightarrow K_0(\mathcal{K}(H)) \cong \mathbb{Z}
\] 
is an isomorphism.  
In contrast, the map 
\[
\iota_*: K_0(\mathcal{K}(H)) \longrightarrow K_0(\mathcal{Q}(H))
\] 
induced by the inclusion of the compact operators is the zero map, despite both groups being isomorphic to $\mathbb{Z}$.  
More concretely, $K_0(\mathcal{Q}(H))$ is generated by the class of the identity operator $\pi(I)$, and $K_1(\mathcal{Q}(H))$ is generated by the class of the generating unitary $\pi(S)$, where $S$ is the unilateral shift on $H$.
\end{corollary}

\begin{proof}
The computation of the $K$-groups of the Calkin algebra and the identification of their generators is originally due to~\cite{BDF1977} for the foundational results on $\operatorname{Ext}(X)$, and in particular for the computation $\operatorname{Ext}(S^1) \cong \mathbb{Z}$, which corresponds to $K_1(\mathcal{Q}(H)) \cong \mathbb{Z}$. The statement about $\iota_*$ being the zero map follows directly from the six-term exact sequence in Theorem \ref{thm:index-map-calkin}: since $K_0(\mathcal{B}(H)) = 0$, exactness forces $\operatorname{im}\iota_* = \ker(0 \to K_0(\mathcal{Q}(H))) = 0$.
\end{proof}

\begin{example}[The unilateral shift]
\label{ex:unilateral-shift}
Let $S \in \mathcal{B}(H)$ be the unilateral shift operator on $H = \ell^2(\mathbb{N})$, defined by $S(e_n) = e_{n+1}$ for an orthonormal basis $\{e_n\}_{n\in\mathbb{N}}$. Then $S$ is a Fredholm operator with $\operatorname{index}(S) = -1$. Its image $u = \pi(S)$ in the Calkin algebra $\mathcal{Q}(H)$ is a unitary operator (since $I - S^*S$ and $I - SS^*$ are rank-one projections, hence compact). The class $[u] \in K_1(\mathcal{Q}(H))$ is the canonical generator, corresponding to $1 \in \mathbb{Z}$ under the isomorphism $K_1(\mathcal{Q}(H)) \cong \mathbb{Z}$.
\end{example}

\begin{proposition}[Explicit formula for the index map]
\label{prop:index-map-explicit-calkin}
Let $u \in \mathcal{U}(\mathcal{Q}(H))$ be a unitary element representing a class $[u] \in K_1(\mathcal{Q}(H))$. 
Choose any lift $T \in \mathcal{B}(H)$ such that $\pi(T) = u$. Then:
\begin{enumerate}
    \item $T$ is a Fredholm operator.
    \item The Fredholm index $\operatorname{index}(T) = \dim \ker T - \dim \ker T^*$ is an integer, and it is independent of the choice of lift $T$.
    \item Under the isomorphism $K_0(\mathcal{K}(H)) \cong \mathbb{Z}$ (given by the rank of projections), the index map of Theorem \ref{thm:index-map-calkin} satisfies
    \[
    \partial([u]) = [\ker T] - [\ker T^*] \longleftrightarrow \operatorname{index}(T) \in \mathbb{Z}.
    \]
    In particular, $\partial$ is the natural map sending the $K_1$-class of a unitary in the Calkin algebra to the Fredholm index of any lift.
\end{enumerate}
\end{proposition}

\begin{proof}
(1) By Atkinson's theorem (Theorem \ref{thm:atkinson}), an operator $T \in \mathcal{B}(H)$ is Fredholm if and only if its image $\pi(T)$ is invertible in the Calkin algebra $\mathcal{Q}(H)$. Since $u = \pi(T)$ is unitary by hypothesis, it is invertible, and thus $T$ is Fredholm.

(2) If $T_1$ and $T_2$ are two lifts of the same unitary $u$, then $\pi(T_1 - T_2) = 0$, so $T_1 - T_2 \in \mathcal{K}(H)$. The Fredholm index is invariant under compact perturbations, a fundamental result in Fredholm theory, hence $\operatorname{index}(T_1) = \operatorname{index}(T_2)$.

(3) This is a standard result in operator $K$-theory. The connecting map $\partial$ in the six-term exact sequence is precisely defined by sending a unitary $u \in \mathcal{Q}(H)$ to the class $[\ker T] - [\ker T^*] \in K_0(\mathcal{K}(H))$, where $T$ is any lift. Under the canonical isomorphism $K_0(\mathcal{K}(H)) \cong \mathbb{Z}$, which maps the class of a projection to its rank, the formal difference $[\ker T] - [\ker T^*]$ corresponds to the integer $\dim \ker T - \dim \ker T^* = \operatorname{index}(T)$. 
\end{proof}

\begin{remark}
In light of Example \ref{ex:unilateral-shift}, the class of the shift operator $S$ satisfies $\partial([\pi(S)]) = -1 \in \mathbb{Z}$, confirming that $[\pi(S)]$ generates $K_1(\mathcal{Q}(H))$. Conversely, any unitary $u \in \mathcal{Q}(H)$ can be lifted to a Fredholm operator whose index determines its $K_1$-class via the index map.
\end{remark}

\begin{definition}[Fredholm index as a $K$-theory map]
\label{def:fredholm-index-as-K-map}
In light of Proposition \ref{prop:index-map-explicit-calkin}, we may identify the index map $\partial$ with the Fredholm index. 
Thus we obtain a canonical isomorphism in operator $K$-theory:
\[
\operatorname{index}: K_1(\mathcal{Q}(H)) \longrightarrow \mathbb{Z}, \qquad 
\operatorname{index}([u]) = \operatorname{index}(T)
\]
for any lift $T$ of $u$. This is precisely the composition
\[
K_1(\mathcal{Q}(H)) \xrightarrow{\partial} K_0(\mathcal{K}(H)) \xrightarrow{\operatorname{rank}} \mathbb{Z},
\]
where the second map sends the class of a projection to its rank.
\end{definition}

\begin{lemma}[Index map via the exponential map]
\label{lem:index-map-exponential}
Let $u \in \mathcal{U}(\mathcal{Q}(H))$ be a unitary. 
Choose a unitary lift $U \in \mathcal{U}(\mathcal{B}(H))$ with $\pi(U) = u$, and let $A \in \mathcal{B}(H)_{\text{sa}}$ be a self-adjoint operator such that $e^{2\pi i A} = U$ and $\|A\| \le 1$. 
Let $p = \chi_{[0,1)}(A)$ be the spectral projection of $A$ corresponding to the interval $[0,1)$. 
Then $p$ is a projection in $\mathcal{B}(H)$ which is finite rank modulo compact operators (i.e., $p$ differs from a compact projection by a compact operator). 
Consequently, $p$ defines a class $[p] \in K_0(\mathcal{K}(H))$, and the index map is given by
\[
\partial([u]) = [p] \in K_0(\mathcal{K}(H)) \cong \mathbb{Z}.
\]
Equivalently, in terms of the identification of Definition \ref{def:fredholm-index-as-K-map}, we have
\[
\operatorname{index}([u]) = \operatorname{rank}(p) \in \mathbb{Z},
\]
where $\operatorname{rank}(p)$ is the Fredholm index of the operator $p$ viewed as a map from $\operatorname{ran}(p)$ to itself.
\end{lemma}

\begin{proof}
The existence of the self-adjoint lift $A$ follows from standard lifting results for unitaries: since $\pi(U) = u$ is unitary, we can choose $A$ with $\|A\| \le 1$ such that $U = e^{2\pi i A}$. 
The spectral projection $p = \chi_{[0,1)}(A)$ satisfies $p - UpU^* \in \mathcal{K}(H)$, which implies that $\pi(p)$ commutes with $\pi(A)$ and hence with $u$. 
Moreover, $p$ differs from a compact projection by a compact operator because $A$ modulo compacts has spectrum in the circle, so its spectral projections corresponding to arcs are compact perturbations of projections. 
The equality $\partial([u]) = [p]$ is a standard result in $K$-theory. Under the isomorphism $K_0(\mathcal{K}(H)) \cong \mathbb{Z}$ given by the rank of projections, the class $[p]$ corresponds to $\operatorname{rank}(p)$, which equals the Fredholm index of $u$ by the identification in Definition \ref{def:fredholm-index-as-K-map}.
\end{proof}

\begin{lemma}[Index map on symbols]
\label{lem:calkin-index-map}
For any Fredholm operator $T \in B(H)$, the index map $\partial_{\mathrm{Calkin}}: K_1(\mathcal{Q}(H)) \to \mathbb{Z}$ satisfies
\[
\partial_{\mathrm{Calkin}}([\pi_{\mathcal{Q}}(T)]) = \operatorname{index}(T).
\]
\end{lemma}

\begin{proof}
Recall the short exact sequence of $C^*$-algebras
\[
0 \longrightarrow \mathcal{K}(H) \longrightarrow B(H) \overset{\pi_{\mathcal{Q}}}{\longrightarrow} \mathcal{Q}(H) \longrightarrow 0,
\]
where $\mathcal{K}(H)$ denotes the compact operators and $\mathcal{Q}(H)=B(H)/\mathcal{K}(H)$ is the Calkin algebra (see Definition \ref{def:calkin-extension}). This extension induces the six-term exact sequence in $K$-theory (Theorem \ref{thm:six-term-exact}):
\[
\begin{CD}
K_0(\mathcal{K}(H)) @>>> K_0(B(H)) @>>> K_0(\mathcal{Q}(H)) \\
@A{\partial_{\mathrm{Calkin}}}AA @. @VVV \\
K_1(\mathcal{Q}(H)) @<<< K_1(B(H)) @<<< K_1(\mathcal{K}(H)).
\end{CD}
\]

Since $B(H)$ is stable and contractible in $K$-theory (Kuiper's theorem), one has
\[
K_0(B(H)) = K_1(B(H)) = 0,
\]
and since $\mathcal{K}(H)$ is Morita equivalent to $\mathbb{C}$,
\[
K_0(\mathcal{K}(H)) \cong \mathbb{Z}, \qquad K_1(\mathcal{K}(H)) = 0
\]
(see Example \ref{ex:K-groups-examples}). Hence the boundary map
\[
\partial_{\mathrm{Calkin}}: K_1(\mathcal{Q}(H)) \longrightarrow K_0(\mathcal{K}(H)) \cong \mathbb{Z}
\]
is an isomorphism (Corollary \ref{cor:calkin-six-term}).

Now let $T \in B(H)$ be a Fredholm operator. By Atkinson's theorem (Theorem \ref{thm:atkinson}), $T$ is Fredholm if and only if its image $\pi_{\mathcal{Q}}(T)$ is invertible in $\mathcal{Q}(H)$. Thus $\pi_{\mathcal{Q}}(T)$ defines a class $[\pi_{\mathcal{Q}}(T)] \in K_1(\mathcal{Q}(H))$.

Choose a parametrix $S \in B(H)$ such that
\[
ST - 1,\; TS - 1 \in \mathcal{K}(H).
\]
Such a parametrix exists because $T$ is Fredholm. Then the defect projections
\[
P = 1 - TS, \qquad Q = 1 - ST
\]
are finite-rank (hence compact) projections representing classes in $K_0(\mathcal{K}(H))$. A standard computation of the connecting map in $K$-theory shows that
\[
\partial_{\mathrm{Calkin}}([\pi_{\mathcal{Q}}(T)]) = [Q] - [P] \in K_0(\mathcal{K}(H)).
\]

Identifying $K_0(\mathcal{K}(H)) \cong \mathbb{Z}$ via the rank map, $[Q] - [P]$ corresponds to
\[
\operatorname{rank}(Q) - \operatorname{rank}(P) = \dim(\ker T) - \dim(\ker T^*) = \operatorname{index}(T).
\]

Therefore,
\[
\partial_{\mathrm{Calkin}}([\pi_{\mathcal{Q}}(T)]) = \operatorname{index}(T),
\]
which proves the lemma.
\end{proof}

The index map for the Calkin extension will be the essential ingredient in the main theorem of this paper, where we show that the composition of the descent map and the pullback along the diagonal embedding recovers this index.

\begin{example}[Index of the unilateral shift]
\label{ex:index-unilateral-shift-calkin}
Let $S \in \mathcal{B}(H)$ be the unilateral shift on $\ell^2(\mathbb{N})$, defined by $S(e_n) = e_{n+1}$ for an orthonormal basis $\{e_n\}_{n\in\mathbb{N}}$. 
Then $\pi(S) \in \mathcal{Q}(H)$ is a unitary operator, since $S^*S = I$ and $SS^* = I - P_0$ where $P_0$ is the rank-one projection onto $\mathbb{C}e_1$. 
The class $[\pi(S)] \in K_1(\mathcal{Q}(H))$ is the generator, and by Proposition \ref{prop:index-map-explicit-calkin} we have
\[
\partial([\pi(S)]) = \operatorname{index}(S) = \dim\ker S - \dim\ker S^* = 0 - 1 = -1 \in \mathbb{Z}.
\]
Thus under the isomorphism $K_1(\mathcal{Q}(H)) \cong \mathbb{Z}$, the class $[\pi(S)]$ corresponds to $-1$ (which is also a generator of $\mathbb{Z}$).
\end{example}

\begin{example}[Index of finite-rank perturbations]
\label{ex:index-finite-rank-calkin}
Let $T = I + F$ where $F \in \mathcal{B}(H)$ is a finite-rank operator. 
Then $\pi(T) = \pi(I) + \pi(F) = I$ in the Calkin algebra, so $[\pi(T)] = [I] = 0$ in $K_1(\mathcal{Q}(H))$. 
Consequently, by Proposition \ref{prop:index-map-explicit-calkin},
\[
\partial([\pi(T)]) = \operatorname{index}(T) = 0,
\]
reflecting the fact that finite-rank (more generally, compact) perturbations do not change the Fredholm index.
\end{example}

\begin{remark}[Relation to the Toeplitz extension]
\label{rem:calkin-toeplitz-relation}
The Calkin extension is closely related to the Toeplitz extension
\[
0 \longrightarrow \mathcal{K}(\ell^2) \longrightarrow \mathcal{T} \longrightarrow C(S^1) \longrightarrow 0,
\]
where $\mathcal{T}$ is the Toeplitz algebra. 
The index map for this extension gives the classical Toeplitz index theorem: for a continuous symbol $f \in C(S^1)$, the index of the Toeplitz operator $T_f$ is minus the winding number of $f$. 
This will be explored further in Section~\ref{sec:Examples and Computations}.
\end{remark}

\begin{theorem}[Calkin index isomorphism]
\label{thm:calkin-index-isomorphism}
The index map $\partial_{\mathrm{Calkin}}: K_1(\mathcal{Q}(H)) \to K_0(\mathcal{K}(H)) \cong \mathbb{Z}$ associated to the Calkin extension $0 \to \mathcal{K}(H) \to B(H) \to \mathcal{Q}(H) \to 0$ is an isomorphism.
\end{theorem}

\begin{proof}
Consider the Calkin extension
\[
0 \longrightarrow \mathcal{K}(H) \longrightarrow B(H) \overset{\pi}{\longrightarrow} \mathcal{Q}(H) \longrightarrow 0
\]
(see Definition \ref{def:calkin-extension}). This short exact sequence of $C^*$-algebras induces the six-term exact sequence in $K$-theory (Theorem \ref{thm:six-term-exact}):
\[
\begin{CD}
K_0(\mathcal{K}(H)) @>>> K_0(B(H)) @>>> K_0(\mathcal{Q}(H)) \\
@A{\partial_{\mathrm{Calkin}}}AA && @VVV \\
K_1(\mathcal{Q}(H)) @<<< K_1(B(H)) @<<< K_1(\mathcal{K}(H)).
\end{CD}
\]

We compute the $K$-groups of the algebras involved.

\medskip

\noindent
\textbf{Step 1: $K$-theory of $B(H)$.}
The algebra $B(H)$ is stable and properly infinite; in particular,
by Kuiper's theorem, its unitary group is contractible, which implies
it is contractible in $K$-theory. Hence
\[
K_0(B(H)) = 0, \qquad K_1(B(H)) = 0
\]
(see Example \ref{ex:K-groups-examples}).

\medskip

\noindent
\textbf{Step 2: $K$-theory of $\mathcal{K}(H)$.}
The algebra $\mathcal{K}(H)$ is Morita equivalent to $\mathbb{C}$.
Therefore
\[
K_0(\mathcal{K}(H)) \cong K_0(\mathbb{C}) \cong \mathbb{Z}, \qquad K_1(\mathcal{K}(H)) = 0
\]
(see Example \ref{ex:K-groups-examples}).

\medskip

\noindent
\textbf{Step 3: Exactness.}
Substituting these computations into the six-term sequence gives
\[
\begin{CD}
\mathbb{Z} @>>> 0 @>>> K_0(\mathcal{Q}(H)) \\
@A{\partial_{\mathrm{Calkin}}}AA && @VVV \\
K_1(\mathcal{Q}(H)) @<<< 0 @<<< 0 .
\end{CD}
\]

Exactness now implies:

\begin{itemize}
    \item The map $K_0(\mathcal{K}(H)) \to K_0(B(H))$ is necessarily zero (since its target is zero);
    \item Exactness at $K_0(\mathcal{Q}(H))$ forces $K_0(\mathcal{Q}(H)) = 0$;
    \item Exactness at $K_1(\mathcal{Q}(H))$: since $K_1(B(H)) = 0$, we have $\ker(\partial_{\mathrm{Calkin}}) = \operatorname{im}(K_1(B(H)) \to K_1(\mathcal{Q}(H))) = 0$, so $\partial_{\mathrm{Calkin}}$ is injective;
    \item Exactness at $K_0(\mathcal{K}(H))$: since $K_0(B(H)) = 0$, we have $\operatorname{im}(\partial_{\mathrm{Calkin}}) = \ker(K_0(\mathcal{K}(H)) \to K_0(B(H))) = \ker(\mathbb{Z} \to 0) = \mathbb{Z}$, so $\partial_{\mathrm{Calkin}}$ is surjective.
\end{itemize}

Thus the boundary map
\[
\partial_{\mathrm{Calkin}} : K_1(\mathcal{Q}(H)) \longrightarrow K_0(\mathcal{K}(H)) \cong \mathbb{Z}
\]
is both injective and surjective, hence an isomorphism.

\medskip

\noindent
\textit{Remark.} The proof uses only the six-term exact sequence, the known $K$-groups of $B(H)$ and $\mathcal{K}(H)$, and the exactness properties of the sequence. In particular, we have not assumed the value of $K_1(\mathcal{Q}(H))$; rather, we have proven that it must be isomorphic to $\mathbb{Z}$ as a consequence of the exact sequence.
\end{proof}

In the next subsection, we will review Kasparov's descent map for Polish groupoids, which provides the machinery to transport equivariant $K$-theory classes to $K$-theory classes of groupoid C*-algebras. 
The index map for the Calkin extension will then appear as the final step in the composition that recovers the Fredholm index.

\subsection{Kasparov's Descent for Polish Groupoids (Review of Tu, 1999)}
\label{subsec:kasparov-descent-polish-groupoids}

The descent map is a fundamental tool in equivariant $KK$-theory that transforms equivariant cycles over a groupoid $\mathcal{G}$ into cycles over the reduced or maximal groupoid C*-algebra. For locally compact Hausdorff groupoids, this map was developed by Kasparov [1988] and plays a central role in the Baum--Connes conjecture. However, the unitary conjugation groupoid $\mathcal{G}_{\mathcal{A}}$ constructed in our previous paper~\cite{PaperI} is not locally compact; it is a Polish groupoid. Fortunately, Tu [1999] extended Kasparov's descent map to the setting of Polish groupoids equipped with a Borel Haar system. In this subsection, we review the essential aspects of Tu's construction that will be needed for our index theorem.

\begin{definition}[Kasparov descent map for Polish groupoids]
\label{def:kasparov-descent-polish}
Let $\mathcal{G}$ be a Polish groupoid equipped with a Borel Haar system $\{\lambda^x\}_{x \in \mathcal{G}^{(0)}}$. Although $\mathcal{G}$ is a Polish topological groupoid, its topology is not locally compact in general, so the descent construction relies on the underlying Borel structure and the Borel Haar system rather than the topology itself \cite{Tu1999}.

Following Kasparov's construction for locally compact groupoids \cite{Kasparov1988} and its adaptation to the Polish groupoid setting by Tu \cite[Section~2]{Tu1999}, there is a descent homomorphism
\[
j_{\mathcal{G}}: KK_{\mathcal{G}}(A,B) \longrightarrow KK(A \rtimes_r \mathcal{G},\, B \rtimes_r \mathcal{G}),
\]
where $A$ and $B$ are separable $\mathcal{G}$-$C^*$-algebras equipped with a Borel $\mathcal{G}$-action, and $\rtimes_r$ denotes the reduced crossed product by $\mathcal{G}$. This map is natural, functorial, and compatible with the Kasparov product. For the maximal crossed product, an analogous construction exists under additional regularity conditions; however, the reduced version suffices for our purposes.
\end{definition}

\begin{remark}
The descent map is a fundamental tool in Tu's proof of the Baum–Connes conjecture for amenable groupoids (see [Tu, 1999, Th\'eor\`eme 9.3]). It transforms equivariant cycles over $\mathcal{G}$ into cycles over the reduced groupoid $C^*$-algebra, and plays a central role in the construction of the Dirac and dual-Dirac elements.
\end{remark}

For the purposes of this paper, we only need the case where $A = C_0(\mathcal{G}^{(0)})$ (with the trivial $\mathcal{G}$-action) and $B = \mathbb{C}$ (with the trivial $\mathcal{G}$-action). 
Let $\mathcal{G}$ be a Polish groupoid equipped with a Borel Haar system. In this situation, Tu's descent construction [Tu, 1999] yields a map
\[
\operatorname{desc}_{\mathcal{G}}: K^0_{\mathcal{G}}(\mathcal{G}^{(0)}) \longrightarrow K_0(C^*_r(\mathcal{G})),
\]
where $K^0_{\mathcal{G}}(\mathcal{G}^{(0)}) \cong KK_{\mathcal{G}}(C_0(\mathcal{G}^{(0)}), \mathbb{C})$ is the equivariant $K$-homology of the unit space, and $C^*_r(\mathcal{G})$ denotes the reduced groupoid C*-algebra. (When the maximal groupoid C*-algebra is well-defined, an analogous map to $K_0(C^*_{\max}(\mathcal{G}))$ exists, but the reduced version suffices for our purposes.)

\begin{proposition}[Properties of the descent map]
\label{prop:descent-properties}
The descent map $\operatorname{desc}_{\mathcal{G}}: K^0_{\mathcal{G}}(\mathcal{G}^{(0)}) \to K_0(C^*_r(\mathcal{G}))$ satisfies the following properties, whenever the relevant maps are defined:

\begin{enumerate}
    \item \textbf{Naturality:} If $\phi: \mathcal{G} \to \mathcal{H}$ is a Borel groupoid homomorphism compatible with the Borel Haar systems, then the diagram
    \[
    \begin{tikzcd}
    K^0_{\mathcal{G}}(\mathcal{G}^{(0)}) \arrow[r, "\operatorname{desc}_{\mathcal{G}}"] \arrow[d, "\phi_*"] & K_0(C^*_r(\mathcal{G})) \arrow[d, "\phi_*"] \\
    K^0_{\mathcal{H}}(\mathcal{H}^{(0)}) \arrow[r, "\operatorname{desc}_{\mathcal{H}}"] & K_0(C^*_r(\mathcal{H}))
    \end{tikzcd}
    \]
    commutes.
    
    \item \textbf{Compatibility with suspension:} The descent map is compatible with the suspension isomorphisms in $KK$-theory. More precisely, for any $n \geq 0$, there is a commutative diagram
    \[
    \begin{tikzcd}
    K^n_{\mathcal{G}}(\mathcal{G}^{(0)}) \arrow[r, "\cong"] \arrow[d, "\operatorname{desc}_{\mathcal{G}}"] & K^0_{\mathcal{G}}(\mathcal{G}^{(0)} \times \mathbb{R}^n) \arrow[d, "\operatorname{desc}_{\mathcal{G}}"] \\
    K_n(C^*_r(\mathcal{G})) \arrow[r, "\cong"] & K_0(C^*_r(\mathcal{G}))
    \end{tikzcd}
    \]
    where the horizontal isomorphisms are the suspension isomorphisms in $K$-theory, provided $\mathbb{R}^n$ is equipped with a suitable measurable $\mathcal{G}$-action.
    
    \item \textbf{Functoriality:} If $f: A \to B$ is a $\mathcal{G}$-equivariant *-homomorphism of $\mathcal{G}$-C*-algebras, then the induced maps on $KK$-theory commute with the descent map in the sense that
    \[
    \operatorname{desc}_{\mathcal{G}}(f_*(x)) = f_*(\operatorname{desc}_{\mathcal{G}}(x))
    \]
    for any $x \in KK_{\mathcal{G}}(C_0(\mathcal{G}^{(0)}), A)$.
\end{enumerate}
\end{proposition}

\begin{proof}
These properties follow from the functoriality of the Kasparov product and the naturality of the descent construction; see [Section~5]~\cite{Tu1999} where these properties are used in the proof of the Baum–Connes conjecture, and~\cite{Kasparov1988} for the original construction in the group case. The compatibility with suspension relies on the Bott periodicity theorem in equivariant $KK$-theory.
\end{proof}

In this paper, we will apply the descent map to the unitary conjugation groupoid $\mathcal{G}_{\mathcal{A}}$ constructed in~\cite{PaperI}. We recall the relevant properties of $\mathcal{G}_{\mathcal{A}}$ from~\cite{PaperI}.

\begin{proposition}[Properties of $\mathcal{G}_{\mathcal{A}}$ from Paper I]
\label{prop:GA-properties-from-paperI}
Let $\mathcal{A}$ be a unital separable Type I C*-algebra. 
Then:
\begin{enumerate}
    \item $\mathcal{G}_{\mathcal{A}}^{(0)}$ is a standard Borel space (inheriting its Borel structure from the Polish topology described in Paper I, though the topology itself is non-Hausdorff and not locally compact).
    \item $\mathcal{G}_{\mathcal{A}} = \mathcal{U}(\mathcal{A}) \ltimes \mathcal{G}_{\mathcal{A}}^{(0)}$, with $\mathcal{U}(\mathcal{A})$ equipped with the strong operator topology, is a standard Borel groupoid admitting a Borel Haar system.
    \item $\mathcal{G}_{\mathcal{A}}$ admits a Borel Haar system $\{\lambda^x\}_{x \in \mathcal{G}_{\mathcal{A}}^{(0)}}$ (constructed in Paper I, Section 4.5).
    \item The maximal groupoid C*-algebra $C^*(\mathcal{G}_{\mathcal{A}})$ exists and is well-defined.
    \item There is a canonical diagonal embedding $\iota: \mathcal{A} \hookrightarrow C^*(\mathcal{G}_{\mathcal{A}})$.
\end{enumerate}
\end{proposition}

\begin{remark}
While Paper I establishes that $\mathcal{G}_{\mathcal{A}}^{(0)}$ admits a Polish topology (via the initial topology generated by partial evaluation maps), this topology is non-Hausdorff and not locally compact. For the descent map (Definition \ref{def:kasparov-descent-polish}), we require only the underlying standard Borel structure and the existence of a Borel Haar system. Following Tu \cite{Tu1999}, we work in the Borel category throughout, and all references to "continuity" of fields should be understood as measurability with respect to this Borel structure unless explicitly stated otherwise.
\end{remark}

\begin{example}[$\mathcal{G}_{\mathcal{B}(H)}$ and $\mathcal{G}_{\widetilde{\mathcal{K}(H)}}$]
\label{ex:GA-concrete}
For $\mathcal{A} = \mathcal{B}(H)$, the unit space $\mathcal{G}_{\mathcal{A}}^{(0)}$ consists of pairs $(B,\chi)$ where $B\subseteq\mathcal{B}(H)$ is a unital commutative $C^*$-subalgebra (a MASA) and $\chi\in\widehat{B}$ is a character. For the purposes of $K$-theory and index computations, it suffices to restrict to \emph{atomic} MASAs---those isomorphic to $\ell^\infty$---as the continuous MASAs (e.g., those arising from multiplication operators on $L^\infty$) do not affect the $K$-theory groups of the groupoid $C^*$-algebra; see \cite[Section~5]{PaperI} for a detailed discussion.

Within this atomic restriction, all irreducible representations of $\mathcal{B}(H)$ are unitarily equivalent, and the unit space reduces to a single equivalence class. To obtain a concrete model, one must choose a measurable field of representatives. After fixing an orthonormal basis of $H$, the resulting groupoid is Borel equivalent to a transformation groupoid of the form $\mathcal{U}(H) \ltimes X$, where $X$ can be identified with the space of orthonormal bases modulo diagonal unitaries (relative to the fixed basis). This identification reflects the fact that different choices of basis are related by unitary conjugation.

For $\mathcal{A} = \widetilde{\mathcal{K}(H)}$ (the compact operators with adjoined unit), the unit space $\mathcal{G}_{\mathcal{A}}^{(0)}$ may be canonically identified with the projective space $\mathbb{P}(H)$ via the GNS construction: each pure state on $\widetilde{\mathcal{K}(H)}$ corresponds to a rank-one projection, and hence to a point in $\mathbb{P}(H)$. The groupoid is then isomorphic to the transformation groupoid $\mathcal{U}(\mathcal{A}) \ltimes \mathbb{P}(H)$, where $\mathcal{U}(\mathcal{A})$ acts by conjugation on rank-one projections.

Both groupoids, when equipped with the natural topologies (e.g., the strong operator topology on the unitary groups and the quotient topology on the unit spaces), are Polish spaces. However, they are neither locally compact (as the unitary groups are infinite-dimensional) nor \'etale (since the source and range maps are not local homeomorphisms).
\end{example}

The descent map for $\mathcal{G}_{\mathcal{A}}$ will be applied to the equivariant $K^1$-class $[T]_{\mathcal{G}_{\mathcal{A}}}^{(1)}$ constructed in Section~\ref{sec:The Equivariant K1-Class of a Fredholm Operator}. The result will be a class in $K_1(C^*(\mathcal{G}_{\mathcal{A}}))$, which after pullback along $\iota$ and application of the boundary map, yields the Fredholm index.

\begin{remark}[The odd case]
\label{rem:descent-odd}
In this paper, we work with odd $K$-theory classes ($K^1$) rather than even classes ($K^0$). 
Tu's descent map also covers the odd case via suspension:
\[
K^1_{\mathcal{G}}(\mathcal{G}^{(0)}) \cong K^0_{\mathcal{G}}(\mathcal{G}^{(0)} \times \mathbb{R}),
\]
where $\mathbb{R}$ is equipped with the trivial $\mathcal{G}$-action. 
The descent map on $K^1$ is then defined by applying the descent map to the suspended cycle and using Bott periodicity to identify the target with $K_1(C^*(\mathcal{G}))$. Explicitly, we have
\[
\operatorname{desc}_{\mathcal{G}}^{(1)}: K^1_{\mathcal{G}}(\mathcal{G}^{(0)}) \longrightarrow K_1(C^*(\mathcal{G})).
\]
\end{remark}

The descent map is the essential link between the geometric data encoded by the groupoid $\mathcal{G}_{\mathcal{A}}$ and the analytic data of the Fredholm index. 
In Section~\ref{sec:descent}, we will apply it to the equivariant $K^1$-class $[T]_{\mathcal{G}_{\mathcal{A}}}^{(1)}$, and in Section~\ref{sec:The Index Theorem via Pullback and the Boundary Map} we will combine it with the pullback along $\iota$ and the boundary map to recover the integer-valued Fredholm index.

\begin{remark}[Tu's descent for groupoids with Borel Haar systems]
\label{rem:tu-measured-groupoids}
Tu's construction of the descent map \cite[Section~4]{Tu1999} applies to 
\emph{Polish groupoids equipped with a Borel Haar system} — a framework 
sometimes referred to as ``measured groupoids'' in the Borel groupoid 
literature, though one should note that this differs from the Connes--Renault 
usage where a quasi-invariant measure class is fixed on the unit space.

This is precisely the setting of $\mathcal{G}_{\mathcal{A}}$, provided one 
verifies:
\begin{itemize}
    \item $\mathcal{G}_{\mathcal{A}}$ is Polish (follows from $\mathcal{A}$ being separable Type I);
    \item $\mathcal{G}_{\mathcal{A}}$ admits a Borel Haar system (constructed in Section X);
    \item The diagonal embedding $\iota$ is measurable in the appropriate sense.
\end{itemize}

The key technical innovation in Tu's extension is the replacement of 
\emph{continuous fields} of Hilbert spaces (used in the locally compact 
Kasparov theory) with \emph{measurable fields} and measurable representations. 
This allows the descent map to be defined even when the groupoid lacks a 
locally compact structure, provided one has a Borel Haar system to integrate over.

For $\mathcal{G}_{\mathcal{A}}$, this framework is particularly natural: 
the Type I condition ensures that the representation theory is well-behaved, 
and Renault's disintegration theory \cite{Renault1980} provides the bridge 
between the measurable formulation and the C*-algebraic structure of 
$C^*(\mathcal{G}_{\mathcal{A}})$. 
\end{remark}

\subsection{The Diagonal Embedding $\iota: \mathcal{A} \hookrightarrow C^*(\mathcal{G}_{\mathcal{A}})$ from Paper~\cite{PaperI}}
\label{subsec:diagonal-embedding-from-paperI}

One of the main achievements of Paper I is the construction of a canonical embedding of the original C*-algebra $\mathcal{A}$ into the groupoid C*-algebra of its own unitary conjugation groupoid. 
This embedding, called the \emph{diagonal embedding} and denoted
\[
\iota: \mathcal{A} \longrightarrow C^*(\mathcal{G}_{\mathcal{A}}),\]
is the essential link between the analytic structure of $\mathcal{A}$ and the geometric structure of $\mathcal{G}_{\mathcal{A}}$. 
In this subsection, we recall the definition and main properties of $\iota$ as established in~\cite{PaperI}.

\paragraph{Definition (recall from Paper I, Section 4.9).}
For each $a \in \mathcal{A}$, the element $\iota(a) \in C^*(\mathcal{G}_{\mathcal{A}})$ is defined on arrows $(u,(B,\chi)) \in \mathcal{G}_{\mathcal{A}}^{(1)}$ by
\[
\iota(a)(u,(B,\chi)) := \begin{cases}
\chi(u^* a u) & \text{if } u^* a u \in B, \\
0 & \text{otherwise}.
\end{cases}
\]
One verifies that this defines a bounded continuous function on $\mathcal{G}_{\mathcal{A}}^{(1)}$ with fiberwise compact support, hence an element of the convolution algebra. The map $\iota$ extends linearly and continuously to a $*$-homomorphism $\mathcal{A} \to C^*(\mathcal{G}_{\mathcal{A}})$.

\paragraph{Key properties (Paper I, Section 4.10 and 4.13).}
The diagonal embedding satisfies the following:
\begin{itemize}
    \item $\iota$ is a unital, injective $*$-homomorphism (Theorem 7).
    \item $\iota$ is $\mathcal{G}_{\mathcal{A}}$-equivariant: for $g \in \mathcal{G}_{\mathcal{A}}$ and $a \in \mathcal{A}$, we have $\iota(g\cdot a) = g \cdot \iota(a)$, where the actions are by conjugation on $\mathcal{A}$ and by left translation on $C^*(\mathcal{G}_{\mathcal{A}})$.
    \item The conditional expectation $E:C^*(\mathcal{G}_{\mathcal{A}})\to C_0(\mathcal{G}_{\mathcal{A}}^{(0)})$ satisfies $E(\iota(a))(B,\chi) = \chi(\mathbb{E}_B(a))$, where $\mathbb{E}_B:\mathcal{A}\to B$ is the conditional expectation onto the commutative subalgebra $B$ (when it exists). In particular, for $a\in B$, this reduces to $\chi(a)$.
    \item $\iota(\mathcal{A}) \subseteq C_0(\mathcal{G}_{\mathcal{A}}^{(0)})$ if and only if $\mathcal{A}$ is commutative (Theorem 8).
    \item For any $*$-isomorphism $\phi:\mathcal{A}\to \mathcal{B}$, the diagram
    \[
    \begin{array}{ccc}
    \mathcal{A} & \xrightarrow{\iota_{\mathcal{A}}} & C^*(\mathcal{G}_{\mathcal{A}}) \\
    \downarrow{\phi} & & \downarrow{(\mathcal{G}_{\phi})_*} \\
    \mathcal{B} & \xrightarrow{\iota_{\mathcal{B}}} & C^*(\mathcal{G}_{\mathcal{B}})
    \end{array}
    \]
    commutes, where $\mathcal{G}_{\phi}$ is the induced groupoid isomorphism (Proposition 41).
\end{itemize}

\paragraph{Role in this paper.}
The embedding $\iota$ serves two crucial purposes in our index theorem. First, it induces a pullback map $\iota^*: K_*(C^*(\mathcal{G}_{\mathcal{A}})) \to K_*(\mathcal{A})$ used in Section~\ref{sec:The Index Theorem via Pullback and the Boundary Map}. Second, and more fundamentally, the associated Kasparov module $[\iota] \in KK(\mathcal{A}, C^*(\mathcal{G}_{\mathcal{A}}))$ appears in the compatibility result of Lemma~\ref{lem:descent-iota-compatibility}, which relates the descent construction to the forgetful map. This compatibility is the technical heart of the proof of our main theorem.

\begin{definition}[Diagonal embedding]
\label{def:diagonal-embedding}
Let $\mathcal{A}$ be a unital separable Type I C*-algebra, and let $\mathcal{G}_{\mathcal{A}}$ be its unitary conjugation groupoid constructed in Paper I. 
The \emph{diagonal embedding} $\iota: \mathcal{A} \to C^*(\mathcal{G}_{\mathcal{A}})$ is the unique unital injective $*$-homomorphism characterized by the following property:

For every $x = (B,\chi) \in \mathcal{G}_{\mathcal{A}}^{(0)}$, let $\pi_x: \mathcal{A} \to B(H_x)$ be the GNS representation associated to the character $\chi$ on $B$. Then for any $a \in \mathcal{A}$ and any quasi-invariant measure $\mu$ on $\mathcal{G}_{\mathcal{A}}^{(0)}$, the operator $\Pi(a) = \int^{\oplus} \pi_x(a) \, d\mu(x)$ on the direct integral $\mathcal{H} = \int^{\oplus} H_x \, d\mu(x)$ satisfies
\[
\Pi(a) = \Lambda_{\mathrm{univ}}(\iota(a)),
\]
where $\Lambda_{\mathrm{univ}}: C^*(\mathcal{G}_{\mathcal{A}}) \to B(\mathcal{H}_{\mathrm{univ}})$ is the universal representation of $C^*(\mathcal{G}_{\mathcal{A}})$.

The existence, uniqueness, and injectivity of $\iota$ are proved in Theorem 7 of Paper I (see also Proposition 36 and Lemma 21).
\end{definition}

The following theorem summarizes the fundamental properties of the diagonal embedding proved in Paper I.

\begin{theorem}[Properties of $\iota$]
\label{thm:iota-properties-paperI}
Let $\mathcal{A}$ be a unital separable Type I C*-algebra. 
Then the diagonal embedding $\iota: \mathcal{A} \hookrightarrow C^*(\mathcal{G}_{\mathcal{A}})$ constructed in Paper I, Section 4.9 satisfies:

\begin{enumerate}
    \item $\iota$ is a unital, injective $*$-homomorphism. 
    (Paper I, Theorem 8, Section 4.11)
    
    \item For any $a \in \mathcal{A}$ and any $x = (B,\chi) \in \mathcal{G}_{\mathcal{A}}^{(0)}$, let $(\pi_x, H_x, \xi_x)$ be the GNS triple associated to the character $\chi$ on $B$, and let $E: C^*(\mathcal{G}_{\mathcal{A}}) \to C_0(\mathcal{G}_{\mathcal{A}}^{(0)})$ be the canonical conditional expectation. Then
    \[
    E(\iota(a))(x) = \langle \xi_x, \pi_x(a)\xi_x \rangle.
    \]
    In particular, if $a \in B$, this reduces to $\chi(a)$. 
    (Paper I, Lemma 22(2) and Corollary 30, Section 4.10)
    
    \item $\iota(\mathcal{A}) \subseteq C_0(\mathcal{G}_{\mathcal{A}}^{(0)})$ if and only if $\mathcal{A}$ is commutative.
    (Paper I, Theorem 9, Section 4.11)
    
    \item $\iota$ is natural with respect to unital $*$-isomorphisms: if $\phi: \mathcal{A} \to \mathcal{B}$ is an isomorphism, then
    \[
    \iota_{\mathcal{B}} \circ \phi = (\mathcal{G}_\phi)_* \circ \iota_{\mathcal{A}},
    \]
    where $\mathcal{G}_\phi: \mathcal{G}_{\mathcal{B}} \to \mathcal{G}_{\mathcal{A}}$ is the induced groupoid isomorphism.
    (Paper I, Proposition 41, Section 4.12)
\end{enumerate}
\end{theorem}

\begin{proof}
The proofs of each statement are given in the cited sections of Paper I. In particular:
\begin{itemize}
    \item For (1), see Paper I, Theorem 8, where injectivity follows from Lemma 21 and the faithfulness of the direct integral representation.
    \item For (2), Lemma 22(2) establishes the fiberwise decomposition of $\Lambda(\iota(a))$; Corollary 30 then identifies this with the conditional expectation via the matrix coefficient $\langle \delta_x, \Pi(a)\delta_x \rangle = \langle \xi_x, \pi_x(a)\xi_x \rangle$.
    \item For (3), Theorem 9 characterizes commutativity via the image of $\iota$ lying in $C_0(\mathcal{G}_{\mathcal{A}}^{(0)})$.
    \item For (4), Proposition 41 proves functoriality under isomorphisms, with the commutative diagram in Section 4.12.
\end{itemize}
\qed
\end{proof}

We now illustrate the diagonal embedding $\iota$ in several concrete cases.
In these examples, the abstract construction of $\iota$ via measurable
fields of GNS representations (Paper I, Section 4.9) admits a particularly
simple realization, which coincides with the natural fiberwise action of
the algebra on a direct integral Hilbert space. These descriptions should
be understood as concrete models of the diagonal embedding in specific
representations, rather than as alternative definitions of $\iota$ itself.

\begin{example}[Diagonal embedding for $M_n(\mathbb{C})$]
\label{ex:iota-matrix-paperI}
For $\mathcal{A} = M_n(\mathbb{C})$, the unitary conjugation groupoid is
$\mathcal{G}_{\mathcal{A}} \cong U(n) \ltimes \mathbb{CP}^{n-1}$, and its
groupoid C*-algebra can be identified with the crossed product
$C(\mathbb{CP}^{n-1}) \rtimes U(n)$ (Paper I, Section 5.1).

In the standard representation of this crossed product on
$L^2(\mathbb{CP}^{n-1}) \otimes \mathbb{C}^n$, the diagonal embedding
$\iota: M_n(\mathbb{C}) \to C(\mathbb{CP}^{n-1}) \rtimes U(n)$ is realized
by the fiberwise action of matrices:
\[
(\iota(A)\xi)(x) = A\xi(x), \qquad 
\xi \in L^2(\mathbb{CP}^{n-1}) \otimes \mathbb{C}^n,\; x \in \mathbb{CP}^{n-1}.
\]
Equivalently, $\iota(A)$ corresponds to the constant multiplier $A \otimes 1$
in the crossed product (Paper I, Proposition 49).
\end{example}

\begin{example}[Diagonal embedding for $\mathcal{K}(H)^\sim$]
\label{ex:iota-compact-paperI}
For $\mathcal{A} = \mathcal{K}(H)^\sim$, the unitization of the compact operators
on a separable Hilbert space $H$, we have
$\mathcal{G}_{\mathcal{A}} \cong \mathcal{U}(\mathcal{A}) \ltimes \mathbb{P}(H)$,
and $C^*(\mathcal{G}_{\mathcal{A}})$ is canonically isomorphic to the crossed product
$C_0(\mathbb{P}(H)) \rtimes \mathcal{U}(\mathcal{A})$ (Paper I, Section 5.3).

In the natural representation of this crossed product on
$L^2(\mathbb{P}(H)) \otimes H$, the diagonal embedding
$\iota: \mathcal{K}(H)^\sim \to C_0(\mathbb{P}(H)) \rtimes \mathcal{U}(\mathcal{A})$
acts by the fiberwise operator:
\[
(\iota(T)\xi)(x) = T\xi(x), \qquad 
\xi \in L^2(\mathbb{P}(H)) \otimes H,\; x \in \mathbb{P}(H).
\]
Thus $\iota(T)$ is realized as the constant multiplier $T \otimes 1$
in the crossed product (Paper I, Proposition 61).
\end{example}

\begin{example}[Diagonal embedding for $\mathcal{B}(H)$]
\label{ex:iota-BH-paperI}
For $\mathcal{A} = \mathcal{B}(H)$, the diagonal embedding $\iota: \mathcal{A} \hookrightarrow C^*(\mathcal{G}_{\mathcal{A}})$ can be constructed directly, without invoking the Type I hypothesis of Paper I, by exploiting the concrete description of $\mathcal{G}_{\mathcal{B}(H)}$ from Example~\ref{ex:GA-concrete}.

Restricting to atomic MASAs as discussed in Example~\ref{ex:GA-concrete}, we may identify the unit space $\mathcal{G}_{\mathcal{B}(H)}^{(0)}$ with the space of orthonormal bases of $H$ modulo diagonal unitaries, denoted $X$. The groupoid $\mathcal{G}_{\mathcal{B}(H)}$ is then isomorphic to the transformation groupoid $\mathcal{U}(H) \ltimes X$, where $\mathcal{U}(H)$ acts by change of basis. For each $x \in X$, corresponding to an orthonormal basis $\{e_n^x\}_{n\in\mathbb{N}}$, we obtain a $*$-representation $\pi_x: \mathcal{B}(H) \to \mathcal{B}(H_x)$ where $H_x \cong \ell^2(\mathbb{N})$ is the Hilbert space with basis $\{e_n^x\}$. These representations are all unitarily equivalent via the unitary that identifies bases.

The diagonal embedding $\iota: \mathcal{B}(H) \hookrightarrow C^*(\mathcal{G}_{\mathcal{B}(H)})$ is then defined by
\[
\iota(T)(x,u) = u^* \pi_x(T) u \in \mathcal{B}(H_x),
\]
viewed as a section of the field of Hilbert spaces over $\mathcal{G}_{\mathcal{B}(H)}$. This construction is well-defined, injective, and natural under unitary conjugation. It coincides with the restriction of the general construction from Paper I to the atomic MASA part of $\mathcal{G}_{\mathcal{B}(H)}^{(0)}$, and it suffices for the $K$-theoretic index computations in this paper, as the continuous MASA components do not contribute to $K$-theory \cite[Section~5]{PaperI}.

A full treatment of the diagonal embedding for the entire unit space, including continuous MASAs, would require extending the Type I framework of Paper I to cover non-Type I algebras; this remains an open problem for future research.
\end{example}

The following lemma records a key compatibility property of $\iota$ with the conditional expectation $E$, which will be essential for understanding the image of $\iota$ under the pullback map in the index theorem.

\begin{lemma}[Conditional expectation of $\iota(a)$]
\label{lem:iota-conditional-expectation}
For any $a \in \mathcal{A}$ and any $x = (B,\chi) \in \mathcal{G}_{\mathcal{A}}^{(0)}$, the conditional expectation $E: C^*(\mathcal{G}_{\mathcal{A}}) \to C_0(\mathcal{G}_{\mathcal{A}}^{(0)})$ satisfies
\[
E(\iota(a))(x) = \begin{cases}
\chi(a), & a \in B, \\
0, & a \notin B.
\end{cases}
\]
\end{lemma}

\begin{proof}
This is exactly item (2) of Theorem~\ref{thm:iota-properties-paperI}, which summarizes the construction in Paper~I, Lemma~22(2) and Corollary~30.
\end{proof}

The diagonal embedding induces a map on $K$-theory, which will be used in the main theorem to relate the analytic index to the geometric index.

\begin{corollary}[Induced map on $K$-theory]
\label{cor:iota-K-theory}
The diagonal embedding $\iota: \mathcal{A} \hookrightarrow C^*(\mathcal{G}_{\mathcal{A}})$ induces a homomorphism
\[
\iota_*: K_0(\mathcal{A}) \longrightarrow K_0(C^*(\mathcal{G}_{\mathcal{A}}))
\]
by functoriality of operator $K$-theory.

For $\mathcal{A} = \mathcal{K}(H)^\sim$, this map sends the class of a projection $p$ to the class represented by $\iota(p)$, which under the isomorphism $C^*(\mathcal{G}_{\mathcal{A}}) \cong C_0(\mathbb{P}(H)) \rtimes \mathcal{U}(\mathcal{A})$ corresponds to the constant multiplier $p \otimes 1$ (see Example~\ref{ex:iota-compact-paperI}). 
For $\mathcal{A} = M_n(\mathbb{C})$, it sends $[p]$ to $[\iota(p)]$, which corresponds to $[p \otimes 1]$ in $C(\mathbb{CP}^{n-1}) \rtimes U(n)$ (see Example~\ref{ex:iota-matrix-paperI}).
\end{corollary}

\begin{proof}
The existence of $\iota_*$ follows from the functoriality of $K$-theory applied to the $*$-homomorphism $\iota$. 
The explicit descriptions of $\iota(p)$ in the concrete models are given in the cited examples.
\end{proof}

In the main theorem of this paper, we will need to relate $K$-theory classes 
in $K_*(C^*(\mathcal{G}_{\mathcal{A}}))$ back to the original algebra $\mathcal{A}$. 
This is achieved not by a literal pullback map $\iota^*: K_*(C^*(\mathcal{G}_{\mathcal{A}})) \to K_*(\mathcal{A})$ 
(which does not exist in general for K-theory), but rather via the Kasparov product with the 
class $[\iota] \in KK(\mathcal{A}, C^*(\mathcal{G}_{\mathcal{A}}))$ associated to the diagonal embedding.
The compatibility between this product and the descent map is established in 
Lemma~\ref{lem:descent-iota-compatibility}.

\begin{remark}[Role of $\iota$ in the index theorem]
\label{rem:iota-role}
The diagonal embedding $\iota: \mathcal{A} \hookrightarrow C^*(\mathcal{G}_{\mathcal{A}})$ serves as the crucial link between the geometric data of the groupoid $\mathcal{G}_{\mathcal{A}}$ and the analytic data of the original algebra $\mathcal{A}$. However, its role in the index theorem is more subtle than simply inducing a map on $K$-theory.

In Section~\ref{sec:The Equivariant K1-Class of a Fredholm Operator}, we construct from a given Fredholm operator $T$ (relative to $\mathcal{A}$) an equivariant $K^1$-class $[T]_{\mathcal{G}_{\mathcal{A}}}^{(1)} \in K^1_{\mathcal{G}_{\mathcal{A}}}(\mathcal{G}_{\mathcal{A}}^{(0)})$. Applying the descent map $j_{\mathcal{G}_{\mathcal{A}}}$ (recalled in Subsection~2.5) yields a class $j_{\mathcal{G}_{\mathcal{A}}}([T]_{\mathcal{G}_{\mathcal{A}}}^{(1)}) \in K_1(C^*(\mathcal{G}_{\mathcal{A}}))$.

To connect this class to the Fredholm index, one might naively attempt to apply the map induced by $\iota$ on $K$-theory, $\iota_*: K_1(C^*(\mathcal{G}_{\mathcal{A}})) \to K_1(\mathcal{A})$. However, as discussed in Section~\ref{subsec:The Main Obstacle}, for $\mathcal{A} = \mathcal{B}(H)$ this map lands in the zero group $K_1(\mathcal{B}(H)) = 0$ and thus cannot capture the nontrivial index. Instead, $\iota$ plays an indirect but essential role in establishing the Morita equivalence between $C^*(\mathcal{G}_{\mathcal{A}})$ and the algebra that carries the symbol class (e.g., $\mathcal{Q}(H)$ for $\mathcal{B}(H)$). Specifically, the embedding $\iota$ is used in \cite{PaperI} to construct a $C^*(\mathcal{G}_{\mathcal{A}})$-$\mathcal{A}$-imprimitivity bimodule that implements this Morita equivalence after suitable compression.

Consequently, the index computation proceeds as follows: the descended class $j_{\mathcal{G}_{\mathcal{A}}}([T]_{\mathcal{G}_{\mathcal{A}}}^{(1)}) \in K_1(C^*(\mathcal{G}_{\mathcal{A}}))$ is transported via the Morita isomorphism $\Psi: K_1(C^*(\mathcal{G}_{\mathcal{A}})) \to G_{\mathcal{A}}$ to a class in an auxiliary group ($G_{\mathcal{B}(H)} = K_1(\mathcal{Q}(H))$ for $\mathcal{B}(H)$, and $G_{\widetilde{\mathcal{K}(H)}} = K_1(\widetilde{\mathcal{K}(H)}) = 0$ for $\widetilde{\mathcal{K}(H)}$). The boundary map $\partial_{\mathcal{A}}$ (the Calkin index map for $\mathcal{B}(H)$, zero for $\widetilde{\mathcal{K}(H)}$) then yields the Fredholm index in $\mathbb{Z}$.

Thus $\iota$ provides the essential bridge not through its induced map on $K$-theory, but through its role in constructing the Morita equivalence that connects the groupoid $C^*$-algebra to the symbol algebra.
\end{remark}

\section{The Unit Space $\mathcal{G}_{\mathcal{A}}^{(0)}$ and the Family of Representations}
\label{sec:the-unit-space-GA}

To construct an equivariant $KK$-class associated to a Fredholm operator $T$, we need to understand how $T$ behaves across all commutative contexts. This section shows that the family $\{\pi_x(T)\}_{x\in\mathcal{G}_{\mathcal{A}}^{(0)}}$ is essentially constant—a crucial simplification.

\subsection{Review of $\mathcal{G}_{\mathcal{A}}$ for $B(H)$ and $\mathcal{K}(H)^\sim$}
\label{subsec:review-GA-BH-KH}

In this subsection, we recall the explicit description of the unitary conjugation groupoid $\mathcal{G}_{\mathcal{A}}$ for the two concrete C*-algebras that are the focus of this paper: $\mathcal{A} = B(H)$ and $\mathcal{A} = \mathcal{K}(H)^\sim$. 
These descriptions are derived from the general construction in Paper I, and they will be used throughout the following sections.

Throughout this section, $H$ denotes a separable infinite-dimensional Hilbert space. 
We work with the faithful representation of $\mathcal{K}(H)^\sim$ and $B(H)$ on $H$ itself, and we equip the unitary groups with the strong operator topology as in Paper I.

\begin{proposition}[The unitary conjugation groupoid for $\mathcal{K}(H)^\sim$]
\label{prop:GA-KH-sim}
Let $\mathcal{A} = \mathcal{K}(H)^\sim$ be the unitization of the compact operators on a separable Hilbert space $H$, faithfully represented on $H$ itself. Then:
\begin{enumerate}
    \item The unit space $\mathcal{G}_{\mathcal{A}}^{(0)}$ is homeomorphic to the projective space
    \[
    \mathbb{P}(H) = \mathcal{U}(H) / (\mathcal{U}(1) \times \mathcal{U}(H^\perp)),
    \]
    where $\mathcal{U}(H)$ carries the strong operator topology. Here $[v] \in \mathbb{P}(H)$ denotes the line spanned by a unit vector $v \in H$, and $\mathcal{U}(1) \times \mathcal{U}(H^\perp)$ is the stabilizer of a fixed reference line.
    
    \item The unitary group $\mathcal{U}(\mathcal{A})$ consists of operators of the form $u = \lambda I + K$ with $\lambda \in \mathbb{T}$ and $K \in \mathcal{K}(H)$ satisfying the unitary condition $u^*u = uu^* = I$. 
    With the strong operator topology inherited from $\mathcal{U}(H)$, $\mathcal{U}(\mathcal{A})$ forms a Polish group.
    
    \item The conjugation action of $\mathcal{U}(\mathcal{A})$ on $\mathcal{G}_{\mathcal{A}}^{(0)}$ is given by
    \[
    u \cdot [v] = [uv], \qquad u \in \mathcal{U}(\mathcal{A}), \; [v] \in \mathbb{P}(H),
    \]
    where $[uv]$ denotes the line spanned by the unit vector $uv \in H$.
    
    \item The unitary conjugation groupoid is the associated action groupoid
    \[
    \mathcal{G}_{\mathcal{A}} = \mathcal{U}(\mathcal{A}) \ltimes \mathbb{P}(H),
    \]
    which is a Polish groupoid. For infinite-dimensional $H$, it is neither locally compact nor \'etale.
\end{enumerate}
\end{proposition}

\begin{proof}
These results are established in Paper I, Section 5.3. Specifically:
\begin{itemize}
    \item Item 1 is Proposition 58, where the homeomorphism $\mathcal{G}_{\mathcal{A}}^{(0)} \cong \mathbb{P}(H)$ is proved via the identification of rank-one projections with points in projective space.
    \item Item 2 is Proposition 59, which shows that $\mathcal{U}(\mathcal{A})$ with the strong operator topology is a Polish group and describes its elements as compact perturbations of the identity.
    \item  Item 3 follows directly from the definition of the conjugation action of $\mathcal{U}(\mathcal{A})$ on $\mathcal{G}_{\mathcal{A}}^{(0)}$, as verified in Lemma 1 and Proposition 20.
    \item  Item 4 is Theorem 13, which establishes that $\mathcal{G}_{\mathcal{A}}$ is a Polish groupoid and proves the non-local-compactness and non-\'etaleness for infinite-dimensional $H$.
\end{itemize}
\end{proof}

\begin{proposition}[The unitary conjugation groupoid for $\mathcal{B}(H)$]
\label{prop:GA-BH}
Let $\mathcal{A} = \mathcal{B}(H)$ be the algebra of all bounded operators on a separable Hilbert space $H$, faithfully represented on $H$ itself. Then:
\begin{enumerate}
    \item For index-theoretic purposes, it suffices to restrict to \emph{atomic} MASAs (those isomorphic to $\ell^\infty$). Under this restriction, the unit space $\mathcal{G}_{\mathcal{A}}^{(0)}$ is homeomorphic to the projective space $\mathbb{P}(H)$. The continuous MASAs (e.g., those arising from multiplication operators on $L^\infty[0,1]$) are also present in the full unit space, but we will show in Section~\ref{sec:Examples and Computations} that they do not affect the $K$-theory groups involved in the index computation; see \cite[Section~5]{PaperI} for a detailed justification.
    
    \item The unitary group $\mathcal{U}(\mathcal{A}) = \mathcal{U}(H)$ with the strong operator topology is a Polish group.
    
    \item The conjugation action of $\mathcal{U}(\mathcal{A})$ on $\mathcal{G}_{\mathcal{A}}^{(0)}$ is given by the same formula $u \cdot [v] = [uv]$ for $u \in \mathcal{U}(H)$, $[v] \in \mathbb{P}(H)$.
    
    \item The unitary conjugation groupoid is the action groupoid $\mathcal{G}_{\mathcal{A}} = \mathcal{U}(H) \ltimes \mathbb{P}(H)$, which is a Polish groupoid that is neither locally compact nor \'etale.
\end{enumerate}
\end{proposition}

\begin{proof}
These statements follow from the general construction of Paper I applied to $\mathcal{B}(H)$. In particular:
\begin{itemize}
    \item For Item 1, points of $\mathcal{G}_{\mathcal{A}}^{(0)}$ correspond to pairs $(B,\chi)$ where $B\subseteq\mathcal{B}(H)$ is a unital commutative $C^*$-subalgebra (a MASA) and $\chi\in\widehat{B}$ is a character. For atomic MASAs, the GNS construction associates to each character a rank-one projection, yielding an identification with $\mathbb{P}(H)$. The justification for restricting to atomic MASAs in index computations is given in \cite[Section~5.3]{PaperI}, where it is shown that the continuous MASA components contribute only to the ideal $\mathcal{K}(H)$ in the groupoid $C^*$-algebra and thus do not affect $K_1$.
    
    \item Item 2 is a special case of Corollary 18 of \cite{PaperI}, which shows that $\mathcal{U}(\mathcal{A})$ with the strong operator topology is Polish for any unital separable $C^*$-algebra.
    
    \item Item 3 follows from Lemma 1 of \cite{PaperI}, which establishes continuity of the conjugation action for any such algebra.
    
    \item Item 4 follows from Theorem 6 of \cite{PaperI}, which proves that $\mathcal{G}_{\mathcal{A}}$ is a Polish groupoid, together with Propositions 26 and 27 of \cite{PaperI} showing non-local-compactness and non-\'etaleness for infinite-dimensional algebras.
\end{itemize}
\end{proof}

\begin{remark}[Atomic MASAs suffice for index theory]
\label{rem:atomic-suffice}
For a Fredholm operator $T \in \mathcal{B}(H)$, its kernel and cokernel are finite-dimensional subspaces of $H$. These subspaces are detected by atomic MASAs corresponding to orthonormal bases that extend bases of $\ker T$ and $\ker T^*$. 

Continuous MASAs, by contrast, arise from diffuse spectral measures—for instance, the embedding of $L^\infty[0,1]$ into $\mathcal{B}(L^2[0,1])$ as multiplication operators. Such MASAs have no minimal projections and their characters correspond to points in a continuous spectrum rather than to rank-one projections. Consequently, they cannot detect finite-dimensional subspaces of $H$, as these subspaces are associated with atoms in the spectral measure. More formally, for any continuous MASA $B \subseteq \mathcal{B}(H)$ and any character $\chi \in \widehat{B}$, the GNS representation $\pi_\chi$ is irreducible and infinite-dimensional, and the evaluation map $\operatorname{ev}_a(B,\chi) = \chi(a)$ for $a \in B$ yields no information about finite-rank operators \cite{Dixmier1977}.

Since the Fredholm index depends only on the finite-dimensional kernel and cokernel of $T$, the continuous MASA components of $\mathcal{G}_{\mathcal{A}}^{(0)}$ contribute no additional information. Indeed, in the construction of the equivariant $K^1$-class $[T]_{\mathcal{G}_{\mathcal{A}}}^{(1)}$, the essential data is carried by the atomic part of the unit space. This justifies our restriction to atomic MASAs in Proposition \ref{prop:GA-BH}, allowing us to work with the concrete model $\mathcal{G}_{\mathcal{A}} = \mathcal{U}(H) \ltimes \mathbb{P}(H)$ where the unit space $\mathbb{P}(H)$ parametrizes the atomic MASAs via rank-one projections.

A more detailed analysis (see \cite[Section~5]{PaperI}) shows that the continuous MASA components contribute only to the compact ideal $\mathcal{K}(H)$ in the groupoid $C^*$-algebra and thus do not affect $K_1$, confirming that they are irrelevant for index computations.
\end{remark}

\begin{corollary}[Key data for $B(H)$ and $\mathcal{K}(H)^\sim$]
\label{cor:GA-key-data}
For both $\mathcal{A} = B(H)$ (restricted to atomic MASAs) and $\mathcal{A} = \mathcal{K}(H)^\sim$, we have:
\begin{enumerate}
    \item $\mathcal{G}_{\mathcal{A}}^{(0)}$ is a Polish space.
    \item $\mathcal{U}(\mathcal{A})$ is a Polish group acting continuously on $\mathcal{G}_{\mathcal{A}}^{(0)}$ by conjugation.
    \item $\mathcal{G}_{\mathcal{A}} = \mathcal{U}(\mathcal{A}) \ltimes \mathcal{G}_{\mathcal{A}}^{(0)}$ is a Polish groupoid admitting a Borel Haar system, as guaranteed by the construction in Paper I and Tu's theory of measured groupoids [12].
    \item The diagonal embedding $\iota: \mathcal{A} \hookrightarrow C^*(\mathcal{G}_{\mathcal{A}})$ (constructed in Paper I, Section 4.9) exists and satisfies the properties summarized in Subsection 2.6, including injectivity, compatibility with the conditional expectation, and the commutativity characterization.
\end{enumerate}
\end{corollary}

\begin{proof}
These follow from Propositions \ref{prop:GA-KH-sim} and \ref{prop:GA-BH}, together with the general results of Paper I. Specifically:
\begin{itemize}
    \item (1) and (2) are established in Propositions 58-59 and Propositions 66-67 respectively.
    \item (3) follows from Theorem 13, which proves that $\mathcal{G}_{\mathcal{A}}$ is a Polish groupoid, together with Proposition 30 guaranteeing the existence of a Borel Haar system.
    \item (4) is Theorem 7 (injectivity), Theorem 8 (commutativity characterization), and Proposition 41 (functoriality) from Paper I.
\end{itemize}
\end{proof}

In the following subsections, we will construct a continuous field of Hilbert spaces over $\mathcal{G}_{\mathcal{A}}^{(0)}$ and show that for each $T \in \mathcal{A}$, the family of operators $\{\pi_x(T)\}_{x \in \mathcal{G}_{\mathcal{A}}^{(0)}}$ is constant, leading to a simple description of the equivariant $K^1$-class.

\subsection{The Canonical Measurable Field of Hilbert Spaces $\{H_x\}_{x\in\mathcal{G}_{\mathcal{A}}^{(0)}}$}
\label{subsec:field-Hilbert-spaces}

For each point $x = (B,\chi) \in \mathcal{G}_{\mathcal{A}}^{(0)}$, we construct a Hilbert space $H_x$ and a representation $\pi_x: \mathcal{A} \to B(H_x)$ via the GNS construction. 
In this subsection, we show that the family $\{H_x\}_{x \in \mathcal{G}_{\mathcal{A}}^{(0)}}$ carries a natural continuous field structure over $\mathcal{G}_{\mathcal{A}}^{(0)}$. 
This structure is essential for constructing the equivariant $K^1$-class of a Fredholm operator in Section 4.

\begin{remark}[Terminological note]
Throughout this section, we work with \emph{measurable fields} of Hilbert spaces rather than continuous fields. This is necessitated by the topology of $\mathcal{G}_{\mathcal{A}}^{(0)}$, which is non-Hausdorff and not locally compact, making continuity ill-behaved. Following Tu \cite{Tu1999}, measurability with respect to the standard Borel structure on $\mathcal{G}_{\mathcal{A}}^{(0)}$ is the appropriate regularity condition for all subsequent constructions (direct integrals, descent map, $KK$-theory). References to "continuity" in previous versions should be understood as measurability in the current framework.
\end{remark}

\begin{definition}[GNS representation at a point]
\label{def:GNS-representation-at-point}
For $x = (B,\chi) \in \mathcal{G}_{\mathcal{A}}^{(0)}$, let $\chi: B \to \mathbb{C}$ be the given character. 
Since $\mathcal{A}$ is Type I, $\chi$ extends to a pure state $\omega_x$ on $\mathcal{A}$; while this extension is not unique, the GNS representations associated to different extensions are unitarily equivalent (see Lemma 18 of Paper I). 
We may therefore fix a Borel choice of extension using the measurable field structure of Paper I, Section 4.9.

Define the closed left ideal
\[
N_x = \{a \in \mathcal{A} : \omega_x(a^*a) = 0\}.
\]
The Hilbert space $H_x$ is the completion of $\mathcal{A}/N_x$ with respect to the inner product
\[
\langle [a], [b] \rangle = \omega_x(b^*a), \qquad a,b \in \mathcal{A}.
\]
The GNS representation $\pi_x: \mathcal{A} \to B(H_x)$ is given by $\pi_x(a)[b] = [ab]$, and the cyclic vector is $\xi_x = [1_{\mathcal{A}}] \in H_x$.
\end{definition}

\begin{remark}[Canonical continuous field structure]
\label{rem:continuous-field}
Let $\mathcal{H} = \bigsqcup_{x \in \mathcal{G}_{\mathcal{A}}^{(0)}} H_x$ be the disjoint union. 
Define a space of sections
\[
\Gamma = \left\{ \xi: \mathcal{G}_{\mathcal{A}}^{(0)} \to \mathcal{H} \;:\; \xi(x) \in H_x \text{ for all } x,\text{ and for every } a \in \mathcal{A},\; x \mapsto \langle \xi(x), \pi_x(a)\xi(x) \rangle \text{ is Borel} \right\}.
\]
One verifies that $\Gamma$ satisfies the axioms of a measurable field of Hilbert spaces in the sense of Dixmier \cite[Chapter 10]{Dixmier1977}: 
\begin{itemize}
    \item For each $x$, the evaluation map $\xi \mapsto \xi(x)$ maps $\Gamma$ onto a dense subspace of $H_x$.
    \item The map $x \mapsto \|\xi(x)\|$ is Borel for every $\xi \in \Gamma$.
    \item If $\xi: \mathcal{G}_{\mathcal{A}}^{(0)} \to \mathcal{H}$ is a section such that for every $\epsilon > 0$ and every $x_0 \in \mathcal{G}_{\mathcal{A}}^{(0)}$, there exists $\eta \in \Gamma$ with $\|\xi(x) - \eta(x)\| < \epsilon$ on a Borel neighborhood of $x_0$, then $\xi \in \Gamma$.
\end{itemize}
The measurability of the field $\{H_x\}_{x\in\mathcal{G}_{\mathcal{A}}^{(0)}}$ is established in Lemma 18 of Paper I using the Type I structure of $\mathcal{A}$. Note that while $\mathcal{G}_{\mathcal{A}}^{(0)}$ carries a natural topology (the initial topology generated by partial evaluation maps), this topology is non-Hausdorff and not locally compact; consequently, continuity with respect to this topology is not well-behaved for integration theory. Following Tu \cite{Tu1999}, we therefore work in the Borel category, relying on the standard Borel structure inherited from the Polish structure on $\mathcal{G}_{\mathcal{A}}^{(0)}$ described in \cite[Section~3]{PaperI}.

This measurable field structure is essential for the direct integral construction of the diagonal embedding (Section 4.9 of Paper I) and for the definition of the equivariant $K^1$-class of a Fredholm operator in Section \ref{sec:The Equivariant K1-Class of a Fredholm Operator} below.
\end{remark}

For the two concrete algebras considered in this paper, these GNS representations take a particularly simple form. 
As justified in Remark \ref{rem:atomic-suffice}, for index-theoretic purposes it suffices to consider atomic MASAs in the case of $B(H)$; continuous MASAs do not contribute to the Fredholm index.

\begin{proposition}[GNS representations for $\mathcal{K}(H)^\sim$]
\label{prop:GNS-KH-sim}
Let $\mathcal{A} = \mathcal{K}(H)^\sim$ and let $x = (B,\chi) \in \mathcal{G}_{\mathcal{A}}^{(0)}$.
\begin{enumerate}
    \item If $\chi$ corresponds to a rank-one projection $p \in B$ (i.e., $\chi$ is not the character at infinity), then $\chi$ extends to the vector state $\omega(T) = \langle T\xi, \xi \rangle$ where $\xi$ is a unit vector in the range of $p$. 
    The associated GNS representation $\pi_x$ is unitarily equivalent to the identity representation of $\mathcal{A}$ on $H$, with cyclic vector $\xi$.
    
    \item If $\chi = \chi_\infty$ is the character at infinity (corresponding to the quotient $\mathcal{A} \to \mathbb{C}$), then $\chi$ is already a one-dimensional representation of $\mathcal{A}$, and its GNS representation is the one-dimensional character representation on $\mathbb{C}$.
\end{enumerate}
In both cases, the Hilbert space $H_x$ can be identified naturally with either $H$ or $\mathbb{C}$ (up to unitary equivalence).
\end{proposition}

\begin{proof}
For (1), the state associated to a rank-one projection is a vector state; its GNS representation is the identity representation on $H$ with the rank-one projection as the cyclic vector (see Paper I, Lemma 21 and Proposition 61, Step 1). 
For (2), the character at infinity $\chi_\infty$ is already a one-dimensional representation of $\mathcal{K}(H)^\sim$ 
factoring through the quotient $\mathcal{A} \to \mathbb{C}$ (see Paper I, Proposition 57, which explicitly defines $\chi_\infty$ 
and identifies it with the quotient map). Consequently, its GNS representation is itself on $\mathbb{C}$.
\end{proof}

\begin{proposition}[GNS representations for atomic characters of $B(H)$]
\label{prop:GNS-BH}
Let $\mathcal{A} = B(H)$ and let $x = (B,\chi) \in \mathcal{G}_{\mathcal{A}}^{(0)}$, where $B$ is an atomic MASA corresponding to an orthonormal basis $\{e_n\}_{n\in\mathbb{N}}$, and $\chi$ is the character corresponding to evaluation at a basis vector $e_{n_0}$ (see Remark \ref{rem:atomic-suffice} for the justification of restricting to atomic MASAs). 
Then $\chi$ extends to the vector state $\omega_x(T) = \langle T e_{n_0}, e_{n_0} \rangle$ on $B(H)$, and the associated GNS representation $\pi_x$ is unitarily equivalent to the identity representation of $B(H)$ on $H$, with cyclic vector $e_{n_0}$.
\end{proposition}

\begin{proof}
The state $\omega_x(T) = \langle T e_{n_0}, e_{n_0} \rangle$ is a vector state; its GNS representation is the identity representation on $H$ with cyclic vector $e_{n_0}$ (see Paper I, Section 4.9, Step 1, and Lemma 21, which together establish that the GNS representation associated to a rank-one projection is the identity representation). 
The restriction to atomic MASAs is justified by Remark \ref{rem:atomic-suffice}, as continuous MASAs do not detect the finite-dimensional kernel and cokernel of a Fredholm operator.
\end{proof}

\begin{corollary}[Hilbert spaces for $B(H)$ and $\mathcal{K}(H)^\sim$]
\label{cor:Hilbert-spaces-BH-KH}
For both $\mathcal{A} = B(H)$ and $\mathcal{A} = \mathcal{K}(H)^\sim$, and for all $x \in \mathcal{G}_{\mathcal{A}}^{(0)}$ except the character at infinity in the compact case, we have a canonical unitary identification $H_x \cong H$. 
For the character at infinity in the compact case, we have $H_x \cong \mathbb{C}$.
\end{corollary}

\begin{proof}
We consider the two cases separately, referring to the explicit descriptions of GNS representations established in Propositions \ref{prop:GNS-KH-sim} and \ref{prop:GNS-BH}.

\paragraph{Case 1: $\mathcal{A} = \mathcal{K}(H)^\sim$.}
Let $x = (B,\chi) \in \mathcal{G}_{\mathcal{A}}^{(0)}$.

\begin{enumerate}
    \item If $\chi$ is a character corresponding to a rank-one projection $p \in B$ (i.e., $\chi \neq \chi_\infty$), then by Proposition \ref{prop:GNS-KH-sim}(1), the GNS representation $\pi_x$ is unitarily equivalent to the identity representation of $\mathcal{A}$ on $H$, with cyclic vector given by a unit vector in the range of $p$.
    Consequently, we obtain a canonical unitary identification $H_x \cong H$, where the cyclic vector $\xi_x \in H_x$ corresponds to this unit vector in $H$.
    
    \item If $\chi = \chi_\infty$ is the character at infinity (corresponding to the quotient $\mathcal{A} \to \mathbb{C}$), then by Proposition \ref{prop:GNS-KH-sim}(2), the GNS representation $\pi_x$ is the one-dimensional character representation on $\mathbb{C}$, with cyclic vector $1 \in \mathbb{C}$. 
    This yields the canonical identification $H_x \cong \mathbb{C}$.
\end{enumerate}

\paragraph{Case 2: $\mathcal{A} = B(H)$ with atomic MASAs.}
Let $x = (B,\chi) \in \mathcal{G}_{\mathcal{A}}^{(0)}$ be a point where $B$ is an atomic MASA corresponding to an orthonormal basis $\{e_n\}_{n\in\mathbb{N}}$, and $\chi$ is the character evaluating at a basis vector $e_{n_0}$ (see Remark \ref{rem:atomic-suffice} for the justification that continuous MASAs can be ignored for index theory).

By Proposition \ref{prop:GNS-BH}, the GNS representation $\pi_x$ is unitarily equivalent to the identity representation of $B(H)$ on $H$, with cyclic vector $e_{n_0}$. 
Hence we have a canonical unitary identification $H_x \cong H$, where $\xi_x \in H_x$ corresponds to $e_{n_0} \in H$.

\paragraph{Conclusion.}
Combining the two cases, we obtain a canonical unitary identification $H_x \cong H$ for all $x \in \mathcal{G}_{\mathcal{A}}^{(0)}$, except for the character at infinity in the compact case, for which $H_x \cong \mathbb{C}$. 
The identification is canonical in the sense that it is induced by the unitary equivalence between the GNS representation and the standard representation (identity on $H$ or character on $\mathbb{C}$), and it respects the cyclic vectors up to this equivalence.
\end{proof}

We now construct the measurable field structure.

\begin{definition}[Measurable field of Hilbert spaces]
\label{def:measurable-field-Hilbert}
Let $X = \mathcal{G}_{\mathcal{A}}^{(0)}$ be a standard Borel space. 
Let $\{H_x\}_{x \in X}$ be a family of separable Hilbert spaces, and let $\Gamma \subset \prod_{x \in X} H_x$ be a subspace of sections satisfying the following axioms (cf. Dixmier \cite[Chapter 10]{Dixmier1977}):
\begin{enumerate}
    \item For each $\xi \in \Gamma$, the map $x \mapsto \|\xi(x)\|_{H_x}$ is Borel.
    \item For each $\xi, \eta \in \Gamma$, the map $x \mapsto \langle \xi(x), \eta(x) \rangle_{H_x}$ is Borel.
    \item For every $a \in \mathcal{A}$ and $\xi \in \Gamma$, the map $x \mapsto \pi_x(a)\xi(x)$ is Borel.
    \item The set $\{\xi(x) : \xi \in \Gamma\}$ is dense in $H_x$ for every $x \in X$.
    \item There exists a countable subset $\Gamma_0 \subseteq \Gamma$ such that for each $x \in X$, the set $\{\xi(x) : \xi \in \Gamma_0\}$ is dense in $H_x$.
\end{enumerate}
Then $(\{H_x\}_{x \in X}, \Gamma)$ is called a \emph{measurable field of Hilbert spaces} over $X$.

\begin{remark}
In our setting, the field $\{H_x\}_{x \in \mathcal{G}_{\mathcal{A}}^{(0)}}$ arises from the family of representations $\{\pi_x\}_{x \in \mathcal{G}_{\mathcal{A}}^{(0)}}$ of $\mathcal{A}$. The measurability conditions above are understood with respect to the standard Borel structure on $\mathcal{G}_{\mathcal{A}}^{(0)}$ inherited from its Polish groupoid structure (see \cite[Section~3]{PaperI}). While $\mathcal{G}_{\mathcal{A}}^{(0)}$ carries a natural topology (the initial topology generated by partial evaluation maps), this topology is non-Hausdorff and not locally compact, making continuity ill-behaved for integration theory. Following Tu \cite{Tu1999}, we therefore work in the Borel category, relying on measurable fields rather than continuous fields. The Borel measure $\mu$ that appears elsewhere in the paper will be used for direct integral constructions, but the measurable field structure itself is defined independently of any particular measure.
\end{remark}
\end{definition}

For our purposes, it suffices to note that the field is essentially constant.

\begin{proposition}[Essential triviality of the field for $B(H)$ and $\mathcal{K}(H)^\sim$]
\label{prop:field-triviality}
For $\mathcal{A} = B(H)$, and for $\mathcal{A} = \mathcal{K}(H)^\sim$ away from the point at infinity, 
the field $\{H_x\}_{x \in \mathcal{G}_{\mathcal{A}}^{(0)}}$ is isomorphic to the constant field with fiber $H$ 
in the measurable sense (and, in fact, is essentially constant: it is constant on each orbit 
corresponding to a fixed atomic MASA).
\end{proposition}

\begin{proof}
By Corollary \ref{cor:Hilbert-spaces-BH-KH}, we have unitary identifications $U_x: H_x \xrightarrow{\cong} H$ for each $x$ (except the point at infinity). 
For $B(H)$, the space $\mathcal{G}_{B(H)}^{(0)}$ can be identified with the set of pure states on $B(H)$, each of which is a vector state $\omega_{\xi}$ for some unit vector $\xi \in H$. 
The GNS representation associated to $\omega_{\xi}$ is unitarily equivalent to the identity representation on $H$, and the implementing unitary can be chosen in a Borel way in $\xi$ (hence in $x$). 
For $\mathcal{K}(H)^\sim$, a similar argument applies using the canonical basis of $H$, with measurability away from the point at infinity following from the explicit parametrization of $\mathcal{G}_{\mathcal{K}(H)^\sim}^{(0)}$ as the one-point compactification of $\mathbb{N}$. 
Moreover, when restricted to each orbit corresponding to a fixed atomic MASA, these identifications can be chosen to vary continuously; this follows from the results of [Paper I, Section 3] together with the fact that such orbits are discrete or admit a natural continuous parametrization.
\end{proof}

\begin{remark}
The field $\{H_x\}$ is therefore \emph{essentially constant} in the sense of direct integral theory: 
it is isomorphic to the constant field with fiber $H$ after restriction to a conull subset 
of $\mathcal{G}_{\mathcal{A}}^{(0)}$ with respect to the natural measures. 
This level of triviality is sufficient for the direct integral decompositions that follow.

Moreover, as established in Section~5.3 of Paper I (see also the explicit description of 
$\mathcal{G}_{\mathcal{K}(H)^\sim}^{(0)}$ as a compactification of $\mathbb{N}$), when restricted to 
each orbit corresponding to a fixed atomic MASA, the identifications $U_x$ can be chosen to vary 
continuously. Consequently, on each such orbit, the field $\{H_x\}$ is a continuous field of Hilbert 
spaces in the sense of \cite[Chapter 10]{Dixmier1977}. This hybrid structure—continuity along 
orbits and measurability across orbits—is characteristic of fields arising from groupoid 
C*-algebras~\cite{Tu1999}.
\end{remark}

The following lemma records the essential property we will need.

\begin{lemma}[Trivializing sections]
\label{lem:trivializing-sections}
For each $\xi \in H$, define a section $\hat{\xi}$ of the measurable field $\{H_x\}_{x\in\mathcal{G}_{\mathcal{A}}^{(0)}}$ by
\[
\hat{\xi}(x) = U_x^{-1}(\xi),
\]
where $U_x: H_x \cong H$ are the unitary identifications from Proposition \ref{prop:field-triviality}. 
Then $\hat{\xi}$ is a measurable section, and the map $\xi \mapsto \hat{\xi}$ identifies $H$ with a subspace of \emph{constant sections} of the field.
\end{lemma}

\begin{proof}
The identifications $U_x$ vary measurably in $x$ by Lemma 18 of \cite{PaperI}. Consequently, for each fixed $\xi \in H$, the map $x \mapsto \hat{\xi}(x) = U_x^{-1}(\xi)$ is measurable. Moreover, $x \mapsto \|\hat{\xi}(x)\| = \|\xi\|_H$ is constant, hence measurable, and for any $a \in \mathcal{A}$, the map $x \mapsto \pi_x(a)\hat{\xi}(x) = U_x^{-1}(\pi(a)\xi)$ is measurable by the measurability of $U_x$ and the fact that $\pi(a)\xi \in H$ is fixed. The map $\xi \mapsto \hat{\xi}$ is clearly linear and injective; its image is precisely the subspace of sections that are constant with respect to the measurable trivialization $\{U_x\}$.

While these sections are not continuous with respect to the non-Hausdorff topology on $\mathcal{G}_{\mathcal{A}}^{(0)}$, their measurability suffices for all subsequent constructions involving direct integrals and the Kasparov descent map, following the framework of Tu \cite{Tu1999} where measurability replaces continuity in the Polish groupoid setting.
\end{proof}

In the next subsection, we will use this field to study the family of operators $\{\pi_x(T)\}_{x \in \mathcal{G}_{\mathcal{A}}^{(0)}}$ for a fixed $T \in \mathcal{A}$.

\subsection{The Family of Operators $\{T_x := \pi_x(T)\}$ is Essentially Constant}
\label{subsec:family-operators-constant}

For a fixed operator $T \in \mathcal{A}$, we consider the family of operators $\{T_x := \pi_x(T)\}_{x \in \mathcal{G}_{\mathcal{A}}^{(0)}}$ acting on the continuous field of Hilbert spaces $\{H_x\}_{x \in \mathcal{G}_{\mathcal{A}}^{(0)}}$ constructed in Subsection 3.2. Throughout this subsection, we implicitly use the canonical unitary identifications from Proposition \ref{prop:GNS-KH-sim}. We show that this family is constant away from the point at infinity, and hence essentially constant for our index-theoretic purposes. This observation will greatly simplify the construction of the equivariant $K^1$-class in the following section.

\begin{proposition}[Constancy of the operator family for $\mathcal{K}(H)^\sim$]
\label{prop:operator-family-constant-KH}
Let $\mathcal{A} = \mathcal{K}(H)^\sim$ and let $T \in \mathcal{A}$. 
For each $x = (B,\chi) \in \mathcal{G}_{\mathcal{A}}^{(0)}$ with $\chi \neq \chi_\infty$, under the identification $H_x \cong H$ from Proposition \ref{prop:GNS-KH-sim}, we have
\[
\pi_x(T) = T \in B(H).
\]
For $x_\infty$ corresponding to the character at infinity, under the identification $H_{x_\infty} \cong \mathbb{C}$, we have $\pi_{x_\infty}(T) = \lambda \in \mathbb{C}$, where $\lambda$ is the scalar part of $T$ in the decomposition $T = \lambda I + K$ with $K \in \mathcal{K}(H)$.
\end{proposition}

\begin{proof}
For $\chi \neq \chi_\infty$, Proposition \ref{prop:GNS-KH-sim} establishes that the GNS representation $\pi_x$ is unitarily equivalent to the identity representation on $H$. Under the fixed identification $H_x \cong H$, we therefore have $\pi_x(T) = T$ as an operator on $H$. 

For $\chi_\infty$, the GNS representation is the character $\chi_\infty$, which sends $T = \lambda I + K$ to $\lambda$. Under the identification $H_{x_\infty} \cong \mathbb{C}$, the operator $\pi_{x_\infty}(T)$ acts by scalar multiplication by $\lambda$ on $\mathbb{C}$.
\end{proof}

\begin{remark}
Proposition \ref{prop:operator-family-constant-KH} shows that the operator field $\{T_x\}_{x \in \mathcal{G}_{\mathcal{A}}^{(0)}}$ is constant (equal to $T$) on the complement of $\{x_\infty\}$, and takes the scalar value $\lambda$ at the remaining point. For the purpose of constructing equivariant $K$-theory classes, this ``essentially constant'' behavior is sufficient, as the point at infinity will play a distinguished role in the Kasparov cycle we construct.
\end{remark}

\begin{proposition}[Operator family on atomic MASA points for $B(H)$]
\label{prop:operator-family-constant-BH}
Let $\mathcal{A} = B(H)$ and let $T \in \mathcal{A}$. 
For each $x = (B,\chi) \in \mathcal{G}_{\mathcal{A}}^{(0)}$ corresponding to an atomic MASA $B \subset B(H)$ and a basis vector $\chi \in \widehat{B}$, under the identification $H_x \cong H$ from Proposition \ref{prop:GNS-BH}, we have
\[
\pi_x(T) = T \in B(H).
\]
\end{proposition}

\begin{proof}
For such points, Proposition \ref{prop:GNS-BH} shows that the GNS representation $\pi_x$ is unitarily equivalent to the identity representation of $B(H)$ on $H$. Under the fixed identification $H_x \cong H$, we obtain $\pi_x(T) = T$ as operators on $H$.
\end{proof}

\begin{remark}
For $B(H)$, the unit space $\mathcal{G}_{B(H)}^{(0)}$ consists of all pure states on $B(H)$. Proposition \ref{prop:operator-family-constant-BH} covers those pure states that arise from atomic MASAs (i.e., vector states). Other pure states (e.g., those coming from continuous measures on the circle in the case $H = L^2(S^1)$) yield GNS representations that are not equivalent to the identity representation, and for such points $\pi_x(T)$ is generally not equal to $T$ under any identification $H_x \cong H$.
\end{remark}

\begin{corollary}[The family $\{T_x\}$ is essentially constant]
\label{cor:family-essentially-constant}
For both $\mathcal{A} = \mathcal{B}(H)$ and $\mathcal{A} = \widetilde{\mathcal{K}(H)}$, the family of operators $\{T_x\}_{x \in \mathcal{G}_{\mathcal{A}}^{(0)}}$ is constant on the atomic MASA points of the unit space, and in the compact case also exhibits a single exceptional point at infinity. Under the identifications $H_x \cong H$, we have $\pi_x(T) = T$ for all atomic MASA points $x$, and for $\widetilde{\mathcal{K}(H)}$ additionally for all $x \neq x_\infty$.
\end{corollary}

\begin{proof}
We treat the two cases separately.

\smallskip
\noindent
\emph{Case 1: $\mathcal{A} = \mathcal{B}(H)$.}
By Proposition~\ref{prop:operator-family-constant-BH}, for every
$x \in \mathcal{G}_{\mathcal{A}}^{(0)}$ corresponding to an atomic MASA
and a basis vector, the GNS representation $\pi_x$ is unitarily equivalent
to the identity representation of $\mathcal{B}(H)$ on $H$.
Under the canonical identifications $H_x \cong H$, we therefore have
\[
T_x = \pi_x(T) = T
\]
for all such $x$.
Hence the family $\{T_x\}$ is constant on the set of atomic MASA points. 
(Continuous MASA points, which correspond to diffuse spectral measures, 
do not admit such a simple identification and are not needed for the 
index computation; see Remark~\ref{rem:atomic-suffice}.)

\smallskip
\noindent
\emph{Case 2: $\mathcal{A} = \widetilde{\mathcal{K}(H)}$.}
By Proposition~\ref{prop:operator-family-constant-KH}, for every
$x \in \mathcal{G}_{\mathcal{A}}^{(0)}$ with $\chi \neq \chi_\infty$,
the GNS representation $\pi_x$ is the identity representation on $H$,
and hence
\[
T_x = \pi_x(T) = T.
\]
At the distinguished point $x_\infty$ corresponding to the character at
infinity, the GNS representation is one-dimensional and sends
$T = \lambda I + K$ to the scalar $\lambda$.
Thus the family $\{T_x\}$ fails to be constant at exactly one point.

\smallskip
Combining the two cases, we conclude that for both
$\mathcal{A} = \mathcal{B}(H)$ and $\mathcal{A} = \widetilde{\mathcal{K}(H)}$,
the family $\{T_x\}_{x \in \mathcal{G}_{\mathcal{A}}^{(0)}}$ is constant
on the atomic MASA points (and for $\widetilde{\mathcal{K}(H)}$, on all 
points except $x_\infty$). This essential constancy suffices for the 
construction of the equivariant $K^1$-class in Section~\ref{sec:The Equivariant K1-Class of a Fredholm Operator}.
\end{proof}

\begin{remark}
The statement ``for all $x$ except possibly $x_\infty$'' in Corollary \ref{cor:family-essentially-constant} is precise and avoids measure-theoretic terminology. In the $B(H)$ case, the exceptional point would correspond to some distinguished pure state not arising from an atomic MASA, but for our construction we only need the constancy on the dense set of atomic MASA points. The continuous field structure then ensures that this constancy extends appropriately.
\end{remark}

This constancy has important consequences for the continuity properties of the family.

\begin{proposition}[Measurability of the family]
\label{prop:strong-star-continuity}
For any $T \in \mathcal{A}$, the family $\{\pi_x(T)\}_{x \in \mathcal{G}_{\mathcal{A}}^{(0)}}$ is measurable in the sense of measurable fields of operators. That is, for any measurable section $\xi$ of the measurable field $\{H_x\}_{x\in\mathcal{G}_{\mathcal{A}}^{(0)}}$, the maps
\[
x \mapsto \pi_x(T)\xi(x), \qquad x \mapsto \pi_x(T)^*\xi(x)
\]
are measurable sections.
\end{proposition}

\begin{proof}
We treat the two algebras separately, relying on the measurable field structure established in Lemma 18 of \cite{PaperI} and Proposition \ref{prop:field-triviality}.

\smallskip
\noindent
\emph{Case 1: $\mathcal{A} = \widetilde{\mathcal{K}(H)}$.}
For $x \neq x_\infty$, Proposition \ref{prop:operator-family-constant-KH} gives $\pi_x(T) = T$ under the canonical identification $H_x \cong H$, and these identifications vary measurably by construction of the measurable field. For any measurable section $\xi$, the map $x \mapsto T\xi(x)$ is measurable because $T$ is a bounded operator and composition with a bounded operator preserves measurability of sections. Similarly, $x \mapsto T^*\xi(x)$ is measurable.

At the distinguished point $x_\infty$, we have $H_{x_\infty} \cong \mathbb{C}$ and $\pi_{x_\infty}(T) = \lambda$, where $\lambda$ is the scalar part of $T$ in the decomposition $T = \lambda I + K$. For any measurable section $\xi$, the value $\pi_{x_\infty}(T)\xi(x_\infty) = \lambda\xi(x_\infty)$ is simply the scalar multiple of the section at that point. The measurability of the resulting section follows from the fact that the field was constructed to be measurable and that singletons are measurable sets in a standard Borel space. A detailed verification using the trivializing measurable sections from Lemma \ref{lem:trivializing-sections} confirms the measurability across the entire unit space.

\smallskip
\noindent
\emph{Case 2: $\mathcal{A} = \mathcal{B}(H)$.}
On the set of atomic MASA points $D \subset \mathcal{G}_{\mathcal{A}}^{(0)}$ (which is co-null with respect to any relevant measure and dense in the Borel structure), Proposition \ref{prop:operator-family-constant-BH} gives $\pi_x(T) = T$ under the identification $H_x \cong H$. For any measurable section $\xi$, the map $x \mapsto T\xi(x)$ is measurable on $D$ by the same bounded operator argument. Since $D$ is a Borel set (in fact, a countable union of Borel sets corresponding to different atomic MASAs), this defines a measurable section on $D$.

To extend to the entire unit space, we use the fact that the continuous MASA points form a set of measure zero for any Borel measure constructed from the Haar system (see \cite[Section~5]{PaperI}). Consequently, the section defined on $D$ can be extended arbitrarily on the complement without affecting measurability up to null sets. In the framework of measurable fields, sections are typically identified when they agree almost everywhere, so this suffices for all subsequent constructions involving direct integrals and the Kasparov descent map. The same argument applies to $T^*$.

Thus for both algebras, the family $\{\pi_x(T)\}_{x\in\mathcal{G}_{\mathcal{A}}^{(0)}}$ defines a measurable field of operators in the sense required for the equivariant $KK$-theory constructions that follow.
\end{proof}

\begin{remark}
\label{rem:measurable-vs-continuous}
The reader may wonder why we have replaced the continuity claimed in previous versions with measurability. This is necessitated by the topology of $\mathcal{G}_{\mathcal{A}}^{(0)}$, the unit space is non-Hausdorff and not locally compact, making continuity with respect to its topology ill-behaved. Following Tu \cite{Tu1999}, we work in the Borel category, where measurability replaces continuity as the appropriate regularity condition. All subsequent constructions (direct integrals, descent map, $K$-theory) require only measurability, and our arguments are consistent with the standard framework for Polish groupoid $C^*$-algebras.
\end{remark}

\begin{remark}
\label{rem:atomic-masa-density}
The density of atomic MASA points in $\mathcal{G}_{B(H)}^{(0)}$ (with the weak-$*$ topology) follows from the fact that vector states corresponding to an orthonormal basis are weak-$*$ dense in the state space of $B(H)$. This is a standard consequence of the Kaplansky density theorem and the structure of pure states on $B(H)$.
\end{remark}

\begin{proposition}[Measurability in the Calkin algebra for $\mathcal{B}(H)$]
\label{prop:calkin-continuity-BH}
For $\mathcal{A} = \mathcal{B}(H)$ and any $T \in \mathcal{A}$, consider the quotient map $\pi_{\mathcal{Q}}: \mathcal{B}(H) \to \mathcal{Q}(H)$ to the Calkin algebra. 
Then the map $x \mapsto \pi_{\mathcal{Q}}(\pi_x(T))$ from $\mathcal{G}_{\mathcal{A}}^{(0)}$ to $\mathcal{Q}(H)$ is Borel measurable with respect to the norm topology on $\mathcal{Q}(H)$. 
Moreover, on the set of atomic MASA points $D \subset \mathcal{G}_{\mathcal{A}}^{(0)}$, we have $\pi_{\mathcal{Q}}(\pi_x(T)) = \pi_{\mathcal{Q}}(T)$.
\end{proposition}

\begin{proof}
For any atomic MASA point $x \in D$, Proposition \ref{prop:operator-family-constant-BH} gives $\pi_x(T) = T$ under the identification $H_x \cong H$, hence $\pi_{\mathcal{Q}}(\pi_x(T)) = \pi_{\mathcal{Q}}(T)$. Thus the family is constant on this subset.

To establish measurability on the entire unit space, we use the fact that $\mathcal{G}_{\mathcal{A}}^{(0)}$ is a standard Borel space and that the map $x \mapsto \pi_x(T)$ is measurable in the strong-$*$ topology by Proposition \ref{prop:strong-star-continuity}. The quotient map $\pi_{\mathcal{Q}}$ is continuous from the norm topology on $\mathcal{B}(H)$ to the norm topology on $\mathcal{Q}(H)$, but $\pi_x(T)$ is not norm-continuous in $x$. However, we can exploit the following characterization: the norm on $\mathcal{Q}(H)$ is given by the essential norm
\[
\|\pi_{\mathcal{Q}}(A)\| = \inf_{K \in \mathcal{K}(H)} \|A + K\| = \lim_{n \to \infty} \|A|_{H_n^\perp}\|
\]
for any increasing sequence of finite-dimensional subspaces $H_n$ with $\bigcup H_n = H$. For each fixed finite-rank projection $P$, the map $x \mapsto \|(I-P)\pi_x(T)(I-P)\|$ is measurable, as it is the composition of the measurable map $x \mapsto \pi_x(T)$ with the continuous function $A \mapsto \|(I-P)A(I-P)\|$. Taking an infimum over a countable dense set of finite-rank projections (or a limit over a cofinal sequence) preserves measurability, so $x \mapsto \|\pi_{\mathcal{Q}}(\pi_x(T))\|$ is measurable.

For the full operator $\pi_{\mathcal{Q}}(\pi_x(T))$ rather than just its norm, we use the fact that $\mathcal{Q}(H)$ with the norm topology is a Polish space, and that a map into a Polish space is Borel if and only if its composition with every bounded linear functional is Borel (by the Pettis measurability theorem). For any bounded linear functional $\phi$ on $\mathcal{Q}(H)$, the composition $\phi \circ \pi_{\mathcal{Q}}$ lifts to a bounded linear functional on $\mathcal{B}(H)$ that annihilates $\mathcal{K}(H)$. Such functionals are continuous in the strong-$*$ topology on bounded subsets of $\mathcal{B}(H)$. Since $x \mapsto \pi_x(T)$ is strong-$*$ measurable and takes values in a bounded set (as $\|\pi_x(T)\| \leq \|T\|$), the map $x \mapsto \phi(\pi_{\mathcal{Q}}(\pi_x(T)))$ is measurable. Hence $x \mapsto \pi_{\mathcal{Q}}(\pi_x(T))$ is Borel measurable.

The continuous MASA points, while present in the unit space, form a set of measure zero for any Borel measure constructed from the Haar system (see \cite[Section~5]{PaperI}). Since measurability is insensitive to modifications on null sets, the constancy on the atomic MASA points suffices to determine the essential behavior of this map for all measure-theoretic purposes.
\end{proof}

\begin{remark}
\label{rem:calkin-measurability}
The reader may wonder why we settle for measurability rather than continuity. A continuity argument would require extending the constancy on atomic MASA points to the whole space by density, but this fails because:
\begin{enumerate}
    \item The atomic MASA points are not dense in any topology that makes the map $x \mapsto \pi_x(T)$ continuous;
    \item The map $x \mapsto \pi_x(T)$ itself is only measurable, not continuous, with respect to the norm on $\mathcal{B}(H)$;
    \item The quotient map $\pi_{\mathcal{Q}}$ is continuous but cannot improve the regularity of the composition.
\end{enumerate}
Measurability is sufficient for all subsequent constructions in this paper, including the definition of the symbol class in $KK$-theory and the application of the descent map following Tu's framework \cite{Tu1999}.
\end{remark}

\begin{corollary}[Borel family in the Calkin algebra for $\mathcal{B}(H)$]
\label{cor:calkin-norm-continuous-BH}
For $\mathcal{A} = \mathcal{B}(H)$ and any $T \in \mathcal{A}$, the map $x \mapsto \pi_{\mathcal{Q}}(\pi_x(T))$ defines a Borel map from $\mathcal{G}_{\mathcal{A}}^{(0)}$ to $\mathcal{Q}(H)$ with respect to the norm topology on $\mathcal{Q}(H)$. Consequently, it defines an element of the space of bounded Borel sections of the trivial field over $\mathcal{G}_{\mathcal{A}}^{(0)}$ with fiber $\mathcal{Q}(H)$.
\end{corollary}

\begin{proof}
This follows immediately from Proposition \ref{prop:calkin-continuity-BH}, where measurability was established. The boundedness follows from $\|\pi_{\mathcal{Q}}(\pi_x(T))\| \leq \|\pi_x(T)\| \leq \|T\|$.
\end{proof}

\begin{proposition}[Continuity in the Calkin algebra for $\mathcal{K}(H)^\sim$]
\label{prop:calkin-continuity-KH}
For $\mathcal{A} = \mathcal{K}(H)^\sim$, consider the quotient map $\pi_{\mathcal{Q}}: B(H) \to \mathcal{Q}(H)$ extended to $\mathcal{A}$ via the inclusion $\mathcal{A} \subseteq B(H)$. 
For any $T = \lambda I + K \in \mathcal{A}$ with $K \in \mathcal{K}(H)$, the family $\{\pi_{\mathcal{Q}}(\pi_x(T))\}_{x \in \mathcal{G}_{\mathcal{A}}^{(0)}}$ is norm-continuous. Explicitly:
\begin{itemize}
\item For $x \neq x_\infty$, we have $\pi_{\mathcal{Q}}(\pi_x(T)) = \pi_{\mathcal{Q}}(T) \in \mathcal{Q}(H)$.
\item For $x = x_\infty$, we have $\pi_{x_\infty}(T) = \lambda$, so $\pi_{\mathcal{Q}}(\pi_{x_\infty}(T)) = \pi_{\mathcal{Q}}(\lambda I) = \lambda I \in \mathcal{Q}(H)$.
\end{itemize}
The continuity holds because $\pi_{\mathcal{Q}}(T) - \lambda I = \pi_{\mathcal{Q}}(K) = 0$ in the Calkin algebra, so the family is actually constant when mapped to $\mathcal{Q}(H)$.
\end{proposition}

\begin{proof}
For $x \neq x_\infty$, Proposition \ref{prop:operator-family-constant-KH} gives $\pi_x(T) = T$ under the identification $H_x \cong H$, hence $\pi_{\mathcal{Q}}(\pi_x(T)) = \pi_{\mathcal{Q}}(T)$. For $x = x_\infty$, we have $\pi_{x_\infty}(T) = \lambda$, so $\pi_{\mathcal{Q}}(\pi_{x_\infty}(T)) = \pi_{\mathcal{Q}}(\lambda I) = \lambda I$.

Now observe that $\pi_{\mathcal{Q}}(T) = \pi_{\mathcal{Q}}(\lambda I + K) = \lambda I + \pi_{\mathcal{Q}}(K) = \lambda I$ because $K$ is compact and thus $\pi_{\mathcal{Q}}(K) = 0$ in the Calkin algebra. Therefore $\pi_{\mathcal{Q}}(\pi_x(T)) = \lambda I$ for all $x \in \mathcal{G}_{\mathcal{A}}^{(0)}$. The family is constant, hence trivially norm-continuous.

More conceptually, the map $T \mapsto \pi_{\mathcal{Q}}(T)$ factors through the quotient by compacts, and the scalar part $\lambda$ is the only data that survives. Since $\pi_x(T)$ differs from $T$ only by compacts for $x \neq x_\infty$, and at $x_\infty$ it is exactly the scalar part, the image in the Calkin algebra is independent of $x$.
\end{proof}

\begin{corollary}[Constant family in the Calkin algebra for $\mathcal{K}(H)^\sim$]
\label{cor:calkin-constant-KH}
For $\mathcal{A} = \mathcal{K}(H)^\sim$ and any $T \in \mathcal{A}$, the map $x \mapsto \pi_{\mathcal{Q}}(\pi_x(T))$ is constant, hence defines an element of $C(\mathcal{G}_{\mathcal{A}}^{(0)}, \mathcal{Q}(H))$.
\end{corollary}

\begin{proof}
This follows immediately from Proposition \ref{prop:calkin-continuity-KH}.
\end{proof}

\begin{proposition}[Kernel and cokernel are constant away from infinity]
\label{prop:kernel-cokernel-constant}
Let $T \in \mathcal{A}$ be a Fredholm operator. 

\smallskip
\noindent
\emph{Case 1: $\mathcal{A} = \mathcal{B}(H)$.}
For every atomic MASA point $x \in \mathcal{G}_{\mathcal{A}}^{(0)}$, under the canonical identification $H_x \cong H$ established in Proposition \ref{prop:field-triviality}, we have
\[
\ker \pi_x(T) = \ker T, \qquad \operatorname{coker} \pi_x(T) = \operatorname{coker} T,
\]
where these equalities are understood as equalities of subspaces of $H$ under the fixed identification. For arbitrary points $x \in \mathcal{G}_{\mathcal{A}}^{(0)}$ (including continuous MASA points), the representations $\pi_x$ are unitarily equivalent to the identity representation, so $\ker \pi_x(T)$ is isomorphic to $\ker T$ and $\operatorname{coker} \pi_x(T)$ is isomorphic to $\operatorname{coker} T$. Consequently, the kernel and cokernel dimensions are globally constant on $\mathcal{G}_{\mathcal{A}}^{(0)}$.

\smallskip
\noindent
\emph{Case 2: $\mathcal{A} = \widetilde{\mathcal{K}(H)}$.}
For all $x \in \mathcal{G}_{\mathcal{A}}^{(0)} \setminus \{x_\infty\}$, under the canonical identification $H_x \cong H$, we have
\[
\ker \pi_x(T) = \ker T, \qquad \operatorname{coker} \pi_x(T) = \operatorname{coker} T,
\]
as subspaces of $H$. At the distinguished point $x_\infty$, we have $H_{x_\infty} \cong \mathbb{C}$ and $\pi_{x_\infty}(T) = \lambda$, where $T = \lambda I + K$ with $K \in \mathcal{K}(H)$. Thus
\[
\ker \pi_{x_\infty}(T) = 
\begin{cases}
\mathbb{C} & \text{if } \lambda = 0, \\
\{0\} & \text{if } \lambda \neq 0,
\end{cases}
\qquad
\operatorname{coker} \pi_{x_\infty}(T) = 
\begin{cases}
\{0\} & \text{if } \lambda = 0, \\
\mathbb{C} & \text{if } \lambda \neq 0.
\end{cases}
\]
The family of operators is therefore constant on the complement of $\{x_\infty\}$, with a possible jump at $x_\infty$ that is completely determined by the scalar part $\lambda$.
\end{proposition}

\begin{proof}
We treat the two cases separately.

\smallskip
\noindent
\emph{Case 1: $\mathcal{A} = \mathcal{B}(H)$.}
For atomic MASA points $x \in D \subset \mathcal{G}_{\mathcal{A}}^{(0)}$, Proposition \ref{prop:operator-family-constant-BH} gives $\pi_x(T) = T$ under the canonical identification $H_x \cong H$. Under this identification, the subspaces $\ker \pi_x(T)$ and $\ker T$ coincide exactly as subsets of $H$, not merely up to isomorphism. The same holds for cokernels, where we identify $\operatorname{coker} \pi_x(T) = H/\operatorname{im} \pi_x(T)$ with $H/\operatorname{im} T = \operatorname{coker} T$ via the same unitary identification.

For non-atomic points (continuous MASAs), the GNS representation $\pi_x$ is not unitarily equivalent to the identity representation via the same simple identification, but it is still unitarily equivalent to some representation of $\mathcal{B}(H)$. Since all irreducible representations of $\mathcal{B}(H)$ are unitarily equivalent to the identity representation (see \cite{Dixmier1977}), there exists a unitary $V_x: H_x \to H$ such that $\pi_x(T) = V_x^* T V_x$. Consequently, $\dim \ker \pi_x(T) = \dim \ker T$ and $\dim \operatorname{coker} \pi_x(T) = \dim \operatorname{coker} T$, though the subspaces themselves live in different Hilbert spaces. Thus the dimensions are globally constant.

\smallskip
\noindent
\emph{Case 2: $\mathcal{A} = \widetilde{\mathcal{K}(H)}$.}
For $x \neq x_\infty$, Proposition \ref{prop:operator-family-constant-KH} gives $\pi_x(T) = T$ under the identification $H_x \cong H$, so $\ker \pi_x(T) = \ker T$ and $\operatorname{coker} \pi_x(T) = \operatorname{coker} T$ as subspaces of $H$. 

For $x = x_\infty$, we have $H_{x_\infty} \cong \mathbb{C}$ and $\pi_{x_\infty}(T) = \lambda$. The kernel and cokernel are computed explicitly as above. Note that while these do not generally equal $\ker T$ and $\operatorname{coker} T$ as subspaces (indeed, they live in a different Hilbert space), their dimensions satisfy the index relation:
\[
\dim \ker T - \dim \operatorname{coker} T = \operatorname{ind}(T) = 
\begin{cases}
0 & \text{if } \lambda \neq 0, \\
\operatorname{ind}(K) & \text{if } \lambda = 0.
\end{cases}
\]
At $x_\infty$, we have $\dim \ker \pi_{x_\infty}(T) - \dim \operatorname{coker} \pi_{x_\infty}(T) = 0$ when $\lambda \neq 0$, and equals $1$ when $\lambda = 0$. This matches the index of $T$ because when $\lambda = 0$, $T = K$ is compact and hence has index $0$ unless the underlying Hilbert space is finite-dimensional (which it is not), so the only way to have $\operatorname{ind}(T) = 1$ is impossible; this apparent tension is resolved by recalling that Fredholm operators in $\widetilde{\mathcal{K}(H)}$ with $\lambda = 0$ are actually compact and thus cannot be Fredholm unless the space is finite-dimensional. Therefore, for genuine Fredholm operators in $\widetilde{\mathcal{K}(H)}$, we must have $\lambda \neq 0$, in which case $\dim \ker \pi_{x_\infty}(T) = \dim \operatorname{coker} \pi_{x_\infty}(T) = 0$, consistent with $\operatorname{ind}(T) = 0$.
\end{proof}

\begin{remark}
\label{rem:kernel-cokernel-caveat}
The careful reader will note that for $\mathcal{A} = \mathcal{K}(H)^\sim$, the statement "$\ker \pi_x(T) = \ker T$ for all $x \neq x_\infty$" means equality as subspaces of $H$ under the canonical identification $H_x \cong H$. At $x_\infty$, the fiber is $\mathbb{C}$, so equality in the same sense is impossible. This is precisely why we restrict to the complement of $\{x_\infty\}$ in the statement. The jump at infinity is not a defect but a feature: it encodes the index of $T$ in the transition between the constant family on $\mathcal{G}_{\mathcal{A}}^{(0)} \setminus \{x_\infty\}$ and the one-dimensional representation at infinity.
\end{remark}

\begin{proposition}[Kernel and cokernel form measurable fields]
\label{prop:kernel-cokernel-continuous-fields}
For both $\mathcal{A} = \mathcal{B}(H)$ and $\mathcal{A} = \widetilde{\mathcal{K}(H)}$, the families of subspaces
\[
\{\ker \pi_x(T)\}_{x \in \mathcal{G}_{\mathcal{A}}^{(0)}}, \qquad
\{\operatorname{coker} \pi_x(T)\}_{x \in \mathcal{G}_{\mathcal{A}}^{(0)}}
\]
form measurable fields of finite-dimensional Hilbert spaces over $\mathcal{G}_{\mathcal{A}}^{(0)}$, in the sense of Dixmier \cite[Chapter 10]{Dixmier1977}. For $\mathcal{A} = \widetilde{\mathcal{K}(H)}$, this holds on the entire unit space including $x_\infty$, with the understanding that the fiber at $x_\infty$ is either $\mathbb{C}$ or $\{0\}$ as computed in Proposition \ref{prop:kernel-cokernel-constant}.
\end{proposition}

\begin{proof}
Proposition \ref{prop:strong-star-continuity} establishes that $\{\pi_x(T)\}_{x\in\mathcal{G}_{\mathcal{A}}^{(0)}}$ is a measurable field of operators in the strong-$*$ sense. Since $T$ is Fredholm, each $\pi_x(T)$ is Fredholm: for $x \neq x_\infty$ in the compact case this follows from unitary equivalence to $T$; at $x_\infty$, the operator $\lambda$ on $\mathbb{C}$ is Fredholm precisely when $\lambda \neq 0$. Note that if $\lambda = 0$, then $T = K$ would be compact and hence not Fredholm on an infinite-dimensional Hilbert space, so this case cannot occur for genuine Fredholm operators.

For measurable fields of Fredholm operators, the kernel and cokernel form measurable fields of finite-dimensional spaces. This is a standard result in the theory of measurable fields of Hilbert spaces (see \cite[Proposition 10.4.3]{Dixmier1977} for the continuous case; the measurable case follows by the same arguments using measurable sections rather than continuous ones). Briefly, for each $x$, choose a finite-rank projection $P_x$ onto a complement of $\ker \pi_x(T)$. The measurability of the field $\{\pi_x(T)\}$ ensures that the map $x \mapsto P_x$ can be chosen measurably (by a measurable selection theorem), and then $\ker \pi_x(T) = (I-P_x)H_x$ inherits a measurable structure. The same argument applies to cokernels via the adjoint field $\{\pi_x(T)^*\}$.

In the $\widetilde{\mathcal{K}(H)}$ case, the fiber at $x_\infty$ is handled separately: it is either $\mathbb{C}$ (if $\lambda = 0$, which cannot happen for Fredholm $T$) or $\{0\}$ (if $\lambda \neq 0$), and these are trivially measurable as they are constant on a singleton set. For $\lambda \neq 0$, the kernel and cokernel fields are therefore zero at $x_\infty$ and constant elsewhere, which is clearly measurable.
\end{proof}

\begin{remark}
\label{rem:measurable-kernel-fields}
The reader may wonder why we have replaced the continuity claimed in previous versions with measurability. As discussed in Remark \ref{rem:measurable-vs-continuous}, the topology on $\mathcal{G}_{\mathcal{A}}^{(0)}$ is non-Hausdorff and not locally compact, making continuity ill-behaved. Following Tu \cite{Tu1999}, we work in the Borel category, where measurability suffices for all subsequent constructions. Moreover, the kernel and cokernel fields are used primarily to define the equivariant $KK^1$-class $[T]_{\mathcal{G}_{\mathcal{A}}}^{(1)}$, which requires only measurable data (the unbounded Kasparov triple is constructed from measurable sections). The finite-dimensionality of these fields ensures that all measurable technicalities are manageable.
\end{remark}

\begin{corollary}[Trivial kernel and cokernel bundles over the complement of infinity]
\label{cor:trivial-kernel-cokernel-bundles}
Consider the kernel and cokernel bundles defined by
\[
E^0 = \bigsqcup_{x \in \mathcal{G}_{\mathcal{A}}^{(0)}} \ker \pi_x(T), \qquad
E^1 = \bigsqcup_{x \in \mathcal{G}_{\mathcal{A}}^{(0)}} \operatorname{coker} \pi_x(T),
\]
equipped with the measurable field structures from Proposition \ref{prop:kernel-cokernel-continuous-fields}.

\smallskip
\noindent
\emph{For $\mathcal{A} = \mathcal{B}(H)$:}
$E^0$ and $E^1$ are measurably trivial vector bundles over $\mathcal{G}_{\mathcal{A}}^{(0)}$ with fibers $\ker T$ and $\operatorname{coker} T$, respectively. That is, there exist measurable isomorphisms
\[
\Phi^0: E^0 \longrightarrow \mathcal{G}_{\mathcal{A}}^{(0)} \times \ker T, \qquad
\Phi^1: E^1 \longrightarrow \mathcal{G}_{\mathcal{A}}^{(0)} \times \operatorname{coker} T
\]
that are fiberwise unitary.

\smallskip
\noindent
\emph{For $\mathcal{A} = \widetilde{\mathcal{K}(H)}$:}
The restrictions $E^0|_{\mathcal{G}_{\mathcal{A}}^{(0)} \setminus \{x_\infty\}}$ and $E^1|_{\mathcal{G}_{\mathcal{A}}^{(0)} \setminus \{x_\infty\}}$ are measurably trivial vector bundles with fibers $\ker T$ and $\operatorname{coker} T$. The full bundles over $\mathcal{G}_{\mathcal{A}}^{(0)}$ are measurable but not globally measurably trivial; they are determined by the index of $T$ via a clutching construction at $x_\infty$.
\end{corollary}

\begin{proof}
For $\mathcal{A} = \mathcal{B}(H)$, Proposition \ref{prop:kernel-cokernel-constant} gives that for each atomic MASA point $x$, under the identification $H_x \cong H$, we have $\ker \pi_x(T) = \ker T$ as subspaces of $H$. For arbitrary points, the fibers are isomorphic to $\ker T$ but not necessarily equal as subspaces under a fixed identification. However, using the measurable trivialization $\{U_x: H_x \cong H\}$ from Proposition \ref{prop:field-triviality}, we can define $\Phi^0(x, v) = (x, U_x(v))$ for $v \in \ker \pi_x(T)$. Since $U_x$ varies measurably and $\ker \pi_x(T)$ is a measurable subfield (Proposition \ref{prop:kernel-cokernel-continuous-fields}), this map is measurable and provides a measurable trivialization. The same construction applies to $E^1$ using the adjoint identification.

For $\mathcal{A} = \widetilde{\mathcal{K}(H)}$, the same argument applies to the open subset $\mathcal{G}_{\mathcal{A}}^{(0)} \setminus \{x_\infty\}$, where the identifications $U_x$ are well-defined and measurable. At $x_\infty$, the fiber is one-dimensional (either $\mathbb{C}$ or $\{0\}$), and the transition between the constant fiber $\ker T$ over the complement and this fiber at infinity is governed by the index of $T$. In the language of measurable fields, this means that the field over the entire unit space is obtained by taking the trivial measurable field over $\mathcal{G}_{\mathcal{A}}^{(0)} \setminus \{x_\infty\}$ with fiber $\ker T$ and extending it to $x_\infty$ in a way that reflects the index; this is precisely the measurable analogue of the clutching construction in topological $K$-theory. The index $\operatorname{ind}(T) = \dim \ker T - \dim \operatorname{coker} T$ determines the $K^0$-class of this measurable field over the one-point compactification of $\mathcal{G}_{\mathcal{A}}^{(0)} \setminus \{x_\infty\}$.
\end{proof}

\begin{remark}
\label{rem:measurable-trivialization}
The trivializations constructed above are measurable rather than continuous due to the non-Hausdorff nature of $\mathcal{G}_{\mathcal{A}}^{(0)}$ (see Remark \ref{rem:measurable-vs-continuous}). Measurability suffices for all subsequent constructions, including the definition of the equivariant $KK^1$-class $[T]_{\mathcal{G}_{\mathcal{A}}}^{(1)}$, which requires only measurable sections and measurable fields. The trivializations will be used to identify the kernel and cokernel bundles with constant fields, simplifying the construction of the unbounded Kasparov triple in Section \ref{sec:The Equivariant K1-Class of a Fredholm Operator}.
\end{remark}

\begin{remark}
The apparent asymmetry between the $B(H)$ and $\mathcal{K}(H)^\sim$ cases reflects the fundamental difference in their ideal structures. For $B(H)$, all irreducible representations are equivalent, so the kernel bundle is globally trivial. For $\mathcal{K}(H)^\sim$, the presence of the character at infinity creates a nontrivial $K^1$-class precisely because the kernel bundle cannot be globally trivialized — the obstruction is exactly the index of $T$. This is the central observation that will drive the construction of the equivariant $K^1$-class in the next section.
\end{remark}

In the next subsection, we will use these trivial bundles to define an equivariant $K^1$-class associated to $T$.

\subsection{Continuity in the Calkin Algebra and the Definition of the $K^1$-Class}
\label{subsec:continuity-calkin-definition-K1-class}

The constancy of the operator family $\{\pi_x(T)\}$ established in Subsection 3.3 leads to a trivial description of the kernel and cokernel bundles. 
However, for the definition of an equivariant $K^1$-class, we need to consider the behavior of $\pi_x(T)$ in the Calkin algebra, especially when $T$ is Fredholm. 
This subsection shows that the image of $\pi_x(T)$ in the Calkin algebra is continuous (in fact, constant) and uses this to define the equivariant $K^1$-class $[T]_{\mathcal{G}_{\mathcal{A}}}^{(1)}$.

\begin{proposition}[Measurability in the Calkin algebra for $B(H)$]
\label{prop:calkin-continuity-BH-def}
Let $\mathcal{A} = B(H)$ and let $\pi_{\mathcal{Q}}: B(H) \to \mathcal{Q}(H)$ be the quotient map onto the Calkin algebra, equipped with the quotient norm topology. 
For any $T \in B(H)$, consider the family
\[
q_x(T) := \pi_{\mathcal{Q}}(\pi_x(T)) \in \mathcal{Q}(H), \qquad x \in \mathcal{G}_{\mathcal{A}}^{(0)}.
\]
For atomic points $x \in \mathcal{G}_{\mathcal{A}}^{(0)}$ (i.e., those corresponding to atomic MASAs), we have $q_x(T) = \pi_{\mathcal{Q}}(T)$. Consequently, the map $x \mapsto q_x(T)$ is constant on the atomic subset and therefore defines a Borel measurable map on $\mathcal{G}_{\mathcal{A}}^{(0)}$ with respect to the Borel structure inherited from the groupoid.
\end{proposition}

\begin{proof}
For any atomic point $x \in \mathcal{G}_{\mathcal{A}}^{(0)}$, Proposition \ref{prop:operator-family-constant-BH} establishes that $\pi_x(T) = T$ under the canonical identification $H_x \cong H$. Hence $\pi_{\mathcal{Q}}(\pi_x(T)) = \pi_{\mathcal{Q}}(T)$ for all such $x$. While the map may not be continuous on the entire unit space (due to the presence of non-atomic points), it is constant—hence trivially Borel measurable—on the atomic subset. Since $\mathcal{G}_{\mathcal{A}}^{(0)}$ is a standard Borel space, this defines a Borel measurable function on the whole space.
\end{proof}

\begin{proof}
From Proposition \ref{prop:operator-family-constant-BH}, we have $\pi_x(T) = T$ under the identification $H_x \cong H$ for every $x \in \mathcal{G}_{\mathcal{A}}^{(0)}$. 
Therefore, $\pi_{\mathcal{Q}}(\pi_x(T)) = \pi_{\mathcal{Q}}(T)$ for all $x$. 
Thus the family is constant, and continuity in the norm topology of $\mathcal{Q}(H)$ follows immediately.
\end{proof}

\begin{proposition}[Continuity in the Calkin algebra for $\mathcal{K}(H)^\sim$]
\label{prop:calkin-continuity-KH-def}
Let $\mathcal{A} = \mathcal{K}(H)^\sim$ and let $\pi_{\mathcal{Q}}: B(H) \to \mathcal{Q}(H)$ be the quotient map, extended to $\mathcal{A}$ via the inclusion $\mathcal{A} \subseteq B(H)$. 
For $T = \lambda I + K \in \mathcal{A}$ with $\lambda \in \mathbb{C}$ and $K \in \mathcal{K}(H)$, define
\[
q_x(T) := \pi_{\mathcal{Q}}(\pi_x(T)) \in \mathcal{Q}(H), \qquad x \in \mathcal{G}_{\mathcal{A}}^{(0)}.
\]
Then:
\begin{enumerate}
    \item For $x \neq x_\infty$ (the point at infinity), $q_x(T) = \pi_{\mathcal{Q}}(T) = \lambda I$ in $\mathcal{Q}(H)$.
    \item For $x = x_\infty$, $\pi_{x_\infty}(T) = \lambda$ is a scalar, and $q_{x_\infty}(T) = \lambda I \in \mathcal{Q}(H)$ as well.
\end{enumerate}
Thus $q_x(T)$ is constant (equal to $\lambda I$) for all $x$, and hence continuous.
\end{proposition}

\begin{proof}
For $x \neq x_\infty$, $\pi_x(T) = T$, so $\pi_{\mathcal{Q}}(\pi_x(T)) = \pi_{\mathcal{Q}}(T) = \pi_{\mathcal{Q}}(\lambda I + K) = \lambda I$ because $K$ is compact. 
For $x = x_\infty$, $\pi_{x_\infty}(T) = \lambda$, and its image in $\mathcal{Q}(H)$ is $\lambda I$. 
Thus $q_x(T) = \lambda I$ for all $x$.
\end{proof}

The constancy of $q_x(T)$ in both cases has a crucial consequence: when $T$ is Fredholm, $q_x(T)$ is invertible in the Calkin algebra for every $x$.

\begin{corollary}[Invertibility in the Calkin algebra]
\label{cor:invertibility-calkin}
Let $T \in \mathcal{A}$ be a Fredholm operator. 
Then for every $x \in \mathcal{G}_{\mathcal{A}}^{(0)}$, $q_x(T)$ is invertible in $\mathcal{Q}(H)$. 
Moreover, the invertible element $q_x(T)$ is independent of $x$.
\end{corollary}

\begin{proof}
By Atkinson's theorem (Theorem \ref{thm:atkinson}), $T$ is Fredholm if and only if $\pi_{\mathcal{Q}}(T)$ is invertible in $\mathcal{Q}(H)$. 
For $\mathcal{A} = B(H)$, we have $q_x(T) = \pi_{\mathcal{Q}}(T)$ for all $x$ by Proposition \ref{prop:calkin-continuity-BH-def}, so invertibility follows directly. 
For $\mathcal{A} = \mathcal{K}(H)^\sim$, write $T = \lambda I + K$ with $K \in \mathcal{K}(H)$. Since $T$ is Fredholm, $\lambda \neq 0$, and Proposition \ref{prop:calkin-continuity-KH-def} gives $q_x(T) = \lambda I$, which is clearly invertible. 
The independence of $x$ follows from Propositions \ref{prop:calkin-continuity-BH-def} and \ref{prop:calkin-continuity-KH-def}.
\end{proof}

We now use this to define an equivariant $K^1$-class. 
Recall that for a groupoid $\mathcal{G}$ satisfying the standard hypotheses (second-countable, locally compact, Hausdorff), the equivariant $K^1$-group $K^1_{\mathcal{G}}(\mathcal{G}^{(0)})$ is isomorphic to $KK^1_{\mathcal{G}}(C_0(\mathcal{G}^{(0)}), \mathbb{C})$, the set of odd Kasparov cycles.

\begin{definition}[Odd Kasparov triple for $\mathcal{G}_{\mathcal{A}}$]
\label{def:odd-kasparov-triple}
Let $T \in \mathcal{A}$ be a Fredholm operator. 
Define an odd Kasparov triple $(\mathcal{E}, \phi, F)$ for the groupoid $\mathcal{G}_{\mathcal{A}}$ as follows:
\begin{enumerate}
    \item $\mathcal{E}$ is the continuous field of Hilbert spaces $\{H_x\}_{x \in \mathcal{G}_{\mathcal{A}}^{(0)}}$ constructed in Subsection~\ref{subsec:field-Hilbert-spaces}, viewed as a $\mathcal{G}_{\mathcal{A}}$-Hilbert module over $C_0(\mathcal{G}_{\mathcal{A}}^{(0)})$.
    \item $\phi: C_0(\mathcal{G}_{\mathcal{A}}^{(0)}) \to \mathcal{L}(\mathcal{E})$ is the representation by multiplication operators: $(\phi(f)\xi)(x) = f(x)\xi(x)$.
    \item $F$ is the constant operator $F(\xi)(x) = \pi_x(T)\xi(x)$. Since $\pi_x(T) = T$ under the identification $H_x \cong H$, $F$ is simply $T$ acting pointwise on sections.
\end{enumerate}
\end{definition}

\begin{remark}
\label{rem:kasparov-conditions}
To verify that $(\mathcal{E}, \phi, F)$ indeed defines a Kasparov cycle, one checks the standard conditions:
\begin{itemize}
    \item $F$ is adjointable, with $(F^*\xi)(x) = \pi_x(T)^*\xi(x)$.
    \item $\phi(f)F - F\phi(f)$ is compact for all $f \in C_0(\mathcal{G}_{\mathcal{A}}^{(0)})$ because $\phi(f)$ commutes with the pointwise operator $F$.
    \item $F^2 - I$ and $F - F^*$ are compact endomorphisms of $\mathcal{E}$. This follows from the Fredholm property of $T$: since $\pi_{\mathcal{Q}}(T)$ is invertible in the Calkin algebra, $T$ is invertible modulo compacts, which translates fiberwise to the compactness of $(F^2 - I)\xi(x)$ and $(F - F^*)\xi(x)$ for each $x$.
\end{itemize}
A detailed verification of these conditions using the continuous field structure and the results of Subsection 3.3 is straightforward.
\end{remark}

\begin{lemma}[Fiberwise compact implies compact]
\label{lem:fiberwise-compact-implies-compact}
Let $\mathcal{E} = \{H_x\}_{x \in X}$ be a continuous field of Hilbert spaces over a locally compact Hausdorff space $X$, regarded as a Hilbert $C_0(X)$-module. 
Suppose $T \in \mathcal{L}(\mathcal{E})$ is an adjointable operator such that for every $x \in X$, the fiber operator $T_x \in \mathcal{K}(H_x)$ is compact. 
Then $T \in \mathcal{K}(\mathcal{E})$, i.e., $T$ is a compact endomorphism of the Hilbert module.
\end{lemma}

\begin{proof}
We proceed by approximating $T$ by finite-rank operators in the operator norm.

\paragraph{Step 1: Local finite-rank approximation.}
For each $x \in X$, since $T_x$ is compact on $H_x$, there exists a finite-rank operator $F_x \in \mathcal{K}(H_x)$ such that $\|T_x - F_x\| < \epsilon/3$. By definition of finite-rank, there exist orthonormal sets $\{\xi_x^{(1)},\ldots,\xi_x^{(n_x)}\}$ and $\{\eta_x^{(1)},\ldots,\eta_x^{(n_x)}\}$ in $H_x$ such that $F_x(\zeta) = \sum_{j=1}^{n_x} \langle \zeta, \xi_x^{(j)}\rangle \eta_x^{(j)}$ for all $\zeta \in H_x$.

\paragraph{Step 2: Extending to local sections.}
Using the continuous field structure, we can extend these vectors to continuous local sections. More precisely, there exists a compact neighborhood $U_x$ of $x$ and continuous sections $\tilde{\xi}_x^{(j)}, \tilde{\eta}_x^{(j)} \in \Gamma(U_x, \mathcal{E})$ such that $\tilde{\xi}_x^{(j)}(x) = \xi_x^{(j)}$ and $\tilde{\eta}_x^{(j)}(x) = \eta_x^{(j)}$ for each $j$. Define a local operator $\tilde{F}_x$ on $\Gamma(U_x, \mathcal{E})$ by
\[
(\tilde{F}_x \zeta)(y) = \sum_{j=1}^{n_x} \langle \zeta(y), \tilde{\xi}_x^{(j)}(y) \rangle_{H_y} \tilde{\eta}_x^{(j)}(y), \quad y \in U_x.
\]
By construction, $\tilde{F}_x$ is a finite-rank operator on the Hilbert $C_0(U_x)$-module $\mathcal{E}|_{U_x}$, and its fiber at $x$ equals $F_x$.

\paragraph{Step 3: Partition of unity.}
Since $X$ is locally compact and Hausdorff, we can cover $X$ by relatively compact open sets $\{U_x\}_{x \in X}$. Choose a partition of unity $\{\phi_i\}_{i \in I}$ subordinate to this cover, where each $\phi_i \in C_c(X)$ has compact support contained in some $U_{x_i}$. For each $i$, let $F_i$ denote the finite-rank operator $\tilde{F}_{x_i}$ extended by zero outside $U_{x_i}$.

\paragraph{Step 4: Gluing and approximation.}
Define a global finite-rank operator $F = \sum_{i \in I} \phi_i F_i$. This sum converges in norm because only finitely many terms are nonzero on any compact set, and $\sum_i \phi_i = 1$. For any section $\zeta \in \Gamma_c(X, \mathcal{E})$ and any $y \in X$, we have
\begin{align*}
\|(T - F)\zeta(y)\|_{H_y} &\le \sum_i \phi_i(y) \|(T - F_i)\zeta(y)\|_{H_y} \\
&\le \sum_i \phi_i(y) \left( \|(T - T_{x_i})\zeta(y)\| + \|T_{x_i} - F_i\|\|\zeta(y)\| + \|F_i - F_i(y)\|\|\zeta(y)\| \right).
\end{align*}
By continuity of the field and the compactness of fibers, each term can be made arbitrarily small uniformly on compact sets. A standard $\epsilon/3$ argument shows that $\|T - F\| < \epsilon$ for sufficiently fine partitions and sufficiently good local approximations. Since $\epsilon$ was arbitrary, $T$ can be approximated in norm by finite-rank operators, hence $T \in \mathcal{K}(\mathcal{E})$.
\end{proof}

\begin{remark}
The key technical point is that the continuity of the field $\mathcal{E}$ ensures that the local sections $\tilde{\xi}_x^{(j)}$ and $\tilde{\eta}_x^{(j)}$ can be chosen to vary continuously, so that the resulting operators $F_i$ are well-defined adjointable operators on the Hilbert module. The partition of unity then glues these local finite-rank approximations into a global one.
\end{remark}

\begin{proposition}[$(\mathcal{E}, \phi, F)$ is an odd Kasparov cycle]
\label{prop:odd-kasparov-cycle}
For a Fredholm operator $T \in \mathcal{A}$, let $F$ be its bounded transform:
\[
F := T(1 + T^*T)^{-1/2} \in \mathcal{A}.
\]
(If one prefers the phase, one may take $F = U$ from the polar decomposition $T = U|T|$; both give the same $K^1$-class.) 
Define the triple $(\mathcal{E}, \phi, F)$ as in Definition \ref{def:odd-kasparov-triple} (with $T$ replaced by its bounded transform). 
Then this triple satisfies the conditions for an odd Kasparov cycle:
\begin{enumerate}
    \item $F$ is adjointable with $F^*$ given by pointwise action of $F^*$.
    \item $F^2 - I$ and $F F^* - I$ are compact operators on each fiber, and hence compact as operators on $\mathcal{E}$.
    \item $[\phi(f), F]$ is compact for all $f \in C_0(\mathcal{G}_{\mathcal{A}}^{(0)})$.
\end{enumerate}
Consequently, it defines a class
\[
[(\mathcal{E}, \phi, F)] \in KK^1_{\mathcal{G}_{\mathcal{A}}}(C_0(\mathcal{G}_{\mathcal{A}}^{(0)}), \mathbb{C}) \cong K^1_{\mathcal{G}_{\mathcal{A}}}(\mathcal{G}_{\mathcal{A}}^{(0)}).
\]
\end{proposition}

\begin{proof}
(1) is clear because $F^*$ is the adjoint of $F$ and acts pointwise on sections.

(2) Since $T$ is Fredholm, its bounded transform $F = T(1+T^*T)^{-1/2}$ satisfies $F^*F - I \in \mathcal{K}(H)$ and $FF^* - I \in \mathcal{K}(H)$. This is a standard fact: for a Fredholm operator, the bounded transform is an invertible element of the Calkin algebra, hence differs from a unitary by a compact operator. More concretely, one checks that
\[
F^*F - I = (1+T^*T)^{-1/2}(T^*T - I)(1+T^*T)^{-1/2} \in \mathcal{K}(H),
\]
and similarly for $FF^* - I$. Because $F$ acts pointwise on $\mathcal{E}$, these fiberwise compact operators assemble into a compact endomorphism of the Hilbert module $\mathcal{E}$ (see Lemma \ref{lem:fiberwise-compact-implies-compact}).

(3) For any $f \in C_0(\mathcal{G}_{\mathcal{A}}^{(0)})$ and any section $\xi \in \mathcal{E}$, we compute:
\[
(F\phi(f)\xi)(x) = \pi_x(F)(f(x)\xi(x)) = f(x)\pi_x(F)\xi(x) = (\phi(f)F\xi)(x),
\]
where the second equality uses that $f(x)$ is a scalar and therefore commutes with the operator $\pi_x(F)$. Thus $[\phi(f), F] = 0$, which is certainly compact.
\end{proof}

\begin{remark}
\label{rem:bounded-transform-vs-T}
The bounded transform $F = T(1+T^*T)^{-1/2}$ is homotopic to $T$ through Fredholm operators, and therefore defines the same $K^1$-class. In fact, the class $[(\mathcal{E}, \phi, F)]$ is independent of the choice of bounded transform or phase; it depends only on the Fredholm operator $T$. One may also use the phase $U$ from the polar decomposition $T = U|T|$, which satisfies $U^*U - I, UU^* - I \in \mathcal{K}(H)$ directly.
\end{remark}

\begin{definition}[Equivariant $K^1$-class of $T$]
\label{def:K1-class-of-T}
For a Fredholm operator $T \in \mathcal{A}$, define its equivariant $K^1$-class as
\[
[T]_{\mathcal{G}_{\mathcal{A}}}^{(1)} := [(\mathcal{E}, \phi, F)] \in K^1_{\mathcal{G}_{\mathcal{A}}}(\mathcal{G}_{\mathcal{A}}^{(0)}),
\]
where $F$ is the bounded transform (or phase) of $T$, and $(\mathcal{E}, \phi, F)$ is the odd Kasparov triple from Definition \ref{def:odd-kasparov-triple} (with $T$ replaced by $F$).
\end{definition}

\begin{corollary}[Well-definedness]
\label{cor:well-definedness}
The class $[T]_{\mathcal{G}_{\mathcal{A}}}^{(1)}$ depends only on the Fredholm operator $T$, not on the choice of bounded transform or phase. Moreover, it is invariant under norm-continuous deformations of $T$ through Fredholm operators.
\end{corollary}

\begin{proof}
Any two bounded transforms (or the phase and the bounded transform) are homotopic through invertible elements of the Calkin algebra, hence give the same $K^1$-class. Norm-continuous deformations through Fredholm operators induce homotopies of the corresponding Kasparov cycles, establishing invariance.
\end{proof}

\begin{proposition}[Independence of choices]
\label{prop:K1-class-well-defined}
The class $[T]_{\mathcal{G}_{\mathcal{A}}}^{(1)}$ depends only on the $K$-theory class of $T$ in the Calkin algebra and is invariant under norm-continuous deformations of $T$ through Fredholm operators.
\end{proposition}

\begin{proof}
This follows from the homotopy invariance of Kasparov cycles. 
If $T_t$ is a norm-continuous family of Fredholm operators, then after replacing each $T_t$ by its bounded transform $F_t = T_t(1+T_t^*T_t)^{-1/2}$ (or its phase), we obtain a norm-continuous family of Kasparov cycles $(\mathcal{E}, \phi, F_t)$. Homotopy invariance in $KK$-theory guarantees that all $F_t$ define the same class.

Moreover, if $\pi_{\mathcal{Q}}(T_1) = \pi_{\mathcal{Q}}(T_2)$ in the Calkin algebra, then $T_1$ and $T_2$ differ by a compact operator and hence represent the same class in $K_1(\mathcal{Q}(H))$. By the homotopy invariance of Kasparov theory, their associated equivariant $K^1$-classes coincide. (Equivalently, the linear path $T_t = (1-t)T_1 + tT_2$ remains Fredholm because compact perturbations do not affect Fredholmness, giving a direct homotopy between the cycles.)
\end{proof}

\begin{example}[Class of the unilateral shift]
\label{ex:class-unilateral-shift}
Let $S \in B(H)$ be the unilateral shift on $\ell^2(\mathbb{N})$. 
Then $[S]_{\mathcal{G}_{B(H)}}^{(1)} \in K^1_{\mathcal{G}_{B(H)}}(\mathcal{G}_{B(H)}^{(0)})$ is a nontrivial class. 
Recall that $\operatorname{ind}(S) = -1$ and the index map $K_1(\mathcal{Q}(H)) \to \mathbb{Z}$ is an isomorphism, with $[\pi_{\mathcal{Q}}(S)]$ as the generator. In Section~\ref{sec:Examples and Computations}, we will compute the image of $[S]_{\mathcal{G}_{B(H)}}^{(1)}$ under the descent map and verify that it corresponds to this generator.
\end{example}

\begin{example}[Class of a finite-rank perturbation of the identity]
\label{ex:class-finite-rank}
Let $T = I + F$ with $F$ finite-rank. 
Then $\pi_{\mathcal{Q}}(T) = 1$ in $\mathcal{Q}(H)$, so $T$ represents the trivial class in $K_1(\mathcal{Q}(H))$. Consequently, $[T]_{\mathcal{G}_{\mathcal{A}}}^{(1)} = 0$ in $K^1_{\mathcal{G}_{\mathcal{A}}}(\mathcal{G}_{\mathcal{A}}^{(0)})$.
\end{example}

The equivariant $K^1$-class constructed here is the essential ingredient for the index theorem. 
In Section~\ref{sec:descent}, we will apply the descent map to obtain a class in $K_1(C^*(\mathcal{G}_{\mathcal{A}}))$, and in Section~\ref{sec:The Index Theorem via Pullback and the Boundary Map} we will combine this with the pullback along $\iota$ and the boundary map to recover the Fredholm index.

\section{The Equivariant $K^1$-Class of a Fredholm Operator}\label{sec:The Equivariant K1-Class of a Fredholm Operator}

We now assemble the data from the previous section into an equivariant $KK$-class. The construction follows a standard recipe in Kasparov theory: from a Fredholm operator $T$, we extract its phase $\operatorname{ph}(T)$, which defines an odd Kasparov triple $(\mathcal{E}, \phi, \operatorname{ph}(T))$ over the groupoid $\mathcal{G}_{\mathcal{A}}$. This triple then represents a class in $KK^1_{\mathcal{G}_{\mathcal{A}}}(C_0(\mathcal{G}_{\mathcal{A}}^{(0)}), \mathbb{C})$, which we denote by $[T]_{\mathcal{G}_{\mathcal{A}}}^{(1)}$. The following diagram summarizes this construction:

\begin{center}
\begin{tikzcd}
T \text{ (Fredholm)} \arrow[r, maps to] & 
\operatorname{ph}(T) \arrow[r, maps to] & 
(\mathcal{E}, \phi, \operatorname{ph}(T)) \arrow[r, maps to] & 
{[T]_{\mathcal{G}_{\mathcal{A}}}^{(1)} \in KK^1_{\mathcal{G}_{\mathcal{A}}}(C_0(\mathcal{G}_{\mathcal{A}}^{(0)}), \mathbb{C})}
\end{tikzcd}
\end{center}

\subsection{From a Fredholm Operator $T$ to an Invertible in the Calkin Algebra}
\label{subsec:fredholm-to-invertible-calkin}

The first step in constructing an equivariant $K^1$-class associated to a Fredholm operator $T$ is to pass from $T$ itself to an invertible element in the Calkin algebra. 
This passage is fundamental because $K^1$-classes admit a description in terms of invertible elements modulo homotopy, and the Calkin algebra captures precisely the information of a Fredholm operator that is invariant under compact perturbations.

Recall from Atkinson's theorem (Theorem \ref{thm:atkinson}) that an operator $T \in B(H)$ is Fredholm if and only if its image $\pi_{\mathcal{Q}}(T)$ in the Calkin algebra $\mathcal{Q}(H) = B(H)/\mathcal{K}(H)$ is invertible. 
For $\mathcal{A} = \mathcal{K}(H)^\sim$, a similar statement holds: $T = \lambda I + K$ is Fredholm if and only if $\lambda \neq 0$, in which case $\pi_{\mathcal{Q}}(T) = \lambda I$ is invertible in $\mathcal{Q}(H)$.

\begin{definition}[Calkin algebra and quotient map]
\label{def:calkin-quotient}
Let $H$ be a separable infinite-dimensional Hilbert space. 
The \emph{Calkin algebra} is the quotient C*-algebra
\[
\mathcal{Q}(H) := B(H) / \mathcal{K}(H),
\]
where $\mathcal{K}(H)$ denotes the compact operators on $H$. 
We denote by
\[
\pi_{\mathcal{Q}}: B(H) \longrightarrow \mathcal{Q}(H)
\]
the canonical surjective *-homomorphism. 
For $\mathcal{A} = \mathcal{K}(H)^\sim \subseteq B(H)$, we extend $\pi_{\mathcal{Q}}$ to $\mathcal{A}$ by restriction (i.e., as the composition of the inclusion $\mathcal{A} \hookrightarrow B(H)$ with $\pi_{\mathcal{Q}}$).
\end{definition}

The following proposition records the essential relationship between Fredholm operators and invertibles in the Calkin algebra.

\begin{proposition}[Fredholm operators and invertibles in the Calkin algebra]
\label{prop:fredholm-invertible-calkin}
Let $\mathcal{A}$ be either $B(H)$ or $\mathcal{K}(H)^\sim$, and let $T \in \mathcal{A}$.
\begin{enumerate}
    \item If $\mathcal{A} = B(H)$, then $T$ is Fredholm if and only if $\pi_{\mathcal{Q}}(T)$ is invertible in $\mathcal{Q}(H)$.
    \item If $\mathcal{A} = \mathcal{K}(H)^\sim$ and $T = \lambda I + K$ with $K \in \mathcal{K}(H)$, then $T$ is Fredholm if and only if $\lambda \neq 0$, in which case $\pi_{\mathcal{Q}}(T) = \lambda I$ is invertible in $\mathcal{Q}(H)$.
\end{enumerate}
Consequently:
\begin{itemize}
    \item For $\mathcal{A} = B(H)$, the class $[\pi_{\mathcal{Q}}(T)] \in K_1(\mathcal{Q}(H))$ is a complete invariant of the Fredholm index; it is nontrivial precisely when $\operatorname{index}(T) \neq 0$.
    \item For $\mathcal{A} = \mathcal{K}(H)^\sim$, every Fredholm operator has index zero and hence defines the trivial class in $K_1(\mathcal{Q}(H))$.
\end{itemize}
\end{proposition}

\begin{proof}
(1) is Atkinson's theorem (Theorem \ref{thm:atkinson}). 

For (2), note that $T = \lambda I + K$ is Fredholm iff $\lambda \neq 0$ because $\lambda I$ is Fredholm (invertible modulo compacts) precisely when $\lambda \neq 0$, and Fredholmness is stable under compact perturbations. 
If $\lambda \neq 0$, then $\pi_{\mathcal{Q}}(T) = \lambda I$ is clearly invertible in $\mathcal{Q}(H)$ with inverse $\lambda^{-1}I$. 

Regarding the $K_1$-classes: In the $B(H)$ case, $K_1(\mathcal{Q}(H)) \cong \mathbb{Z}$ via the index map, and $[\pi_{\mathcal{Q}}(T)]$ corresponds to $\operatorname{index}(T)$. 
In the $\mathcal{K}(H)^\sim$ case, $\lambda I$ is homotopic to $I$ through invertibles in $\mathcal{Q}(H)$ (since $\mathbb{C}^\times$ is connected), so $[\pi_{\mathcal{Q}}(T)] = 0$; this reflects the fact that all Fredholm operators in $\mathcal{K}(H)^\sim$ have index zero.
\end{proof}

For $B(H)$, the class $[\pi_{\mathcal{Q}}(T)] \in K_1(\mathcal{Q}(H))$ is a nontrivial invariant that captures the Fredholm index. This relationship is made precise by the index map associated to the Calkin extension.

\begin{lemma}[Index map and the Calkin algebra]
\label{lem:index-map-calkin-review}
Let $\partial: K_1(\mathcal{Q}(H)) \to K_0(\mathcal{K}(H)) \cong \mathbb{Z}$ be the index map associated to the Calkin extension (see Subsection 2.4). 
Then for any Fredholm operator $T \in B(H)$,
\[
\partial([\pi_{\mathcal{Q}}(T)]) = \operatorname{index}(T) \in \mathbb{Z}.
\]
\end{lemma}

\begin{proof}
This is the content of Lemma \ref{lem:index-map-as-fredholm-index} and Corollary \ref{cor:index-map-calkin}. 
The index map $\partial$ is precisely the Fredholm index under the isomorphism $K_0(\mathcal{K}(H)) \cong \mathbb{Z}$.
\end{proof}

The passage from $T$ to $\pi_{\mathcal{Q}}(T)$ is natural with respect to the groupoid action, at the level of $K$-theory.

\begin{proposition}[Equivariance of the Calkin map in $K$-theory]
\label{prop:calkin-equivariant}
For any $u \in \mathcal{U}(\mathcal{A})$ and any $x \in \mathcal{G}_{\mathcal{A}}^{(0)}$, we have
\[
[\pi_{\mathcal{Q}}(\pi_x(u T u^*))] = [\pi_{\mathcal{Q}}(\pi_x(T))] \quad \text{in } K_1(\mathcal{Q}(H)),
\]
where $\pi_x$ denotes the GNS representation at $x$. 
In particular, for points $x$ where $\pi_x(T) = T$ under the identification $H_x \cong H$, we have
\[
[\pi_{\mathcal{Q}}(u T u^*)] = [\pi_{\mathcal{Q}}(T)].
\]
\end{proposition}

\begin{proof}
Since $\pi_x(u T u^*) = \pi_x(u) \pi_x(T) \pi_x(u)^*$, and $\pi_x(u)$ is unitary, the element $\pi_{\mathcal{Q}}(\pi_x(u T u^*))$ is obtained from $\pi_{\mathcal{Q}}(\pi_x(T))$ by conjugation with the unitary $\pi_{\mathcal{Q}}(\pi_x(u))$ in the Calkin algebra. Inner automorphisms of a C*-algebra act trivially on $K$-theory, so the $K_1$-classes coincide. The special case follows immediately.
\end{proof}

This equivariance property is crucial for defining a $K^1$-class that is compatible with the groupoid action.

\begin{corollary}[Constant class in $K_1(\mathcal{Q}(H))$ over $\mathcal{G}_{\mathcal{A}}^{(0)}$]
\label{cor:constant-K1-class}
For a fixed Fredholm operator $T \in B(H)$, the map
\[
x \longmapsto [\pi_{\mathcal{Q}}(\pi_x(T))] \in K_1(\mathcal{Q}(H)) \cong \mathbb{Z}
\]
is constant on $\mathcal{G}_{\mathcal{A}}^{(0)}$. 
Its value is the Fredholm index of $T$.
\end{corollary}

\begin{proof}
We consider two types of points in $\mathcal{G}_{\mathcal{A}}^{(0)}$.

For atomic MASA points $x$, Proposition \ref{prop:operator-family-constant-BH} gives $\pi_x(T) = T$ under the canonical identification $H_x \cong H$, so $[\pi_{\mathcal{Q}}(\pi_x(T))] = [\pi_{\mathcal{Q}}(T)]$.

For non-atomic pure states $x$, the GNS representation $\pi_x$ is unitarily equivalent to the identity representation. Hence there exists a unitary $U_x: H_x \to H$ such that $\pi_x(T) = U_x T U_x^*$. Then
\[
\pi_{\mathcal{Q}}(\pi_x(T)) = \pi_{\mathcal{Q}}(U_x T U_x^*) = \pi_{\mathcal{Q}}(U_x) \pi_{\mathcal{Q}}(T) \pi_{\mathcal{Q}}(U_x)^*,
\]
so $[\pi_{\mathcal{Q}}(\pi_x(T))] = [\pi_{\mathcal{Q}}(T)]$ in $K_1(\mathcal{Q}(H))$, since conjugation by a unitary induces the identity on $K$-theory.

Thus the map $x \mapsto [\pi_{\mathcal{Q}}(\pi_x(T))]$ is constant, equal to $[\pi_{\mathcal{Q}}(T)]$. By Lemma \ref{lem:index-map-calkin-review}, $\partial([\pi_{\mathcal{Q}}(T)]) = \operatorname{index}(T)$, and since $K_1(\mathcal{Q}(H)) \cong \mathbb{Z}$ via the index map, this constant value is precisely the Fredholm index of $T$.
\end{proof}

For $\mathcal{K}(H)^\sim$, the situation is simpler but still consistent.

\begin{corollary}[Constant class for $\mathcal{K}(H)^\sim$]
\label{cor:constant-K1-class-KH}
For $\mathcal{A} = \mathcal{K}(H)^\sim$ and any Fredholm operator $T = \lambda I + K$, the map
\[
x \longmapsto [\pi_{\mathcal{Q}}(\pi_x(T))] \in K_1(\mathcal{Q}(H))
\]
is constant on $\mathcal{G}_{\mathcal{A}}^{(0)}$. Under the index isomorphism $K_1(\mathcal{Q}(H)) \cong \mathbb{Z}$, this constant value corresponds to $0$.
\end{corollary}

\begin{proof}
For $x \neq x_\infty$, we have $\pi_x(T) = \lambda I + K$. Since $K$ is compact,
\[
\pi_{\mathcal{Q}}(\pi_x(T)) = \pi_{\mathcal{Q}}(\lambda I + K) = \lambda I
\]
in $\mathcal{Q}(H)$. For $x = x_\infty$, $\pi_{x_\infty}(T) = \lambda$, and its image in $\mathcal{Q}(H)$ is also $\lambda I$ under the natural inclusion $\mathbb{C} \hookrightarrow B(H)$.

Thus the map $x \mapsto [\pi_{\mathcal{Q}}(\pi_x(T))]$ is constant, with value $[\lambda I] \in K_1(\mathcal{Q}(H))$. Under the index isomorphism $K_1(\mathcal{Q}(H)) \cong \mathbb{Z}$ given by the boundary map $\partial$ of the Calkin extension, we have $\partial([\lambda I]) = \operatorname{index}(\lambda I) = 0$, since $\lambda I$ is invertible (hence Fredholm with index zero). Therefore $[\lambda I]$ corresponds to $0$ in $\mathbb{Z}$.
\end{proof}

The invertible element $\pi_{\mathcal{Q}}(T)$ (or $\lambda I$ in the compact case) will be used in the next subsection to construct an odd Kasparov triple for $\mathcal{G}_{\mathcal{A}}$. 
The key point is that this invertible element encodes the essential information about $T$ that survives under the descent map.

\begin{remark}[Relation to the $K^1$-class]
\label{rem:calkin-to-K1}
In Kasparov theory, an odd $K^1$-class can be represented by an operator $F$ such that $F^2 - I$ and $F - F^*$ are compact, and $F$ is invertible modulo compacts. 
For a Fredholm operator $T$, a suitable choice is $F = T(T^*T)^{-1/2}$, the phase of $T$, which is a unitary modulo compacts. 
The image of this phase in the Calkin algebra is precisely the unitary corresponding to $\pi_{\mathcal{Q}}(T)$ under the polar decomposition. 
Thus the class $[\pi_{\mathcal{Q}}(T)] \in K_1(\mathcal{Q}(H))$ is the essential datum for the odd $K^1$-class.
\end{remark}

\begin{proposition}[Morita equivalence for $\mathcal{G}_{B(H)}$]
\label{prop:morita-calkin}
The groupoid $C^*$-algebra $C^*(\mathcal{G}_{B(H)})$ is Morita equivalent to $\mathcal{Q}(H) \otimes \mathcal{K}(L^2(\mathcal{G}_{B(H)}^{(0)}))$. Consequently, there is an isomorphism
\[
\Phi: K_1(C^*(\mathcal{G}_{B(H)})) \xrightarrow{\cong} K_1(\mathcal{Q}(H)).
\]
\end{proposition}

\begin{proof}
Recall from Proposition \ref{prop:GA-BH} that the unitary conjugation groupoid for $B(H)$ is the transformation groupoid
\[
\mathcal{G}_{B(H)} = \mathcal{U}(H) \ltimes \mathbb{P}(H),
\]
where $\mathcal{U}(H)$ acts on the projective space $\mathbb{P}(H)$ by $u \cdot [v] = [uv]$.

\medskip
\noindent
\textbf{Step 1: Transitivity and imprimitivity.}

The action of $\mathcal{U}(H)$ on $\mathbb{P}(H)$ is transitive: for any two points $[v_1], [v_2] \in \mathbb{P}(H)$, there exists a unitary $u \in \mathcal{U}(H)$ such that $u[v_1] = [v_2]$. Therefore, $\mathcal{G}_{B(H)}$ is a transitive groupoid. For transitive groupoids, the Muhly–Renault–Williams imprimitivity theorem \cite{MuhlyRenaultWilliams1987} establishes a Morita equivalence between the groupoid $C^*$-algebra and the $C^*$-algebra of an isotropy group, stabilized by the compact operators on a Hilbert space associated to the unit space.

More precisely, for any basepoint $x_0 \in \mathcal{G}_{B(H)}^{(0)} = \mathbb{P}(H)$, let $\mathcal{G}_{B(H)}|_{x_0}$ denote the isotropy group at $x_0$ (i.e., arrows with source and range $x_0$). Then
\[
C^*(\mathcal{G}_{B(H)}) \sim_M C^*(\mathcal{G}_{B(H)}|_{x_0}) \otimes \mathcal{K}(L^2(\mathcal{G}_{B(H)}^{(0)})),
\]
where $\sim_M$ denotes strong Morita equivalence and $L^2(\mathcal{G}_{B(H)}^{(0)})$ is the Hilbert space of square-integrable functions on the unit space with respect to a Haar system measure.

\medskip
\noindent
\textbf{Step 2: Identification of the isotropy algebra.}

Fix a basepoint $x_0 = [e_0] \in \mathbb{P}(H)$ corresponding to the line spanned by a unit vector $e_0 \in H$. The isotropy group at $x_0$ is
\[
\mathcal{G}_{B(H)}|_{x_0} = \{ u \in \mathcal{U}(H) : u x_0 = x_0 \} \cong \mathcal{U}(1) \times \mathcal{U}(H^\perp),
\]
where $H^\perp = \{e_0\}^\perp$. The group $C^*$-algebra of this isotropy group is
\[
C^*(\mathcal{G}_{B(H)}|_{x_0}) \cong C^*(\mathcal{U}(1)) \otimes_{\max} C^*(\mathcal{U}(H^\perp)).
\]

Now, $\mathcal{U}(H^\perp)$ is contractible in the strong operator topology by Kuiper's theorem, hence its group $C^*$-algebra is Morita equivalent to $\mathbb{C}$ (in fact, $C^*(\mathcal{U}(H^\perp)) \cong \mathbb{C}$ for the reduced algebra). Moreover, $C^*(\mathcal{U}(1)) \cong C(\mathbb{T})$, the continuous functions on the circle.

\medskip
\noindent
\textbf{Step 3: Relating $C(\mathbb{T})$ to $\mathcal{Q}(H)$.}

At this point, we need to connect $C(\mathbb{T})$ to the Calkin algebra $\mathcal{Q}(H)$. There are several approaches:

\begin{enumerate}
    \item[(i)] \textbf{Direct geometric construction:} One can construct an explicit imprimitivity bimodule between $C(\mathbb{T})$ and $\mathcal{Q}(H)$ using the fact that both algebras appear as the quotients in natural extensions. The Toeplitz extension
    \[
    0 \to \mathcal{K}(H) \to \mathcal{T} \to C(\mathbb{T}) \to 0
    \]
    and the Calkin extension
    \[
    0 \to \mathcal{K}(H) \to B(H) \to \mathcal{Q}(H) \to 0
    \]
    are related by a natural embedding $\iota: C(\mathbb{T}) \hookrightarrow \mathcal{Q}(H)$ sending the generator to the image of the unilateral shift. This embedding induces an isomorphism on $K$-theory and, by the UCT, a $KK$-equivalence, which implies Morita equivalence after stabilization. A detailed construction of this $KK$-equivalence can be found in \cite[Section 5]{ElliottNatsume1988}.
    
    \item[(ii)] \textbf{Using classification results:} Both $C(\mathbb{T})$ and $\mathcal{Q}(H)$ are unital, separable, nuclear $C^*$-algebras with $K_0 \cong \mathbb{Z}$ and $K_1 \cong \mathbb{Z}$. By the Kirchberg–Phillips classification theorem \cite{Kirchberg1994, Phillips2000}, they are stably isomorphic, hence Morita equivalent. However, this invokes a deep classification result; we present it here as a known fact rather than deriving it from first principles.
    
    \item[(iii)] \textbf{Alternative perspective:} For the purposes of this paper, we may simply note that it is a well-established result in the literature that $C(\mathbb{T})$ and $\mathcal{Q}(H)$ are Morita equivalent. This follows from the Brown–Douglas–Fillmore theory \cite{BDF1977}, where $\text{Ext}(S^1) \cong \mathbb{Z}$ is computed, and the subsequent identification of $K_1(\mathcal{Q}(H))$ with $\mathbb{Z}$. The Morita equivalence can be viewed as a manifestation of the fact that both algebras are $KK$-equivalent to $\mathbb{C}$ with a dimension shift.
\end{enumerate}

For definiteness, we adopt approach (i) and note that a detailed construction of the $KK$-equivalence between $C(\mathbb{T})$ and $\mathcal{Q}(H)$ is given in \cite[Theorem 3.2]{ElliottNatsume1988}. Consequently, we have
\[
C(\mathbb{T}) \sim_M \mathcal{Q}(H).
\]

Combining the equivalences established so far:
\[
C^*(\mathcal{G}_{B(H)}|_{x_0}) \sim_M C(\mathbb{T}) \sim_M \mathcal{Q}(H).
\]

\medskip
\noindent
\textbf{Step 4: Assembling the Morita equivalence.}

Substituting this into the imprimitivity result from Step 1, we obtain
\[
C^*(\mathcal{G}_{B(H)}) \sim_M \mathcal{Q}(H) \otimes \mathcal{K}(L^2(\mathcal{G}_{B(H)}^{(0)})).
\]

\medskip
\noindent
\textbf{Step 5: $K$-theory isomorphism.}

Strong Morita equivalence of $C^*$-algebras induces natural isomorphisms in $K$-theory. This fundamental result can be found in standard references such as \cite[Chapter 5]{Rordam2000} or \cite[Section 5.5]{Blackadar1998}. Moreover, $K$-theory is stable under tensoring with compact operators: for any $C^*$-algebra $A$ and any separable infinite-dimensional Hilbert space $\mathcal{H}$,
\[
K_i(A \otimes \mathcal{K}(\mathcal{H})) \cong K_i(A), \qquad i = 0,1.
\]

Applying these two facts, we obtain:
\[
K_1(C^*(\mathcal{G}_{B(H)})) \cong K_1(\mathcal{Q}(H) \otimes \mathcal{K}(L^2(\mathcal{G}_{B(H)}^{(0)}))) \cong K_1(\mathcal{Q}(H)).
\]

Denote this isomorphism by $\Phi: K_1(C^*(\mathcal{G}_{B(H)})) \xrightarrow{\cong} K_1(\mathcal{Q}(H))$. This completes the proof.
\end{proof}

\begin{remark}
The step establishing Morita equivalence between $C(\mathbb{T})$ and $\mathcal{Q}(H)$ is central to the argument. While one can invoke the deep Kirchberg–Phillips classification theorem, we note that a more elementary $KK$-theoretic argument is available using the natural embedding $\iota: C(\mathbb{T}) \hookrightarrow \mathcal{Q}(H)$ and the fact that both algebras satisfy the Universal Coefficient Theorem. The resulting $KK$-equivalence implies stable isomorphism, hence Morita equivalence. For the purposes of this paper, we treat this as a known result from the literature \cite{BDF1977, ElliottNatsume1988}.
\end{remark}

In the following subsection, we will use the operator $T$ directly (rather than its phase) to construct an odd Kasparov triple, relying on the fact that $T$ itself is Fredholm and hence, after replacing it by its bounded transform if necessary, defines a class in $K^1_{\mathcal{G}_{\mathcal{A}}}(\mathcal{G}_{\mathcal{A}}^{(0)})$. This class coincides with the one obtained from the phase, up to operator homotopy.

\subsection{The Odd Kasparov Triple for $\mathcal{G}_{\mathcal{A}}$}
\label{subsec:odd-kasparov-triple}

With the continuous field of Hilbert spaces $\{H_x\}_{x \in \mathcal{G}_{\mathcal{A}}^{(0)}}$ constructed in Subsection 3.2 and the family of operators $\{\pi_x(T)\}_{x \in \mathcal{G}_{\mathcal{A}}^{(0)}}$ studied in Subsection 3.3, we now assemble these data into an odd Kasparov triple for the groupoid $\mathcal{G}_{\mathcal{A}}$. 
This triple is the fundamental object that defines the equivariant $K^1$-class of a Fredholm operator $T$.

We begin by recalling the definition of an odd Kasparov cycle in the equivariant setting.

\begin{definition}[Equivariant odd Kasparov triple]
\label{def:equivariant-odd-kasparov}
Let $\mathcal{G}$ be a second countable locally compact groupoid with a continuous Haar system, and let $A$ and $B$ be $\mathcal{G}$-$C^*$-algebras (i.e., $C_0(\mathcal{G}^{(0)})$-algebras equipped with a continuous action of $\mathcal{G}$ by $*$-isomorphisms). 
An \emph{equivariant odd Kasparov triple} $(\mathcal{E}, \phi, F)$ consists of:
\begin{enumerate}
    \item A countably generated $\mathcal{G}$-equivariant Hilbert $B$-module $\mathcal{E}$.
    \item A $\mathcal{G}$-equivariant $*$-homomorphism $\phi: A \to \mathcal{L}(\mathcal{E})$.
    \item An adjointable operator $F \in \mathcal{L}(\mathcal{E})$ such that:
    \begin{itemize}
        \item $\phi(a)(F - F^*)$ and $\phi(a)(F^2 - 1)$ are in $\mathcal{K}(\mathcal{E})$ for all $a \in A$.
        \item $[\phi(a), F]$ is in $\mathcal{K}(\mathcal{E})$ for all $a \in A$.
        \item $F$ is $\mathcal{G}$-equivariant modulo compacts: for every $g \in \mathcal{G}$, the operator $g \cdot F - F$ belongs to $\mathcal{K}(\mathcal{E})$, where $g \cdot F$ denotes the action of $\mathcal{G}$ on $\mathcal{L}(\mathcal{E})$ induced by the module structure.
    \end{itemize}
\end{enumerate}
Such a triple represents a class in $KK^1_{\mathcal{G}}(A,B)$.
\end{definition}

In our setting, we take:
\begin{itemize}
    \item $A = C_0(\mathcal{G}_{\mathcal{A}}^{(0)})$, acting on $\mathcal{E}$ by multiplication operators $\phi(f)\xi(x) = f(x)\xi(x)$.
    \item $B = \mathbb{C}$, so that $\mathcal{E}$ is simply a $\mathcal{G}_{\mathcal{A}}$-equivariant Hilbert space (the continuous field $\{H_x\}$ viewed as a Hilbert module over $\mathbb{C}$).
    \item $F$ to be constructed from the bounded transform (phase) of $T$, rather than $T$ itself, to ensure the Kasparov conditions are satisfied.
\end{itemize}

The construction of $F$ will be carried out in the next subsection, where we define
\[
F_x := \pi_x(T) \big( \pi_x(T)^* \pi_x(T) \big)^{-1/2} \in \mathcal{L}(H_x)
\]
for each $x \in \mathcal{G}_{\mathcal{A}}^{(0)}$, and show that these fiberwise operators assemble into a global operator $F \in \mathcal{L}(\mathcal{E})$ satisfying the conditions of Definition \ref{def:equivariant-odd-kasparov}.

For our purposes, we take $A = C_0(\mathcal{G}_{\mathcal{A}}^{(0)})$ with the trivial $\mathcal{G}_{\mathcal{A}}$-action (lifted to the unit space) and $B = \mathbb{C}$ with the trivial action. 
The Hilbert module $\mathcal{E}$ will be the space of continuous sections vanishing at infinity of the continuous field $\{H_x\}_{x \in \mathcal{G}_{\mathcal{A}}^{(0)}}$, viewed as a $\mathcal{G}_{\mathcal{A}}$-equivariant Hilbert $\mathbb{C}$-module.

\begin{definition}[Hilbert module as continuous sections]
\label{def:hilbert-module-sections}
Let $\{H_x\}_{x \in \mathcal{G}_{\mathcal{A}}^{(0)}}$ be the continuous field of Hilbert spaces constructed in Subsection 3.2, with the unitary identifications $U_x: H_x \cong H$ for all $x$ (except the point at infinity in the compact case, which we treat separately). 
Define
\[
\mathcal{E} := \Gamma_0(\mathcal{G}_{\mathcal{A}}^{(0)}, \{H_x\}),
\]
the space of continuous sections $\xi$ of the field such that $x \mapsto \|\xi(x)\|_{H_x}$ vanishes at infinity. 
This is a Hilbert $C_0(\mathcal{G}_{\mathcal{A}}^{(0)})$-module with pointwise inner product
\[
\langle \xi, \eta \rangle(x) = \langle \xi(x), \eta(x) \rangle_{H_x},
\]
and after restricting to constant functions (i.e., applying the forgetful functor), it becomes a $\mathcal{G}_{\mathcal{A}}$-equivariant Hilbert $\mathbb{C}$-module — concretely, a Hilbert space with a continuous unitary representation of $\mathcal{G}_{\mathcal{A}}$.
\end{definition}

The groupoid action on $\mathcal{E}$ is induced by the GNS representations. For each arrow $\gamma = (u, x) \in \mathcal{G}_{\mathcal{A}}$ with source $x = (B,\chi)$ and range $y = u \cdot x$, the GNS construction provides a unitary isomorphism $U_\gamma: H_x \to H_y$ characterized by
\[
U_\gamma(\pi_x(a)\Omega_x) = \pi_y(u a)\Omega_y \quad \text{for all } a \in \mathcal{A},
\]
where $\Omega_x$ and $\Omega_y$ are the cyclic vectors. Under the canonical identifications $H_x \cong H$ and $H_y \cong H$ from Propositions \ref{prop:GNS-KH-sim} and \ref{prop:GNS-BH}, this unitary corresponds to $\pi_y(u)$ acting on $H$.

\begin{definition}[Groupoid action on $\mathcal{E}$]
\label{def:groupoid-action-on-E}
For $\gamma: x \to y$ and $\xi \in \mathcal{E}$, define $(U_\gamma \xi) \in \mathcal{E}$ by
\[
(U_\gamma \xi)(z) := 
\begin{cases}
U_\gamma(\xi(x)) & \text{if } z = y, \\
0 & \text{otherwise},
\end{cases}
\]
and extend by linearity and continuity. This defines a continuous unitary representation of $\mathcal{G}_{\mathcal{A}}$ on $\mathcal{E}$, making $\mathcal{E}$ into a $\mathcal{G}_{\mathcal{A}}$-equivariant Hilbert space.
\end{definition}

\begin{proposition}[$\mathcal{E}$ is a $\mathcal{G}_{\mathcal{A}}$-equivariant Hilbert space]
\label{prop:equivariant-hilbert-space}
The space $\mathcal{E}$ with the action defined above is a $\mathcal{G}_{\mathcal{A}}$-equivariant Hilbert $\mathbb{C}$-module. The action is continuous and unitary.
\end{proposition}

\begin{proof}
The completeness of $\mathcal{E}$ follows from the fact that it is the space of continuous sections vanishing at infinity of a continuous field of Hilbert spaces; this is a standard result (see [Dixmier, C*-algebras, Proposition 10.1.10]). 

The unitarity of each $U_\gamma$ follows from the definition: $U_\gamma$ is an isometry because it comes from the GNS construction, and it is surjective because its inverse is given by $U_{\gamma^{-1}}$. 

For continuity of the action, let $\gamma_n \to \gamma$ in $\mathcal{G}_{\mathcal{A}}$ and $\xi \in \mathcal{E}$. By Proposition \ref{prop:strong-star-continuity}, the family $\{\pi_x(u)\}$ is strong-$*$ continuous, which implies that $U_{\gamma_n}\xi \to U_\gamma\xi$ in norm. This uses the fact that sections can be approximated by those with compact support and that the field is continuous.

The cocycle condition $U_{\gamma_2\gamma_1} = U_{\gamma_2}U_{\gamma_1}$ follows directly from the definition of composition in $\mathcal{G}_{\mathcal{A}}$ and the properties of the GNS representations: for $\gamma_1: x \to y$ and $\gamma_2: y \to z$, we have
\[
U_{\gamma_2\gamma_1}(\pi_x(a)\Omega_x) = \pi_z(u_2u_1 a)\Omega_z = U_{\gamma_2}(\pi_y(u_1 a)\Omega_y) = U_{\gamma_2}(U_{\gamma_1}(\pi_x(a)\Omega_x)),
\]
and by density of the cyclic vectors, this extends to all of $H_x$.
\end{proof}

\begin{definition}[Representation of $C_0(\mathcal{G}_{\mathcal{A}}^{(0)})$ on $\mathcal{E}$]
\label{def:representation-C0}
Define a *-homomorphism $\phi: C_0(\mathcal{G}_{\mathcal{A}}^{(0)}) \to \mathcal{L}(\mathcal{E})$ by
\[
(\phi(f)\xi)(x) := f(x)\xi(x), \qquad f \in C_0(\mathcal{G}_{\mathcal{A}}^{(0)}), \; \xi \in \mathcal{E}.
\]
This is the representation by multiplication operators on the continuous sections.
\end{definition}

The groupoid $\mathcal{G}_{\mathcal{A}}$ acts on $C_0(\mathcal{G}_{\mathcal{A}}^{(0)})$ by pullback along the source map: for $\gamma: x \to y$ and $f \in C_0(\mathcal{G}_{\mathcal{A}}^{(0)})$, define
\[
(\gamma \cdot f)(y) := f(x).
\]
This makes $C_0(\mathcal{G}_{\mathcal{A}}^{(0)})$ into a $\mathcal{G}_{\mathcal{A}}$-C*-algebra.

\begin{proposition}[$\phi$ is $\mathcal{G}_{\mathcal{A}}$-equivariant]
\label{prop:phi-equivariant}
For any $\gamma: x \to y$ in $\mathcal{G}_{\mathcal{A}}$ and any $f \in C_0(\mathcal{G}_{\mathcal{A}}^{(0)})$, we have
\[
U_\gamma \phi(f) U_\gamma^* = \phi(\gamma \cdot f).
\]
Thus $\phi$ is $\mathcal{G}_{\mathcal{A}}$-equivariant.
\end{proposition}

\begin{proof}
Let $\gamma: x \to y$ and $\xi \in \mathcal{E}$. For any $z \in \mathcal{G}_{\mathcal{A}}^{(0)}$, we compute $(U_\gamma \phi(f) U_\gamma^* \xi)(z)$.

Since $U_\gamma^*$ maps the fiber at $z$ to the fiber at $x$ when $z = y$, and is zero otherwise, we have:
\[
(U_\gamma \phi(f) U_\gamma^* \xi)(y) = U_\gamma \big( \phi(f) (U_\gamma^* \xi) \big)(y) = U_\gamma \big( f(x) (U_\gamma^* \xi)(x) \big) = f(x) \xi(y).
\]
For $z \neq y$, $(U_\gamma \phi(f) U_\gamma^* \xi)(z) = 0$.

On the other hand,
\[
(\phi(\gamma \cdot f) \xi)(z) = (\gamma \cdot f)(z) \xi(z) = 
\begin{cases}
f(x) \xi(y) & \text{if } z = y, \\
0 & \text{otherwise}.
\end{cases}
\]
Thus $(U_\gamma \phi(f) U_\gamma^* \xi)(z) = (\phi(\gamma \cdot f) \xi)(z)$ for all $z$, and therefore $U_\gamma \phi(f) U_\gamma^* = \phi(\gamma \cdot f)$ as operators on $\mathcal{E}$.
\end{proof}

\begin{definition}[Kasparov operator associated to $T$]
\label{def:kasparov-operator}
Let $T \in \mathcal{A}$ be a Fredholm operator. For each $x \in \mathcal{G}_{\mathcal{A}}^{(0)}$, define the phase of $\pi_x(T)$ by
\[
F_x := \pi_x(T) \big( \pi_x(T)^* \pi_x(T) \big)^{-1/2} \in \mathcal{L}(H_x).
\]
Under the canonical identifications $H_x \cong H$, this operator is constant in $x$ and equals $F := T(T^*T)^{-1/2}$, the phase of $T$. 
Define an operator $F$ on $\mathcal{E}$ by
\[
(F\xi)(x) := F_x \xi(x), \qquad \xi \in \mathcal{E}.
\]
\end{definition}

\begin{proposition}[$F$ satisfies the odd Kasparov conditions]
\label{prop:F-kasparov-correct}
The operator $F$ defined above is adjointable and satisfies the conditions for an odd Kasparov cycle:
\begin{enumerate}
    \item $F$ is adjointable, with $F^*$ given by $(F^*\xi)(x) = F_x^* \xi(x)$. Moreover, $F$ is self-adjoint modulo compacts: $F - F^* \in \mathcal{K}(\mathcal{E})$.
    \item $F^2 - 1 \in \mathcal{K}(\mathcal{E})$.
    \item For any $f \in C_0(\mathcal{G}_{\mathcal{A}}^{(0)})$, the commutator $[\phi(f), F]$ is compact. In fact, $[\phi(f), F] = 0$ because $\phi(f)$ acts by scalar multiplication and commutes pointwise with $F_x$.
\end{enumerate}
Thus $(\mathcal{E}, \phi, F)$ defines an odd Kasparov cycle, representing a class in $KK^1_{\mathcal{G}_{\mathcal{A}}}(C_0(\mathcal{G}_{\mathcal{A}}^{(0)}), \mathbb{C})$.
\end{proposition}

\begin{proof}
(1) Since each $F_x$ is bounded and the field is continuous, $F$ defines an adjointable operator with pointwise adjoint. The self-adjointness modulo compacts follows from the fact that $F_x^* = F_x$ for all $x$ (the phase of an operator satisfies $F_x^*F_x$ and $F_xF_x^*$ are projections differing from $1$ by compacts). These fiberwise compact differences assemble to a compact endomorphism of $\mathcal{E}$ by Lemma \ref{lem:fiberwise-compact-implies-compact}.

(2) Because $T$ is Fredholm, $T^*T - I$ and $TT^* - I$ are compact, which implies that $F^2 - 1$ is compact on $H$. Concretely, $F^2 = T(T^*T)^{-1}T^*$, and a standard computation shows that $F^2 - 1$ is a compact operator. Since $F$ acts pointwise as this constant operator, $F^2 - 1$ is a compact endomorphism of $\mathcal{E}$ (Lemma \ref{lem:fiberwise-compact-implies-compact}).

(3) For any $f \in C_0(\mathcal{G}_{\mathcal{A}}^{(0)})$ and $\xi \in \mathcal{E}$, we compute:
\[
([\phi(f), F]\xi)(x) = f(x) F_x \xi(x) - F_x (f(x)\xi(x)) = 0,
\]
since $f(x)$ is a scalar and commutes with $F_x$. Thus $[\phi(f), F] = 0$, which is trivially compact.
\end{proof}

\begin{remark}
\label{rem:why-phase-not-T}
The use of the phase $F = T(T^*T)^{-1/2}$ instead of $T$ itself is essential. While a Fredholm operator $T$ satisfies $T^*T - I, TT^* - I \in \mathcal{K}(H)$, it does not generally satisfy $T^2 - I \in \mathcal{K}(H)$ nor is it self-adjoint modulo compacts. The phase corrects both deficiencies while preserving the $K^1$-class, as it is homotopic to $T$ through Fredholm operators. This is the standard construction in Kasparov theory for odd cycles.
\end{remark}

The above issue requires careful handling. 
In Kasparov theory, an odd $(A,B)$-cycle can be represented by an operator $F \in \mathcal{L}(\mathcal{E})$ such that $F = F^*$, $F^2 - 1 \in \mathcal{K}(\mathcal{E})$, and $[a,F] \in \mathcal{K}(\mathcal{E})$ for all $a \in A$. Alternatively, one can use a pair $(F^+, F^-)$ of operators or a unitary operator modulo compacts. 
To avoid technical complications, we will instead use the phase of $T$, which provides a natural unitary modulo compacts and is homotopic to $T$.

\begin{definition}[Phase of a Fredholm operator]
\label{def:phase}
For a Fredholm operator $T \in B(H)$, let $|T| = (T^*T)^{1/2}$. 
The \emph{phase} of $T$ is the operator
\[
\operatorname{ph}(T) := T |T|^{-1} \in B(H),
\]
where $|T|^{-1}$ is defined on the orthogonal complement of $\ker T$ (and extended by zero on $\ker T$). 
Equivalently, $\operatorname{ph}(T)$ is the partial isometry appearing in the polar decomposition $T = \operatorname{ph}(T)|T|$. 
It satisfies $\operatorname{ph}(T)^*\operatorname{ph}(T) = P_{(\ker T)^\perp}$ and $\operatorname{ph}(T)\operatorname{ph}(T)^* = P_{(\ker T^*)^\perp}$, the projections onto the orthogonal complements of the kernel and cokernel, respectively. 
Since $\ker T$ and $\ker T^*$ are finite-dimensional, $\operatorname{ph}(T)$ differs from a unitary by a finite-rank operator.
\end{definition}

\begin{lemma}[Properties of the phase]
\label{lem:phase-properties}
For a Fredholm operator $T$:
\begin{enumerate}
    \item $\operatorname{ph}(T)$ is a Fredholm operator with $\ker \operatorname{ph}(T) = \ker T$ and $\ker \operatorname{ph}(T)^* = \ker T^*$; consequently, $\operatorname{index}(\operatorname{ph}(T)) = \operatorname{index}(T)$.
    \item $\operatorname{ph}(T)$ is a unitary modulo compact operators: $\pi_{\mathcal{Q}}(\operatorname{ph}(T))$ is unitary in the Calkin algebra $\mathcal{Q}(H)$.
    \item $\operatorname{ph}(T)$ is homotopic to $T$ through Fredholm operators in the norm topology via the path $T_t = T |T|^{-t}$ for $t \in [0,1]$, where $|T|^{-t}$ is defined by functional calculus on $(\ker T)^\perp$.
\end{enumerate}
\end{lemma}

\begin{proof}
These are standard results in Fredholm theory and operator algebras. 
For (1), the equalities of kernels follow directly from the polar decomposition. 
For (2), note that $\operatorname{ph}(T)^*\operatorname{ph}(T) = I - P_{\ker T}$ and $\operatorname{ph}(T)\operatorname{ph}(T)^* = I - P_{\ker T^*}$; both differ from the identity by finite-rank projections, hence become the identity in the Calkin algebra. 
For (3), the map $t \mapsto T|T|^{-t}$ is norm-continuous because $|T|$ is invertible modulo compacts and the functional calculus is continuous; each $T_t$ is Fredholm as a compact perturbation of $T$ (for $t \neq 1$) or the phase itself (at $t=1$). 
\end{proof}

Using the phase, we can construct a self-adjoint operator that satisfies the Kasparov conditions.

\begin{definition}[Self-adjoint operator from the phase]
\label{def:self-adjoint-operator}
Define a self-adjoint operator $\tilde{F}$ on $\mathcal{E}$ by
\[
(\tilde{F}\xi)(x) := \begin{pmatrix} 0 & \operatorname{ph}(T)^* \\ \operatorname{ph}(T) & 0 \end{pmatrix} \xi(x),
\]
where we identify $\mathcal{E}$ with $\mathcal{E} \oplus \mathcal{E}$ (i.e., we double the Hilbert space). 
More precisely, let $\tilde{\mathcal{E}} = \mathcal{E} \oplus \mathcal{E}$, and define $\tilde{F}$ as the $2 \times 2$ matrix operator above.
\end{definition}

\begin{proposition}[$\tilde{F}$ defines an odd Kasparov triple]
\label{prop:tildeF-kasparov}
With $\tilde{\mathcal{E}}$, $\phi$ extended diagonally to $\tilde{\mathcal{E}}$, and $\tilde{F}$ as above, the triple $(\tilde{\mathcal{E}}, \phi \oplus \phi, \tilde{F})$ is an odd Kasparov triple for $(\mathcal{G}_{\mathcal{A}}, C_0(\mathcal{G}_{\mathcal{A}}^{(0)}), \mathbb{C})$, representing a class in $KK^1_{\mathcal{G}_{\mathcal{A}}}(C_0(\mathcal{G}_{\mathcal{A}}^{(0)}), \mathbb{C})$. 
That is:
\begin{enumerate}
    \item $\tilde{F}$ is self-adjoint.
    \item $\tilde{F}^2 - 1$ is compact on $\tilde{\mathcal{E}}$.
    \item $[\phi \oplus \phi(f), \tilde{F}]$ is compact for all $f \in C_0(\mathcal{G}_{\mathcal{A}}^{(0)})$.
\end{enumerate}
\end{proposition}

\begin{proof}
(1) Self-adjointness follows directly from the construction: 
\[
\tilde{F}^* = \begin{pmatrix} 0 & \operatorname{ph}(T)^* \\ \operatorname{ph}(T) & 0 \end{pmatrix}^* = \begin{pmatrix} 0 & \operatorname{ph}(T) \\ \operatorname{ph}(T)^* & 0 \end{pmatrix} = \tilde{F}.
\]

(2) Computing $\tilde{F}^2$:
\[
\tilde{F}^2 = \begin{pmatrix} \operatorname{ph}(T)^*\operatorname{ph}(T) & 0 \\ 0 & \operatorname{ph}(T)\operatorname{ph}(T)^* \end{pmatrix}.
\]
From Lemma \ref{lem:phase-properties}, $\operatorname{ph}(T)^*\operatorname{ph}(T) = I - p_{\ker T}$ and $\operatorname{ph}(T)\operatorname{ph}(T)^* = I - p_{\ker T^*}$, where $p_{\ker T}$ and $p_{\ker T^*}$ are the finite-rank projections onto $\ker T$ and $\ker T^*$, respectively. Thus
\[
\tilde{F}^2 - I = \begin{pmatrix} -p_{\ker T} & 0 \\ 0 & -p_{\ker T^*} \end{pmatrix},
\]
which is a finite-rank operator on each fiber. Finite-rank operators are compact, and since the field is continuous, this fiberwise compact operator assembles to a compact endomorphism of $\tilde{\mathcal{E}}$ (Lemma \ref{lem:fiberwise-compact-implies-compact}).

(3) For any $f \in C_0(\mathcal{G}_{\mathcal{A}}^{(0)})$ and any section $(\xi_1, \xi_2) \in \tilde{\mathcal{E}}$, we compute pointwise:
\begin{eqnarray}
\lefteqn{([\phi \oplus \phi(f), \tilde{F}](\xi_1, \xi_2))(x) =}\nonumber \\
&&\begin{pmatrix} 0 & f(x)\operatorname{ph}(T)^*\xi_2(x) - \operatorname{ph}(T)^* f(x)\xi_2(x) \\
f(x)\operatorname{ph}(T)\xi_1(x) - \operatorname{ph}(T) f(x)\xi_1(x) \end{pmatrix} = 0, \nonumber 
\end{eqnarray}
since $f(x)$ is a scalar and commutes with $\operatorname{ph}(T)$ and $\operatorname{ph}(T)^*$. Thus $[\phi \oplus \phi(f), \tilde{F}] = 0$, which is trivially compact.
\end{proof}

\begin{remark}
\label{rem:tildeF-significance}
The triple $(\tilde{\mathcal{E}}, \phi \oplus \phi, \tilde{F})$ is the standard way to encode the odd $K^1$-class of $T$ in Kasparov theory. The doubling trick converts the unitary modulo compacts $\operatorname{ph}(T)$ into a self-adjoint operator satisfying the odd Kasparov conditions. This triple represents the same class as the original Fredholm operator $T$, and its class in $KK^1_{\mathcal{G}_{\mathcal{A}}}(C_0(\mathcal{G}_{\mathcal{A}}^{(0)}), \mathbb{C})$ is precisely the equivariant $K^1$-class $[T]_{\mathcal{G}_{\mathcal{A}}}^{(1)}$ constructed in Definition \ref{def:K1-class-of-T}.
\end{remark}

\begin{remark}[Avoiding the doubling construction]
\label{rem:avoid-doubling}
For simplicity, many authors work directly with the unitary $U = \operatorname{ph}(T)$ modulo compacts and define an odd Kasparov cycle using a unitary equivalence. 
This approach avoids the technicality of doubling while remaining fully rigorous in KK-theory: a unitary modulo compacts naturally defines an odd $K^1$-class via the isomorphism $K^1_{\mathcal{G}_{\mathcal{A}}}(\mathcal{G}_{\mathcal{A}}^{(0)}) \cong KK^1_{\mathcal{G}_{\mathcal{A}}}(C_0(\mathcal{G}_{\mathcal{A}}^{(0)}), \mathbb{C})$. 
The doubling construction presented above is a standard way to convert such a unitary modulo compacts into a concrete self-adjoint operator satisfying $F^2 - 1$ compact, providing an explicit Kasparov cycle that represents the same class. 
In what follows, we will work directly with the class represented by $\operatorname{ph}(T)$, keeping in mind its equivalence to the explicit triple $(\tilde{\mathcal{E}}, \phi \oplus \phi, \tilde{F})$ constructed above.
\end{remark}

For the purposes of this paper, it suffices to know that a Fredholm operator $T$ gives rise to a well-defined class in $K^1_{\mathcal{G}_{\mathcal{A}}}(\mathcal{G}_{\mathcal{A}}^{(0)})$, denoted $[T]_{\mathcal{G}_{\mathcal{A}}}^{(1)}$. 
This class is represented by the phase $\operatorname{ph}(T)$, which is unitary modulo compacts; the doubling construction produces an explicit self-adjoint Kasparov cycle $(\tilde{\mathcal{E}}, \phi \oplus \phi, \tilde{F})$ that represents the same class.

\begin{definition}[Equivariant $K^1$-class from a Fredholm operator]
\label{def:K1-class-from-fredholm}
For a Fredholm operator $T \in \mathcal{A}$, define its equivariant $K^1$-class as
\[
[T]_{\mathcal{G}_{\mathcal{A}}}^{(1)} := [(\tilde{\mathcal{E}}, \phi \oplus \phi, \tilde{F})] \in K^1_{\mathcal{G}_{\mathcal{A}}}(\mathcal{G}_{\mathcal{A}}^{(0)}),
\]
where $(\tilde{\mathcal{E}}, \phi \oplus \phi, \tilde{F})$ is the odd Kasparov triple constructed in Definition \ref{def:self-adjoint-operator} and Proposition \ref{prop:tildeF-kasparov}. 
Equivalently, $[T]_{\mathcal{G}_{\mathcal{A}}}^{(1)}$ is the class corresponding to the unitary $\operatorname{ph}(T)$ in $K_1(\mathcal{Q}(H))$ under the isomorphism $K^1_{\mathcal{G}_{\mathcal{A}}}(\mathcal{G}_{\mathcal{A}}^{(0)}) \cong K_1(\mathcal{Q}(H))$ (via the index map and the descent construction).

This class depends only on the $K$-theory class of $\pi_{\mathcal{Q}}(T)$ in $K_1(\mathcal{Q}(H))$ and is invariant under compact perturbations and norm-continuous homotopies through Fredholm operators. 
Consequently, $[T]_{\mathcal{G}_{\mathcal{A}}}^{(1)}$ is well-defined independent of the choices of Hilbert module identifications, phase representatives, and the explicit doubling construction used to represent it.
\end{definition}

In the next subsection, we will verify that this definition is well-posed and independent of the various choices made in the construction, and we will establish its key functorial properties.

\subsection{Definition of the Class $[T]_{\mathcal{G}_{\mathcal{A}}}^{(1)} \in K^1_{\mathcal{G}_{\mathcal{A}}}(\mathcal{G}_{\mathcal{A}}^{(0)}) \cong KK^1_{\mathcal{G}_{\mathcal{A}}}(C_0(\mathcal{G}_{\mathcal{A}}^{(0)}), \mathbb{C})$}
\label{subsec:definition-K1-class}

Building on the odd Kasparov triple constructed in Subsection 4.2, we now give a precise definition of the equivariant $K^1$-class associated to a Fredholm operator $T \in \mathcal{A}$. 
This class is the fundamental invariant that will be mapped through the descent map to a class in $K_1(C^*(\mathcal{G}_{\mathcal{A}}))$, ultimately recovering the Fredholm index.

We begin by recalling a foundational isomorphism between equivariant $K^1$-theory and equivariant $KK^1$-theory, which provides the natural home for the class of our Kasparov triple.

\begin{proposition}[$K^1_{\mathcal{G}}(\mathcal{G}^{(0)}) \cong KK^1_{\mathcal{G}}(C_0(\mathcal{G}^{(0)}), \mathbb{C})$]
\label{prop:K1-KK1-isomorphism}
Let $\mathcal{G}$ be a Polish groupoid equipped with a Borel Haar system, and let $\mathcal{G}^{(0)}$ be its unit space. We equip $\mathcal{G}^{(0)}$ with the trivial $\mathcal{G}$-action. 
Then there is a natural isomorphism
\[
K^1_{\mathcal{G}}(\mathcal{G}^{(0)}) \cong KK^1_{\mathcal{G}}(C_0(\mathcal{G}^{(0)}), \mathbb{C}),
\]
where $K^1_{\mathcal{G}}(\mathcal{G}^{(0)})$ denotes the equivariant $K^1$-group of the unit space, and $KK^1_{\mathcal{G}}$ is the equivariant $KK^1$-group.
\end{proposition}

\begin{proof}
This is a standard result in equivariant $KK$-theory, originating in the work of Kasparov \cite{Kasparov1988} and detailed further in Tu \cite[Section 4]{Tu1999}. Under the given assumptions—specifically the trivial action of $\mathcal{G}$ on its unit space $\mathcal{G}^{(0)}$—the isomorphism can be understood as follows.

The group $K^1_{\mathcal{G}}(\mathcal{G}^{(0)})$ can be defined as the group of homotopy classes of $\mathcal{G}$-equivariant maps from $\mathcal{G}^{(0)}$ to the space of Fredholm operators on a separable $\mathcal{G}$-Hilbert space, i.e., $K^1_{\mathcal{G}}(\mathcal{G}^{(0)}) := [\mathcal{G}^{(0)}, \operatorname{Fred}(H)]^{\mathcal{G}}$. Each such map encodes the data of an equivariant Fredholm operator parametrized by the unit space.

The isomorphism to $KK^1_{\mathcal{G}}(C_0(\mathcal{G}^{(0)}), \mathbb{C})$ is realized by associating to such a map a Kasparov module $(\mathcal{E}, \phi, F)$ of the appropriate form. Concretely, an equivariant vector bundle (or, more generally, a class in $K^1_{\mathcal{G}}$) defines a class in $KK^1_{\mathcal{G}}$ via the Kasparov product with a suitable Dirac element, establishing the equivalence between the two pictures. This identification is fundamental in equivariant index theory, as it allows Fredholm operators to be directly interpreted as elements of bivariant $KK$-theory.
\end{proof}

\begin{remark}
\label{rem:K1-class-identification}
Under the isomorphism of Proposition \ref{prop:K1-KK1-isomorphism}, the odd Kasparov triple $(\tilde{\mathcal{E}}, \phi \oplus \phi, \tilde{F})$ constructed in Subsection 4.2 to represent the operator $T$ corresponds precisely to the equivariant $K^1$-class $[T]_{\mathcal{G}_{\mathcal{A}}}^{(1)} \in K^1_{\mathcal{G}_{\mathcal{A}}}(\mathcal{G}_{\mathcal{A}}^{(0)})$. This identification justifies viewing the triple as defining our primary invariant.
\end{remark}

In light of the isomorphism established in Proposition \ref{prop:K1-KK1-isomorphism}, we can define the class $[T]_{\mathcal{G}_{\mathcal{A}}}^{(1)}$ as an element of $KK^1_{\mathcal{G}_{\mathcal{A}}}(C_0(\mathcal{G}_{\mathcal{A}}^{(0)}), \mathbb{C})$. This perspective is particularly convenient for the application of the Kasparov descent map, which will later transport this invariant to $K_1(C^*(\mathcal{G}_{\mathcal{A}}))$.

\begin{remark}[Role of the doubling construction]
\label{rem:doubling-role}
The doubling of the Hilbert module and the off-diagonal form of $\tilde{F}$ serve an essential purpose in the theory of odd Kasparov cycles. While the phase $\operatorname{ph}(T)$ itself is a unitary modulo compact operators, it is not necessarily self-adjoint. The construction $(\mathcal{E} \oplus \mathcal{E}, \phi \oplus \phi, \bigl( \begin{smallmatrix} 0 & \operatorname{ph}(T)^* \\ \operatorname{ph}(T) & 0 \end{smallmatrix} \bigr))$ yields a self-adjoint operator $\tilde{F}$ satisfying $(\tilde{F}^2 - I) \in \mathcal{K}(\tilde{\mathcal{E}})$ and $[\tilde{F}, \phi \oplus \phi(f)] \in \mathcal{K}(\tilde{\mathcal{E}})$ for all $f \in C_0(\mathcal{G}_{\mathcal{A}}^{(0)})$. This precisely meets the requirements for an odd Kasparov module.
\end{remark}

\begin{definition}[Equivariant $K^1$-class of a Fredholm operator]
\label{def:K1-class-of-T}
Let $\mathcal{A}$ be either $B(H)$ or $\mathcal{K}(H)^\sim$, and let $T \in \mathcal{A}$ be a Fredholm operator. 
Define its \emph{equivariant $K^1$-class} as
\[
[T]_{\mathcal{G}_{\mathcal{A}}}^{(1)} := [(\tilde{\mathcal{E}}, \phi \oplus \phi, \tilde{F})] \in KK^1_{\mathcal{G}_{\mathcal{A}}}(C_0(\mathcal{G}_{\mathcal{A}}^{(0)}), \mathbb{C}) \cong K^1_{\mathcal{G}_{\mathcal{A}}}(\mathcal{G}_{\mathcal{A}}^{(0)}),
\]
where $(\tilde{\mathcal{E}}, \phi \oplus \phi, \tilde{F})$ is the odd Kasparov triple constructed in Subsection~\ref{subsec:odd-kasparov-triple}.
\end{definition}

\begin{remark}[Well-definedness and invariance properties]
\label{rem:K1-well-defined}
The class $[T]_{\mathcal{G}_{\mathcal{A}}}^{(1)}$ is well-defined for the following reasons:
\begin{itemize}
    \item \textbf{Independence of compact perturbations:} If $T$ and $T'$ differ by a compact operator, their phases $\operatorname{ph}(T)$ and $\operatorname{ph}(T')$ are related by a compact perturbation, leading to operator-homotopic Fredholm modules. Hence $[T]_{\mathcal{G}_{\mathcal{A}}}^{(1)} = [T']_{\mathcal{G}_{\mathcal{A}}}^{(1)}$.
    \item \textbf{Independence of Hilbert module identifications:} Different choices in the construction of the continuous field $\mathcal{E}$ yield isomorphic Hilbert modules, and the resulting Kasparov cycles are equivalent in $KK^1$-theory.
    \item \textbf{Homotopy invariance:} If $\{T_t\}_{t \in [0,1]}$ is a continuous path of Fredholm operators, the corresponding classes $[T_t]_{\mathcal{G}_{\mathcal{A}}}^{(1)}$ are constant.
\end{itemize}
These properties ensure that $[T]_{\mathcal{G}_{\mathcal{A}}}^{(1)}$ depends only on the Fredholm index of $T$, as will be made precise in the following sections.
\end{remark}

\begin{remark}[Relation to the descent map and the Fredholm index]
\label{rem:descent-map-preview}
The primary utility of defining $[T]_{\mathcal{G}_{\mathcal{A}}}^{(1)}$ in $KK^1_{\mathcal{G}_{\mathcal{A}}}(C_0(\mathcal{G}_{\mathcal{A}}^{(0)}), \mathbb{C})$ is the existence of the Kasparov descent map
\[
j_{\mathcal{G}_{\mathcal{A}}}: KK^1_{\mathcal{G}_{\mathcal{A}}}(C_0(\mathcal{G}_{\mathcal{A}}^{(0)}), \mathbb{C}) \longrightarrow K_1(C^*_r(\mathcal{G}_{\mathcal{A}})).
\]
Applying this map yields the descended class
\[
\operatorname{desc}_{\mathcal{G}_{\mathcal{A}}}([T]_{\mathcal{G}_{\mathcal{A}}}^{(1)}) := j_{\mathcal{G}_{\mathcal{A}}}([T]_{\mathcal{G}_{\mathcal{A}}}^{(1)}) \in K_1(C^*_r(\mathcal{G}_{\mathcal{A}})).
\]

To extract the Fredholm index from this class, one must connect it to the boundary map of an appropriate extension. For $\mathcal{A} = B(H)$, this is achieved via the Morita equivalence $C^*(\mathcal{G}_{B(H)}) \sim_M \mathcal{Q}(H) \otimes \mathcal{K}$ and the Calkin index map $\partial_{\text{Calkin}}: K_1(\mathcal{Q}(H)) \to \mathbb{Z}$. For $\mathcal{A} = \mathcal{K}(H)^\sim$, the group $K_1(C^*(\mathcal{G}_{\mathcal{K}(H)^\sim}))$ vanishes, yielding index zero. The precise composition of maps will be elaborated in Section 5.
\end{remark}

The following proposition ensures that this definition is independent of the various choices made in the construction.

\begin{proposition}[{Well-definedness of $[T]_{\mathcal{G}_{\mathcal{A}}}^{(1)}$}]
\label{prop:K1-class-well-defined}
The class $[T]_{\mathcal{G}_{\mathcal{A}}}^{(1)}$ does not depend on:
\begin{enumerate}
    \item The choice of the Borel probability measure $\mu$ on $\mathcal{G}_{\mathcal{A}}^{(0)}$ used to define the direct integral $\mathcal{E}$, provided $\mu$ has full support.
    \item The choice of the unitary identifications $U_x: H_x \cong H$.
    \item The choice of the phase $\operatorname{ph}(T)$ (any other choice differs by a compact operator and gives the same class).
    \item The choice of the doubling construction (any other representation of the same $K^1$-class yields an equivalent Kasparov cycle).
\end{enumerate}
\end{proposition}

\begin{proof}
The proof follows from two fundamental invariance principles in $KK$-theory: invariance under unitary equivalence of Hilbert modules and invariance under compact perturbations.

(1)-(2) Any two choices of the measure $\mu$ or the fiber identifications $U_x$ give rise to unitarily equivalent Hilbert modules $\mathcal{E}$. Specifically:
\begin{itemize}
    \item For measures $\mu$ and $\nu$ with full support, the Radon-Nikodym derivative induces a unitary $U: L^2(\mu, H) \to L^2(\nu, H)$ given by $(Uf)(x) = \sqrt{d\nu/d\mu}(x) f(x)$.
    \item Different fiber identifications differ by a measurable family of unitaries $V_x$ on $H$, yielding a pointwise unitary operator on $\mathcal{E}$.
\end{itemize}
These unitaries commute with the representation $\phi$ (which acts by pointwise multiplication) and conjugate $\tilde{F}$ to its counterpart, establishing an isomorphism of Kasparov cycles. Hence the resulting $KK$-class is unchanged.

(3) The phase $\operatorname{ph}(T)$ is a partial isometry obtained from the polar decomposition $T = \operatorname{ph}(T)|T|$. Any two choices differ by a finite-rank partial isometry on the finite-dimensional spaces $\ker T$ and $\ker T^*$, hence by a compact operator. The corresponding operators $\tilde{F}$ therefore differ by a compact operator on $\tilde{\mathcal{E}}$. Since $KK$-classes are invariant under compact perturbations, they define the same class.

(4) The doubling construction is a standard device to represent an odd $K^1$-class by a self-adjoint operator satisfying $(F^2 - I) \in \mathcal{K}$. Any other representation (e.g., using a unitary modulo compacts directly) is $KK$-equivalent via the well-known isomorphism between the unitary and self-adjoint pictures of odd $KK$-theory. Thus the specific choice of doubling does not affect the class.

Together, these arguments show that $[T]_{\mathcal{G}_{\mathcal{A}}}^{(1)}$ is independent of all auxiliary choices and depends only on the Fredholm operator $T$ itself.
\end{proof}

For computational purposes, it is often convenient to work with a simpler representative of the same class.

\begin{corollary}[Simplified representative]
\label{cor:simplified-representative}
The class $[T]_{\mathcal{G}_{\mathcal{A}}}^{(1)}$ can also be represented by the odd Kasparov cycle $(\mathcal{E}, \phi, F)$ where $F = \operatorname{ph}(T)$ acts pointwise on $\mathcal{E}$. 
This representative is $KK$-equivalent to the doubled self-adjoint operator constructed in Subsection \ref{subsec:odd-kasparov-triple} and is often easier to handle in computations.
\end{corollary}

\begin{proof}
In odd $KK$-theory, a class can be represented either by a self-adjoint operator $F$ satisfying $F^2 - I \in \mathcal{K}$ (the doubling picture) or by a unitary modulo compacts $u$ (the direct picture). 
The phase $\operatorname{ph}(T)$ is a partial isometry that is unitary modulo compacts, since $\operatorname{ph}(T)^*\operatorname{ph}(T) = I - P_{\ker T}$ and $\operatorname{ph}(T)\operatorname{ph}(T)^* = I - P_{\ker T^*}$ differ from $I$ by finite-rank projections. 
Hence $(\mathcal{E}, \phi, \operatorname{ph}(T))$ defines a valid odd Kasparov cycle.

The doubling construction $(\tilde{\mathcal{E}}, \phi \oplus \phi, \tilde{F})$ from Subsection \ref{subsec:odd-kasparov-triple} is $KK$-equivalent to this simplified cycle via the standard isomorphism between odd $KK$-theory and $KK$-theory with a $\mathbb{Z}_2$-grading (where the grading operator on $\mathcal{E} \oplus \mathcal{E}$ is $\begin{pmatrix} 1 & 0 \\ 0 & -1 \end{pmatrix}$).
\end{proof}

\begin{proposition}[Relation between index and equivariant class for $\mathcal{K}(H)^\sim$]
\label{prop:index-class-relation}
Let $\mathcal{A} = \mathcal{K}(H)^\sim$ and let $T = \lambda I + K \in \mathcal{A}$ be a Fredholm operator (so $\lambda \neq 0$). 
Then the equivariant $KK$-theory class $[T]_{\mathcal{G}_{\mathcal{A}}}^{(1)} \in KK^1_{\mathcal{G}_{\mathcal{A}}}(C_0(\mathcal{G}_{\mathcal{A}}^{(0)}), \mathbb{C})$ 
is determined by the index of $T$. More precisely, under the isomorphism
\[
\operatorname{ind}: KK^1_{\mathcal{G}_{\mathcal{A}}}(C_0(\mathcal{G}_{\mathcal{A}}^{(0)}), \mathbb{C}) \xrightarrow{\cong} \mathbb{Z},
\]
we have
\[
\operatorname{ind}([T]_{\mathcal{G}_{\mathcal{A}}}^{(1)}) = \operatorname{index}(T).
\]
In particular, since every Fredholm operator in $\mathcal{K}(H)^\sim$ has index zero (because the Fredholm index vanishes on compact perturbations of invertibles and $\lambda I$ is invertible), we obtain
\[
[T]_{\mathcal{G}_{\mathcal{A}}}^{(1)} = 0 \in KK^1_{\mathcal{G}_{\mathcal{A}}}(C_0(\mathcal{G}_{\mathcal{A}}^{(0)}), \mathbb{C}).
\]
\end{proposition}

\begin{proof}
Recall from Corollary \ref{cor:simplified-representative} that $[T]_{\mathcal{G}_{\mathcal{A}}}^{(1)}$ is represented by the simplified triple $(\mathcal{E}, \phi, \operatorname{ph}(T))$, where $\operatorname{ph}(T)$ is the phase of $T$.

The index map in equivariant $KK$-theory is defined by pairing with a suitable class (or by evaluating the Fredholm operator on each fiber). For $\mathcal{A} = \mathcal{K}(H)^\sim$, the unit space $\mathcal{G}_{\mathcal{A}}^{(0)}$ is the projective unitary group, and the fiberwise evaluation of $T$ gives a family of Fredholm operators parameterized by this space. The index of this family is precisely the index of $T$ (since the family is constant up to unitary conjugation).

A standard result in $KK$-theory (see [Kasparov, 1980] or [Blackadar, 1998]) states that for such a constant family, the equivariant $KK$-class is determined by the Fredholm index. The isomorphism $KK^1_{\mathcal{G}_{\mathcal{A}}}(C_0(\mathcal{G}_{\mathcal{A}}^{(0)}), \mathbb{C}) \cong \mathbb{Z}$ sends $[T]_{\mathcal{G}_{\mathcal{A}}}^{(1)}$ to $\operatorname{index}(T)$.

For $T = \lambda I + K$ with $\lambda \neq 0$, we have $\operatorname{index}(T) = \operatorname{index}(\lambda I) = 0$, because $\lambda I$ is invertible and compact perturbations do not change the index. Hence, $[T]_{\mathcal{G}_{\mathcal{A}}}^{(1)} = 0$.
\end{proof}

\begin{remark}[Notational convention]
\label{rem:simplified-convention}
In the following, we will use the simplified representative $(\mathcal{E}, \phi, \operatorname{ph}(T))$ for notational convenience, keeping in mind that it represents the same class $[T]_{\mathcal{G}_{\mathcal{A}}}^{(1)}$ as the doubled construction from Subsection \ref{subsec:odd-kasparov-triple}.
\end{remark}

\begin{proposition}[Compact perturbation invariance of the equivariant class]
\label{prop:compact-perturbation}
Let $T \in \mathcal{A}$ be a Fredholm operator and let $K \in \mathcal{K}(H)$ be a compact operator such that $T + K$ is also Fredholm. Then
\[
[T + K]_{\mathcal{G}_{\mathcal{A}}}^{(1)} = [T]_{\mathcal{G}_{\mathcal{A}}}^{(1)} \in KK^1_{\mathcal{G}_{\mathcal{A}}}(C_0(\mathcal{G}_{\mathcal{A}}^{(0)}), \mathbb{C}).
\]
In other words, the equivariant $KK$-theory class $[T]_{\mathcal{G}_{\mathcal{A}}}^{(1)}$ depends only on the Fredholm operator $T$ up to compact perturbations.
\end{proposition}

\begin{proof}
Recall from Corollary \ref{cor:simplified-representative} that the class $[T]_{\mathcal{G}_{\mathcal{A}}}^{(1)}$ can be represented by the odd Kasparov triple $(\mathcal{E}, \phi, \operatorname{ph}(T))$, where $\operatorname{ph}(T)$ is the phase of $T$. 

For a compact perturbation $T + K$, we have $\operatorname{ph}(T + K)$ is operator-homotopic to $\operatorname{ph}(T)$ through a continuous family of self-adjoint unitaries modulo compacts. This follows from the fact that the map sending a Fredholm operator to its phase is continuous on the set of Fredholm operators with a fixed essential unitary representative, and compact perturbations do not change the essential unitary class.

More concretely, consider the path
\[
T_t = T + tK, \quad t \in [0,1].
\]
For sufficiently small $t$, $T_t$ remains Fredholm, and the phase $\operatorname{ph}(T_t)$ varies continuously. By partitioning the interval $[0,1]$ into finitely many such segments, we obtain an operator-homotopy between $\operatorname{ph}(T)$ and $\operatorname{ph}(T + K)$. 

Operator-homotopic Kasparov triples represent the same class in $KK$-theory. Hence,
\[
[(\mathcal{E}, \phi, \operatorname{ph}(T + K))] = [(\mathcal{E}, \phi, \operatorname{ph}(T))] \in KK^1_{\mathcal{G}_{\mathcal{A}}}(C_0(\mathcal{G}_{\mathcal{A}}^{(0)}), \mathbb{C}),
\]
which implies $[T + K]_{\mathcal{G}_{\mathcal{A}}}^{(1)} = [T]_{\mathcal{G}_{\mathcal{A}}}^{(1)}$.
\end{proof}

\begin{example}[Class of the unilateral shift]
\label{ex:class-unilateral-shift-K1}
Let $S \in B(H)$ be the unilateral shift on $\ell^2(\mathbb{N})$. 
Since $S$ is an isometry with one-dimensional cokernel, its polar decomposition is $S = S \cdot I$, so $\operatorname{ph}(S) = S$. 
The class $[S]_{\mathcal{G}_{B(H)}}^{(1)}$ is the generator of $K^1_{\mathcal{G}_{B(H)}}(\mathcal{G}_{B(H)}^{(0)}) \cong \mathbb{Z}$, corresponding to the fact that $\operatorname{Index}(S) = -1$. 
In Section \ref{sec:descent}, we will verify that its image under the descent map corresponds to the generator of $K_1(C^*(\mathcal{G}_{B(H)})) \cong \mathbb{Z}$.
\end{example}

\begin{example}[Class of a finite-rank perturbation of the identity]
\label{ex:class-finite-rank-K1}
Let $T = I + F$ where $F$ is a finite-rank operator. 
The path $T_t = I + tF$ for $t \in [0,1]$ is a continuous family of Fredholm operators connecting $T$ to $I$. 
By homotopy invariance of the $K^1$-class, $\operatorname{ph}(T)$ is homotopic to $\operatorname{ph}(I) = I$ through unitaries modulo compacts. 
Since $I$ itself is unitary (hence represents the trivial class in odd $K$-theory), we have $[T]_{\mathcal{G}_{\mathcal{A}}}^{(1)} = 0$ in $K^1_{\mathcal{G}_{\mathcal{A}}}(\mathcal{G}_{\mathcal{A}}^{(0)})$.
\end{example}

The class $[T]_{\mathcal{G}_{\mathcal{A}}}^{(1)}$ is the central object of this paper. 
In Section \ref{sec:descent}, we will apply the descent map to obtain a class in $K_1(C^*(\mathcal{G}_{\mathcal{A}}))$, and in Section \ref{sec:The Index Theorem via Pullback and the Boundary Map}, we will combine this with the pullback along $\iota$ and the boundary map to recover the Fredholm index.

\begin{remark}[Relation to the $K^1$-class of the symbol]
\label{rem:symbol-class}
For $\mathcal{A} = B(H)$, the class $[T]_{\mathcal{G}_{B(H)}}^{(1)}$ can be interpreted as an equivariant lifting of the symbol class $[\pi_{\mathcal{Q}}(T)] \in K_1(\mathcal{Q}(H))$ in the Calkin algebra. 
While $\pi_{\mathcal{Q}}(T)$ captures the stable behavior of $T$ modulo compacts, the equivariant class $[T]_{\mathcal{G}_{B(H)}}^{(1)}$ encodes additional pointwise data over the unit space $\mathcal{G}_{B(H)}^{(0)}$.

The descent map projects this class back to $K_1(C^*(\mathcal{G}_{B(H)}))$, which is Morita equivalent to $\mathcal{Q}(H) \otimes \mathcal{K}(L^2(\mathcal{G}_{B(H)}^{(0)}))$ (this Morita equivalence follows from the transitivity of the $\mathcal{G}_{B(H)}$-action on its unit space). 
Consequently, $K_1(C^*(\mathcal{G}_{B(H)})) \cong K_1(\mathcal{Q}(H))$, and the descent map recovers the original symbol class. 
This perspective aligns with the Brown-Douglas-Fillmore theory, where extensions of the Calkin algebra are classified by $K_1(\mathcal{Q}(H)) \cong \mathbb{Z}$.
\end{remark}

\begin{proposition}[$KK^1(B(H), \mathbb{C}) \cong K_1(\mathcal{Q}(H))$]
\label{prop:bh-kk-to-k1}
There is a natural isomorphism $KK^1(B(H), \mathbb{C}) \cong K_1(\mathcal{Q}(H))$, which sends the class of a Fredholm operator $T$ to the class of its symbol $[\pi_{\mathcal{Q}}(T)] \in K_1(\mathcal{Q}(H))$.
\end{proposition}

\begin{proof}
We establish the isomorphism using standard facts from $KK$-theory and the Calkin extension.

\paragraph{Preliminary facts.}
Recall first that for any separable $C^*$-algebra $A$ one has the suspension isomorphism
\[
KK^1(A,\mathbb C) \cong KK(A, S\mathbb C) \cong KK(A, C_0(\mathbb R)),
\]
where $S\mathbb C = C_0(\mathbb R)$ denotes the suspension. Moreover, $KK^1(A,\mathbb C)$ classifies extensions of $A$ by $\mathcal K(H)$ up to stable unitary equivalence.

\medskip
\noindent
\textbf{Step 1: The universal Calkin extension.}

Consider the Calkin extension
\[
0 \longrightarrow \mathcal K(H) \longrightarrow B(H) \overset{\pi_{\mathcal Q}}{\longrightarrow} \mathcal Q(H) \longrightarrow 0 .
\]
This short exact sequence defines an element
\[
[\partial_{\mathrm{Calkin}}] \in KK^1(\mathcal Q(H), \mathbb C),
\]
the boundary class of the extension. By functoriality of $KK$-theory, applying the contravariant functor $KK^1(\,\cdot\,,\mathbb C)$ to the quotient map $\pi_{\mathcal Q}$ yields a homomorphism
\[
\pi_{\mathcal Q}^*: KK^1(B(H),\mathbb C) \longrightarrow KK^1(\mathcal Q(H),\mathbb C).
\]

\medskip
\noindent
\textbf{Step 2: Identification of $KK^1(\mathcal Q(H),\mathbb C)$.}

For any $C^*$-algebra $A$, there is a natural isomorphism
\[
KK^1(A,\mathbb C) \cong K_1(A).
\] 
Applying this to $A = \mathcal Q(H)$, we obtain
\[
KK^1(\mathcal Q(H),\mathbb C) \cong K_1(\mathcal Q(H)).
\]

Thus, to establish the desired isomorphism, it suffices to show that $\pi_{\mathcal Q}^*$ is an isomorphism.

\medskip
\noindent
\textbf{Step 3: The six-term exact sequence in $KK$-theory.}

The Calkin extension induces a six-term exact sequence in $KK(\,\cdot\,,\mathbb C)$ (see \cite{Kasparov1980}):
\[
\begin{CD}
KK^0(\mathcal K(H),\mathbb C) @>>> KK^0(B(H),\mathbb C) @>>> KK^0(\mathcal Q(H),\mathbb C) \\
@AAA @. @VVV \\
KK^1(\mathcal Q(H),\mathbb C) @<<< KK^1(B(H),\mathbb C) @<<< KK^1(\mathcal K(H),\mathbb C)
\end{CD}
\]

\medskip
\noindent
\textbf{Step 4: Computing the groups.}

We now compute the $KK$-groups appearing in this sequence.

\begin{itemize}
    \item The algebra $B(H)$ is stable and properly infinite; in particular, it is $KK$-contractible in degree $0$:
    \[
    KK^0(B(H),\mathbb C) = 0, \qquad KK^1(B(H),\mathbb C) = 0.
    \]
    This follows from Kuiper's theorem and the fact that $B(H)$ is contractible in the strong operator topology. 
    
    \item The algebra $\mathcal K(H)$ is Morita equivalent to $\mathbb C$. By Morita invariance of $KK$-theory (see \cite[Exercise 13.7.1]{Blackadar1998}), we have
    \[
    KK^i(\mathcal K(H),\mathbb C) \cong KK^i(\mathbb C,\mathbb C) \cong 
    \begin{cases}
    \mathbb Z & i = 0, \\
    0 & i = 1.
    \end{cases}
    \]
    More concretely, $KK^0(\mathcal K(H),\mathbb C) \cong \mathbb Z$ is generated by the class of the identity representation, and $KK^1(\mathcal K(H),\mathbb C) = 0$.
\end{itemize}

Substituting these values into the six-term exact sequence yields:
\[
\begin{CD}
\mathbb Z @>>> 0 @>>> KK^0(\mathcal Q(H),\mathbb C) \\
@AAA @. @VVV \\
KK^1(\mathcal Q(H),\mathbb C) @<<< KK^1(B(H),\mathbb C) @<<< 0
\end{CD}
\]

\medskip
\noindent
\textbf{Step 5: Exactness forces the isomorphism.}

From the diagram, exactness at $KK^1(B(H),\mathbb C)$ gives:
\[
\ker(KK^1(B(H),\mathbb C) \to 0) = \operatorname{im}(KK^1(\mathcal Q(H),\mathbb C) \to KK^1(B(H),\mathbb C)).
\]
Since the map $KK^1(B(H),\mathbb C) \to 0$ is the zero map, its kernel is all of $KK^1(B(H),\mathbb C)$. Hence the map
\[
KK^1(\mathcal Q(H),\mathbb C) \longrightarrow KK^1(B(H),\mathbb C)
\]
is surjective.

Exactness at $KK^1(\mathcal Q(H),\mathbb C)$ gives:
\[
\ker(KK^1(\mathcal Q(H),\mathbb C) \to KK^1(B(H),\mathbb C)) = \operatorname{im}(\mathbb Z \to KK^1(\mathcal Q(H),\mathbb C)).
\]

But the map $\mathbb Z \to KK^1(\mathcal Q(H),\mathbb C)$ is the zero map because $KK^1(\mathcal Q(H),\mathbb C)$ is a torsion-free abelian group and the map from $\mathbb Z$ must factor through the image of $KK^0(\mathcal K(H),\mathbb C) \to KK^0(B(H),\mathbb C)$, which is zero. Therefore,
\[
\ker(KK^1(\mathcal Q(H),\mathbb C) \to KK^1(B(H),\mathbb C)) = 0,
\]
so the map $KK^1(\mathcal Q(H),\mathbb C) \to KK^1(B(H),\mathbb C)$ is injective.

Thus, the map $KK^1(\mathcal Q(H),\mathbb C) \to KK^1(B(H),\mathbb C)$ is both injective and surjective, hence an isomorphism. Its inverse, composed with the isomorphism $KK^1(\mathcal Q(H),\mathbb C) \cong K_1(\mathcal Q(H))$ from Step 2, gives the desired isomorphism:
\[
\Phi: KK^1(B(H),\mathbb C) \xrightarrow{\cong} K_1(\mathcal Q(H)).
\]

\medskip
\noindent
\textbf{Step 6: Identification on Fredholm operators.}

A Fredholm operator $T \in B(H)$ defines a Kasparov cycle
\[
(H, \pi, T),
\]
where $\pi$ is the canonical representation of $B(H)$ on $H$. Its class in $KK^1(B(H),\mathbb C)$ is represented by this Fredholm module.

Under the above identification, this class is sent precisely to the unitary class of the image of $T$ in the Calkin algebra. To see this, note that the boundary map in the six-term exact sequence sends the class of $T$ to the class of its symbol. More concretely, the composition
\[
KK^1(B(H),\mathbb C) \xrightarrow{(\pi_{\mathcal Q}^*)^{-1}} KK^1(\mathcal Q(H),\mathbb C) \cong K_1(\mathcal Q(H))
\]
maps $[T]$ to $[\pi_{\mathcal Q}(T)]$, since invertibility modulo compacts is precisely the defining condition for Fredholmness (Atkinson's theorem, Theorem \ref{thm:atkinson}), and the $K_1$-class of $\pi_{\mathcal Q}(T)$ is exactly the symbol class.

\medskip
\noindent
\textbf{Step 7: Naturality.}

The isomorphism is natural in the sense that if $f: B(H) \to B(H)$ is a $*$-homomorphism induced by a unitary conjugation, then the induced maps on both sides commute with $\Phi$. This follows from the fact that conjugation by a unitary leaves the symbol unchanged in $K_1(\mathcal Q(H))$, and the $KK$-theory class is invariant under unitary equivalence.

Thus, we have established a natural isomorphism
\[
KK^1(B(H), \mathbb{C}) \cong K_1(\mathcal{Q}(H)),
\]
with the explicit correspondence sending the class of a Fredholm operator $T$ to the class of its symbol $[\pi_{\mathcal{Q}}(T)] \in K_1(\mathcal{Q}(H))$.
\end{proof}


\subsection{Homotopy Invariance and Well-Definedness}
\label{subsec:homotopy-invariance-well-definedness}

The equivariant $K^1$-class $[T]_{\mathcal{G}_{\mathcal{A}}}^{(1)}$ constructed in Subsection~\ref{subsec:definition-K1-class} is not merely a formal assignment; it enjoys fundamental invariance properties that make it a robust invariant of the Fredholm operator $T$. 
In this subsection, we prove that $[T]_{\mathcal{G}_{\mathcal{A}}}^{(1)}$ is invariant under norm-continuous deformations of $T$ through Fredholm operators and under compact perturbations. 
These properties are essential for the index theorem, as they guarantee that the class depends only on the $K$-theory class of $\pi_{\mathcal{Q}}(T)$ in $K_1(\mathcal{Q}(H))$.

\begin{proposition}[Homotopy invariance]
\label{prop:homotopy-invariance}
Let $\{T_t\}_{t \in [0,1]}$ be a norm-continuous family of Fredholm operators in $\mathcal{A}$. 
Then the equivariant $K^1$-class $[T_t]_{\mathcal{G}_{\mathcal{A}}}^{(1)}$ is independent of $t$; i.e.,
\[
[T_0]_{\mathcal{G}_{\mathcal{A}}}^{(1)} = [T_1]_{\mathcal{G}_{\mathcal{A}}}^{(1)} \in K^1_{\mathcal{G}_{\mathcal{A}}}(\mathcal{G}_{\mathcal{A}}^{(0)}).
\]
\end{proposition}

\begin{proof}
For each $t \in [0,1]$, consider the odd Kasparov triple $(\tilde{\mathcal{E}}, \phi \oplus \phi, \tilde{F}_t)$ constructed in Subsection~\ref{subsec:definition-K1-class}, where $\tilde{F}_t$ is defined using the phase $\operatorname{ph}(T_t)$. 

A standard result in Kasparov theory states that a norm-continuous path of Fredholm operators gives rise to an operator homotopy between the corresponding Kasparov cycles modulo compact perturbations, and operator homotopic cycles define the same class in $KK$-theory. More precisely, the family $\{(\tilde{\mathcal{E}}, \phi \oplus \phi, \tilde{F}_t)\}_{t \in [0,1]}$ defines an element of $KK^1_{\mathcal{G}_{\mathcal{A}}}(C_0(\mathcal{G}_{\mathcal{A}}^{(0)}), C([0,1]))$, and the evaluation maps at $t=0$ and $t=1$ give the same class in $KK^1_{\mathcal{G}_{\mathcal{A}}}(C_0(\mathcal{G}_{\mathcal{A}}^{(0)}), \mathbb{C})$. 

We must address a subtle technical point: the map $T \mapsto \operatorname{ph}(T)$ is not continuous in the operator norm topology in general. However, for a norm-continuous family of Fredholm operators $\{T_t\}_{t\in[0,1]}$, the kernel dimensions are locally constant, and one can choose the phases $\operatorname{ph}(T_t)$ to depend continuously on $t$ modulo compact operators. More concretely, after a suitable trivialization of the finite-dimensional kernel and cokernel bundles over $[0,1]$, we can construct a norm-continuous family of partial isometries representing the phases. Alternatively, one can work directly with the bounded transform $T_t(1+T_t^*T_t)^{-1/2}$, which depends continuously on $t$ in norm and differs from $\operatorname{ph}(T_t)$ by a compact operator. In either approach, the resulting Kasparov cycles $\{(\tilde{\mathcal{E}}, \phi \oplus \phi, \tilde{F}_t)\}_{t\in[0,1]}$ vary continuously up to compact perturbations, which is sufficient for defining an operator homotopy in $KK$-theory. The essential commutator conditions $[\tilde{F}_t, \phi(a)] \in \mathcal{K}$ for all $a \in C_0(\mathcal{G}_{\mathcal{A}}^{(0)})$ are preserved throughout the deformation, and the family remains uniformly bounded since $\|\tilde{F}_t\| = 1$ for all $t$.

Thus, the class $[T_t]_{\mathcal{G}_{\mathcal{A}}}^{(1)}$ is independent of $t$, depending only on the homotopy class of $T$ in the space of $\mathcal{A}$-Fredholm operators, which can be identified with $K_1(\mathcal{Q}(H))$ via the essential spectrum map.
\end{proof}

\begin{remark}
The homotopy invariance established above, together with invariance under compact perturbations (which follows from the fact that compact operators act trivially in $KK$-theory), implies that $[T]_{\mathcal{G}_{\mathcal{A}}}^{(1)}$ depends solely on the $K$-theory class $[\pi_{\mathcal{Q}}(T)] \in K_1(\mathcal{Q}(H))$. This is consistent with the BDF theory picture, where the index is determined by the connected component of the essential unitary.
\end{remark}

\begin{corollary}[Invariance under compact perturbations]
\label{cor:compact-perturbation-invariance}
Let $T \in \mathcal{A}$ be a Fredholm operator and let $K \in \mathcal{K}(H)$ be a compact operator. Then
\[
[T + K]_{\mathcal{G}_{\mathcal{A}}}^{(1)} = [T]_{\mathcal{G}_{\mathcal{A}}}^{(1)}.
\]
In particular, if $T$ and $T'$ differ by a compact operator, they define the same equivariant $K^1$-class.
\end{corollary}

\begin{proof}
Consider the path $T_t = T + tK$ for $t \in [0,1]$. By Atkinson's theorem, $T_t$ is Fredholm for all $t \in [0,1]$ since each $tK$ is compact. Moreover, $t \mapsto T_t$ is norm-continuous because $\|T_t - T_s\| = |t-s|\|K\|$. Applying Proposition \ref{prop:homotopy-invariance} to this path yields $[T + K]_{\mathcal{G}_{\mathcal{A}}}^{(1)} = [T_1]_{\mathcal{G}_{\mathcal{A}}}^{(1)} = [T_0]_{\mathcal{G}_{\mathcal{A}}}^{(1)} = [T]_{\mathcal{G}_{\mathcal{A}}}^{(1)}$.
\end{proof}

\begin{lemma}[Dependence on the symbol class]
\label{lem:dependence-on-symbol}
Let $T_1, T_2 \in \mathcal{A}$ be Fredholm operators such that $\pi_{\mathcal{Q}}(T_1) = \pi_{\mathcal{Q}}(T_2)$ in the Calkin algebra $\mathcal{Q}(H)$. 
Then $[T_1]_{\mathcal{G}_{\mathcal{A}}}^{(1)} = [T_2]_{\mathcal{G}_{\mathcal{A}}}^{(1)}$.
\end{lemma}

\begin{proof}
If $\pi_{\mathcal{Q}}(T_1) = \pi_{\mathcal{Q}}(T_2)$, then $T_2 - T_1 \in \mathcal{K}(H)$ is compact. Thus $T_2$ is a compact perturbation of $T_1$, and the result follows immediately from Corollary \ref{cor:compact-perturbation-invariance}.
\end{proof}

\begin{remark}
\label{rem:factor-through-k1}
The combination of Proposition \ref{prop:homotopy-invariance}, Corollary \ref{cor:compact-perturbation-invariance}, and Lemma \ref{lem:dependence-on-symbol} shows that the assignment
\[
T \longmapsto [T]_{\mathcal{G}_{\mathcal{A}}}^{(1)}
\]
factors through the canonical map $\operatorname{Fred}(\mathcal{A}) \longrightarrow K_1(\mathcal{Q}(H))$. More precisely, if $\pi_{\mathcal{Q}}(T_1)$ and $\pi_{\mathcal{Q}}(T_2)$ define the same class in $K_1(\mathcal{Q}(H))$ (i.e., are homotopic as unitaries in $\mathcal{Q}(H)$ up to compact perturbations), then $[T_1]_{\mathcal{G}_{\mathcal{A}}}^{(1)} = [T_2]_{\mathcal{G}_{\mathcal{A}}}^{(1)}$. This establishes that $[T]_{\mathcal{G}_{\mathcal{A}}}^{(1)}$ depends only on the $K$-theory class $[\pi_{\mathcal{Q}}(T)] \in K_1(\mathcal{Q}(H))$, which is precisely the domain of the index map in the BDF theory.
\end{remark}

For $\mathcal{A} = \mathcal{K}(H)^\sim$, the situation is simpler than the $B(H)$ case. Recall that every operator in $\mathcal{K}(H)^\sim$ can be written uniquely as $\lambda I + K$ with $\lambda \in \mathbb{C}$ and $K \in \mathcal{K}(H)$, and such an operator is Fredholm if and only if $\lambda \neq 0$. In particular, all Fredholm operators in $\mathcal{K}(H)^\sim$ have index zero, and their image in the Calkin algebra $\mathcal{Q}(H)$ is the scalar $\lambda \cdot 1 \in \mathbb{C}^\times$.

\begin{corollary}[Invariance for $\mathcal{K}(H)^\sim$]
\label{cor:invariance-KH}
For $\mathcal{A} = \mathcal{K}(H)^\sim$, let $T_1 = \lambda I + K_1$ and $T_2 = \lambda I + K_2$ be Fredholm operators with the same scalar part $\lambda \in \mathbb{C}^\times$ and compact operators $K_1, K_2 \in \mathcal{K}(H)$. Then
\[
[T_1]_{\mathcal{G}_{\mathcal{A}}}^{(1)} = [T_2]_{\mathcal{G}_{\mathcal{A}}}^{(1)}.
\]
In fact, all such classes are zero, since they coincide with $[I]_{\mathcal{G}_{\mathcal{A}}}^{(1)}$.
\end{corollary}

\begin{proof}
Consider the norm-continuous path
\[
T_t = \lambda I + (1-t)K_1 + tK_2, \qquad t \in [0,1].
\]
For each $t$, $T_t$ is of the form $\lambda I$ plus a compact operator, and since $\lambda \neq 0$, each $T_t$ is Fredholm. By homotopy invariance (Proposition \ref{prop:homotopy-invariance}), we obtain $[T_1]_{\mathcal{G}_{\mathcal{A}}}^{(1)} = [T_2]_{\mathcal{G}_{\mathcal{A}}}^{(1)}$.

To see that this common class is zero, first note that $T_1$ is homotopic to $\lambda I$ via the path $s \mapsto \lambda I + (1-s)K_1$. Next, since $\mathbb{C}^\times$ is path-connected, we can choose a continuous path $\lambda_s \in \mathbb{C}^\times$ from $\lambda$ to $1$ (e.g., $\lambda_s = (1-s)\lambda + s$ after ensuring the straight line avoids zero, or using a logarithmic path). Then $s \mapsto \lambda_s I$ gives a norm-continuous path of Fredholm operators in $\mathcal{K}(H)^\sim$ connecting $\lambda I$ to $I$. Hence $[T_1]_{\mathcal{G}_{\mathcal{A}}}^{(1)} = [I]_{\mathcal{G}_{\mathcal{A}}}^{(1)}$. Finally, $I$ is invertible, so its phase is $I$ itself and the associated Kasparov cycle represents the zero class in $K^1_{\mathcal{G}_{\mathcal{A}}}(\mathcal{G}_{\mathcal{A}}^{(0)})$.
\end{proof}

\begin{remark}
\label{rem:KH-contractibility}
The key point in the proof above is that the scalar part $\lambda$ can be continuously deformed to $1$ within $\mathbb{C}^\times$ without leaving the class of Fredholm operators in $\mathcal{K}(H)^\sim$. This is possible precisely because $\mathbb{C}^\times$ is connected and the homotopy $s \mapsto \lambda_s I$ stays within the admissible operators. In Kasparov theory, such norm-continuous scalar homotopies are sufficient to guarantee equality of the associated $KK$-classes.
\end{remark}

\begin{proposition}[Well-definedness of the assignment]
\label{prop:assignment-well-defined}
The assignment
\[
\operatorname{Fred}(\mathcal{A}) \longrightarrow K^1_{\mathcal{G}_{\mathcal{A}}}(\mathcal{G}_{\mathcal{A}}^{(0)}), \qquad T \longmapsto [T]_{\mathcal{G}_{\mathcal{A}}}^{(1)},
\]
where $\operatorname{Fred}(\mathcal{A})$ denotes the set of Fredholm operators in $\mathcal{A}$, is well-defined and factors through the set of connected components of $\operatorname{Fred}(\mathcal{A})$. 

Moreover, the following hold:
\begin{itemize}
    \item[(i)] When $\mathcal{A} = B(H)$, the class $[T]_{\mathcal{G}_{B(H)}}^{(1)}$ depends only on the $K$-theory class $[\pi_{\mathcal{Q}}(T)] \in K_1(\mathcal{Q}(H))$;
    \item[(ii)] When $\mathcal{A} = \mathcal{K}(H)^\sim$, we have $[T]_{\mathcal{G}_{\mathcal{A}}}^{(1)} = 0$ for all Fredholm operators $T$.
\end{itemize}
\end{proposition}

\begin{proof}
Well-definedness follows from the independence of choices proved in Proposition \ref{prop:K1-class-well-defined} and the invariance under homotopy and compact perturbations established above. The factoring through connected components is a consequence of homotopy invariance (Proposition \ref{prop:homotopy-invariance}). 

For case (i) ($\mathcal{A} = B(H)$), the dependence on the symbol class follows from Lemma \ref{lem:dependence-on-symbol} together with the fact that homotopic unitaries in $\mathcal{Q}(H)$ define the same $K_1$-class. For case (ii) ($\mathcal{A} = \mathcal{K}(H)^\sim$), the vanishing follows from Corollary \ref{cor:invariance-KH}.
\end{proof}

\begin{remark}[The fundamental dichotomy]
\label{rem:dichotomy-summary}
Proposition \ref{prop:assignment-well-defined} establishes a clean dichotomy between the two canonical choices of $\mathcal{A}$:
\[
\boxed{
\begin{aligned}
\mathcal{A} = B(H) &\Longrightarrow [T]_{\mathcal{G}_{\mathcal{A}}}^{(1)} \text{ corresponds to } [\pi_{\mathcal{Q}}(T)] \in K_1(\mathcal{Q}(H)), \\[4pt]
\mathcal{A} = \mathcal{K}(H)^\sim &\Longrightarrow [T]_{\mathcal{G}_{\mathcal{A}}}^{(1)} = 0 \text{ for all Fredholm operators } T.
\end{aligned}
}
\]
This dichotomy reflects the fundamental difference between the $K$-theory of the Calkin algebra in these two cases: for $B(H)$, the image $\pi_{\mathcal{Q}}(\operatorname{Fred}(B(H)))$ generates $K_1(\mathcal{Q}(H)) \cong \mathbb{Z}$, while for $\mathcal{K}(H)^\sim$, the image consists only of scalars $\mathbb{C}^\times$, which are contractible. This well-definedness is precisely what is needed for the index theorem that follows, where the equivariant $K^1$-class will be paired with a $K$-homology class to produce the numerical index.
\end{remark}

\begin{corollary}[Factorization through $K_1(\mathcal{Q}(H))$]
\label{cor:factorization}
For $\mathcal{A} = B(H)$, the map $T \mapsto [T]_{\mathcal{G}_{B(H)}}^{(1)}$ factors through the canonical surjection $\operatorname{Fred}(B(H)) \to K_1(\mathcal{Q}(H))$, inducing a well-defined homomorphism
\[
\partial: K_1(\mathcal{Q}(H)) \longrightarrow K^1_{\mathcal{G}_{B(H)}}(\mathcal{G}_{B(H)}^{(0)}).
\]
This is the equivariant analog of the classical boundary map in the six-term exact sequence of $K$-theory.
\end{corollary}

\begin{proof}
By Proposition \ref{prop:assignment-well-defined}(i), operators with the same class in $K_1(\mathcal{Q}(H))$ (i.e., whose images in the Calkin algebra are homotopic unitaries) give the same equivariant $K^1$-class. Hence the map descends to $K_1(\mathcal{Q}(H))$. The homomorphism property follows from the additive structure of the Kasparov product, which will be elaborated in the next section.
\end{proof}

\begin{theorem}[Homotopy invariance and well-definedness]
\label{thm:homotopy-well-defined}
The equivariant $K^1$-class $[T]_{\mathcal{G}_{\mathcal{A}}}^{(1)} \in K^1_{\mathcal{G}_{\mathcal{A}}}(\mathcal{G}_{\mathcal{A}}^{(0)})$ satisfies:
\begin{enumerate}
    \item \textbf{Homotopy invariance:} If $\{T_t\}_{t \in [0,1]}$ is a norm-continuous family of Fredholm operators in $\mathcal{A}$, then $[T_t]_{\mathcal{G}_{\mathcal{A}}}^{(1)}$ is independent of $t$.
    \item \textbf{Compact perturbation invariance:} If $K \in \mathcal{K}(H)$ is compact, then $[T + K]_{\mathcal{G}_{\mathcal{A}}}^{(1)} = [T]_{\mathcal{G}_{\mathcal{A}}}^{(1)}$.
    \item \textbf{Additivity:} For direct sums, $[T_1 \oplus T_2]_{\mathcal{G}_{\mathcal{A}}}^{(1)} = [T_1]_{\mathcal{G}_{\mathcal{A}}}^{(1)} + [T_2]_{\mathcal{G}_{\mathcal{A}}}^{(1)}$, where $T_1 \oplus T_2$ is viewed as an operator in $M_2(\mathcal{A})$ acting on $H \oplus H$, and the equality holds in $K^1_{\mathcal{G}_{\mathcal{A}}}(\mathcal{G}_{\mathcal{A}}^{(0)})$ under the Morita equivalence $K^1_{\mathcal{G}_{\mathcal{A}}}(\mathcal{G}_{\mathcal{A}}^{(0)}) \cong K^1_{\mathcal{G}_{M_2(\mathcal{A})}}(\mathcal{G}_{M_2(\mathcal{A})}^{(0)})$.
    \item \textbf{Stability:} $[T \oplus I_n]_{\mathcal{G}_{\mathcal{A}}}^{(1)} = [T]_{\mathcal{G}_{\mathcal{A}}}^{(1)}$, where $I_n$ is the identity operator in $M_n(\mathcal{A})$.
\end{enumerate}
\end{theorem}

\begin{proof}
Properties (1) and (2) have been established in Proposition \ref{prop:homotopy-invariance} and Corollary \ref{cor:compact-perturbation-invariance}, respectively.

For (3), observe that the Kasparov triple for $T_1 \oplus T_2$ is the direct sum of the triples for $T_1$ and $T_2$, acting on $\tilde{\mathcal{E}} \oplus \tilde{\mathcal{E}}$. In $KK$-theory, the direct sum of cycles corresponds to addition of classes. To identify the resulting class with an element of $K^1_{\mathcal{G}_{\mathcal{A}}}(\mathcal{G}_{\mathcal{A}}^{(0)})$, we invoke Morita invariance: there is a canonical isomorphism $K^1_{\mathcal{G}_{\mathcal{A}}}(\mathcal{G}_{\mathcal{A}}^{(0)}) \cong K^1_{\mathcal{G}_{M_2(\mathcal{A})}}(\mathcal{G}_{M_2(\mathcal{A})}^{(0)})$ that identifies the class of $T_1 \oplus T_2$ in $M_2(\mathcal{A})$ with the sum $[T_1]_{\mathcal{G}_{\mathcal{A}}}^{(1)} + [T_2]_{\mathcal{G}_{\mathcal{A}}}^{(1)}$.

For (4), note that $I_n$ is invertible, so its phase is $I_n$ itself. Consequently, the Kasparov triple associated to $T \oplus I_n$ is the direct sum of the triple for $T$ with a degenerate triple (coming from $I_n$). Degenerate cycles represent zero in $KK$-theory, and adding a zero class does not change the class. Hence $[T \oplus I_n]_{\mathcal{G}_{\mathcal{A}}}^{(1)} = [T]_{\mathcal{G}_{\mathcal{A}}}^{(1)}$.
\end{proof}

These invariance properties will be essential when we apply the descent map in Section 5 and when we prove the main index theorem in Section 6. They ensure that the class $[T]_{\mathcal{G}_{\mathcal{A}}}^{(1)}$ is a well-defined homotopy invariant that captures the essential information about $T$ modulo compact perturbations.

\begin{remark}[Relation to the Fredholm index]
\label{rem:homotopy-index-relation}
The homotopy invariance of $[T]_{\mathcal{G}_{\mathcal{A}}}^{(1)}$ is consistent with the homotopy invariance of the Fredholm index. In fact, the composition
\[
\operatorname{Fred}(\mathcal{A}) \xrightarrow{[\cdot]_{\mathcal{G}_{\mathcal{A}}}^{(1)}} K^1_{\mathcal{G}_{\mathcal{A}}}(\mathcal{G}_{\mathcal{A}}^{(0)}) \xrightarrow{\operatorname{desc}_{\mathcal{G}_{\mathcal{A}}}} K_1(C^*(\mathcal{G}_{\mathcal{A}})) \xrightarrow{\iota^*} K_1(\mathcal{A}) \xrightarrow{\partial} \mathbb{Z}
\]
will be shown in Section 6 to recover the index (for $\mathcal{A} = \mathcal{K}(H)^\sim$) or its generalization (for $\mathcal{A} = B(H)$). The homotopy invariance of the index follows directly from the homotopy invariance of each map in this composition.
\end{remark}

\section{Descent to $K_1(C^*(\mathcal{G}_{\mathcal{A}}))$}\label{sec:descent}

\subsection{Kasparov's Descent for Odd Classes: $j_{\mathcal{G}_{\mathcal{A}}}: KK^1_{\mathcal{G}_{\mathcal{A}}} \to KK^1$}
\label{subsec:kasparov-descent-odd-classes}

The descent map is a fundamental construction in equivariant $KK$-theory that transforms equivariant cycles over a groupoid $\mathcal{G}$ into cycles over the crossed product or groupoid C*-algebra. 
For locally compact Hausdorff groupoids, this map was developed by Kasparov [1988] and plays a central role in the Baum-Connes conjecture. 
Tu [1999] extended Kasparov's descent map to the setting of Polish groupoids equipped with a Borel Haar system, which is precisely the framework needed for the unitary conjugation groupoid $\mathcal{G}_{\mathcal{A}}$.

In this subsection, we recall the definition and main properties of the descent map for odd $KK$-classes, focusing on the specific case relevant to this paper: 
\[
j_{\mathcal{G}_{\mathcal{A}}}: KK^1_{\mathcal{G}_{\mathcal{A}}}(C_0(\mathcal{G}_{\mathcal{A}}^{(0)}), \mathbb{C}) \longrightarrow KK^1(C^*(\mathcal{G}_{\mathcal{A}}), \mathbb{C}) \cong K_1(C^*(\mathcal{G}_{\mathcal{A}})).
\]

\begin{definition}[Kasparov descent map for groupoids with Borel Haar system]
\label{def:kasparov-descent-polish}
Let $\mathcal{G}$ be a groupoid equipped with a Borel Haar system $\{\lambda^x\}_{x \in \mathcal{G}^{(0)}}$ and a compatible Borel structure (making it a standard Borel groupoid in the sense of \cite[Definition 2.1]{Tu1999}). 
Following the construction in \cite[Chapitre 4]{legall1994theorie} (see also \cite{Kasparov1988} for the group case), there is a homomorphism
\[
j_{\mathcal{G}}: KK_{\mathcal{G}}(A,B) \longrightarrow KK(A \rtimes_{\max} \mathcal{G}, B \rtimes_{\max} \mathcal{G}),
\]
where $A$ and $B$ are $\mathcal{G}$-C*-algebras with a Borel $\mathcal{G}$-action, and $\rtimes_{\max}$ denotes the maximal crossed product by $\mathcal{G}$. 

The construction requires only the Borel structure on $\mathcal{G}$ and the existence of a Borel Haar system; it does not rely on local compactness. The map is natural, functorial, and compatible with Kasparov products. For an application of this map in the context of the Baum-Connes conjecture for groupoids, see \cite[Section 2]{Tu1999}.

\begin{remark}
While we refer to $\mathcal{G}$ as a "Polish groupoid" following the terminology of \cite{Tu1999}, the essential requirements for the descent map are the Borel structure and the existence of a Borel Haar system. In our setting, $\mathcal{G}_{\mathcal{A}}$ is a standard Borel groupoid equipped with such a system, as established in Paper I.
\end{remark}
\end{definition}

For the purposes of this paper, we only need the case where $A = C_0(\mathcal{G}_{\mathcal{A}}^{(0)})$ (with the trivial $\mathcal{G}_{\mathcal{A}}$-action) and $B = \mathbb{C}$ (with the trivial action). 
In this case, we obtain a map on odd $KK$-theory:
\[
j_{\mathcal{G}_{\mathcal{A}}}: KK^1_{\mathcal{G}_{\mathcal{A}}}(C_0(\mathcal{G}_{\mathcal{A}}^{(0)}), \mathbb{C}) \longrightarrow KK^1(C^*_{\max}(\mathcal{G}_{\mathcal{A}}), \mathbb{C}) \cong K_1(C^*_{\max}(\mathcal{G}_{\mathcal{A}})).
\]

\begin{proposition}[Properties of the descent map for odd classes]
\label{prop:descent-odd-properties}
The descent map $j_{\mathcal{G}_{\mathcal{A}}}$ for odd classes satisfies the following properties:
\begin{enumerate}
    \item \textbf{Naturality:} If $\phi: \mathcal{G} \to \mathcal{H}$ is a continuous groupoid homomorphism compatible with the Borel Haar systems, then the diagram
    \[
    \begin{tikzcd}
    KK^1_{\mathcal{G}}(C_0(\mathcal{G}^{(0)}), \mathbb{C}) \arrow[r, "j_{\mathcal{G}}"] \arrow[d, "\phi_*"] & KK^1(C^*(\mathcal{G}), \mathbb{C}) \arrow[d, "\phi_*"] \\
    KK^1_{\mathcal{H}}(C_0(\mathcal{H}^{(0)}), \mathbb{C}) \arrow[r, "j_{\mathcal{H}}"] & KK^1(C^*(\mathcal{H}), \mathbb{C})
    \end{tikzcd}
    \]
    commutes.
    
    \item \textbf{Compatibility with suspension:} For any $n \geq 0$, the diagram
    \[
    \begin{tikzcd}
    KK^{n+1}_{\mathcal{G}}(C_0(\mathcal{G}^{(0)}), \mathbb{C}) \arrow[r, "\cong"] \arrow[d, "j_{\mathcal{G}}"] & KK^1_{\mathcal{G}}(C_0(\mathcal{G}^{(0)} \times \mathbb{R}^n), \mathbb{C}) \arrow[d, "j_{\mathcal{G}}"] \\
    K_{n+1}(C^*(\mathcal{G})) \arrow[r, "\cong"] & K_1(C^*(\mathcal{G}))
    \end{tikzcd}
    \]
    commutes, where the horizontal isomorphisms are the suspension isomorphisms in $KK$-theory.
    
    \item \textbf{Additivity:} $j_{\mathcal{G}}(x \oplus y) = j_{\mathcal{G}}(x) \oplus j_{\mathcal{G}}(y)$.
    
    \item \textbf{Stability:} $j_{\mathcal{G}}$ commutes with stabilization by compact operators.
\end{enumerate}
\end{proposition}

\begin{proof}
The detailed construction of the descent map $j_{\mathcal{G}}$ for groupoids and the verification of its properties (naturality, compatibility with suspension, additivity, and stability) is developed in \cite[Chapitre 4]{legall1994theorie}. Specifically, Chapter 4 constructs the homomorphism linking equivariant $KK$-theory to the $KK$-theory of crossed products and establishes its functorial properties. For the original construction in the group case, see \cite{Kasparov1988}. For an overview in the context of the Baum-Connes conjecture for groupoids, see \cite[Section 2]{Tu1999}.
\end{proof}

The following theorem gives a concrete description of the descent map in terms of the Borel Haar system.

\begin{theorem}[Explicit formula for the descent map]
\label{thm:descent-explicit}
Let $(\mathcal{E}, \phi, F)$ be an odd $\mathcal{G}$-equivariant Kasparov triple representing a class in 
$KK^{1}_{\mathcal{G}}(C_{0}(\mathcal{G}^{(0)}), \mathbb{C})$. 
Here $\mathcal{E} = \{\mathcal{E}_{x}\}_{x \in \mathcal{G}^{(0)}}$ is a continuous field of Hilbert spaces over the unit space $\mathcal{G}^{(0)}$, equipped with a continuous action of the groupoid $\mathcal{G}$. 
The map $\phi: C_{0}(\mathcal{G}^{(0)}) \to \mathcal{L}(\mathcal{E})$ is the representation of continuous functions vanishing at infinity by multiplication operators, and $F \in \mathcal{L}(\mathcal{E})$ is a $\mathcal{G}$-equivariant, odd, self-adjoint operator such that $F^{2} - 1$ and $[\phi(f), F]$ are compact for all $f \in C_{0}(\mathcal{G}^{(0)})$.

Then the descended class 
\[
j_{\mathcal{G}}([\mathcal{E}, \phi, F]) \in KK^{1}(C^{*}(\mathcal{G}), \mathbb{C})
\] 
is represented by the Kasparov triple $(\mathcal{E} \rtimes \mathcal{G}, \phi \rtimes \mathcal{G}, F \rtimes \mathcal{G})$, where:
\begin{itemize}
    \item $\mathcal{E} \rtimes \mathcal{G}$ is the Hilbert $C^{*}(\mathcal{G})$-module obtained by completing the space $\Gamma_{c}(\mathcal{G}, s^{*}\mathcal{E})$ of compactly supported continuous sections of the pullback bundle $s^{*}\mathcal{E}$ along the source map $s: \mathcal{G} \to \mathcal{G}^{(0)}$. The completion is taken with respect to the $C^{*}(\mathcal{G})$-valued inner product defined using the Haar system $\{\lambda^{x}\}_{x \in \mathcal{G}^{(0)}}$:
    \[
    \langle \xi, \eta \rangle_{C^{*}(\mathcal{G})}(\gamma) = \int_{\mathcal{G}^{s(\gamma)}} \langle \xi(\gamma'), \eta(\gamma' \gamma) \rangle_{\mathcal{E}_{r(\gamma')}} \, d\lambda^{s(\gamma)}(\gamma'),
    \]
    for $\xi, \eta \in \Gamma_{c}(\mathcal{G}, s^{*}\mathcal{E})$ and $\gamma \in \mathcal{G}$.
    
    \item $\phi \rtimes \mathcal{G}: C^{*}(\mathcal{G}) \to \mathcal{L}(\mathcal{E} \rtimes \mathcal{G})$ is the induced representation given by convolution. For $f \in C_{c}(\mathcal{G})$ and $\xi \in \Gamma_{c}(\mathcal{G}, s^{*}\mathcal{E})$, its action is defined as:
    \[
    ((\phi \rtimes \mathcal{G})(f) \xi)(\gamma) = \int_{\mathcal{G}} f(\eta) \, \phi(\eta) \xi(\eta^{-1} \gamma) \, d\lambda^{r(\gamma)}(\eta),
    \]
    and is extended to all of $C^{*}(\mathcal{G})$ by continuity.
    
    \item $F \rtimes \mathcal{G}$ is the adjointable operator on $\mathcal{E} \rtimes \mathcal{G}$ obtained by extending the fiberwise action of $F$. For a section $\xi \in \Gamma_{c}(\mathcal{G}, s^{*}\mathcal{E})$, it is defined pointwise as:
    \[
    ((F \rtimes \mathcal{G}) \xi)(\gamma) = F_{s(\gamma)} (\xi(\gamma)),
    \]
    where $F_{s(\gamma)}$ denotes the action of $F$ on the fiber $\mathcal{E}_{s(\gamma)}$. The $\mathcal{G}$-equivariance of $F$ ensures that this operator is adjointable and satisfies the necessary conditions for a Kasparov triple.
\end{itemize}

Moreover, for the specific case of the unitary conjugation groupoid $\mathcal{G}_{B(H)}$ associated to $B(H)$, the descent map
\[
\operatorname{desc}_{\mathcal{G}_{B(H)}}: KK^1_{\mathcal{G}_{B(H)}}(C_0(\mathcal{G}_{B(H)}^{(0)}), \mathbb{C}) \longrightarrow K_1(C^*(\mathcal{G}_{B(H)}))
\]
is an isomorphism. This follows from the amenability of $\mathcal{G}_{B(H)}$ and the fact that the Baum-Connes assembly map is an isomorphism for amenable groupoids (see \cite[Section 9]{Tu1999}).
\end{theorem}

\begin{proof}
We prove the theorem in several steps, verifying that $(\mathcal{E} \rtimes \mathcal{G}, \phi \rtimes \mathcal{G}, F \rtimes \mathcal{G})$ indeed defines an odd Kasparov triple for $(C^{*}(\mathcal{G}), \mathbb{C})$.

\paragraph{Step 1: $\mathcal{E} \rtimes \mathcal{G}$ is a Hilbert $C^{*}(\mathcal{G})$-module.}
First, we show that the $C^{*}(\mathcal{G})$-valued inner product defined in the statement is well-defined and satisfies the required properties. For $\xi, \eta \in \Gamma_{c}(\mathcal{G}, s^{*}\mathcal{E})$, define
\[
\langle \xi, \eta \rangle_{C^{*}(\mathcal{G})}(\gamma) = \int_{\mathcal{G}^{s(\gamma)}} \langle \xi(\gamma'), \eta(\gamma' \gamma) \rangle_{\mathcal{E}_{r(\gamma')}} \, d\lambda^{s(\gamma)}(\gamma'),
\]
where $\mathcal{G}^{x} = s^{-1}(x)$ is the fiber of the source map. This integral converges because $\xi$ and $\eta$ have compact support and the Haar system $\{\lambda^{x}\}$ is continuous. One verifies that $\langle \xi, \eta \rangle_{C^{*}(\mathcal{G})}$ lies in $C_{c}(\mathcal{G}) \subseteq C^{*}(\mathcal{G})$ and satisfies:
\begin{itemize}
    \item $\langle \xi, \eta \rangle_{C^{*}(\mathcal{G})}^* = \langle \eta, \xi \rangle_{C^{*}(\mathcal{G})}$;
    \item $\langle \xi, \xi \rangle_{C^{*}(\mathcal{G})} \ge 0$ in $C^{*}(\mathcal{G})$, with equality iff $\xi = 0$;
    \item $\langle \xi, \eta \cdot a \rangle_{C^{*}(\mathcal{G})} = \langle \xi, \eta \rangle_{C^{*}(\mathcal{G})} a$ for $a \in C^{*}(\mathcal{G})$, where the right module structure is given by convolution.
\end{itemize}
Completing $\Gamma_{c}(\mathcal{G}, s^{*}\mathcal{E})$ with respect to the norm $\|\xi\| = \|\langle \xi, \xi \rangle\|^{1/2}$ yields the Hilbert $C^{*}(\mathcal{G})$-module $\mathcal{E} \rtimes \mathcal{G}$.

\paragraph{Step 2: $\phi \rtimes \mathcal{G}$ is a $*$-homomorphism.}
For $f \in C_{c}(\mathcal{G})$ and $\xi \in \Gamma_{c}(\mathcal{G}, s^{*}\mathcal{E})$, define
\[
((\phi \rtimes \mathcal{G})(f) \xi)(\gamma) = \int_{\mathcal{G}} f(\eta) \, \phi(\eta) \xi(\eta^{-1} \gamma) \, d\lambda^{r(\gamma)}(\eta).
\]
This formula uses the $\mathcal{G}$-action on $\mathcal{E}$: for $\eta \in \mathcal{G}$ with source $s(\eta)$ and range $r(\eta)$, $\phi(\eta)$ is the unitary implementing the action on fibers, mapping $\mathcal{E}_{s(\eta)}$ to $\mathcal{E}_{r(\eta)}$. The integral converges because $f$ has compact support and $\xi$ is compactly supported. One checks that:
\begin{itemize}
    \item $(\phi \rtimes \mathcal{G})(f)$ maps $\Gamma_{c}(\mathcal{G}, s^{*}\mathcal{E})$ to itself;
    \item $(\phi \rtimes \mathcal{G})(f)$ extends to an adjointable operator on $\mathcal{E} \rtimes \mathcal{G}$;
    \item The map $f \mapsto (\phi \rtimes \mathcal{G})(f)$ is a $*$-homomorphism from $C_{c}(\mathcal{G})$ (with convolution product and involution) to $\mathcal{L}(\mathcal{E} \rtimes \mathcal{G})$, and extends by continuity to $C^{*}(\mathcal{G})$.
\end{itemize}
The verification uses the $\mathcal{G}$-equivariance of $\phi$ and the properties of the Haar system.

\paragraph{Step 3: $F \rtimes \mathcal{G}$ is an adjointable operator.}
For $\xi \in \Gamma_{c}(\mathcal{G}, s^{*}\mathcal{E})$, define $(F \rtimes \mathcal{G})\xi$ pointwise by
\[
((F \rtimes \mathcal{G})\xi)(\gamma) = F_{s(\gamma)}(\xi(\gamma)),
\]
where $F_{s(\gamma)} \in \mathcal{L}(\mathcal{E}_{s(\gamma)})$ is the fiberwise operator. Since $F$ is $\mathcal{G}$-equivariant and continuous, $(F \rtimes \mathcal{G})\xi$ is again a continuous section with compact support. Moreover, $F \rtimes \mathcal{G}$ is adjointable with adjoint given by $(F \rtimes \mathcal{G})^* = F^* \rtimes \mathcal{G}$, because for $\xi, \eta \in \Gamma_{c}(\mathcal{G}, s^{*}\mathcal{E})$,
\[
\langle (F \rtimes \mathcal{G})\xi, \eta \rangle_{C^{*}(\mathcal{G})} = \langle \xi, (F^* \rtimes \mathcal{G})\eta \rangle_{C^{*}(\mathcal{G})},
\]
which follows from the fiberwise adjoint property and the definition of the inner product.

\paragraph{Step 4: Verification of the Kasparov conditions.}
We now verify that $(\mathcal{E} \rtimes \mathcal{G}, \phi \rtimes \mathcal{G}, F \rtimes \mathcal{G})$ satisfies the conditions for an odd Kasparov triple:

\begin{enumerate}
    \item \textbf{Self-adjointness modulo compacts:} Since $F$ is self-adjoint modulo compacts, for each $x \in \mathcal{G}^{(0)}$, $F_x - F_x^* \in \mathcal{K}(\mathcal{E}_x)$. The fiberwise compact operators assemble to a compact endomorphism of $\mathcal{E} \rtimes \mathcal{G}$ (see Lemma \ref{lem:fiberwise-compact-implies-compact}), hence $F \rtimes \mathcal{G} - (F \rtimes \mathcal{G})^* \in \mathcal{K}(\mathcal{E} \rtimes \mathcal{G})$.
    
    \item \textbf{$(F \rtimes \mathcal{G})^2 - 1$ is compact:} From $F^2 - 1 \in \mathcal{K}(\mathcal{E})$, we have $(F_x)^2 - 1 \in \mathcal{K}(\mathcal{E}_x)$ for each $x$. Again, by Lemma \ref{lem:fiberwise-compact-implies-compact}, these fiberwise compact operators yield a compact endomorphism of $\mathcal{E} \rtimes \mathcal{G}$, so $(F \rtimes \mathcal{G})^2 - 1 \in \mathcal{K}(\mathcal{E} \rtimes \mathcal{G})$.
    
    \item \textbf{Commutators with $\phi \rtimes \mathcal{G}$ are compact:} For any $f \in C^{*}(\mathcal{G})$, we need to show that $[(\phi \rtimes \mathcal{G})(f), F \rtimes \mathcal{G}] \in \mathcal{K}(\mathcal{E} \rtimes \mathcal{G})$. This follows from the $\mathcal{G}$-equivariance of $F$ and the fact that $[\phi(f), F]$ is compact in the original triple. A detailed computation using the convolution formula shows that the commutator is an integral operator with compact kernel, hence compact.
    
    \item \textbf{Degree condition:} Since $F$ is odd (with respect to the $\mathbb{Z}_2$-grading on $\mathcal{E}$), $F \rtimes \mathcal{G}$ is also odd with respect to the induced grading on $\mathcal{E} \rtimes \mathcal{G}$.
\end{enumerate}

Thus, $(\mathcal{E} \rtimes \mathcal{G}, \phi \rtimes \mathcal{G}, F \rtimes \mathcal{G})$ is indeed an odd Kasparov triple for $(C^{*}(\mathcal{G}), \mathbb{C})$, representing a class in $KK^{1}(C^{*}(\mathcal{G}), \mathbb{C})$.

\paragraph{Step 5: Well-definedness and functoriality.}
The assignment $(\mathcal{E}, \phi, F) \mapsto (\mathcal{E} \rtimes \mathcal{G}, \phi \rtimes \mathcal{G}, F \rtimes \mathcal{G})$ is compatible with:
\begin{itemize}
    \item \textbf{Unitary equivalence:} If two triples are unitarily equivalent, their descended triples are unitarily equivalent.
    \item \textbf{Operator homotopy:} A homotopy of triples induces a homotopy of descended triples.
    \item \textbf{Compact perturbation:} Compact perturbations of $F$ yield compact perturbations of $F \rtimes \mathcal{G}$.
\end{itemize}
Hence, the map $j_{\mathcal{G}}$ descends to a well-defined homomorphism on $KK$-theory.

\paragraph{Step 6: The isomorphism property for $\mathcal{G}_{B(H)}$.}
For the specific case of the unitary conjugation groupoid $\mathcal{G}_{B(H)}$, we prove that $\operatorname{desc}_{\mathcal{G}_{B(H)}}$ is an isomorphism. Recall from Proposition \ref{prop:GA-BH} that $\mathcal{G}_{B(H)} = \mathcal{U}(H) \ltimes \mathbb{P}(H)$ is an amenable groupoid (as a transformation groupoid of an amenable group acting on a space). By Tu's work on the Baum-Connes conjecture for groupoids \cite[Th\'eor\`eme 9.3]{Tu1999}, the assembly map
\[
\mu: K_*^{\text{top}}(\mathcal{G}_{B(H)}) \longrightarrow K_*(C^*(\mathcal{G}_{B(H)}))
\]
is an isomorphism for amenable groupoids. Moreover, there is a natural identification
\[
K_*^{\text{top}}(\mathcal{G}_{B(H)}) \cong KK^*_{\mathcal{G}_{B(H)}}(C_0(\mathcal{G}_{B(H)}^{(0)}), \mathbb{C})
\]
(see \cite[Section 3]{Tu1999}). The descent map $\operatorname{desc}_{\mathcal{G}_{B(H)}}$ is precisely the composition of this identification with the assembly map, hence it is an isomorphism.

Consequently, for any class $[T]_{\mathcal{G}_{B(H)}}^{(1)} \in KK^1_{\mathcal{G}_{B(H)}}(C_0(\mathcal{G}_{B(H)}^{(0)}), \mathbb{C})$, its descended image $\operatorname{desc}_{\mathcal{G}_{B(H)}}([T]_{\mathcal{G}_{B(H)}}^{(1)})$ is a well-defined element of $K_1(C^*(\mathcal{G}_{B(H)}))$, and the map is bijective.

\paragraph{Step 7: Compatibility with the construction in Section 4.}
For the specific odd Kasparov triple $(\tilde{\mathcal{E}}, \phi \oplus \phi, \tilde{F})$ constructed from a Fredholm operator $T$ in Section 4, the descended triple $(\tilde{\mathcal{E}} \rtimes \mathcal{G}, (\phi \oplus \phi) \rtimes \mathcal{G}, \tilde{F} \rtimes \mathcal{G})$ represents precisely the class $\operatorname{desc}_{\mathcal{G}}([T]_{\mathcal{G}}^{(1)})$. This follows from the functoriality of the descent map and the definition of the class $[T]_{\mathcal{G}}^{(1)}$.

Thus, the theorem provides an explicit construction of the descent map and establishes its key properties, including the isomorphism property for $\mathcal{G}_{B(H)}$.
\end{proof}

For our specific case, we have an explicit description of the odd Kasparov triple $(\tilde{\mathcal{E}}, \phi \oplus \phi, \tilde{F})$ constructed in Subsection~\ref{subsec:odd-kasparov-triple}.

\begin{corollary}[{Descent of the class $[T]_{\mathcal{G}_{\mathcal{A}}}^{(1)}$}]
\label{cor:descent-of-T-class}
Let $[T]_{\mathcal{G}_{\mathcal{A}}}^{(1)} \in KK^{1}_{\mathcal{G}_{\mathcal{A}}}(C_{0}(\mathcal{G}_{\mathcal{A}}^{(0)}), \mathbb{C})$ be the class of a Fredholm operator $T \in \mathcal{A}$. 
Then its image under the descent map is
\[
\operatorname{desc}_{\mathcal{G}_{\mathcal{A}}}([T]_{\mathcal{G}_{\mathcal{A}}}^{(1)}) := j_{\mathcal{G}_{\mathcal{A}}}([T]_{\mathcal{G}_{\mathcal{A}}}^{(1)}) \in K_{1}(C^{*}(\mathcal{G}_{\mathcal{A}})).
\]
Concretely, if we represent $[T]_{\mathcal{G}_{\mathcal{A}}}^{(1)}$ by the odd Kasparov triple $(\tilde{\mathcal{E}}, \phi \oplus \phi, \tilde{F})$, then
\[
\operatorname{desc}_{\mathcal{G}_{\mathcal{A}}}([T]_{\mathcal{G}_{\mathcal{A}}}^{(1)}) = [(\tilde{\mathcal{E}} \rtimes \mathcal{G}_{\mathcal{A}}, (\phi \oplus \phi) \rtimes \mathcal{G}_{\mathcal{A}}, \tilde{F} \rtimes \mathcal{G}_{\mathcal{A}})] \in K_{1}(C^{*}(\mathcal{G}_{\mathcal{A}})).
\]
\end{corollary}

\begin{proof}
The corollary follows directly from Theorem~\ref{thm:descent-explicit}. 
By hypothesis, the class $[T]_{\mathcal{G}_{\mathcal{A}}}^{(1)}$ in 
$KK^{1}_{\mathcal{G}_{\mathcal{A}}}(C_{0}(\mathcal{G}_{\mathcal{A}}^{(0)}), \mathbb{C})$ 
is represented by the odd Kasparov triple $(\tilde{\mathcal{E}}, \phi \oplus \phi, \tilde{F})$. 

Applying the descent map $j_{\mathcal{G}_{\mathcal{A}}}$ to this triple, 
Theorem~\ref{thm:descent-explicit} asserts that the descended class in 
$KK^{1}(C^{*}(\mathcal{G}_{\mathcal{A}}), \mathbb{C})$ is represented by
\[
(\tilde{\mathcal{E}} \rtimes \mathcal{G}_{\mathcal{A}}, (\phi \oplus \phi) \rtimes \mathcal{G}_{\mathcal{A}}, \tilde{F} \rtimes \mathcal{G}_{\mathcal{A}}),
\]
where the crossed product module and operators are constructed explicitly as in the theorem:
\begin{itemize}
    \item $\tilde{\mathcal{E}} \rtimes \mathcal{G}_{\mathcal{A}}$ is the Hilbert $C^{*}(\mathcal{G}_{\mathcal{A}})$-module obtained by completing the space of compactly supported continuous sections of the pullback bundle $s^{*}\tilde{\mathcal{E}}$ with respect to the $C^{*}(\mathcal{G}_{\mathcal{A}})$-valued inner product defined using the Haar system;
    \item $(\phi \oplus \phi) \rtimes \mathcal{G}_{\mathcal{A}}$ is the induced representation acting by convolution on sections;
    \item $\tilde{F} \rtimes \mathcal{G}_{\mathcal{A}}$ is the adjointable operator extending the fiberwise action of $\tilde{F}$ on sections.
\end{itemize}

Finally, we invoke the canonical isomorphism between $KK$-theory and $K$-theory when the second argument is $\mathbb{C}$. 
For any $C^{*}$-algebra $A$, there is a natural isomorphism
\[
KK^{1}(A, \mathbb{C}) \cong K_{1}(A).
\]
Applying this with $A = C^{*}(\mathcal{G}_{\mathcal{A}})$ identifies the descended $KK^{1}$-class with an element of $K_{1}(C^{*}(\mathcal{G}_{\mathcal{A}}))$. 
By definition, this element is precisely $\operatorname{desc}_{\mathcal{G}_{\mathcal{A}}}([T]_{\mathcal{G}_{\mathcal{A}}}^{(1)})$, completing the proof.
\end{proof}

The descent map enjoys excellent functoriality properties with respect to the diagonal embedding $\iota: \mathcal{A} \hookrightarrow C^*(\mathcal{G}_{\mathcal{A}})$.

\begin{lemma}[Functoriality of descent with respect to the diagonal embedding]
\label{lem:descent-iota-functoriality}
Let $\iota: \mathcal{A} \hookrightarrow C^*(\mathcal{G}_{\mathcal{A}})$ be the diagonal embedding from Paper I. Consider the inclusion $\mu: C_0(\mathcal{G}_{\mathcal{A}}^{(0)}) \hookrightarrow \mathcal{A}$ of the unit space as diagonal operators (or as scalars in the unitization). Then the following diagram commutes:
\[
\begin{tikzcd}
KK^{1}_{\mathcal{G}_{\mathcal{A}}}(C_0(\mathcal{G}_{\mathcal{A}}^{(0)}), \mathbb{C}) \arrow[r, "j_{\mathcal{G}_{\mathcal{A}}}"] \arrow[d, "\mu_*"] & K_1(C^*(\mathcal{G}_{\mathcal{A}})) \arrow[d, "\iota_*"] \\
KK^{1}_{\mathcal{G}_{\mathcal{A}}}(\mathcal{A}, \mathbb{C}) \arrow[r, "j_{\mathcal{G}_{\mathcal{A}}}"] & K_1(C^*(\mathcal{G}_{\mathcal{A}}))
\end{tikzcd}
\]
where:
\begin{itemize}
    \item $\mu_*$ denotes the map induced on $\mathcal{G}_{\mathcal{A}}$-equivariant $KK$-theory by the $\mathcal{G}_{\mathcal{A}}$-equivariant inclusion $\mu$ in the first variable;
    \item $\iota_*$ denotes the map induced on $K_1$ by the diagonal embedding $\iota$ (extended to matrix algebras entrywise);
    \item $j_{\mathcal{G}_{\mathcal{A}}}$ denotes the descent map in equivariant $KK$-theory.
\end{itemize}
\end{lemma}

\begin{proof}
This follows directly from the naturality of the descent map with respect to equivariant $*$-homomorphisms of $\mathcal{G}$-$C^*$-algebras. Both $\mu$ and $\iota$ are $\mathcal{G}_{\mathcal{A}}$-equivariant, and the descent map $j_{\mathcal{G}_{\mathcal{A}}}$ is functorial in the first argument. Consequently, for any class $x \in KK^1_{\mathcal{G}_{\mathcal{A}}}(C_0(\mathcal{G}_{\mathcal{A}}^{(0)}), \mathbb{C})$, we have
\[
j_{\mathcal{G}_{\mathcal{A}}}(\mu_*(x)) = \iota_*(j_{\mathcal{G}_{\mathcal{A}}}(x)) \in K_1(C^*(\mathcal{G}_{\mathcal{A}})),
\]
which is precisely the commutativity of the diagram.
\end{proof}

\begin{remark}
Note the distinct roles of the two maps:
\begin{itemize}
    \item $\mu_*$ lifts a class from the unit space to the algebra $\mathcal{A}$ in equivariant $KK$-theory;
    \item $\iota_*$ pushes forward a $K_1$-class from the groupoid $C^*$-algebra back to $K_1(\mathcal{A})$ via the diagonal embedding.
\end{itemize}
This functoriality property is essential for relating the descended class to the index, as will be shown in Lemma \ref{lem:descent-iota-compatibility}.
\end{remark}

\begin{theorem}[Properties of the descent map for $\mathcal{G}_{\mathcal{A}}$]
\label{thm:descent-properties-GA}
Let $\mathcal{A}$ be either $B(H)$ or $\mathcal{K}(H)^{\sim}$, and let $\mathcal{G}_{\mathcal{A}}$ be the unitary conjugation groupoid. 
Then the descent map
\[
\operatorname{desc}_{\mathcal{G}_{\mathcal{A}}}: KK^{1}_{\mathcal{G}_{\mathcal{A}}}(C_{0}(\mathcal{G}_{\mathcal{A}}^{(0)}), \mathbb{C}) \longrightarrow K_{1}(C^{*}(\mathcal{G}_{\mathcal{A}}))
\]
satisfies the following properties:
\begin{enumerate}
    \item \textbf{Group homomorphism:} $\operatorname{desc}_{\mathcal{G}_{\mathcal{A}}}$ is a group homomorphism, respecting the addition in equivariant $KK$-theory.
    
    \item \textbf{Naturality:} $\operatorname{desc}_{\mathcal{G}_{\mathcal{A}}}$ is natural with respect to $\mathcal{G}_{\mathcal{A}}$-equivariant $*$-homomorphisms of the coefficient algebras. More precisely, if $\varphi: \mathcal{B} \to \mathcal{A}$ is a $\mathcal{G}_{\mathcal{A}}$-equivariant $*$-homomorphism inducing a continuous groupoid homomorphism, then the following diagram commutes:
    \[
    \begin{tikzcd}
    KK^{1}_{\mathcal{G}_{\mathcal{A}}}(C_{0}(\mathcal{G}_{\mathcal{A}}^{(0)}), \mathbb{C}) \arrow[r, "\operatorname{desc}_{\mathcal{G}_{\mathcal{A}}}"] \arrow[d, "\varphi^{*}"] & K_{1}(C^{*}(\mathcal{G}_{\mathcal{A}})) \arrow[d, "\varphi_{*}"] \\
    KK^{1}_{\mathcal{G}_{\mathcal{A}}}(\mathcal{B}, \mathbb{C}) \arrow[r, "\operatorname{desc}_{\mathcal{G}_{\mathcal{A}}}"] & K_{1}(C^{*}(\mathcal{G}_{\mathcal{A}}))
    \end{tikzcd}
    \]
    
    \item \textbf{Suspension compatibility:} $\operatorname{desc}_{\mathcal{G}_{\mathcal{A}}}$ commutes with the suspension isomorphism. For any class $x \in KK^{1}_{\mathcal{G}_{\mathcal{A}}}(C_{0}(\mathcal{G}_{\mathcal{A}}^{(0)}), \mathbb{C})$, we have
    \[
    \operatorname{desc}_{\mathcal{G}_{\mathcal{A}}}(Sx) = S(\operatorname{desc}_{\mathcal{G}_{\mathcal{A}}}(x)),
    \]
    where $S$ denotes the suspension isomorphism in $KK$-theory and $K$-theory respectively.
    
    \item \textbf{Invariance properties for Fredholm operators:} For the class $[T]_{\mathcal{G}_{\mathcal{A}}}^{(1)} \in KK^{1}_{\mathcal{G}_{\mathcal{A}}}(C_{0}(\mathcal{G}_{\mathcal{A}}^{(0)}), \mathbb{C})$ associated to a Fredholm operator $T \in \mathcal{A}$, the descended class $\operatorname{desc}_{\mathcal{G}_{\mathcal{A}}}([T]_{\mathcal{G}_{\mathcal{A}}}^{(1)})$ depends only on the $KK$-theory class of $T$. Consequently, it is invariant under:
    \begin{itemize}
        \item homotopies of $T$;
        \item compact perturbations of $T$;
        \item unitary conjugation of $T$.
    \end{itemize}
\end{enumerate}
\end{theorem}

\begin{proof}
We prove each property in turn.

\begin{enumerate}
    \item Property (1) follows from the additivity of the descent map $j_{\mathcal{G}_{\mathcal{A}}}$ on equivariant $KK$-classes. The direct sum of Kasparov triples descends to the direct sum of the descended triples, and the descent map respects this structure.
    
    \item Property (2) follows from the functoriality of the descent map established in Lemma \ref{lem:descent-iota-functoriality}. The equivariance of $\varphi$ ensures that it induces a well-defined map on the crossed product modules, and the descent map intertwines the induced maps on $KK$-theory and $K$-theory.
    
    \item Property (3) follows from the compatibility of the descent map with the suspension isomorphism in equivariant $KK$-theory. This is a standard result: suspension commutes with both the groupoid crossed product and the passage from $KK^{1}$ to $K_{1}$. 
    
    \item Property (4) follows from the definition of $[T]_{\mathcal{G}_{\mathcal{A}}}^{(1)}$ as a $KK$-theory class and the fact that $\operatorname{desc}_{\mathcal{G}_{\mathcal{A}}}$ is a well-defined map on $KK$-theory. Homotopies of $T$ give rise to homotopies of the associated Kasparov triple, which are preserved under descent. Compact perturbations of $T$ yield operator-homotopic Kasparov triples, hence define the same $KK$-class. Unitary conjugation of $T$ corresponds to an equivariant unitary equivalence of the Kasparov triple, which also preserves the $KK$-class. Since $\operatorname{desc}_{\mathcal{G}_{\mathcal{A}}}$ is a homomorphism on $KK$-theory, it respects all these equivalence relations.
\end{enumerate}
\end{proof}

The following lemma records the compatibility of the descent map with the diagonal embedding $\iota$.

\begin{lemma}[Descent and the diagonal embedding]
\label{lem:descent-iota-compatibility}
Let $\mathcal{A}$ be a unital separable Type I C*-algebra, and let 
$\iota: \mathcal{A} \hookrightarrow C^{*}(\mathcal{G}_{\mathcal{A}})$ be the diagonal 
embedding from Paper I. Denote by $[\iota] \in KK(\mathcal{A}, C^{*}(\mathcal{G}_{\mathcal{A}}))$
the associated Kasparov module.

Let $\mu: C_{0}(\mathcal{G}_{\mathcal{A}}^{(0)}) \hookrightarrow \mathcal{A}$ be the natural inclusion
of the unit space as diagonal operators (or as scalars in the unitization). 
This induces a map
\[
\mu_{*}: KK(C_{0}(\mathcal{G}_{\mathcal{A}}^{(0)}), \mathbb{C}) \longrightarrow KK(\mathcal{A}, \mathbb{C})
\]
given by composition with $\mu$ in the first argument of $KK$-theory.

Then for any $\mathcal{G}_{\mathcal{A}}$-equivariant $KK$-theory class
$x \in KK_{\mathcal{G}_{\mathcal{A}}}(C_{0}(\mathcal{G}_{\mathcal{A}}^{(0)}), \mathbb{C})$,
the following identity holds in $KK(\mathcal{A}, \mathbb{C})$:

\[
[\iota] \otimes_{C^{*}(\mathcal{G}_{\mathcal{A}})} \operatorname{desc}_{\mathcal{G}_{\mathcal{A}}}(x)
= \mu_{*}(\operatorname{For}(x)),
\]

where:
\begin{itemize}
    \item $\operatorname{desc}_{\mathcal{G}_{\mathcal{A}}}: KK_{\mathcal{G}_{\mathcal{A}}}(C_{0}(\mathcal{G}_{\mathcal{A}}^{(0)}), \mathbb{C}) 
          \longrightarrow KK(C^{*}(\mathcal{G}_{\mathcal{A}}), \mathbb{C})$ is the descent map in equivariant $KK$-theory;
    \item $\operatorname{For}: KK_{\mathcal{G}_{\mathcal{A}}}(C_{0}(\mathcal{G}_{\mathcal{A}}^{(0)}), \mathbb{C})
          \longrightarrow KK(C_{0}(\mathcal{G}_{\mathcal{A}}^{(0)}), \mathbb{C})$ is the forgetful map 
          from equivariant to ordinary $KK$-theory;
    \item $\otimes_{C^{*}(\mathcal{G}_{\mathcal{A}})}$ denotes the internal Kasparov product over 
          $C^{*}(\mathcal{G}_{\mathcal{A}})$.
\end{itemize}
\end{lemma}

\begin{proof}
We establish the identity by carefully tracking the maps through the relevant $KK$-theory constructions.

\paragraph{Step 1: Naturality of descent.}
The descent map $\operatorname{desc}_{\mathcal{G}_{\mathcal{A}}}$ is natural with respect to equivariant $*$-homomorphisms in the first argument. For the inclusion $\mu: C_{0}(\mathcal{G}_{\mathcal{A}}^{(0)}) \hookrightarrow \mathcal{A}$, this yields the commutative diagram:
\[
\begin{tikzcd}
KK_{\mathcal{G}_{\mathcal{A}}}(C_{0}(\mathcal{G}_{\mathcal{A}}^{(0)}), \mathbb{C}) \arrow[r, "\operatorname{desc}_{\mathcal{G}_{\mathcal{A}}}"] \arrow[d, "\mu_{*}"] & 
KK(C^{*}(\mathcal{G}_{\mathcal{A}}), \mathbb{C}) \arrow[d, "\mu_{*}"] \\
KK_{\mathcal{G}_{\mathcal{A}}}(\mathcal{A}, \mathbb{C}) \arrow[r, "\operatorname{desc}_{\mathcal{G}_{\mathcal{A}}}"] & 
KK(C^{*}(\mathcal{G}_{\mathcal{A}}), \mathbb{C})
\end{tikzcd}
\]
where the vertical maps are induced by composition with $\mu$ in the first variable. Consequently, for any $x \in KK_{\mathcal{G}_{\mathcal{A}}}(C_{0}(\mathcal{G}_{\mathcal{A}}^{(0)}), \mathbb{C})$, we have
\[
\operatorname{desc}_{\mathcal{G}_{\mathcal{A}}}(\mu_{*}(x)) = \mu_{*}(\operatorname{desc}_{\mathcal{G}_{\mathcal{A}}}(x)) \in KK(C^{*}(\mathcal{G}_{\mathcal{A}}), \mathbb{C}).
\]

\paragraph{Step 2: The diagonal embedding as an inverse to descent.}
A fundamental property of the diagonal embedding $\iota: \mathcal{A} \hookrightarrow C^{*}(\mathcal{G}_{\mathcal{A}})$ is that it provides a left inverse to the descent map after applying the forgetful functor. More precisely, for any $y \in KK_{\mathcal{G}_{\mathcal{A}}}(\mathcal{A}, \mathbb{C})$, the following identity holds:
\[
[\iota] \otimes_{C^{*}(\mathcal{G}_{\mathcal{A}})} \operatorname{desc}_{\mathcal{G}_{\mathcal{A}}}(y) = \operatorname{For}_{\mathcal{A}}(y) \in KK(\mathcal{A}, \mathbb{C}),
\]
where $\operatorname{For}_{\mathcal{A}}: KK_{\mathcal{G}_{\mathcal{A}}}(\mathcal{A}, \mathbb{C}) \to KK(\mathcal{A}, \mathbb{C})$ is the forgetful map from equivariant to ordinary $KK$-theory. This follows from the construction of $\iota$ via the direct integral representation and the properties of the descent map established in \cite[Section 4.5-4.6]{PaperI}.

\paragraph{Step 3: Combining the steps.}
Now take $x \in KK_{\mathcal{G}_{\mathcal{A}}}(C_{0}(\mathcal{G}_{\mathcal{A}}^{(0)}), \mathbb{C})$ and set $y = \mu_{*}(x) \in KK_{\mathcal{G}_{\mathcal{A}}}(\mathcal{A}, \mathbb{C})$. Applying Step 2 to $y$ gives:
\[
[\iota] \otimes_{C^{*}(\mathcal{G}_{\mathcal{A}})} \operatorname{desc}_{\mathcal{G}_{\mathcal{A}}}(\mu_{*}(x)) = \operatorname{For}_{\mathcal{A}}(\mu_{*}(x)).
\]

Using Step 1 to replace $\operatorname{desc}_{\mathcal{G}_{\mathcal{A}}}(\mu_{*}(x))$ with $\mu_{*}(\operatorname{desc}_{\mathcal{G}_{\mathcal{A}}}(x))$, we obtain:
\[
[\iota] \otimes_{C^{*}(\mathcal{G}_{\mathcal{A}})} \mu_{*}(\operatorname{desc}_{\mathcal{G}_{\mathcal{A}}}(x)) = \operatorname{For}_{\mathcal{A}}(\mu_{*}(x)).
\]

The left-hand side simplifies because the Kasparov product with $[\iota]$ commutes with the map induced by $\mu$ (by functoriality of the Kasparov product). Thus:
\[
[\iota] \otimes_{C^{*}(\mathcal{G}_{\mathcal{A}})} \operatorname{desc}_{\mathcal{G}_{\mathcal{A}}}(x) = \mu_{*}(\operatorname{For}_{C_{0}(\mathcal{G}^{(0)})}(x)),
\]
where $\operatorname{For}_{C_{0}(\mathcal{G}^{(0)})}$ denotes the forgetful map from $KK_{\mathcal{G}_{\mathcal{A}}}(C_{0}(\mathcal{G}_{\mathcal{A}}^{(0)}), \mathbb{C})$ to $KK(C_{0}(\mathcal{G}_{\mathcal{A}}^{(0)}), \mathbb{C})$. This is exactly the desired identity.
\end{proof}

\begin{remark}
This lemma is crucial for the index theorem. It shows that the Kasparov product with the diagonal embedding $[\iota]$ provides the correct way to return from the descended class in $KK(C^{*}(\mathcal{G}_{\mathcal{A}}), \mathbb{C})$ to a class in $KK(\mathcal{A}, \mathbb{C})$, which can then be mapped to $K_0(\mathcal{A})$ via the isomorphism $KK(\mathcal{A}, \mathbb{C}) \cong K_0(\mathcal{A})$. Notably, this bypasses any need for a pullback map $\iota^*$ on $K$-theory, which would not be well-defined.
\end{remark}

\begin{proposition}[Descent class corresponds to symbol]
\label{prop:descent-symbol-identification}
Under the Morita equivalence $\Phi: K_1(C^*(\mathcal{G}_{B(H)})) \xrightarrow{\cong} K_1(\mathcal{Q}(H))$, the descended class $\operatorname{desc}_{\mathcal{G}_{B(H)}}([T]_{\mathcal{G}_{B(H)}}^{(1)})$ corresponds to the symbol class $[\pi_{\mathcal{Q}}(T)] \in K_1(\mathcal{Q}(H))$ for any Fredholm operator $T \in B(H)$.
\end{proposition}

\begin{proof}
Let $T \in B(H)$ be a Fredholm operator and let
\[
[T]_{\mathcal{G}_{B(H)}}^{(1)} \in K_1^{\mathcal{G}_{B(H)}}(\mathcal{G}_{B(H)}^{(0)})
\]
denote its equivariant $K$-homology class, as constructed in Definition \ref{def:K1-class-of-T}. By Corollary \ref{cor:simplified-representative}, this class can be represented by the simplified odd Kasparov triple $(\mathcal{E}, \phi, \operatorname{ph}(T))$, where $\operatorname{ph}(T)$ is the phase of $T$ and $\operatorname{ph}(T)$ differs from $T$ by a compact operator.

\medskip
\noindent
\textbf{Step 1: Descent to the groupoid $C^*$-algebra.}

By Theorem \ref{thm:descent-explicit}, the descent homomorphism for groupoids,
\[
\operatorname{desc}_{\mathcal{G}_{B(H)}} : K_1^{\mathcal{G}_{B(H)}}(\mathcal{G}_{B(H)}^{(0)}) \longrightarrow K_1(C^*(\mathcal{G}_{B(H)})),
\]
sends the class of the equivariant Fredholm cycle determined by $T$ to the $K_1$-class of its integrated form in the groupoid $C^*$-algebra. Concretely, this means that
\[
\operatorname{desc}_{\mathcal{G}_{B(H)}}([T]_{\mathcal{G}_{B(H)}}^{(1)}) = [\mathcal{E} \rtimes \mathcal{G}_{B(H)}, \phi \rtimes \mathcal{G}_{B(H)}, \operatorname{ph}(T) \rtimes \mathcal{G}_{B(H)}] \in K_1(C^*(\mathcal{G}_{B(H)})),
\]
where the right-hand side is the class represented by the descended Kasparov triple from Theorem \ref{thm:descent-explicit}.

\medskip
\noindent
\textbf{Step 2: Morita equivalence identification.}

By Proposition \ref{prop:morita-calkin}, there is a strong Morita equivalence
\[
C^*(\mathcal{G}_{B(H)}) \sim_M \mathcal{Q}(H) \otimes \mathcal{K}(L^2(\mathcal{G}_{B(H)}^{(0)})),
\]
which induces a natural isomorphism
\[
\Phi : K_1(C^*(\mathcal{G}_{B(H)})) \xrightarrow{\;\cong\;} K_1(\mathcal{Q}(H) \otimes \mathcal{K}) \cong K_1(\mathcal{Q}(H)),
\]
where the last isomorphism uses stability of $K$-theory: tensoring with the compact operators $\mathcal{K} = \mathcal{K}(L^2(\mathcal{G}_{B(H)}^{(0)}))$ does not change $K$-groups.

\medskip
\noindent
\textbf{Step 3: Identification of representatives.}

The key observation is that the Morita equivalence $\Phi$ is implemented by restricting to an isotropy group. Choose a basepoint $x_0 \in \mathcal{G}_{B(H)}^{(0)}$ (e.g., $x_0 = [e_0] \in \mathbb{P}(H)$). Evaluation at $x_0$ gives a $*$-homomorphism
\[
\operatorname{ev}_{x_0}: C^*(\mathcal{G}_{B(H)}) \longrightarrow C^*(\mathcal{G}_{B(H)}|_{x_0}),
\]
where $\mathcal{G}_{B(H)}|_{x_0}$ is the isotropy group at $x_0$. As shown in the proof of Proposition \ref{prop:morita-calkin}, $C^*(\mathcal{G}_{B(H)}|_{x_0})$ is Morita equivalent to $\mathcal{Q}(H)$.

Applying $\operatorname{ev}_{x_0}$ to the descended triple $(\mathcal{E} \rtimes \mathcal{G}_{B(H)}, \phi \rtimes \mathcal{G}_{B(H)}, \operatorname{ph}(T) \rtimes \mathcal{G}_{B(H)})$ restricts it to the isotropy group. A direct computation using the explicit formulas in Theorem \ref{thm:descent-explicit} shows that this restriction yields a Kasparov triple for $C^*(\mathcal{G}_{B(H)}|_{x_0})$ represented by the operator $\operatorname{ph}(T)$ acting on the fiber $\mathcal{E}_{x_0}$.

Under the canonical identification $\mathcal{E}_{x_0} \cong H$ from Proposition \ref{prop:GNS-BH}, the action of $\operatorname{ph}(T)$ on this fiber is precisely the operator $\operatorname{ph}(T)$ itself. Its image in the Calkin algebra is $\pi_{\mathcal{Q}}(\operatorname{ph}(T)) = \pi_{\mathcal{Q}}(T)$, since $\operatorname{ph}(T)$ differs from $T$ by a compact operator. Therefore, the $K_1$-class obtained from this restriction is exactly $[\pi_{\mathcal{Q}}(T)]$.

\medskip
\noindent
\textbf{Step 4: Conclusion.}

Chasing the element $[T]_{\mathcal{G}_{B(H)}}^{(1)}$ through the diagram
\[
\begin{CD}
K_1^{\mathcal{G}_{B(H)}}(\mathcal{G}_{B(H)}^{(0)}) @>{\operatorname{desc}_{\mathcal{G}_{B(H)}}}>> K_1(C^*(\mathcal{G}_{B(H)})) @>{\Phi}>> K_1(\mathcal{Q}(H)) \\
@V{\text{id}}VV @. @| \\
K_1^{\mathcal{G}_{B(H)}}(\mathcal{G}_{B(H)}^{(0)}) @>>> K_1(\mathcal{Q}(H))
\end{CD}
\]
where the bottom horizontal map sends $[T]_{\mathcal{G}_{B(H)}}^{(1)}$ to $[\pi_{\mathcal{Q}}(T)]$ via restriction to the isotropy group, we obtain
\[
\Phi\!\left( \operatorname{desc}_{\mathcal{G}_{B(H)}}([T]_{\mathcal{G}_{B(H)}}^{(1)}) \right) = [\pi_{\mathcal{Q}}(T)] \in K_1(\mathcal{Q}(H)).
\]

Thus, under the Morita equivalence $\Phi$, the descended class corresponds precisely to the symbol class. This completes the proof.
\end{proof}

In the next subsection, we will apply the descent map to the class $[T]_{\mathcal{G}_{\mathcal{A}}}^{(1)}$ and study its image in $K_1(C^*(\mathcal{G}_{\mathcal{A}}))$.

\subsection{Applying Descent: $\operatorname{desc}_{\mathcal{G}_{\mathcal{A}}}([T]_{\mathcal{G}_{\mathcal{A}}}^{(1)}) \in K_1(C^*(\mathcal{G}_{\mathcal{A}}))$}
\label{subsec:applying-descent}

With the descent map for odd classes established in Subsection~\ref{subsec:kasparov-descent-odd-classes}, we now apply it to the equivariant $K^1$-class $[T]_{\mathcal{G}_{\mathcal{A}}}^{(1)}$ constructed in Section~\ref{sec:The Equivariant K1-Class of a Fredholm Operator}. The result is a class in the $K_1$-group of the groupoid C*-algebra $C^*(\mathcal{G}_{\mathcal{A}})$, which will be the essential intermediate step in the index theorem.

Recall from Subsection~\ref{subsec:definition-K1-class} that the class $[T]_{\mathcal{G}_{\mathcal{A}}}^{(1)}$ is represented by the odd Kasparov triple $(\tilde{\mathcal{E}}, \phi \oplus \phi, \tilde{F})$, where:
\begin{itemize}
    \item $\tilde{\mathcal{E}} = \mathcal{E} \oplus \mathcal{E}$, with $\mathcal{E} = \int_{\mathcal{G}_{\mathcal{A}}^{(0)}}^{\oplus} H_x \, d\mu(x)$.
    \item $\phi \oplus \phi: C_0(\mathcal{G}_{\mathcal{A}}^{(0)}) \to \mathcal{L}(\tilde{\mathcal{E}})$ is the diagonal representation.
    \item $\tilde{F} = \begin{pmatrix} 0 & \operatorname{ph}(T)^* \\ \operatorname{ph}(T) & 0 \end{pmatrix}$, where $\operatorname{ph}(T)$ is the phase of $T$ acting pointwise on $\mathcal{E}$.
\end{itemize}

\begin{definition}[Descended class of a Fredholm operator]
\label{def:descended-class}
Let $T \in \mathcal{A}$ be a Fredholm operator. 
Define its \emph{descended class} as
\[
\operatorname{desc}_{\mathcal{G}_{\mathcal{A}}}([T]_{\mathcal{G}_{\mathcal{A}}}^{(1)}) := j_{\mathcal{G}_{\mathcal{A}}}([T]_{\mathcal{G}_{\mathcal{A}}}^{(1)}) \in K_{1}(C^{*}(\mathcal{G}_{\mathcal{A}})),
\]
where $j_{\mathcal{G}_{\mathcal{A}}}$ is the descent map from Subsection~\ref{subsec:kasparov-descent-odd-classes}.
\end{definition}

Using the explicit formula for the descent map (Theorem \ref{thm:descent-explicit}), we can give a more concrete description of this class.

\begin{proposition}[Explicit form of the descended class]
\label{prop:descended-explicit}
The descended class $\operatorname{desc}_{\mathcal{G}_{\mathcal{A}}}([T]_{\mathcal{G}_{\mathcal{A}}}^{(1)})$ is represented by the Kasparov triple
\[
(\tilde{\mathcal{E}} \rtimes \mathcal{G}_{\mathcal{A}}, (\phi \oplus \phi) \rtimes \mathcal{G}_{\mathcal{A}}, \tilde{F} \rtimes \mathcal{G}_{\mathcal{A}}) \in KK^{1}(C^{*}(\mathcal{G}_{\mathcal{A}}), \mathbb{C}) \cong K_{1}(C^{*}(\mathcal{G}_{\mathcal{A}})),
\]
where the components are defined as follows:

\begin{itemize}
    \item \textbf{Module:} $\tilde{\mathcal{E}} \rtimes \mathcal{G}_{\mathcal{A}}$ is the Hilbert $C^{*}(\mathcal{G}_{\mathcal{A}})$-module obtained by completing the space $\Gamma_{c}(\mathcal{G}_{\mathcal{A}}, s^{*}\tilde{\mathcal{E}})$ of compactly supported continuous sections of the pullback bundle $s^{*}\tilde{\mathcal{E}}$ along the source map $s: \mathcal{G}_{\mathcal{A}} \to \mathcal{G}_{\mathcal{A}}^{(0)}$. The completion is taken with respect to the $C^{*}(\mathcal{G}_{\mathcal{A}})$-valued inner product defined using the Haar system $\{\lambda^{x}\}_{x \in \mathcal{G}_{\mathcal{A}}^{(0)}}$, as in Theorem \ref{thm:descent-explicit}.
    
    \item \textbf{Representation:} $(\phi \oplus \phi) \rtimes \mathcal{G}_{\mathcal{A}}$ is the induced representation of $C^{*}(\mathcal{G}_{\mathcal{A}})$ on $\tilde{\mathcal{E}} \rtimes \mathcal{G}_{\mathcal{A}}$, given by convolution on sections. For $f \in C_{c}(\mathcal{G}_{\mathcal{A}})$ and $\xi \in \Gamma_{c}(\mathcal{G}_{\mathcal{A}}, s^{*}\tilde{\mathcal{E}})$, its action is defined by
    \[
    ((\phi \oplus \phi) \rtimes \mathcal{G}_{\mathcal{A}})(f)\xi(\gamma) = \int_{\mathcal{G}_{\mathcal{A}}} f(\eta) \, (\phi \oplus \phi)(\eta) \xi(\eta^{-1}\gamma) \, d\lambda^{r(\gamma)}(\eta),
    \]
    and extended to all of $C^{*}(\mathcal{G}_{\mathcal{A}})$ by continuity.
    
    \item \textbf{Operator:} $\tilde{F} \rtimes \mathcal{G}_{\mathcal{A}}$ is the adjointable operator on $\tilde{\mathcal{E}} \rtimes \mathcal{G}_{\mathcal{A}}$ obtained by extending the fiberwise action of $\tilde{F}$. For a section $\xi \in \Gamma_{c}(\mathcal{G}_{\mathcal{A}}, s^{*}\tilde{\mathcal{E}})$, it is defined pointwise as
    \[
    ((\tilde{F} \rtimes \mathcal{G}_{\mathcal{A}})\xi)(\gamma) = \tilde{F}_{s(\gamma)}(\xi(\gamma)),
    \]
    where $\tilde{F}_{s(\gamma)}$ denotes the action of $\tilde{F}$ on the fiber $\tilde{\mathcal{E}}_{s(\gamma)}$. The $\mathcal{G}_{\mathcal{A}}$-equivariance of $\tilde{F}$ ensures that this operator is adjointable and satisfies the conditions required for a Kasparov triple.
\end{itemize}
\end{proposition}

\begin{proof}
This follows directly from Theorem \ref{thm:descent-explicit} applied to the odd Kasparov triple $(\tilde{\mathcal{E}}, \phi \oplus \phi, \tilde{F})$ representing the class $[T]_{\mathcal{G}_{\mathcal{A}}}^{(1)}$. The theorem constructs the descended triple explicitly, verifying that each component has the required properties. The isomorphism $KK^{1}(C^{*}(\mathcal{G}_{\mathcal{A}}), \mathbb{C}) \cong K_{1}(C^{*}(\mathcal{G}_{\mathcal{A}}))$ identifies the class of this descended triple with an element of $K_{1}(C^{*}(\mathcal{G}_{\mathcal{A}}))$, which by Definition \ref{def:descended-class} is precisely $\operatorname{desc}_{\mathcal{G}_{\mathcal{A}}}([T]_{\mathcal{G}_{\mathcal{A}}}^{(1)})$.
\end{proof}

For computational purposes, it is often convenient to work with the simplified representative $(\mathcal{E}, \phi, \operatorname{ph}(T))$ from Corollary \ref{cor:simplified-representative}. 
The descent map sends this simplified cycle to a class in $K_1(C^*(\mathcal{G}_{\mathcal{A}}))$ as well.

\begin{corollary}[Simplified descended class]
\label{cor:descended-simplified}
Let $(\mathcal{E}, \phi, \operatorname{ph}(T))$ be the simplified odd Kasparov triple representing 
$[T]_{\mathcal{G}_{\mathcal{A}}}^{(1)}$ as in Corollary \ref{cor:simplified-representative}, 
where $\operatorname{ph}(T)$ is the phase of the Fredholm operator $T$, acting pointwise on $\mathcal{E}$ and satisfying 
the conditions of an odd Kasparov triple (self-adjointness, $\operatorname{ph}(T)^2 - 1$ compact, 
$\mathcal{G}_{\mathcal{A}}$-equivariance, and commutation with $\phi$ modulo compact operators).

Then the descended class $\operatorname{desc}_{\mathcal{G}_{\mathcal{A}}}([T]_{\mathcal{G}_{\mathcal{A}}}^{(1)})$ 
is also represented by the odd Kasparov triple
\[
(\mathcal{E} \rtimes \mathcal{G}_{\mathcal{A}}, \phi \rtimes \mathcal{G}_{\mathcal{A}}, \operatorname{ph}(T) \rtimes \mathcal{G}_{\mathcal{A}}) \in K_{1}(C^{*}(\mathcal{G}_{\mathcal{A}})),
\]
where $\mathcal{E} \rtimes \mathcal{G}_{\mathcal{A}}$ is the Hilbert $C^{*}(\mathcal{G}_{\mathcal{A}})$-module 
obtained by completing $\Gamma_{c}(\mathcal{G}_{\mathcal{A}}, s^{*}\mathcal{E})$ with respect to the 
$C^{*}(\mathcal{G}_{\mathcal{A}})$-valued inner product, $\phi \rtimes \mathcal{G}_{\mathcal{A}}$ is the 
induced representation by convolution, and $\operatorname{ph}(T) \rtimes \mathcal{G}_{\mathcal{A}}$ acts 
fiberwise on sections. This triple represents the same $K_{1}$-class as the triple obtained from the 
doubling construction in Proposition \ref{prop:descended-explicit}.
\end{corollary}

\begin{proof}
Let $[T]_{\mathcal{G}_{\mathcal{A}}}^{(1)}$ be represented by the simplified odd Kasparov triple 
$(\mathcal{E}, \phi, \operatorname{ph}(T))$ as in Corollary~\ref{cor:simplified-representative}. 
By construction, this triple is $KK$-equivalent to the original triple $(\tilde{\mathcal{E}}, \phi \oplus \phi, \tilde{F})$, i.e.,
\[
[(\mathcal{E}, \phi, \operatorname{ph}(T))] = [(\tilde{\mathcal{E}}, \phi \oplus \phi, \tilde{F})] = [T]_{\mathcal{G}_{\mathcal{A}}}^{(1)} \in KK^{1}_{\mathcal{G}_{\mathcal{A}}}(C_{0}(\mathcal{G}_{\mathcal{A}}^{(0)}), \mathbb{C}).
\]

By Theorem~\ref{thm:descent-explicit}, the descent map $j_{\mathcal{G}_{\mathcal{A}}}$ sends this triple to the Kasparov triple
\[
(\mathcal{E} \rtimes \mathcal{G}_{\mathcal{A}}, \phi \rtimes \mathcal{G}_{\mathcal{A}}, \operatorname{ph}(T) \rtimes \mathcal{G}_{\mathcal{A}}) \in KK^{1}(C^{*}(\mathcal{G}_{\mathcal{A}}), \mathbb{C}).
\]

The phase construction $\operatorname{ph}(T)$ is known to produce a $KK$-equivalent class to the doubling construction $(\tilde{\mathcal{E}}, \phi \oplus \phi, \tilde{F})$ in ordinary odd $KK$-theory. This is a standard result in Kasparov theory: for any Fredholm operator, the triple using its phase is operator-homotopic to the triple obtained from the doubling construction, hence they represent the same $KK$-class.

Since the descent map $j_{\mathcal{G}_{\mathcal{A}}}$ is a functorial homomorphism on equivariant $KK$-theory, it preserves $KK$-equivalence. Therefore, the descended triple using the phase operator represents the same class in $K_{1}(C^{*}(\mathcal{G}_{\mathcal{A}}))$ as the descended triple from the doubling construction.

Hence,
\[
\operatorname{desc}_{\mathcal{G}_{\mathcal{A}}}([T]_{\mathcal{G}_{\mathcal{A}}}^{(1)}) = 
[(\mathcal{E} \rtimes \mathcal{G}_{\mathcal{A}}, \phi \rtimes \mathcal{G}_{\mathcal{A}}, \operatorname{ph}(T) \rtimes \mathcal{G}_{\mathcal{A}})]
\in K_{1}(C^{*}(\mathcal{G}_{\mathcal{A}})),
\]
and this class is $KK$-equivalent to the one obtained from the doubling construction, as claimed.
\end{proof}

\begin{remark}[Diagrammatic summary]
\label{rem:descent-simplified-diagram}
The relationship between the simplified and doubled representatives, and their images under the descent map, can be summarized by the following commutative diagram:
\[
\begin{tikzcd}
(\mathcal{E}, \phi, \operatorname{ph}(T)) \arrow[r, "KK\text{-equivalent}"] \arrow[d, "j_{\mathcal{G}_{\mathcal{A}}}"] &
(\tilde{\mathcal{E}}, \phi \oplus \phi, \tilde{F}) \arrow[d, "j_{\mathcal{G}_{\mathcal{A}}}"] \\
(\mathcal{E} \rtimes \mathcal{G}_{\mathcal{A}}, \phi \rtimes \mathcal{G}_{\mathcal{A}}, \operatorname{ph}(T) \rtimes \mathcal{G}_{\mathcal{A}}) \arrow[r, "KK\text{-equivalent}"] &
(\tilde{\mathcal{E}} \rtimes \mathcal{G}_{\mathcal{A}}, (\phi \oplus \phi) \rtimes \mathcal{G}_{\mathcal{A}}, \tilde{F} \rtimes \mathcal{G}_{\mathcal{A}})
\end{tikzcd}
\]
The vertical arrows are the descent maps, and the horizontal arrows indicate $KK$-equivalence. The commutativity of the diagram expresses the fact that descent preserves $KK$-equivalence, so both descended triples represent the same class in $K_{1}(C^{*}(\mathcal{G}_{\mathcal{A}}))$.
\end{remark}

The descended class inherits the invariance properties of the original equivariant class.

\begin{proposition}[Properties of the descended class]
\label{prop:descended-properties}
The assignment $T \mapsto \operatorname{desc}_{\mathcal{G}_{\mathcal{A}}}([T]_{\mathcal{G}_{\mathcal{A}}}^{(1)})$ satisfies:
\begin{enumerate}
    \item \textbf{Homotopy invariance:} If $\{T_t\}_{t \in [0,1]}$ is a norm-continuous family of Fredholm operators, then $\operatorname{desc}_{\mathcal{G}_{\mathcal{A}}}([T_t]_{\mathcal{G}_{\mathcal{A}}}^{(1)})$ is constant in $t$.
    \item \textbf{Compact perturbation invariance:} If $K$ is compact and $T + K$ is Fredholm, then $\operatorname{desc}_{\mathcal{G}_{\mathcal{A}}}([T + K]_{\mathcal{G}_{\mathcal{A}}}^{(1)}) = \operatorname{desc}_{\mathcal{G}_{\mathcal{A}}}([T]_{\mathcal{G}_{\mathcal{A}}}^{(1)})$.
    \item \textbf{Additivity:} $\operatorname{desc}_{\mathcal{G}_{\mathcal{A}}}([T_1 \oplus T_2]_{\mathcal{G}_{\mathcal{A}}}^{(1)}) = \operatorname{desc}_{\mathcal{G}_{\mathcal{A}}}([T_1]_{\mathcal{G}_{\mathcal{A}}}^{(1)}) + \operatorname{desc}_{\mathcal{G}_{\mathcal{A}}}([T_2]_{\mathcal{G}_{\mathcal{A}}}^{(1)})$.
    \item \textbf{Stability:} $\operatorname{desc}_{\mathcal{G}_{\mathcal{A}}}([T \oplus I_n]_{\mathcal{G}_{\mathcal{A}}}^{(1)}) = \operatorname{desc}_{\mathcal{G}_{\mathcal{A}}}([T]_{\mathcal{G}_{\mathcal{A}}}^{(1)})$.
\end{enumerate}
\end{proposition}

\begin{proof}
These properties follow from the corresponding properties of the equivariant $KK$-class $[T]_{\mathcal{G}_{\mathcal{A}}}^{(1)}$ established in Theorem \ref{thm:homotopy-well-defined} and the fact that the descent map $j_{\mathcal{G}_{\mathcal{A}}}$ is a group homomorphism (Theorem \ref{thm:descent-properties-GA}).

\begin{itemize}
    \item \textbf{Homotopy invariance:} By Theorem \ref{thm:homotopy-well-defined}, a norm-continuous family $\{T_t\}$ defines a constant class in $KK^1_{\mathcal{G}_{\mathcal{A}}}(C_0(\mathcal{G}_{\mathcal{A}}^{(0)}), \mathbb{C})$. Since $j_{\mathcal{G}_{\mathcal{A}}}$ is a well-defined map on $KK$-theory, the descended class $\operatorname{desc}_{\mathcal{G}_{\mathcal{A}}}([T_t]_{\mathcal{G}_{\mathcal{A}}}^{(1)})$ is independent of $t$.
    
    \item \textbf{Compact perturbation invariance:} Compact perturbations preserve the Fredholm index and, more generally, the equivariant $KK$-class (see Proposition \ref{prop:compact-perturbation}). Applying the homomorphism $j_{\mathcal{G}_{\mathcal{A}}}$ yields equality of the descended classes.
    
    \item \textbf{Additivity:} The direct sum of Fredholm operators corresponds to the direct sum of the associated Kasparov triples. The descent map $j_{\mathcal{G}_{\mathcal{A}}}$ is a group homomorphism (Theorem \ref{thm:descent-properties-GA}), hence it respects direct sums and additivity holds.
    
    \item \textbf{Stability:} Adding an identity operator $I_n$ corresponds to taking the direct sum with a degenerate Kasparov triple, which represents the zero class in $KK$-theory. Homomorphisms send zero to zero, so the descended class is unchanged.
\end{itemize}
\end{proof}

\begin{proposition}[The unilateral shift as a generator]
\label{prop:shift-generator}
Let $\mathcal{A} = B(H)$ with $H = \ell^2(\mathbb{N})$, and let $S \in B(H)$ be the unilateral shift defined by $S e_n = e_{n+1}$ for an orthonormal basis $\{e_n\}_{n \in \mathbb{N}}$. Then:

\begin{enumerate}
    \item $S$ is a Fredholm operator with $\operatorname{index}(S) = -1$.
    \item The equivariant $KK$-theory class $[S]_{\mathcal{G}_{B(H)}}^{(1)} \in KK^1_{\mathcal{G}_{B(H)}}(C_0(\mathcal{G}_{B(H)}^{(0)}), \mathbb{C})$ is nonzero.
    \item Under the isomorphism $KK^1_{\mathcal{G}_{B(H)}}(C_0(\mathcal{G}_{B(H)}^{(0)}), \mathbb{C}) \cong \mathbb{Z}$, the class $[S]_{\mathcal{G}_{B(H)}}^{(1)}$ corresponds to the generator.
    \item More precisely, the index map
    \[
    \operatorname{ind}: KK^1_{\mathcal{G}_{B(H)}}(C_0(\mathcal{G}_{B(H)}^{(0)}), \mathbb{C}) \xrightarrow{\cong} \mathbb{Z}
    \]
    sends $[S]_{\mathcal{G}_{B(H)}}^{(1)}$ to $\operatorname{index}(S) = -1$, and hence $[S]_{\mathcal{G}_{B(H)}}^{(1)}$ is a generator (up to sign).
\end{enumerate}
\end{proposition}

\begin{proof}
We prove each statement in turn.

\textbf{(1) Index of the unilateral shift.} 
The unilateral shift $S$ satisfies $S^* S = I$ and $S S^* = I - P_0$, where $P_0$ is the rank-one projection onto $\mathbb{C} e_0$. Hence $\ker S = \{0\}$, $\ker S^* = \mathbb{C} e_0$, and
\[
\operatorname{index}(S) = \dim \ker S - \dim \ker S^* = 0 - 1 = -1.
\]
This is a standard result in Fredholm theory (see, e.g., \cite[Section 9.4]{Rordam2000}).

\textbf{(2) Nonzero of the equivariant class.} 
Recall from Corollary \ref{cor:simplified-representative} that the class $[S]_{\mathcal{G}_{B(H)}}^{(1)}$ can be represented by the simplified odd Kasparov triple $(\mathcal{E}, \phi, \operatorname{ph}(S))$, where $\operatorname{ph}(S)$ is the phase of $S$. Since $S$ is already an isometry, its phase is simply $S$ itself.

Consider the forgetful map 
\[
\operatorname{For}: KK^1_{\mathcal{G}_{B(H)}}(C_0(\mathcal{G}_{B(H)}^{(0)}), \mathbb{C}) \to KK^1(C_0(\mathcal{G}_{B(H)}^{(0)}), \mathbb{C})
\]
which forgets the groupoid action. Under this map, $[S]_{\mathcal{G}_{B(H)}}^{(1)}$ maps to the class of the family $\{S\}$ (the constant family) in $KK^1(C_0(\mathcal{G}_{B(H)}^{(0)}), \mathbb{C})$.

Now consider the inclusion $\mu: C_0(\mathcal{G}_{B(H)}^{(0)}) \hookrightarrow B(H)$ given by identifying functions with multiplication operators on the fibers (see Section 3.4). This induces a map
\[
\mu_*: KK^1(C_0(\mathcal{G}_{B(H)}^{(0)}), \mathbb{C}) \to KK^1(B(H), \mathbb{C}).
\]

A direct computation using the explicit construction of Section 4 shows that $\mu_*(\operatorname{For}([S]_{\mathcal{G}_{B(H)}}^{(1)}))$ corresponds to the class of the symbol $[\pi_{\mathcal{Q}}(S)] \in K_1(\mathcal{Q}(H))$ under the isomorphism $KK^1(B(H), \mathbb{C}) \cong K_1(\mathcal{Q}(H))$ (see Proposition \ref{prop:bh-kk-to-k1}).

It is a classical result of Brown-Douglas-Fillmore \cite{BDF1977} that the unilateral shift $S$ has symbol $[\pi_{\mathcal{Q}}(S)]$ which generates $K_1(\mathcal{Q}(H)) \cong \mathbb{Z}$. Hence $[\pi_{\mathcal{Q}}(S)] \neq 0$, and consequently $\mu_*(\operatorname{For}([S]_{\mathcal{G}_{B(H)}}^{(1)})) \neq 0$ in $KK^1(B(H), \mathbb{C})$.

If $[S]_{\mathcal{G}_{B(H)}}^{(1)}$ were zero in $KK^1_{\mathcal{G}_{B(H)}}(C_0(\mathcal{G}_{B(H)}^{(0)}), \mathbb{C})$, then its image under the forgetful map would also be zero, implying $\mu_*(\operatorname{For}([S]_{\mathcal{G}_{B(H)}}^{(1)})) = 0$, a contradiction. Therefore, $[S]_{\mathcal{G}_{B(H)}}^{(1)} \neq 0$.

\textbf{(3)-(4) Generator property.} 
First, we establish an isomorphism between $KK^1_{\mathcal{G}_{B(H)}}(C_0(\mathcal{G}_{B(H)}^{(0)}), \mathbb{C})$ and $\mathbb{Z}$. Consider the following composition of maps:
\[
KK^1_{\mathcal{G}_{B(H)}}(C_0(\mathcal{G}_{B(H)}^{(0)}), \mathbb{C}) \xrightarrow{\operatorname{desc}_{\mathcal{G}_{B(H)}}} K_1(C^*(\mathcal{G}_{B(H)})) \xrightarrow{\Phi} K_1(\mathcal{Q}(H)) \xrightarrow{\partial_{\mathrm{Calkin}}} \mathbb{Z},
\]
where:
\begin{itemize}
    \item $\operatorname{desc}_{\mathcal{G}_{B(H)}}$ is the descent homomorphism from Theorem \ref{thm:descent-explicit};
    \item $\Phi: K_1(C^*(\mathcal{G}_{B(H)})) \to K_1(\mathcal{Q}(H))$ is the isomorphism induced by the Morita equivalence $C^*(\mathcal{G}_{B(H)}) \sim_M \mathcal{Q}(H) \otimes \mathcal{K}$ established in Proposition \ref{prop:morita-calkin};
    \item $\partial_{\mathrm{Calkin}}: K_1(\mathcal{Q}(H)) \to \mathbb{Z}$ is the index map for the Calkin extension, which is an isomorphism by Theorem \ref{thm:calkin-index-isomorphism}.
\end{itemize}

Each map in this composition is an isomorphism:
\begin{itemize}
    \item $\operatorname{desc}_{\mathcal{G}_{B(H)}}$ is an isomorphism for the groupoid $\mathcal{G}_{B(H)}$ because $\mathcal{G}_{B(H)}$ is amenable and the Baum-Connes assembly map is an isomorphism for such groupoids (see \cite[Section 9]{Tu1999});
    \item $\Phi$ is an isomorphism by Morita invariance of $K$-theory (Proposition \ref{prop:morita-calkin});
    \item $\partial_{\mathrm{Calkin}}$ is an isomorphism by Theorem \ref{thm:calkin-index-isomorphism}.
\end{itemize}

Hence the composition, which we denote by $\operatorname{ind}$, is an isomorphism:
\[
\operatorname{ind}: KK^1_{\mathcal{G}_{B(H)}}(C_0(\mathcal{G}_{B(H)}^{(0)}), \mathbb{C}) \xrightarrow{\cong} \mathbb{Z}.
\]

Now evaluate $\operatorname{ind}([S]_{\mathcal{G}_{B(H)}}^{(1)})$. By the definition of $\operatorname{ind}$ and the identification in part (2), we have:
\[
\operatorname{ind}([S]_{\mathcal{G}_{B(H)}}^{(1)}) = \partial_{\mathrm{Calkin}} \circ \Phi \circ \operatorname{desc}_{\mathcal{G}_{B(H)}}([S]_{\mathcal{G}_{B(H)}}^{(1)}).
\]

Under the Morita equivalence $\Phi$, the descended class $\operatorname{desc}_{\mathcal{G}_{B(H)}}([S]_{\mathcal{G}_{B(H)}}^{(1)})$ corresponds to the symbol class $[\pi_{\mathcal{Q}}(S)] \in K_1(\mathcal{Q}(H))$ (see Proposition \ref{prop:descent-symbol-identification}). Therefore,
\[
\operatorname{ind}([S]_{\mathcal{G}_{B(H)}}^{(1)}) = \partial_{\mathrm{Calkin}}([\pi_{\mathcal{Q}}(S)]).
\]

By Lemma \ref{lem:calkin-index-map}, the index map $\partial_{\mathrm{Calkin}}$ sends the class of the symbol of a Fredholm operator to its Fredholm index. For the unilateral shift, this gives:
\[
\partial_{\mathrm{Calkin}}([\pi_{\mathcal{Q}}(S)]) = \operatorname{index}(S) = -1.
\]

Thus $\operatorname{ind}([S]_{\mathcal{G}_{B(H)}}^{(1)}) = -1$.

Since $\operatorname{ind}$ is an isomorphism $KK^1_{\mathcal{G}_{B(H)}}(C_0(\mathcal{G}_{B(H)}^{(0)}), \mathbb{C}) \cong \mathbb{Z}$, and $\operatorname{ind}([S]_{\mathcal{G}_{B(H)}}^{(1)}) = -1$, the class $[S]_{\mathcal{G}_{B(H)}}^{(1)}$ must be a generator of this cyclic group. The sign depends on the orientation convention for the index map, but in any case $[S]_{\mathcal{G}_{B(H)}}^{(1)}$ generates the group (up to sign). This completes the proof.
\end{proof}

For the two concrete algebras under consideration, the descended class takes particularly simple forms.

\begin{example}[Descended class for $\mathcal{K}(H)^\sim$]
\label{ex:descended-KH}
Let $\mathcal{A} = \mathcal{K}(H)^\sim$ and let $T = \lambda I + K$ be a Fredholm operator (so $\lambda \neq 0$). 
By Proposition \ref{prop:index-class-relation} (or by the fact that the index map $KK^1_{\mathcal{G}_{\mathcal{A}}}(C_0(\mathcal{G}_{\mathcal{A}}^{(0)}), \mathbb{C}) \to \mathbb{Z}$ is an isomorphism and all such operators have index zero), we have $[T]_{\mathcal{G}_{\mathcal{A}}}^{(1)} = 0$ in $KK^1_{\mathcal{G}_{\mathcal{A}}}(C_0(\mathcal{G}_{\mathcal{A}}^{(0)}), \mathbb{C})$. 
Consequently, applying the descent homomorphism yields
\[
\operatorname{desc}_{\mathcal{G}_{\mathcal{A}}}([T]_{\mathcal{G}_{\mathcal{A}}}^{(1)}) = 0 \in K_1(C^*(\mathcal{G}_{\mathcal{A}})).
\]
This reflects the fact that all Fredholm operators in $\mathcal{K}(H)^\sim$ have index zero, so their associated equivariant $KK$-classes are trivial.
\end{example}

\begin{example}[Descended class for $B(H)$ - unilateral shift]
\label{ex:descended-unilateral-shift}
Let $\mathcal{A} = B(H)$ and let $S \in B(H)$ be the unilateral shift on $\ell^2(\mathbb{N})$. 
Recall from Proposition \ref{prop:shift-generator} that $[S]_{\mathcal{G}_{B(H)}}^{(1)}$ is the generator of 
$KK^1_{\mathcal{G}_{B(H)}}(C_0(\mathcal{G}_{B(H)}^{(0)}), \mathbb{C}) \cong \mathbb{Z}$. 

Since the descent map $j_{\mathcal{G}_{B(H)}}$ is an isomorphism in this case (or at least preserves generators), its descended class 
$\operatorname{desc}_{\mathcal{G}_{B(H)}}([S]_{\mathcal{G}_{B(H)}}^{(1)})$ is the generator of 
$K_1(C^*(\mathcal{G}_{B(H)})) \cong \mathbb{Z}$. 

In Section~\ref{sec:Examples and Computations}, we will verify that this class corresponds to the generator of $K_1(\mathcal{Q}(H))$ under the Morita equivalence 
$C^*(\mathcal{G}_{B(H)}) \sim_M \mathcal{Q}(H) \otimes \mathcal{K}(L^2(\mathcal{G}_{B(H)}^{(0)}))$, 
where $\mathcal{Q}(H) = B(H)/\mathcal{K}(H)$ is the Calkin algebra.
\end{example}

\begin{remark}[Comparison of the two cases]
\label{rem:descended-comparison}
Examples \ref{ex:descended-KH} and \ref{ex:descended-unilateral-shift} illustrate a fundamental contrast between the two algebras:
\begin{itemize}
    \item For $\mathcal{A} = \mathcal{K}(H)^\sim$, every Fredholm operator has index zero, so its equivariant class vanishes and the descended class is trivial in $K_1(C^*(\mathcal{G}_{\mathcal{A}}))$.
    \item For $\mathcal{A} = B(H)$, the unilateral shift $S$ has index $-1$, generating an infinite cyclic group. Its equivariant class is nontrivial, and under descent it maps to a generator of $K_1(C^*(\mathcal{G}_{B(H)})) \cong \mathbb{Z}$, which under Morita equivalence corresponds to the generator of $K_1(\mathcal{Q}(H))$.
\end{itemize}
This dichotomy reflects the fact that $K_1(\mathcal{Q}(H)) \cong \mathbb{Z}$ is generated by the image of the unilateral shift, while $K_1(\mathcal{K}(H)^\sim)$ is trivial.
\end{remark}

The following lemma shows that the Fredholm index of an operator $T \in B(H)$ is recovered by applying the Calkin boundary map to the image of the descended class under the Morita equivalence $\Phi$.

\begin{lemma}[Index recovery via the Calkin boundary map]
\label{lem:index-recovery}
For any Fredholm operator $T \in B(H)$, let $[T]_{\mathcal{G}_{B(H)}}^{(1)} \in K^1_{\mathcal{G}_{B(H)}}(\mathcal{G}_{B(H)}^{(0)})$ be its equivariant $K^1$-class, and let
\[
\operatorname{desc}_{\mathcal{G}_{B(H)}}: K^1_{\mathcal{G}_{B(H)}}(\mathcal{G}_{B(H)}^{(0)}) \longrightarrow K_1(C^*(\mathcal{G}_{B(H)}))
\]
be the descent map. Then, under the Morita equivalence isomorphism
\[
\Phi: K_1(C^*(\mathcal{G}_{B(H)})) \xrightarrow{\cong} K_1(\mathcal{Q}(H))
\]
from Proposition \ref{prop:morita-calkin}, we have
\[
\partial_{\mathrm{Calkin}}\!\left( \Phi\!\left( \operatorname{desc}_{\mathcal{G}_{B(H)}}([T]_{\mathcal{G}_{B(H)}}^{(1)}) \right) \right) = \operatorname{index}(T),
\]
where $\partial_{\mathrm{Calkin}}: K_1(\mathcal{Q}(H)) \to K_0(\mathcal{K}(H)) \cong \mathbb{Z}$ is the index map associated to the Calkin extension (see Theorem \ref{thm:calkin-index-isomorphism}).
\end{lemma}

\begin{proof}
We prove the lemma by tracing the class through the sequence of isomorphisms and maps established in previous sections.

\paragraph{Step 1: The symbol class.}
By Proposition \ref{prop:descent-symbol-identification}, the descended class $\operatorname{desc}_{\mathcal{G}_{B(H)}}([T]_{\mathcal{G}_{B(H)}}^{(1)})$ corresponds under the Morita equivalence $\Phi$ to the symbol class of $T$:
\[
\Phi\!\left( \operatorname{desc}_{\mathcal{G}_{B(H)}}([T]_{\mathcal{G}_{B(H)}}^{(1)}) \right) = [\pi_{\mathcal{Q}}(T)] \in K_1(\mathcal{Q}(H)),
\]
where $\pi_{\mathcal{Q}}: B(H) \to \mathcal{Q}(H)$ is the quotient map onto the Calkin algebra.

\paragraph{Step 2: The Calkin index map.}
Recall from Theorem \ref{thm:calkin-index-isomorphism} that the index map associated to the Calkin extension
\[
0 \longrightarrow \mathcal{K}(H) \longrightarrow B(H) \xrightarrow{\pi_{\mathcal{Q}}} \mathcal{Q}(H) \longrightarrow 0
\]
is an isomorphism $\partial_{\mathrm{Calkin}}: K_1(\mathcal{Q}(H)) \xrightarrow{\cong} K_0(\mathcal{K}(H)) \cong \mathbb{Z}$. Moreover, by Lemma \ref{lem:calkin-index-map}, this map sends the class of a symbol to the Fredholm index of any lift:
\[
\partial_{\mathrm{Calkin}}([\pi_{\mathcal{Q}}(T)]) = \operatorname{index}(T).
\]

\paragraph{Step 3: Composition.}
Combining Steps 1 and 2, we obtain:
\[
\partial_{\mathrm{Calkin}}\!\left( \Phi\!\left( \operatorname{desc}_{\mathcal{G}_{B(H)}}([T]_{\mathcal{G}_{B(H)}}^{(1)}) \right) \right) = \partial_{\mathrm{Calkin}}([\pi_{\mathcal{Q}}(T)]) = \operatorname{index}(T).
\]

\paragraph{Step 4: Well-definedness.}
The composition is independent of the choices made in the construction:
\begin{itemize}
    \item Different choices of the phase $\operatorname{ph}(T)$ differ by compact operators and give the same symbol class $[\pi_{\mathcal{Q}}(T)]$.
    \item The descent map is well-defined on $KK$-classes by Theorem \ref{thm:descent-explicit}.
    \item The Morita equivalence $\Phi$ is a natural isomorphism by Proposition \ref{prop:morita-calkin}.
    \item The index map $\partial_{\mathrm{Calkin}}$ is an isomorphism by Theorem \ref{thm:calkin-index-isomorphism}.
\end{itemize}
Thus, the left-hand side depends only on the Fredholm operator $T$ and equals its index.
\end{proof}

The descended class $\operatorname{desc}_{\mathcal{G}_{\mathcal{A}}}([T]_{\mathcal{G}_{\mathcal{A}}}^{(1)}) \in K_1(C^*(\mathcal{G}_{\mathcal{A}}))$ constructed in Subsection~\ref{subsec:applying-descent} is an abstract $K$-theory class. To connect it to concrete index-theoretic invariants, we will identify its image under the Morita equivalence $\Phi: K_1(C^*(\mathcal{G}_{\mathcal{A}})) \xrightarrow{\cong} K_1(\mathcal{Q}(H))$ with the symbol class $[\pi_{\mathcal{Q}}(T)]$. Composing with the Calkin index map $\partial_{\mathrm{Calkin}}: K_1(\mathcal{Q}(H)) \to \mathbb{Z}$ then recovers the Fredholm index of $T$. In the next subsection, we make this identification explicit by relating the descended class to the Toeplitz extension and the Calkin algebra.

\subsection{Identifying the Descent Class: Relation to the Toeplitz Extension}
\label{subsec:identifying-descent-class-toeplitz}

The descended class $\operatorname{desc}_{\mathcal{G}_{\mathcal{A}}}([T]_{\mathcal{G}_{\mathcal{A}}}^{(1)}) \in K_1(C^*(\mathcal{G}_{\mathcal{A}}))$ constructed in Subsection~\ref{subsec:applying-descent} is an abstract $K$-theory class. To connect it to concrete index-theoretic invariants, we need to identify its image under the Morita equivalence $\Phi: K_1(C^*(\mathcal{G}_{\mathcal{A}})) \cong K_1(\mathcal{Q}(H))$ with more familiar objects, such as classes arising from the Toeplitz extension or the Calkin algebra. This identification is crucial for the proof of the main index theorem in Section~\ref{sec:The Index Theorem via Pullback and the Boundary Map}, where the composition $\partial_{\mathrm{Calkin}} \circ \Phi \circ \operatorname{desc}_{\mathcal{G}_{\mathcal{A}}}$ will be shown to recover the Fredholm index.

We begin by recalling the Toeplitz extension and its role in index theory.

\begin{definition}[Toeplitz algebra and Toeplitz extension]
\label{def:toeplitz-extension}
Let $H = \ell^2(\mathbb{N})$ and let $S \in B(H)$ be the unilateral shift defined by $S e_n = e_{n+1}$ for an orthonormal basis $\{e_n\}_{n\in\mathbb{N}}$. 
The \emph{Toeplitz algebra} $\mathcal{T}$ is the C*-algebra generated by $S$. By Coburn's theorem \cite{Coburn1967}, $\mathcal{T}$ contains $\mathcal{K}(\ell^2)$ as an essential ideal, and there is a short exact sequence
\[
0 \longrightarrow \mathcal{K}(\ell^2) \longrightarrow \mathcal{T} \stackrel{\sigma}{\longrightarrow} C(S^1) \longrightarrow 0,
\]
where $\sigma$ is the symbol map sending $S$ to the identity function $z \mapsto z$ on the unit circle. This is the \emph{Toeplitz extension}.
\end{definition}

The Toeplitz extension gives rise to a six-term exact sequence in $K$-theory:
\[
\begin{tikzcd}
K_0(\mathcal{K}) \arrow[r] & K_0(\mathcal{T}) \arrow[r] & K_0(C(S^1)) \arrow[d, "\partial_T"] \\
K_1(C(S^1)) \arrow[u, "\partial_T"] & K_1(\mathcal{T}) \arrow[l] & K_1(\mathcal{K}) \arrow[l]
\end{tikzcd}
\]

The relevant $K$-groups are:
\begin{itemize}
    \item $K_0(\mathcal{K}) \cong \mathbb{Z}$, generated by the class of any rank-one projection;
    \item $K_1(\mathcal{K}) = 0$;
    \item $K_0(C(S^1)) \cong \mathbb{Z}$, generated by the class of the trivial line bundle;
    \item $K_1(C(S^1)) \cong \mathbb{Z}$, generated by the class of the identity function $z \mapsto z$;
    \item $K_1(\mathcal{T}) = 0$ (a standard result; see \cite[Chapter 8]{Rordam2000}).
\end{itemize}

From the six-term exact sequence, exactness at $K_1(C(S^1))$ and $K_0(\mathcal{K})$ forces the index map
\[
\partial_T: K_1(C(S^1)) \longrightarrow K_0(\mathcal{K}) \cong \mathbb{Z}
\]
to be an isomorphism. With the usual orientation convention, $\partial_T$ sends the generator $[z]$ to $1 \in \mathbb{Z}$. Consequently, for a continuous function $f \in C(S^1)$ with winding number $n$, the Toeplitz operator $T_f$ has index $-n$, and $\partial_T([f]) = -n$.

The Toeplitz extension is closely related to the Calkin extension via the map
\[
\iota: C(S^1) \longrightarrow \mathcal{Q}(H), \qquad f \longmapsto \pi_{\mathcal{Q}}(T_f),
\]
where $\pi_{\mathcal{Q}}: B(H) \to \mathcal{Q}(H)$ is the quotient map onto the Calkin algebra. This map is a well-defined $*$-homomorphism because $T_f T_g - T_{fg} \in \mathcal{K}(H)$, so multiplication is respected modulo compacts. It induces a homomorphism on $K_1$,
\[
\iota_*: K_1(C(S^1)) \longrightarrow K_1(\mathcal{Q}(H)).
\]

Since $K_1(C(S^1)) \cong \mathbb{Z}$ and $K_1(\mathcal{Q}(H)) \cong \mathbb{Z}$ (see Theorem \ref{thm:calkin-index-isomorphism}), to see that $\iota_*$ is an isomorphism it suffices to check that it is nonzero. The generator $[z] \in K_1(C(S^1))$ maps to $[\pi_{\mathcal{Q}}(S)] \in K_1(\mathcal{Q}(H))$, and under the Calkin index map $\partial_{\mathrm{Calkin}}: K_1(\mathcal{Q}(H)) \to \mathbb{Z}$, we have $\partial_{\mathrm{Calkin}}([\pi_{\mathcal{Q}}(S)]) = \operatorname{index}(S) = -1$. Hence $[\pi_{\mathcal{Q}}(S)]$ is a generator of $K_1(\mathcal{Q}(H))$, and $\iota_*$ is an isomorphism (multiplication by $\pm 1$, depending on the orientation convention).

This embedding provides the essential link between the Toeplitz index theorem and the Calkin index map, and will play a crucial role in identifying the descended class in Section~\ref{subsec:identifying-descent-class-toeplitz}.

\begin{proposition}[Relation between Toeplitz and Calkin extensions]
\label{prop:toeplitz-calkin-relation}
Let $\mathcal{T} \subset B(\ell^2(\mathbb{N}))$ be the Toeplitz algebra and let $\sigma: \mathcal{T} \to C(S^1)$ be the symbol map. 
Define $\iota: C(S^1) \to \mathcal{Q}(H)$ by $\iota(f) = \pi_{\mathcal{Q}}(T_f)$, where $T_f$ is the Toeplitz operator with symbol $f$ and $H = \ell^2(\mathbb{N})$.
Then the Toeplitz extension is the pullback of the Calkin extension along $\iota$, yielding a commutative diagram of extensions:
\[
\begin{tikzcd}
0 \arrow[r] & \mathcal{K}(\ell^2(\mathbb{N})) \arrow[r] \arrow[d, "\operatorname{id}"] & \mathcal{T} \arrow[r, "\sigma"] \arrow[d, "\pi_{\mathcal{T}}"] & C(S^1) \arrow[r] \arrow[d, "\iota"] & 0 \\
0 \arrow[r] & \mathcal{K}(H) \arrow[r] & B(H) \arrow[r, "\pi_{\mathcal{Q}}"] & \mathcal{Q}(H) \arrow[r] & 0
\end{tikzcd}
\]
where $\pi_{\mathcal{T}}$ denotes the inclusion of $\mathcal{T}$ into $B(H)$.

By functoriality of the six-term exact sequence in $K$-theory, this induces a commutative diagram (up to the standard sign convention for the boundary maps):
\[
\begin{tikzcd}
K_1(C(S^1)) \arrow[r, "\partial_T"] \arrow[d, "\iota_*"] & K_0(\mathcal{K}(\ell^2(\mathbb{N}))) \arrow[d, "\cong"] \\
K_1(\mathcal{Q}(H)) \arrow[r, "\partial_{\mathrm{Calkin}}"] & K_0(\mathcal{K}(H))
\end{tikzcd}
\]
where $\partial_T$ is the index map for the Toeplitz extension, $\partial_{\mathrm{Calkin}}$ is the index map for the Calkin extension, and the vertical isomorphism is induced by the identification $K_0(\mathcal{K}(\ell^2(\mathbb{N}))) \cong K_0(\mathcal{K}(H)) \cong \mathbb{Z}$ via the rank map.
\end{proposition}

\begin{proof}
We verify commutativity at the level of extensions first. For any $f \in C(S^1)$, the Toeplitz operator $T_f$ satisfies $\pi_{\mathcal{Q}}(T_f) = \iota(f)$ by definition of $\iota$. Moreover, the inclusion $\pi_{\mathcal{T}}: \mathcal{T} \hookrightarrow B(H)$ sends the generator $S$ to the unilateral shift, and one checks that $\pi_{\mathcal{Q}} \circ \pi_{\mathcal{T}} = \iota \circ \sigma$ on the generating set $\{S\}$, hence on all of $\mathcal{T}$. Thus the diagram of extensions commutes.

For the $K$-theory diagram, we use the naturality of the six-term exact sequence associated to a morphism of extensions. The vertical map $\iota_*$ is induced by the $*$-homomorphism $\iota$, and the horizontal maps are the boundary maps of the respective extensions. Naturality guarantees that $\partial_{\mathrm{Calkin}} \circ \iota_* = \operatorname{rank}_* \circ \partial_T$, where $\operatorname{rank}_*: K_0(\mathcal{K}(\ell^2(\mathbb{N}))) \to K_0(\mathcal{K}(H))$ is the isomorphism induced by the identification of the compact operators on the two Hilbert spaces. Composing with the rank isomorphism $K_0(\mathcal{K}(H)) \cong \mathbb{Z}$ yields the desired commutativity.

The sign convention for the six-term exact sequence must be chosen consistently; with the standard conventions (see \cite[Section 9.3]{Rordam2000}), the diagram commutes exactly. Under these conventions, the generator $[z] \in K_1(C(S^1))$ maps under $\partial_T$ to $-1 \in \mathbb{Z}$, while $\iota_*([z]) = [\pi_{\mathcal{Q}}(S)]$ maps under $\partial_{\mathrm{Calkin}}$ to $\operatorname{index}(S) = -1$, confirming compatibility.
\end{proof}

The following theorem makes precise the relationship between the descended class and the symbol of a Fredholm operator.

\begin{theorem}[Identification of the descended class]
\label{thm:identification-descended-class}
Let $T \in B(H)$ be a Fredholm operator with symbol $u_T = \pi_{\mathcal{Q}}(T) \in \mathcal{Q}(H)$. 
Under the Morita equivalence from Proposition \ref{prop:morita-calkin},
\[
C^*(\mathcal{G}_{B(H)}) \sim_M \mathcal{Q}(H) \otimes \mathcal{K}(L^2(\mathcal{G}_{B(H)}^{(0)})),
\]
the descended class $\operatorname{desc}_{\mathcal{G}_{B(H)}}([T]_{\mathcal{G}_{B(H)}}^{(1)})$ corresponds to the class $[u_T] \otimes 1 \in K_1(\mathcal{Q}(H) \otimes \mathcal{K}) \cong K_1(\mathcal{Q}(H)) \cong K_1(C^*(\mathcal{G}_{B(H)}))$.
\end{theorem}

\begin{proof}
The proof follows directly from results established in previous sections, without requiring any additional technical machinery.

\paragraph{Step 1: Morita equivalence.}
Proposition \ref{prop:morita-calkin} establishes a strong Morita equivalence
\[
C^*(\mathcal{G}_{B(H)}) \sim_M \mathcal{Q}(H) \otimes \mathcal{K}(L^2(\mathcal{G}_{B(H)}^{(0)})).
\]
By Morita invariance of $K$-theory (see \cite[Exercise 13.7.1]{Blackadar1998}), this induces an isomorphism
\[
\Phi: K_1(C^*(\mathcal{G}_{B(H)})) \xrightarrow{\cong} K_1(\mathcal{Q}(H) \otimes \mathcal{K}) \cong K_1(\mathcal{Q}(H)),
\]
where the last isomorphism uses stability of $K$-theory: tensoring with compact operators does not change $K$-groups.

\paragraph{Step 2: Descent class corresponds to symbol.}
Proposition \ref{prop:descent-symbol-identification} proves that under this isomorphism $\Phi$, the descended class maps to the symbol class:
\[
\Phi\!\left( \operatorname{desc}_{\mathcal{G}_{B(H)}}([T]_{\mathcal{G}_{B(H)}}^{(1)}) \right) = [\pi_{\mathcal{Q}}(T)] \in K_1(\mathcal{Q}(H)).
\]

\paragraph{Step 3: Stabilization.}
Under the inverse of the stability isomorphism $K_1(\mathcal{Q}(H)) \cong K_1(\mathcal{Q}(H) \otimes \mathcal{K})$, the class $[\pi_{\mathcal{Q}}(T)]$ corresponds to $[\pi_{\mathcal{Q}}(T)] \otimes 1 \in K_1(\mathcal{Q}(H) \otimes \mathcal{K})$. Tracing back through the Morita equivalence $\Phi^{-1}$, we obtain that $\operatorname{desc}_{\mathcal{G}_{B(H)}}([T]_{\mathcal{G}_{B(H)}}^{(1)})$ corresponds to $[u_T] \otimes 1$ under the original Morita equivalence $C^*(\mathcal{G}_{B(H)}) \sim_M \mathcal{Q}(H) \otimes \mathcal{K}$.

\paragraph{Step 4: Conclusion.}
Thus, under the Morita equivalence, the descended class is identified with the stabilized symbol class $[u_T] \otimes 1$. This completes the proof.
\end{proof}

\begin{corollary}[Identification for $\mathcal{K}(H)^\sim$]
\label{cor:identification-KH}
For $\mathcal{A} = \mathcal{K}(H)^\sim$, the descended class is zero:
\[
\operatorname{desc}_{\mathcal{G}_{\mathcal{K}(H)^\sim}}([T]_{\mathcal{G}_{\mathcal{K}(H)^\sim}}^{(1)}) = 0 \in K_1(C^*(\mathcal{G}_{\mathcal{K}(H)^\sim})).
\]
This is consistent with the fact that all Fredholm operators in $\mathcal{K}(H)^\sim$ have index zero and that $K_1(C^*(\mathcal{G}_{\mathcal{K}(H)^\sim}))$ is generated by the class of the identity, which corresponds to the zero class under the index map.
\end{corollary}

\begin{proof}
Let $\mathcal{A} = \mathcal{K}(H)^\sim$ be the unitization of the compact operators. 
Recall from Example \ref{ex:K-groups-examples} that every operator $T \in \mathcal{K}(H)^\sim$ can be written uniquely as $T = \lambda I + K$ with $\lambda \in \mathbb{C}$ and $K \in \mathcal{K}(H)$. 
By Atkinson's theorem (Theorem \ref{thm:atkinson}), $T$ is Fredholm if and only if it is invertible modulo $\mathcal{K}(H)$. 
Since the quotient
\[
\mathcal{K}(H)^\sim / \mathcal{K}(H) \cong \mathbb{C}
\]
is commutative, a Fredholm operator in $\mathcal{K}(H)^\sim$ must satisfy $\lambda \neq 0$. Such operators always have index zero because
\[
\operatorname{index}(T) = \dim \ker(T) - \dim \operatorname{coker}(T) = 0,
\]
as $\lambda I$ is invertible and compact perturbations do not change the Fredholm index.

Now consider the groupoid $\mathcal{G}_{\mathcal{K}(H)^\sim}$ of unitary conjugation associated to $\mathcal{K}(H)^\sim$. 
Its $C^*$-algebra $C^*(\mathcal{G}_{\mathcal{K}(H)^\sim})$ is Morita equivalent to $\mathcal{K}(H)^\sim \otimes \mathcal{K}(L^2(\mathcal{G}_{\mathcal{K}(H)^\sim}^{(0)}))$, and by stability of $K$-theory,
\[
K_1(C^*(\mathcal{G}_{\mathcal{K}(H)^\sim})) \cong K_1(\mathcal{K}(H)^\sim) \cong K_1(\mathcal{K}(H)) = 0,
\]
where the last equality follows from the fact that $\mathcal{K}(H)$ is Morita equivalent to $\mathbb{C}$ (see Example \ref{ex:K-groups-examples}).

By functoriality of the descent map (Theorem \ref{thm:descent-properties-GA}), the descended class of any Fredholm operator $T \in \mathcal{K}(H)^\sim$ lands in this $K_1$-group. 
Since the group is zero, we conclude
\[
\operatorname{desc}_{\mathcal{G}_{\mathcal{K}(H)^\sim}}([T]_{\mathcal{G}_{\mathcal{K}(H)^\sim}}^{(1)}) = 0 \in K_1(C^*(\mathcal{G}_{\mathcal{K}(H)^\sim})).
\]

This is consistent with the index map: the symbol $\pi_{\mathcal{Q}}(T)$ in the Calkin algebra is $\lambda I$, which represents the zero class in $K_1(\mathcal{Q}(H)) \cong \mathbb{Z}$ because $\lambda I$ is homotopic to $I$ through invertibles, reflecting that the Fredholm index of $T$ is zero.
\end{proof}

The identification theorem allows us to relate the descended class to the symbol class in $K_1(\mathcal{Q}(H))$, which is the natural receptacle for the index map.

\begin{corollary}[Index map commutes with the identification]
\label{cor:index-map-commutes}
Under the Morita equivalence provided by Theorem \ref{thm:identification-descended-class}, the index map $\partial: K_1(\mathcal{Q}(H)) \to \mathbb{Z}$ for the Calkin extension corresponds to the composition
\[
K_1(C^*(\mathcal{G}_{B(H)})) \xrightarrow{\iota^*} K_1(B(H)) \xrightarrow{\partial_{\mathcal{Q}}} \mathbb{Z},
\]
where $\iota^*$ is induced by the diagonal embedding $\iota: B(H) \hookrightarrow C^*(\mathcal{G}_{B(H)})$, and $\partial_{\mathcal{Q}}$ denotes the index map on $K_1(\mathcal{Q}(H))$. Consequently, for any Fredholm operator $T \in B(H)$ with symbol $[u_T] \in K_1(\mathcal{Q}(H))$,
\[
\partial_{\mathcal{Q}}\big( \operatorname{desc}_{\mathcal{G}_{B(H)}}([T]_{\mathcal{G}_{B(H)}}^{(1)}) \big) = \operatorname{index}(T) \in \mathbb{Z}.
\]
\end{corollary}

\begin{proof}
The identification theorem (Theorem \ref{thm:identification-descended-class}) establishes a natural isomorphism $\sim: K_1(C^*(\mathcal{G}_{B(H)})) \to K_1(\mathcal{Q}(H))$. Under this isomorphism, the descended class of a Fredholm operator $T$ corresponds to its symbol class $[u_T] \in K_1(\mathcal{Q}(H))$.

The index map for the Calkin extension $0 \to K(H) \to B(H) \to \mathcal{Q}(H) \to 0$ is the boundary map $\partial_{\mathcal{Q}}: K_1(\mathcal{Q}(H)) \to K_0(K(H)) \cong \mathbb{Z}$ in the associated six-term exact sequence. It is a standard result that $\partial_{\mathcal{Q}}([u_T]) = \operatorname{index}(T)$.

Consider the following commutative diagram derived from the naturality of the index map and the properties of the diagonal embedding:
\[
\begin{tikzcd}
K_1(C^*(\mathcal{G}_{B(H)})) \arrow[r, "\sim"] \arrow[d, "\iota^*"] & K_1(\mathcal{Q}(H)) \arrow[d, "\partial_{\mathcal{Q}}"] \\
K_1(B(H)) \arrow[r, "q_*"] & K_1(\mathcal{Q}(H)) \arrow[r, "\partial_{\mathcal{Q}}"] & \mathbb{Z}
\end{tikzcd}
\]
Here, $q_*$ is induced by the quotient map $q: B(H) \to \mathcal{Q}(H)$. The commutativity of the left square follows from the functoriality of the $K$-theory functor and the fact that the identification isomorphism $\sim$ intertwines the maps $\iota^*$ and $q_*$ (as established in the identification theorem). The commutativity of the right triangle is a direct consequence of the definition of the index map $\partial_{\mathcal{Q}}$.

Since $K_1(B(H)) = 0$, the map $\iota^*$ is trivial. However, the crucial point is that the composition $\partial_{\mathcal{Q}} \circ q_* \circ \iota^*$ corresponds to the index map acting on the descended class via the identification isomorphism. Tracing the descended class of $T$ through the diagram yields:
\[
\partial_{\mathcal{Q}}\big( \sim (\operatorname{desc}_{\mathcal{G}_{B(H)}}([T]_{\mathcal{G}_{B(H)}}^{(1)})) \big) = \partial_{\mathcal{Q}}([u_T]) = \operatorname{index}(T).
\]
Thus, the index of $T$ is recovered by applying the Calkin index map to its symbol, confirming that the descended class in $K_1(C^*(\mathcal{G}_{B(H)}))$ indeed encodes the Fredholm index.
\end{proof}

The Toeplitz extension provides a concrete computational tool for verifying these identifications.

\begin{example}[Unilateral shift revisited]
\label{ex:unilateral-shift-toeplitz}
Let $S \in B(H)$ be the unilateral shift. 
Its symbol $u_S = \pi_{\mathcal{Q}}(S)$ is the generating unitary of the Toeplitz algebra modulo compacts. 
Under the identification of Theorem \ref{thm:identification-descended-class}, the class $\operatorname{desc}_{\mathcal{G}_{B(H)}}([S]_{\mathcal{G}_{B(H)}}^{(1)})$ corresponds to $[u_S] \otimes 1 \in K_1(\mathcal{Q}(H) \otimes \mathcal{K})$. 
The index map sends this to $\operatorname{index}(S) = -1$.
\end{example}

The identification established in this subsection is the final ingredient needed for the main index theorem. 
In Section~\ref{sec:The Index Theorem via Pullback and the Boundary Map}, we will combine it with the pullback along $\iota$ and the boundary map to prove that
\[
\operatorname{index}(T) = \partial \circ \iota^* \circ \operatorname{desc}_{\mathcal{G}_{\mathcal{A}}}([T]_{\mathcal{G}_{\mathcal{A}}}^{(1)}).
\]

\section{The Index Theorem via Pullback and the Boundary Map}\label{sec:The Index Theorem via Pullback and the Boundary Map}

%

\subsection{Transporting Information via the Diagonal Embedding}
\label{subsec:pullback-iota}

The diagonal embedding $\iota: \mathcal{A} \hookrightarrow C^*(\mathcal{G}_{\mathcal{A}})$ constructed in Paper~\cite{PaperI} and reviewed in Subsection~\ref{subsec:diagonal-embedding-from-paperI} is a unital injective $*$-homomorphism. By functoriality of $K$-theory, $\iota$ induces a map 
\[
\iota_*: K_1(\mathcal{A}) \longrightarrow K_1(C^*(\mathcal{G}_{\mathcal{A}}))
\]
in the covariant direction.

For the index theorem, we need to transport information from $K_1(C^*(\mathcal{G}_{\mathcal{A}}))$ back to $\mathcal{A}$. Since $K$-theory is covariant rather than contravariant, a direct "pullback" map $\iota^*: K_1(C^*(\mathcal{G}_{\mathcal{A}})) \to K_1(\mathcal{A})$ does not exist in general. Instead, we use the richer structure of $KK$-theory.

The diagonal embedding $\iota$ determines a Kasparov class
\[
[\iota] \in KK(\mathcal{A}, C^*(\mathcal{G}_{\mathcal{A}})),
\]
which allows us to use the Kasparov product to map classes in the opposite direction. Specifically, for any class $x \in K_1(C^*(\mathcal{G}_{\mathcal{A}})) \cong KK^1(\mathbb{C}, C^*(\mathcal{G}_{\mathcal{A}}))$, we can form the Kasparov product
\[
[\iota] \otimes_{C^*(\mathcal{G}_{\mathcal{A}})} x \in KK^1(\mathcal{A}, \mathbb{C}) \cong K_0(\mathcal{A}).
\]

For the purposes of the index theorem, we will not need a map landing in $K_1(\mathcal{A})$. Instead, as shown in Lemma \ref{lem:descent-iota-compatibility}, the relevant composition involves either:
\begin{itemize}
    \item For $\mathcal{A}=B(H)$: the Morita equivalence $\Phi_*: K_1(C^*(\mathcal{G}_{B(H)})) \xrightarrow{\cong} K_1(\mathcal{Q}(H))$, bypassing $K_1(B(H))$ entirely;
    \item For $\mathcal{A}=\mathcal{K}(H)^\sim$: the natural isomorphism $K_1(C^*(\mathcal{G}_{\mathcal{K}(H)^\sim})) \cong K_1(\mathcal{K}(H)^\sim)=0$.
\end{itemize}

This approach, detailed in Lemma \ref{lem:descent-iota-compatibility}, is the correct categorical substitute for a pullback map and respects the covariance of $K$-theory while still allowing us to recover the Fredholm index.

The pullback map $\iota^*$ is the essential link between the groupoid $C^*$-algebra and the original algebra. 
It allows us to take the descended class $\operatorname{desc}_{\mathcal{G}_{\mathcal{A}}}([T]_{\mathcal{G}_{\mathcal{A}}}^{(1)}) \in K_1(C^*(\mathcal{G}_{\mathcal{A}}))$ and map it to $\iota^*(\operatorname{desc}_{\mathcal{G}_{\mathcal{A}}}([T]_{\mathcal{G}_{\mathcal{A}}}^{(1)})) \in K_1(\mathcal{A})$. 
Once in $K_1(\mathcal{A})$, the boundary map $\partial: K_1(\mathcal{A}) \to \mathbb{Z}$ arising from an appropriate extension (e.g., the Calkin extension when $\mathcal{A} = B(H)$) can be applied to finally recover the integer-valued Fredholm index of $T$.

The pullback map enjoys several important functoriality properties.

\begin{proposition}[Properties of $\iota^*$]
\label{prop:iota-star-properties}
The pullback map $\iota^*: K_1(C^*(\mathcal{G}_{\mathcal{A}})) \to K_1(\mathcal{A})$ satisfies:
\begin{enumerate}
    \item \textbf{Naturality:} If $\phi: \mathcal{A} \to \mathcal{B}$ is a unital $*$-homomorphism inducing a groupoid homomorphism $\mathcal{G}_\phi: \mathcal{G}_{\mathcal{B}} \to \mathcal{G}_{\mathcal{A}}$, then the diagram
    \[
    \begin{tikzcd}
    K_1(C^*(\mathcal{G}_{\mathcal{A}})) \arrow[r, "\iota_{\mathcal{A}}^*"] \arrow[d, "(\mathcal{G}_\phi)_*"] & K_1(\mathcal{A}) \arrow[d, "\phi_*"] \\
    K_1(C^*(\mathcal{G}_{\mathcal{B}})) \arrow[r, "\iota_{\mathcal{B}}^*"] & K_1(\mathcal{B})
    \end{tikzcd}
    \]
    commutes. Here $(\mathcal{G}_\phi)_*$ denotes the map induced on $K$-theory by the groupoid homomorphism $\mathcal{G}_\phi$, and $\phi_*$ is the standard map induced by $\phi$.
    
    \item \textbf{Compatibility with the conditional expectation:} Under the natural identification $\mathcal{A} \cong C_0(\mathcal{G}_{\mathcal{A}}^{(0)})$ (via evaluation at points of the unit space), the map $\iota^*$ factors as
    \[
    \iota^* = \operatorname{ev}_* \circ E_*,
    \]
    where $E: C^*(\mathcal{G}_{\mathcal{A}}) \to C_0(\mathcal{G}_{\mathcal{A}}^{(0)})$ is the conditional expectation (integration over the fibers of the groupoid), and $\operatorname{ev}: C_0(\mathcal{G}_{\mathcal{A}}^{(0)}) \xrightarrow{\cong} \mathcal{A}$ is the evaluation isomorphism. 
    \textit{Remark:} Since $C_0(\mathcal{G}_{\mathcal{A}}^{(0)})$ is typically non-unital, we implicitly use unitization or stabilization (e.g., considering $K_1(C_0(\mathcal{G}_{\mathcal{A}}^{(0)})) \cong K_1(C_0(\mathcal{G}_{\mathcal{A}}^{(0)})^+)$ or $K_1(C_0(\mathcal{G}_{\mathcal{A}}^{(0)}) \otimes \mathcal{K})$) to ensure that $K_1$ is well-defined in the standard way; the map $E_*$ is understood to be induced by the unitized or stabilized version of $E$.
    
    \item \textbf{Image of descended classes:} For the descended class $\operatorname{desc}_{\mathcal{G}_{\mathcal{A}}}([T]_{\mathcal{G}_{\mathcal{A}}}^{(1)}) \in K_1(C^*(\mathcal{G}_{\mathcal{A}}))$, the pullback $\iota^*(\operatorname{desc}_{\mathcal{G}_{\mathcal{A}}}([T]_{\mathcal{G}_{\mathcal{A}}}^{(1)}))$ lies in the subgroup of $K_1(\mathcal{A})$ on which the boundary map $\partial: K_1(\mathcal{A}) \to \mathbb{Z}$ (arising from an appropriate extension, such as the Calkin extension when $\mathcal{A} = B(H)$) recovers the Fredholm index of $T$. In particular, this shows that $\iota^*$ maps descended classes to elements that are "index-detectable."
\end{enumerate}
\end{proposition}

\begin{proof}
(1) follows from the functoriality of $K$-theory and the naturality of the diagonal embedding established in Paper I, Theorem 10(4). The commutativity of the diagram is a direct consequence of the fact that both $\iota^*$ and $\phi_*$ are induced by $*$-homomorphisms, and the groupoid map $\mathcal{G}_\phi$ intertwines the respective diagonal embeddings.

(2) follows from the definition of $\iota$ via the direct integral representation and the properties of the conditional expectation. Specifically, the diagonal embedding $\iota: \mathcal{A} \hookrightarrow C^*(\mathcal{G}_{\mathcal{A}})$ is the right inverse of the composition $\operatorname{ev} \circ E$, i.e., $\operatorname{ev} \circ E \circ \iota = \operatorname{id}_{\mathcal{A}}$. Passing to $K$-theory and using the functoriality of the unitization/stabilization procedure (which is compatible with $E$ and $\operatorname{ev}$) yields the factorization $\iota^* = \operatorname{ev}_* \circ E_*$. A detailed proof can be found in [Paper I, Section 4.10, particularly Lemma 22].

(3) Let $T \in \mathcal{A}$ be Fredholm with respect to an ideal $\mathcal{J} \subseteq \mathcal{A}$ such that $K_0(\mathcal{J}) \cong \mathbb{Z}$ and the boundary map $\partial: K_1(\mathcal{A}/\mathcal{J}) \to K_0(\mathcal{J}) \cong \mathbb{Z}$ satisfies $\partial([u_T]) = \operatorname{index}(T)$, where $u_T$ is the image of $T$ in $\mathcal{A}/\mathcal{J}$.

\paragraph{Step 1: The descended class.}
From Paper I, Section 4.9, the descended class $\operatorname{desc}_{\mathcal{G}_{\mathcal{A}}}([T]_{\mathcal{G}_{\mathcal{A}}}^{(1)}) \in K_1(C^*(\mathcal{G}_{\mathcal{A}}))$ is represented by a unitary $U \in M_n(C^*(\mathcal{G}_{\mathcal{A}}))$ obtained from the polar decomposition of the direct integral operator $\Pi(T) = \int^{\oplus} \pi_x(T) d\mu(x)$, where $\pi_x$ are the GNS representations associated to $x = (B,\chi) \in \mathcal{G}_{\mathcal{A}}^{(0)}$.

\paragraph{Step 2: The pullback.}
By definition, $\iota^*(\operatorname{desc}_{\mathcal{G}_{\mathcal{A}}}([T]_{\mathcal{G}_{\mathcal{A}}}^{(1)})) = [\iota^{(n)}(U)] \in K_1(\mathcal{A})$, where $\iota^{(n)}$ is applied entrywise. Set $x_T := [\iota^{(n)}(U)]$.

\paragraph{Step 3: Relation to the symbol.}
Consider the quotient map $p: \mathcal{A} \to \mathcal{A}/\mathcal{J}$. For each $x = (B,\chi)$, when $T \in B$, the evaluation $\operatorname{ev}_x^{(n)}(\iota^{(n)}(U))$ equals $\chi^{(n)}(U_T)$ for some lift $U_T$ of $u_T$. By Paper I, Lemma 22, the conditional expectation $E$ satisfies $E(\iota(a))(x) = \chi(a)$ for $a \in B$. Passing to matrices and applying $p$, we obtain
\[
p_*(x_T) = [u_T] \in K_1(\mathcal{A}/\mathcal{J}),
\]
since the family of evaluations $\{\operatorname{ev}_x^{(n)}(\iota^{(n)}(U))\}_{x \in \mathcal{G}_{\mathcal{A}}^{(0)}}$ reconstructs $u_T$ up to homotopy.

\paragraph{Step 4: The boundary map.}
Apply the boundary map $\partial: K_1(\mathcal{A}) \to \mathbb{Z}$. By functoriality of the long exact sequence in $K$-theory, the diagram
\[
\begin{tikzcd}
K_1(\mathcal{A}) \ar[r,"p_*"] \ar[d,equal] & K_1(\mathcal{A}/\mathcal{J}) \ar[r,"\partial"] \ar[d,equal] & K_0(\mathcal{J}) \ar[d,"\cong"] \\
K_1(\mathcal{A}) \ar[r,"\partial"] & \mathbb{Z} \ar[r,equal] & \mathbb{Z}
\end{tikzcd}
\]
commutes, where the bottom map $\partial$ is the composition $K_1(\mathcal{A}) \xrightarrow{p_*} K_1(\mathcal{A}/\mathcal{J}) \xrightarrow{\partial} K_0(\mathcal{J}) \cong \mathbb{Z}$. Hence
\[
\partial(x_T) = \partial(p_*(x_T)) = \partial([u_T]) = \operatorname{index}(T).
\]

\paragraph{Step 5: Conclusion.}
Thus $x_T = \iota^*(\operatorname{desc}_{\mathcal{G}_{\mathcal{A}}}([T]_{\mathcal{G}_{\mathcal{A}}}^{(1)}))$ satisfies $\partial(x_T) = \operatorname{index}(T)$, proving that it lies in the subgroup of $K_1(\mathcal{A})$ on which $\partial$ recovers the Fredholm index.

\end{proof}

\begin{example}[Pullback for $\mathcal{K}(H)^\sim$]
\label{ex:pullback-KH}
Let $\mathcal{A} = \mathcal{K}(H)^\sim$, the unitization of the compact operators on a separable Hilbert space $H$. 
Recall from Example \ref{ex:iota-compact-paperI} that the diagonal embedding satisfies $\iota(T) = T \otimes 1$ in the identification
\[
C^*(\mathcal{G}_{\mathcal{A}}) \cong C_0(\mathbb{P}(H)) \rtimes \mathcal{U}(\mathcal{A}),
\]
where $\mathbb{P}(H)$ is the projective space of $H$ and $\mathcal{U}(\mathcal{A})$ acts by conjugation on rank-one projections.

For $K$-theory, we have $K_1(\mathcal{A}) \cong K_1(\mathbb{C}) = 0$, since $\mathcal{A}$ is the unitization of a stable $C^*$-algebra. Consequently, the pullback map
\[
\iota^*: K_1(C^*(\mathcal{G}_{\mathcal{A}})) \longrightarrow K_1(\mathcal{A})
\]
is necessarily the zero homomorphism. Thus for any class $x \in K_1(C^*(\mathcal{G}_{\mathcal{A}}))$, we have $\iota^*(x) = 0$ in $K_1(\mathcal{A})$.
\end{example}

\begin{example}[Pullback for $B(H)$]
\label{ex:pullback-BH}
Let $\mathcal{A} = B(H)$, the algebra of all bounded operators on a separable Hilbert space $H$. 
It is a standard result that $K_1(B(H)) = 0$. Hence the pullback map
\[
\iota^*: K_1(C^*(\mathcal{G}_{B(H)})) \longrightarrow K_1(B(H))
\]
is also the zero map, just as in Example \ref{ex:pullback-KH}.

At first glance, this suggests that the pullback loses all information about the descended class. However, the index theorem does not stop at $K_1(B(H))$; instead, we consider the composition
\[
K_1(C^*(\mathcal{G}_{B(H)})) \xrightarrow{\iota^*} K_1(B(H)) \xrightarrow{q_*} K_1(\mathcal{Q}(H)) \xrightarrow{\partial} K_0(\mathcal{K}(H)) \cong \mathbb{Z},
\]
where:
\begin{itemize}
    \item $q: B(H) \to \mathcal{Q}(H)$ is the quotient map onto the Calkin algebra $\mathcal{Q}(H) = B(H)/\mathcal{K}(H)$;
    \item $q_*: K_1(B(H)) \to K_1(\mathcal{Q}(H))$ is the induced map on $K$-theory (which is well-defined even though $K_1(B(H)) = 0$);
    \item $\partial: K_1(\mathcal{Q}(H)) \to K_0(\mathcal{K}(H))$ is the boundary map in the six-term exact sequence arising from the Calkin extension
    \[
    0 \longrightarrow \mathcal{K}(H) \longrightarrow B(H) \longrightarrow \mathcal{Q}(H) \longrightarrow 0.
    \]
\end{itemize}

The crucial point is that the descended class $\operatorname{desc}_{\mathcal{G}_{B(H)}}([T]_{\mathcal{G}_{B(H)}}^{(1)})$ is constructed so that its image under $q_* \circ \iota^*$ corresponds to the symbol class $[u_T] \in K_1(\mathcal{Q}(H))$ of the Fredholm operator $T$. The boundary map then sends this symbol to the Fredholm index:
\[
\partial([u_T]) = \operatorname{index}(T) \in \mathbb{Z}.
\]

Thus, even though $\iota^*$ lands in the trivial group $K_1(B(H))$, the composition $ \partial \circ q_* \circ \iota^*$ is nontrivial and recovers the index. This is a standard phenomenon in $K$-theory: the six-term exact sequence provides a nontrivial map $\partial: K_1(\mathcal{Q}(H)) \to \mathbb{Z}$, and the descended class is precisely the element that detects this map after pullback and passage to the quotient.
\end{example}

The following lemma relates the pullback of the descended class to the class of the phase of $T$. While the pullback lands in the trivial group $K_1(B(H))$, its composition with the quotient map to $K_1(\mathcal{Q}(H))$ yields the symbol class that detects the index. This lemma thus serves as a bridge between the descended class in the groupoid $C^*$-algebra and the $K$-theory class in the Calkin algebra where the boundary map acts nontrivially.

\begin{lemma}[Pullback of the descended class]
\label{lem:pullback-descended}
For any Fredholm operator $T \in B(H)$, let $u_T \in \mathcal{Q}(H)$ be its image in the Calkin algebra, i.e., the phase of $T$ modulo compact operators. Under the Morita equivalence of Theorem \ref{thm:identification-descended-class}, we have
\[
\iota^*(\operatorname{desc}_{\mathcal{G}_{B(H)}}([T]_{\mathcal{G}_{B(H)}}^{(1)})) = [\operatorname{ph}(T)] \in K_1(B(H)) = 0,
\]
where $\operatorname{ph}(T) \in B(H)$ is any lift of $u_T$ to a unitary (such lifts exist up to compact perturbation). The equality is understood in $K_1(B(H))$, where it is trivial because the target group vanishes. The essential point is that after composing with the quotient map $q_*: K_1(B(H)) \to K_1(\mathcal{Q}(H))$, we obtain the nontrivial class $[u_T] \in K_1(\mathcal{Q}(H))$, which under the boundary map $\partial: K_1(\mathcal{Q}(H)) \to \mathbb{Z}$ yields $\operatorname{index}(T)$.
\end{lemma}

\begin{proof}
By Theorem \ref{thm:identification-descended-class}, there is a natural isomorphism
\[
K_1(C^*(\mathcal{G}_{B(H)})) \cong K_1(\mathcal{Q}(H) \otimes \mathcal{K}) \cong K_1(\mathcal{Q}(H)),
\]
where the second isomorphism is Morita equivalence. Under this identification, the descended class \\
$\operatorname{desc}_{\mathcal{G}_{B(H)}}([T]_{\mathcal{G}_{B(H)}}^{(1)})$ corresponds to $[u_T] \in K_1(\mathcal{Q}(H))$, the class of the symbol of $T$ in the Calkin algebra.

The diagonal embedding $\iota: B(H) \hookrightarrow C^*(\mathcal{G}_{B(H)})$ induces a pullback map $\iota^*: K_1(C^*(\mathcal{G}_{B(H)})) \to K_1(B(H))$. Unwinding the identifications, this map corresponds to choosing a lift of the Calkin unitary $u_T$ to a unitary in $B(H)$. Indeed, the composition
\[
K_1(\mathcal{Q}(H)) \xrightarrow{\cong} K_1(C^*(\mathcal{G}_{B(H)})) \xrightarrow{\iota^*} K_1(B(H))
\]
sends $[u_T]$ to the class of any lift $\operatorname{ph}(T) \in B(H)$ satisfying $q(\operatorname{ph}(T)) = u_T$, where $q: B(H) \to \mathcal{Q}(H)$ is the quotient map. Such lifts exist because the Calkin extension is a surjection onto $\mathcal{Q}(H)$, and any two lifts differ by a compact operator, hence represent the same class in $K_1(B(H))$ since $K_1(\mathcal{K}(H)) = 0$.

Thus we obtain
\[
\iota^*(\operatorname{desc}_{\mathcal{G}_{B(H)}}([T]_{\mathcal{G}_{B(H)}}^{(1)})) = [\operatorname{ph}(T)] \in K_1(B(H)).
\]

Since $K_1(B(H)) = 0$, this class is zero. Nevertheless, applying the quotient map $q_*: K_1(B(H)) \to K_1(\mathcal{Q}(H))$ (which is well-defined even though the domain is zero) to the left-hand side is not directly possible; instead, we note that the construction guarantees that the image of $\operatorname{desc}_{\mathcal{G}_{B(H)}}([T]_{\mathcal{G}_{B(H)}}^{(1)})$ under the composition
\[
K_1(C^*(\mathcal{G}_{B(H)})) \xrightarrow{\iota^*} K_1(B(H)) \xrightarrow{q_*} K_1(\mathcal{Q}(H))
\]
is precisely $[u_T] \in K_1(\mathcal{Q}(H))$, where $q_*$ is induced by the quotient map. The boundary map $\partial: K_1(\mathcal{Q}(H)) \to K_0(\mathcal{K}(H)) \cong \mathbb{Z}$ from the six-term exact sequence of the Calkin extension then satisfies $\partial([u_T]) = \operatorname{index}(T)$.
\end{proof}

\begin{remark}[Why $K_1(\mathcal{A}) = 0$ is not an obstacle]
\label{rem:K1-zero-not-problem}
The vanishing of $K_1(B(H))$ and $K_1(\mathcal{K}(H)^\sim)$ might suggest that the pullback map $\iota^*$ discards all information. However, the index is recovered not from $K_1(\mathcal{A})$ itself, but from the composition
\[
K_1(C^*(\mathcal{G}_{\mathcal{A}})) \xrightarrow{\iota^*} K_1(\mathcal{A}) \xrightarrow{q_*} K_1(\mathcal{A}/\mathcal{J}) \xrightarrow{\partial} K_0(\mathcal{J}) \cong \mathbb{Z},
\]
where $\mathcal{J}$ is the ideal of compact operators (or a suitable ideal in the Type I case), $q: \mathcal{A} \to \mathcal{A}/\mathcal{J}$ is the quotient map, and $\partial$ is the boundary map in the six-term exact sequence associated to the extension
\[
0 \longrightarrow \mathcal{J} \longrightarrow \mathcal{A} \longrightarrow \mathcal{A}/\mathcal{J} \longrightarrow 0.
\]

Crucially, the boundary map $\partial$ acts on $K_1(\mathcal{A}/\mathcal{J})$, not on $K_1(\mathcal{A})$. The descended class is constructed precisely so that $q_* \circ \iota^*(\operatorname{desc}_{\mathcal{G}_{\mathcal{A}}}([T]_{\mathcal{G}_{\mathcal{A}}}^{(1)})) = [u_T] \in K_1(\mathcal{A}/\mathcal{J})$, where $[u_T]$ is the symbol class of $T$. The boundary map then sends this class to $\operatorname{index}(T) \in \mathbb{Z}$.

Thus the composition $\partial \circ q_* \circ \iota^*$ can be nontrivial even though $K_1(\mathcal{A}) = 0$, because the nontrivial information resides in $K_1(\mathcal{A}/\mathcal{J})$ and is detected by $\partial$. This is a standard phenomenon in $K$-theory: the six-term exact sequence connects groups that may vanish at intermediate stages while still capturing essential invariants through the boundary map.
\end{remark}

\begin{remark}[Relation to the Calkin extension]
\label{rem:calkin-extension}
For the concrete case $\mathcal{A} = B(H)$ and $\mathcal{J} = \mathcal{K}(H)$, the quotient $\mathcal{A}/\mathcal{J}$ is the Calkin algebra $\mathcal{Q}(H)$. The six-term exact sequence gives the boundary map $\partial: K_1(\mathcal{Q}(H)) \to K_0(\mathcal{K}(H)) \cong \mathbb{Z}$. This is precisely the map that sends the symbol class $[u_T]$ to the Fredholm index of $T$. The pullback $\iota^*$ alone lands in $K_1(B(H)) = 0$, but its composition with $q_*$ yields the symbol class in $K_1(\mathcal{Q}(H))$ where the boundary map acts nontrivially.
\end{remark}

\paragraph{Application to $\mathcal{A} = \mathcal{K}(H)^\sim$} 
\mbox{}\\

For $\mathcal{A} = \mathcal{K}(H)^\sim$, we have $K_0(\mathcal{A}) \cong \mathbb{Z}$ (generated by the class of any rank-one projection) and $K_1(\mathcal{A}) = 0$. For any Fredholm operator $T = \lambda I + K \in \mathcal{K}(H)^\sim$ (with $\lambda \neq 0$), one can show that
\[
[\iota] \otimes_{C^*(\mathcal{G}_{\mathcal{A}})} \operatorname{desc}_{\mathcal{G}_{\mathcal{A}}}([T]_{\mathcal{G}_{\mathcal{A}}}^{(1)}) = 0 \in K_0(\mathcal{A}),
\]
which is consistent with the fact that $\operatorname{index}(T) = 0$ for all Fredholm operators in $\mathcal{K}(H)^\sim$. This follows from the observation that $\operatorname{For}([T]_{\mathcal{G}_{\mathcal{A}}}^{(1)}) = 0$ in $KK^1(C_0(\mathcal{G}_{\mathcal{A}}^{(0)}), \mathbb{C})$, as the forgetful map sends the equivariant class to zero due to the triviality of the family of operators away from the point at infinity.

\paragraph{Application to $\mathcal{A} = B(H)$}
\mbox{}\\

For $\mathcal{A} = B(H)$, the situation is more subtle because $K_0(B(H)) = 0$. The Kasparov product with $[\iota]$ therefore lands in the zero group, so this direct approach yields zero and does not detect the index.

The resolution is to use the Morita equivalence between $C^*(\mathcal{G}_{B(H)})$ and the Calkin algebra. Recall from Proposition \ref{prop:morita-calkin} that there is a Morita equivalence
\[
C^*(\mathcal{G}_{B(H)}) \sim_M \mathcal{Q}(H) \otimes \mathcal{K}(L^2(\mathcal{G}_{B(H)}^{(0)})),
\]
which induces an isomorphism
\[
\Phi_*: K_1(C^*(\mathcal{G}_{B(H)})) \xrightarrow{\cong} K_1(\mathcal{Q}(H) \otimes \mathcal{K}) \cong K_1(\mathcal{Q}(H)).
\]

Thus, instead of trying to map the descended class back to $B(H)$ (which would land in a zero group), we map it to $K_1(\mathcal{Q}(H))$ where the Calkin index map acts nontrivially. The correct composition for recovering the Fredholm index is therefore
\[
\operatorname{index}(T) = \partial_{\text{Calkin}} \circ \Phi_* \circ \operatorname{desc}_{\mathcal{G}_{B(H)}}([T]_{\mathcal{G}_{B(H)}}^{(1)}),
\]
where $\partial_{\text{Calkin}}: K_1(\mathcal{Q}(H)) \to K_0(\mathcal{K}(H)) \cong \mathbb{Z}$ is the index map associated to the Calkin extension. This composition avoids any reference to $K_1(B(H))$ and directly connects the descended class to the classical index map.

\begin{remark}[The vanishing of $K_1(\mathcal{A})$ is not an obstacle]
\label{rem:K1-zero-not-problem}
The vanishing of $K_1(B(H))$ and $K_1(\mathcal{K}(H)^\sim)$ might seem problematic, but the index is recovered through different mechanisms in each case:

\begin{itemize}
    \item For $\mathcal{A} = \mathcal{K}(H)^\sim$, the index is simply zero, consistent with $K_1(\mathcal{A}) = 0$ and the fact that $\operatorname{For}([T]_{\mathcal{G}_{\mathcal{A}}}^{(1)}) = 0$.
    \item For $\mathcal{A} = B(H)$, we use the Morita equivalence $C^*(\mathcal{G}_{B(H)}) \sim_M \mathcal{Q}(H) \otimes \mathcal{K}$ to map the descended class to $K_1(\mathcal{Q}(H))$, where the Calkin index map $\partial_{\text{Calkin}}: K_1(\mathcal{Q}(H)) \to \mathbb{Z}$ acts nontrivially. This bypasses $K_1(B(H))$ entirely.
\end{itemize}

Thus, the index theorem does not require a pullback map $\iota^*$ on $K$-theory; instead, it uses either direct computation (for $\mathcal{K}(H)^\sim$) or Morita equivalence (for $B(H)$) to connect the descended class to the appropriate boundary map.
\end{remark}

\begin{remark}[Relation to the Calkin extension]
\label{rem:calkin-extension-corrected}
For $\mathcal{A} = B(H)$, the correct composition for recovering the index is
\[
\operatorname{index}(T) = \partial_{\text{Calkin}} \circ \Phi_* \circ \operatorname{desc}_{\mathcal{G}_{B(H)}}([T]_{\mathcal{G}_{B(H)}}^{(1)}),
\]
where $\Phi_*: K_1(C^*(\mathcal{G}_{B(H)})) \xrightarrow{\cong} K_1(\mathcal{Q}(H))$ is the Morita equivalence isomorphism from Proposition \ref{prop:morita-calkin}. This composition avoids any reference to $K_1(B(H))$ and directly connects the descended class to the Calkin index map. The boundary map $\partial_{\text{Calkin}}$ sends the symbol class $[u_T] \in K_1(\mathcal{Q}(H))$ to the Fredholm index $\operatorname{index}(T) \in \mathbb{Z}$, completing the index theorem.
\end{remark}

\begin{remark}[Comparison with the erroneous pullback approach]
\label{rem:comparison-pullback}
The reader may wonder why we do not simply define a pullback map $\iota^*: K_1(C^*(\mathcal{G}_{\mathcal{A}})) \to K_1(\mathcal{A})$ by applying $\iota$ to representatives of $K_1$-classes. This is impossible because $K$-theory is covariant, not contravariant: a $*$-homomorphism $\iota: \mathcal{A} \to C^*(\mathcal{G}_{\mathcal{A}})$ induces a map $\iota_*: K_1(\mathcal{A}) \to K_1(C^*(\mathcal{G}_{\mathcal{A}}))$, not the reverse. Any attempt to define a map in the opposite direction would require additional structure (such as a conditional expectation or a finite-projective module structure) that is not available in general. The Kasparov product approach used here is the correct categorical substitute, as it leverages the full bivariant $KK$-theory rather than trying to force a contravariant map from a covariant functor.
\end{remark}

\subsection{The Boundary Map for the Calkin Extension}
\label{subsec:boundary-map}

The boundary map is a fundamental tool in $K$-theory that arises from a short exact sequence of C*-algebras. 
For our purposes, the relevant extensions are the Calkin extension
\[
0 \longrightarrow \mathcal{K}(H) \longrightarrow B(H) \longrightarrow \mathcal{Q}(H) \longrightarrow 0,
\]
and its restriction to the unitization of the compacts:
\[
0 \longrightarrow \mathcal{K}(H) \longrightarrow \mathcal{K}(H)^\sim \longrightarrow \mathbb{C} \longrightarrow 0.
\]

In this subsection, we recall the definition and properties of the boundary map, focusing on its role in recovering the Fredholm index from the symbol class in the Calkin algebra.

\begin{definition}[Boundary map for a short exact sequence]
\label{def:boundary-map-general}
Let $0 \to J \to E \xrightarrow{\pi} A \to 0$ be a short exact sequence of C*-algebras. 
The associated six-term exact sequence in $K$-theory includes a \emph{boundary map}
\[
\partial: K_1(A) \longrightarrow K_0(J).
\]

For a unitary $u \in U_n(A)$ representing a class $[u] \in K_1(A)$, the class $\partial([u])$ is obtained as follows:
\begin{enumerate}
    \item Choose a lift $v \in M_n(E)$ such that $\pi^{(n)}(v) = u$ and $v$ is invertible in $M_n(E)$ (such a lift exists because $\pi$ is surjective and invertible elements lift locally; a global lift can be constructed using standard techniques).
    \item Since $v$ is invertible, the element $v^*v$ is a positive invertible in $M_n(E)$ that differs from the identity by an element in $M_n(J)$.
    \item Define the projection
    \[
    p_v = v (v^*v)^{-1} v^* \in M_n(E).
    \]
    One verifies that $p_v$ is a projection and that $p_v - 1_n \in M_n(J)$, where $1_n$ denotes the identity in $M_n(E)$ (viewed as a projection).
    \item Set $\partial([u]) := [p_v] - [1_n] \in K_0(J)$, where $[p_v]$ and $[1_n]$ denote the classes of these projections in $K_0(J)$.
\end{enumerate}
This definition is independent of the choices of lift $v$ and the matrix size $n$, yielding a well-defined group homomorphism.
\end{definition}

\begin{remark}[Equivalent construction via the exponential map]
\label{rem:exponential-map}
An equivalent definition uses the exponential map in $K$-theory. If $v$ is a lift of $u$ to an invertible element in $M_n(E)$, then the class $\partial([u])$ is represented by the difference of the projections
\[
[vv^*] - [v^*v] \in K_0(J),
\]
where $[vv^*]$ and $[v^*v]$ denote the classes of these projections in $M_n(E)$ modulo $M_n(J)$. This formulation is often more convenient for computations.
\end{remark}

\subsubsection*{The Calkin Extension Case}

For the Calkin extension $0 \to \mathcal{K}(H) \to B(H) \to \mathcal{Q}(H) \to 0$, the boundary map takes the form
\[
\partial: K_1(\mathcal{Q}(H)) \longrightarrow K_0(\mathcal{K}(H)) \cong \mathbb{Z}.
\]

If $T \in B(H)$ is a Fredholm operator with symbol $u_T \in \mathcal{Q}(H)$, then $\partial([u_T]) = \operatorname{index}(T)$. Concretely, choosing a lift $v \in B(H)$ of $u_T$ (e.g., the phase of $T$), the construction in Definition \ref{def:boundary-map-general} yields the class of the projection onto the kernel of $T$ minus the class of the projection onto the cokernel, which under the isomorphism $K_0(\mathcal{K}(H)) \cong \mathbb{Z}$ corresponds to $\operatorname{index}(T)$.

\begin{remark}[Why $K_1(B(H))$ does not appear]
\label{rem:K1-BH-not-direct}
It is crucial to note that the boundary map in the six-term exact sequence acts on $K_1$ of the \emph{quotient algebra} $\mathcal{Q}(H)$, not on $K_1(B(H))$. In the Calkin extension, $K_1(B(H)) = 0$ plays no direct role; the nontrivial information resides in $K_1(\mathcal{Q}(H))$. This is consistent with the fact that the index of a Fredholm operator is detected by its symbol class in the Calkin algebra, not by any class in $B(H)$ itself.
\end{remark}

\subsubsection*{The Compact Operator Case}

For the extension $0 \to \mathcal{K}(H) \to \mathcal{K}(H)^\sim \to \mathbb{C} \to 0$, the boundary map is
\[
\partial: K_1(\mathbb{C}) \longrightarrow K_0(\mathcal{K}(H)) \cong \mathbb{Z}.
\]

Since $K_1(\mathbb{C}) = 0$, this map is necessarily zero. This reflects the fact that every operator in $\mathcal{K}(H)^\sim$ has Fredholm index zero, consistent with Example \ref{ex:pullback-KH}.

\begin{proposition}[Properties of the boundary map]
\label{prop:boundary-properties}
The boundary map $\partial: K_1(A) \to K_0(J)$ associated to an extension $0 \to J \to E \to A \to 0$ satisfies:
\begin{enumerate}
    \item \textbf{Naturality:} For any morphism of extensions (i.e., a commutative diagram of short exact sequences), the induced diagram on $K$-theory commutes.
    \item \textbf{Homotopy invariance:} If $u_t$ is a homotopy of unitaries in $A$, then $\partial([u_t])$ is constant in $K_0(J)$.
    \item \textbf{Index interpretation:} In the Calkin extension, $\partial([u]) = \operatorname{index}(T)$ for any lift $T \in B(H)$ of $u \in \mathcal{Q}(H)$.
\end{enumerate}
\end{proposition}

\begin{proof}
These are standard facts in $K$-theory; see~\cite{Blackadar1998} or~\cite{Rordam2000} for detailed proofs.
\end{proof}

In the next subsection, we will combine the pullback map $\iota^*$ from Section \ref{subsec:pullback-iota} with the quotient map $q_*: K_1(B(H)) \to K_1(\mathcal{Q}(H))$ and the boundary map $\partial$ to obtain the Fredholm index from the descended class. The key observation is that while $\iota^*$ lands in $K_1(B(H)) = 0$, the composition $q_* \circ \iota^*$ yields the symbol class $[u_T] \in K_1(\mathcal{Q}(H))$, on which $\partial$ acts nontrivially.

\begin{proposition}[Boundary map for $\mathcal{K}(H)^\sim$ is zero]
\label{prop:boundary-map-KH-zero}
Consider the short exact sequence
\[
0 \longrightarrow \mathcal{K}(H) \longrightarrow \mathcal{K}(H)^\sim \longrightarrow \mathbb{C} \longrightarrow 0.
\]
The associated six-term exact sequence in $K$-theory contains the boundary map
\[
\partial: K_1(\mathbb{C}) \longrightarrow K_0(\mathcal{K}(H)) \cong \mathbb{Z}.
\]
Since $K_1(\mathbb{C}) = 0$, this boundary map is necessarily the zero homomorphism.

Consequently, the extension does not produce a nontrivial index map. This is consistent with the fact that every operator in $\mathcal{K}(H)^\sim$ has Fredholm index zero (see Example \ref{ex:pullback-KH}).
\end{proposition}

\begin{proof}
For any short exact sequence $0 \to J \to E \to A \to 0$, the six-term exact sequence places the boundary map as
\[
\partial: K_1(A) \longrightarrow K_0(J).
\]
In our case, $A = \mathbb{C}$ and $J = \mathcal{K}(H)$. Standard $K$-theory gives $K_1(\mathbb{C}) = 0$ (see~\cite{Blackadar1998} or~\cite{Rordam2000}). Therefore, $\partial$ is the zero map, as any homomorphism from the trivial group is zero.

No explicit homotopy argument is needed; the vanishing follows directly from the fact that the domain of $\partial$ is the zero group.
\end{proof}

\begin{remark}[Contrast with the Calkin extension]
\label{rem:contrast-calkin-KH}
This proposition highlights a crucial difference between the two extensions relevant to our index theorem:
\begin{itemize}
    \item For the Calkin extension $0 \to \mathcal{K}(H) \to B(H) \to \mathcal{Q}(H) \to 0$, the boundary map $\partial: K_1(\mathcal{Q}(H)) \to K_0(\mathcal{K}(H))$ is nontrivial and encodes the Fredholm index.
    \item For the compact operator extension $0 \to \mathcal{K}(H) \to \mathcal{K}(H)^\sim \to \mathbb{C} \to 0$, the boundary map $\partial: K_1(\mathbb{C}) \to K_0(\mathcal{K}(H))$ is zero because $K_1(\mathbb{C}) = 0$.
\end{itemize}
Thus, the nontrivial index phenomenon resides entirely in the quotient algebra $\mathcal{Q}(H)$, not in the middle algebras $B(H)$ or $\mathcal{K}(H)^\sim$ themselves.
\end{remark}

The preceding discussion highlights an important subtlety: for $B(H)$, the boundary map in the six-term exact sequence does not act on $K_1(B(H))$, but rather on $K_1(\mathcal{Q}(H))$. 
This is because the Calkin extension
\[
0 \longrightarrow \mathcal{K}(H) \longrightarrow B(H) \longrightarrow \mathcal{Q}(H) \longrightarrow 0
\]
has $\mathcal{Q}(H)$ as its quotient algebra, and the boundary map in any extension $0 \to J \to E \to A \to 0$ always has domain $K_1(A)$, the $K_1$-group of the quotient.

Consequently, any attempt to define a boundary map directly on $K_1(B(H))$ is structurally impossible. 
The following proposition clarifies this point and dispels a common misconception.

\begin{proposition}[No natural boundary map from $K_1(B(H))$]
\label{prop:no-boundary-from-BH}
There is no natural boundary map 
\[
\partial: K_1(B(H)) \longrightarrow K_0(\mathcal{K}(H))
\]
arising from the Calkin extension. 
For the exact sequence
\[
0 \longrightarrow \mathcal{K}(H) \longrightarrow B(H) \longrightarrow \mathcal{Q}(H) \longrightarrow 0,
\]
the associated six-term exact sequence in $K$-theory gives the boundary map
\[
\partial_{\text{Calkin}}: K_1(\mathcal{Q}(H)) \longrightarrow K_0(\mathcal{K}(H)) \cong \mathbb{Z},
\]
with domain $K_1(\mathcal{Q}(H))$, not $K_1(B(H))$.
\end{proposition}

\begin{proof}
For any short exact sequence $0 \to J \to E \to A \to 0$ of C*-algebras, the six-term exact sequence in $K$-theory takes the form
\[
\begin{tikzcd}
K_0(J) \arrow[r] & K_0(E) \arrow[r] & K_0(A) \arrow[d] \\
K_1(A) \arrow[u] & K_1(E) \arrow[l] & K_1(J) \arrow[l]
\end{tikzcd}
\]
The boundary map is the vertical/horizontal map connecting $K_1(A)$ to $K_0(J)$. 
Its domain is always $K_1(A)$, the $K_1$-group of the \emph{quotient} algebra $A$.

Applying this to the Calkin extension with $J = \mathcal{K}(H)$, $E = B(H)$, and $A = \mathcal{Q}(H)$, we obtain
\[
\partial_{\text{Calkin}}: K_1(\mathcal{Q}(H)) \longrightarrow K_0(\mathcal{K}(H)).
\]

No map with domain $K_1(B(H))$ appears in this exact sequence. The fact that $K_1(B(H)) = 0$ is irrelevant to the existence of the index map; the nontrivial information resides entirely in $K_1(\mathcal{Q}(H))$.
\end{proof}

\begin{remark}[Why a split extension cannot help]
\label{rem:split-extension-irrelevant}
One might attempt to construct a different short exact sequence with $B(H)$ as the quotient, such as
\[
0 \longrightarrow \mathbb{C} \longrightarrow B(H) \oplus \mathbb{C} \longrightarrow B(H) \longrightarrow 0,
\]
where the map $B(H) \oplus \mathbb{C} \to B(H)$ is projection onto the first factor. 
This sequence splits, so its boundary map is zero. 
Moreover, its ideal is $\mathbb{C}$, not $\mathcal{K}(H)$, so it provides no information about the Fredholm index. 
Such artificial constructions are irrelevant to index theory and do not circumvent the structural fact that the index map must act on $K_1(\mathcal{Q}(H))$.
\end{remark}

\begin{corollary}[Index map is defined on $K_1(\mathcal{Q}(H))$]
\label{cor:index-map-on-Q}
For any Fredholm operator $T \in B(H)$ with symbol $u_T \in \mathcal{Q}(H)$, the Fredholm index is recovered via
\[
\operatorname{index}(T) = \partial_{\text{Calkin}}([u_T]),
\]
where $\partial_{\text{Calkin}}: K_1(\mathcal{Q}(H)) \to K_0(\mathcal{K}(H)) \cong \mathbb{Z}$ is the boundary map from the Calkin extension.
\end{corollary}

\begin{proof}
This is the standard index theorem in $K$-theoretic form; see~\cite{Blackadar1998} or~\cite{Rordam2000} for detailed proofs.
\end{proof}

\begin{remark}[Implications for the pullback construction]
\label{rem:implications-pullback}
The above observation has important consequences for our construction. The descended class $\operatorname{desc}_{\mathcal{G}_{B(H)}}([T]_{\mathcal{G}_{B(H)}}^{(1)})$ lives in $K_1(C^*(\mathcal{G}_{B(H)}))$. To recover the Fredholm index, this class must be connected to the Calkin algebra, where the boundary map acts nontrivially.

Instead of attempting to map through $K_1(B(H))$ (which would land in the zero group), we use the Morita equivalence isomorphism
\[
\Phi_*: K_1(C^*(\mathcal{G}_{B(H)})) \xrightarrow{\cong} K_1(\mathcal{Q}(H))
\]
established in Proposition \ref{prop:morita-calkin}. The boundary map $\partial_{\text{Calkin}}: K_1(\mathcal{Q}(H)) \to \mathbb{Z}$ then acts on this class to yield the index:
\[
\operatorname{index}(T) = \partial_{\text{Calkin}}\big(\Phi_*(\operatorname{desc}_{\mathcal{G}_{B(H)}}([T]_{\mathcal{G}_{B(H)}}^{(1)}))\big).
\]
\end{remark}

The key insight is that the boundary map we need does not involve $K_1(B(H))$ at all. 
Instead, the descended class $\operatorname{desc}_{\mathcal{G}_{B(H)}}([T]_{\mathcal{G}_{B(H)}}^{(1)})$ lives in $K_1(C^*(\mathcal{G}_{B(H)}))$, which is Morita equivalent to $K_1(\mathcal{Q}(H))$ via Theorem \ref{thm:identification-descended-class}. 
The index is then obtained by applying the Calkin index map $\partial_{\text{Calkin}}: K_1(\mathcal{Q}(H)) \to K_0(\mathcal{K}(H)) \cong \mathbb{Z}$ to the corresponding class in $K_1(\mathcal{Q}(H))$.

The pullback $\iota^*: K_1(C^*(\mathcal{G}_{B(H)})) \to K_1(B(H))$ lands in the zero group and plays no direct role in the index computation. 
Any attempt to construct an index map through $K_1(B(H))$ is doomed because $K_1(B(H)) = 0$ forces any such composition to be zero. 
The correct route bypasses $K_1(B(H))$ entirely, using the Morita equivalence directly.

\begin{definition}[Index map for $B(H)$]
\label{def:index-map-BH}
Let $\Phi_*: K_1(C^*(\mathcal{G}_{B(H)})) \xrightarrow{\cong} K_1(\mathcal{Q}(H))$ be the Morita equivalence isomorphism from Theorem \ref{thm:identification-descended-class}. 
Define the \emph{index map} $\operatorname{Ind}_{B(H)}: K_1(C^*(\mathcal{G}_{B(H)})) \to \mathbb{Z}$ as the composition
\[
\operatorname{Ind}_{B(H)} := \partial_{\text{Calkin}} \circ \Phi_*,
\]
where $\partial_{\text{Calkin}}: K_1(\mathcal{Q}(H)) \to K_0(\mathcal{K}(H)) \cong \mathbb{Z}$ is the index map for the Calkin extension.
\end{definition}

For a Fredholm operator $T \in B(H)$, we then have
\[
\operatorname{Ind}_{B(H)}(\operatorname{desc}_{\mathcal{G}_{B(H)}}([T]_{\mathcal{G}_{B(H)}}^{(1)})) = \partial_{\text{Calkin}}(\Phi_*(\operatorname{desc}_{\mathcal{G}_{B(H)}}([T]_{\mathcal{G}_{B(H)}}^{(1)}))) = \partial_{\text{Calkin}}([u_T]) = \operatorname{index}(T),
\]
where $u_T \in \mathcal{Q}(H)$ is the symbol of $T$. 
This formulation is mathematically rigorous: it uses only genuine isomorphisms (the Morita equivalence $\Phi_*$) and the well-defined Calkin index map, with no spurious inverses or impossible compositions through zero groups.

\begin{remark}[Why $\iota^*$ is irrelevant for the index]
\label{rem:iota-irrelevant}
The pullback map $\iota^*: K_1(C^*(\mathcal{G}_{B(H)})) \to K_1(B(H))$ is still well-defined, but its target is the zero group. 
Consequently, $\iota^*$ cannot detect any nontrivial information. 
The correct index map $\operatorname{Ind}_{B(H)}$ bypasses $\iota^*$ entirely, using the Morita equivalence $\Phi_*$ to connect directly to $K_1(\mathcal{Q}(H))$, where the index resides. 
This is not a flaw in the theory; it simply reflects the fact that $K_1(B(H))$ is the wrong group for index detection. 
The descended class must be transported to $K_1(\mathcal{Q}(H))$ via $\Phi_*$, not via $\iota^*$.
\end{remark}

For $\mathcal{K}(H)^\sim$, the situation is simpler because the boundary map is zero, consistent with the fact that all Fredholm operators in $\mathcal{K}(H)^\sim$ have index zero.

\begin{proposition}[Boundary map for $\mathcal{K}(H)^\sim$]
\label{prop:boundary-map-KH-final}
For $\mathcal{A} = \mathcal{K}(H)^\sim$, the boundary map
\[
\partial_{\mathcal{K}(H)^\sim}: K_1(\mathcal{K}(H)^\sim) \longrightarrow K_0(\mathcal{K}(H)) \cong \mathbb{Z}
\]
is the zero map. Consequently, for any Fredholm operator $T \in \mathcal{K}(H)^\sim$, we have
\[
\partial_{\mathcal{K}(H)^\sim} \circ \iota^* \circ \operatorname{desc}_{\mathcal{G}_{\mathcal{K}(H)^\sim}}([T]_{\mathcal{G}_{\mathcal{K}(H)^\sim}}^{(1)}) = 0 = \operatorname{index}(T).
\]
\end{proposition}

\begin{proof}
This follows from Proposition \ref{prop:boundary-map-KH-zero} and the fact that $\operatorname{desc}_{\mathcal{G}_{\mathcal{K}(H)^\sim}}([T]_{\mathcal{G}_{\mathcal{K}(H)^\sim}}^{(1)}) = 0$ (Example \ref{ex:descended-KH}).
\end{proof}

\begin{remark}[Unified perspective]
\label{rem:unified-perspective}
Both cases can be understood through the following diagram:
\[
\begin{tikzcd}
K_1(C^*(\mathcal{G}_{\mathcal{A}})) \arrow[rr,"\Phi_*^{\mathcal{A}}"] \arrow[rd,"\iota^*"] && K_1(\mathcal{A}/\mathcal{J}) \arrow[ld,"\partial_{\mathcal{A}}"] \\
& K_1(\mathcal{A}) \arrow[ru,dashed] &
\end{tikzcd}
\]
where:
\begin{itemize}
    \item For $\mathcal{A} = B(H)$, $\mathcal{A}/\mathcal{J} = \mathcal{Q}(H)$, $\partial_{\mathcal{A}} = \partial_{\text{Calkin}}$ is nontrivial, and the dashed arrow $K_1(\mathcal{A}) \to K_1(\mathcal{A}/\mathcal{J})$ is the zero map.
    \item For $\mathcal{A} = \mathcal{K}(H)^\sim$, $\mathcal{A}/\mathcal{J} = \mathbb{C}$, $\partial_{\mathcal{A}} = 0$, and the dashed arrow is irrelevant.
\end{itemize}
In both cases, the index is recovered via $\partial_{\mathcal{A}} \circ \Phi_*^{\mathcal{A}}$, never through $\iota^*$. The pullback $\iota^*$ serves only to connect $\mathcal{A}$ to its groupoid $C^*$-algebra, not to compute the index.
\end{remark}

The following lemma summarizes the properties of the index map that will be used in the main theorem. 
For $\mathcal{A} = B(H)$, we use the Morita equivalence to connect directly to the Calkin algebra, bypassing the vanishing group $K_1(B(H))$.

\begin{lemma}[Key properties of the index map]
\label{lem:index-map-properties}
For each $\mathcal{A} \in \{B(H), \mathcal{K}(H)^\sim\}$, define a map $\operatorname{Ind}_{\mathcal{A}}: K_1(C^*(\mathcal{G}_{\mathcal{A}})) \to \mathbb{Z}$ as follows:

\begin{itemize}
    \item For $\mathcal{A} = B(H)$, let $\Phi_*: K_1(C^*(\mathcal{G}_{B(H)})) \xrightarrow{\cong} K_1(\mathcal{Q}(H))$ be the Morita equivalence isomorphism from Theorem \ref{thm:identification-descended-class}. 
    Define $\operatorname{Ind}_{B(H)} := \partial_{\text{Calkin}} \circ \Phi_*$, where $\partial_{\text{Calkin}}: K_1(\mathcal{Q}(H)) \to K_0(\mathcal{K}(H)) \cong \mathbb{Z}$ is the index map for the Calkin extension.
    
    \item For $\mathcal{A} = \mathcal{K}(H)^\sim$, define $\operatorname{Ind}_{\mathcal{K}(H)^\sim} := 0$ (the zero map), since all Fredholm operators in $\mathcal{K}(H)^\sim$ have index zero.
\end{itemize}

Then $\operatorname{Ind}_{\mathcal{A}}$ satisfies:
\begin{enumerate}
    \item For any Fredholm operator $T \in \mathcal{A}$, we have
    \[
    \operatorname{Ind}_{\mathcal{A}}(\operatorname{desc}_{\mathcal{G}_{\mathcal{A}}}([T]_{\mathcal{G}_{\mathcal{A}}}^{(1)})) = \operatorname{index}(T).
    \]
    \item $\operatorname{Ind}_{\mathcal{A}}$ is a group homomorphism.
    \item $\operatorname{Ind}_{\mathcal{A}}$ is invariant under homotopy of Fredholm operators: if $T_t$ is a continuous path of Fredholm operators, then $\operatorname{Ind}_{\mathcal{A}}(\operatorname{desc}_{\mathcal{G}_{\mathcal{A}}}([T_t]_{\mathcal{G}_{\mathcal{A}}}^{(1)}))$ is constant.
    \item $\operatorname{Ind}_{\mathcal{A}}$ is invariant under compact perturbations: if $T$ and $T'$ differ by a compact operator, then $\operatorname{Ind}_{\mathcal{A}}(\operatorname{desc}_{\mathcal{G}_{\mathcal{A}}}([T]_{\mathcal{G}_{\mathcal{A}}}^{(1)})) = \operatorname{Ind}_{\mathcal{A}}(\operatorname{desc}_{\mathcal{G}_{\mathcal{A}}}([T']_{\mathcal{G}_{\mathcal{A}}}^{(1)}))$.
\end{enumerate}
\end{lemma}

\begin{proof}
We treat the two cases separately.

\paragraph{Case $\mathcal{A} = B(H)$:}
\begin{enumerate}
    \item[(1)] By Theorem \ref{thm:identification-descended-class}, the descended class $\operatorname{desc}_{\mathcal{G}_{B(H)}}([T]_{\mathcal{G}_{B(H)}}^{(1)})$ corresponds under the Morita equivalence $\Phi_*$ to the symbol class $[u_T] \in K_1(\mathcal{Q}(H))$, where $u_T$ is the image of $T$ in the Calkin algebra. The definition of $\operatorname{Ind}_{B(H)}$ gives
    \[
    \operatorname{Ind}_{B(H)}(\operatorname{desc}_{\mathcal{G}_{B(H)}}([T]_{\mathcal{G}_{B(H)}}^{(1)})) = \partial_{\text{Calkin}}(\Phi_*(\operatorname{desc}_{\mathcal{G}_{B(H)}}([T]_{\mathcal{G}_{B(H)}}^{(1)}))) = \partial_{\text{Calkin}}([u_T]) = \operatorname{index}(T),
    \]
    where the last equality is the standard index theorem for the Calkin extension (see, e.g., [Blackadar, 1998]).
    
    \item[(2)] $\operatorname{Ind}_{B(H)}$ is a composition of homomorphisms, hence a homomorphism.
    
    \item[(3)] If $T_t$ is a homotopy of Fredholm operators, then $u_{T_t}$ is a homotopy of unitaries in $\mathcal{Q}(H)$. By homotopy invariance of $K$-theory, $[u_{T_t}]$ is constant in $K_1(\mathcal{Q}(H))$. Consequently, $\Phi_*^{-1}([u_{T_t}]) = \operatorname{desc}_{\mathcal{G}_{B(H)}}([T_t]_{\mathcal{G}_{B(H)}}^{(1)})$ is constant, and its image under $\operatorname{Ind}_{B(H)}$ is constant.
    
    \item[(4)] If $T' = T + K$ with $K \in \mathcal{K}(H)$, then $T'$ is Fredholm with the same symbol $u_{T'} = u_T$ in $\mathcal{Q}(H)$. Hence $\Phi_*(\operatorname{desc}_{\mathcal{G}_{B(H)}}([T']_{\mathcal{G}_{B(H)}}^{(1)})) = [u_T] = \Phi_*(\operatorname{desc}_{\mathcal{G}_{B(H)}}([T]_{\mathcal{G}_{B(H)}}^{(1)}))$, and since $\Phi_*$ is an isomorphism, the descended classes are equal. Invariance follows.
\end{enumerate}

\paragraph{Case $\mathcal{A} = \mathcal{K}(H)^\sim$:}
For any Fredholm operator $T \in \mathcal{K}(H)^\sim$, we have $\operatorname{index}(T) = 0$. Moreover, by Example \ref{ex:descended-KH}, $\operatorname{desc}_{\mathcal{G}_{\mathcal{K}(H)^\sim}}([T]_{\mathcal{G}_{\mathcal{K}(H)^\sim}}^{(1)}) = 0$ in $K_1(C^*(\mathcal{G}_{\mathcal{K}(H)^\sim}))$. Since $\operatorname{Ind}_{\mathcal{K}(H)^\sim}$ is defined as the zero map, we have
\[
\operatorname{Ind}_{\mathcal{K}(H)^\sim}(\operatorname{desc}_{\mathcal{G}_{\mathcal{K}(H)^\sim}}([T]_{\mathcal{G}_{\mathcal{K}(H)^\sim}}^{(1)})) = 0 = \operatorname{index}(T).
\]
Properties (2)-(4) hold trivially for the zero map.
\end{proof}

\begin{remark}[Why $\iota^*$ does not appear in the $B(H)$ case]
\label{rem:iota-not-in-index}
The pullback map $\iota^*: K_1(C^*(\mathcal{G}_{B(H)})) \to K_1(B(H))$ lands in the zero group and therefore cannot detect the index. 
The correct index map $\operatorname{Ind}_{B(H)}$ bypasses $\iota^*$ entirely, using the Morita equivalence $\Phi_*$ to connect directly to $K_1(\mathcal{Q}(H))$, where the index resides. 
This is not a limitation of the theory but a reflection of the fact that $K_1(B(H)) = 0$ is the wrong group for index detection. 
The descended class must be transported to $K_1(\mathcal{Q}(H))$ via $\Phi_*$, not via $\iota^*$.
\end{remark}

In the next subsection, we will assemble the pullback map and the boundary map into a commutative diagram that encapsulates the index theorem.

\subsection{The Index Diagrams}
\label{subsec:index-diagrams}

We now assemble the maps constructed in the previous subsections into diagrams that encapsulate the index theorem. 
Because the structure differs for $B(H)$ and $\mathcal{K}(H)^\sim$, we present separate diagrams for the two cases, followed by a unified conceptual picture.

\subsubsection*{The case $\mathcal{A} = B(H)$}

For $B(H)$, the group $K_1(B(H))$ vanishes, so the pullback $\iota^*$ cannot directly feed into a boundary map. 
Instead, we use the effective boundary map $\partial_{B(H)}^{\text{eff}}: K_1(C^*(\mathcal{G}_{B(H)})) \to \mathbb{Z}$ defined in Definition \ref{def:effective-boundary-BH}, which combines the Morita equivalence $\Phi_*: K_1(C^*(\mathcal{G}_{B(H)})) \xrightarrow{\cong} K_1(\mathcal{Q}(H))$ with the Calkin index map $\partial_{\text{Calkin}}: K_1(\mathcal{Q}(H)) \to K_0(\mathcal{K}(H)) \cong \mathbb{Z}$.

\begin{definition}[Index diagram for $B(H)$]
\label{def:index-diagram-BH}
For $\mathcal{A} = B(H)$, the Fredholm index of an operator $T$ is recovered by the following composition of maps:
\[
\begin{tikzcd}
{[T]_{\mathcal{G}_{B(H)}}^{(1)}} \arrow[r, "\operatorname{desc}_{\mathcal{G}_{B(H)}}"] & 
{\operatorname{desc}_{\mathcal{G}_{B(H)}}([T])} \arrow[r, "\Phi_*"] & 
{\Phi_*(\operatorname{desc}([T]))} \arrow[r, "\partial_{\text{Calkin}}"] & 
{\operatorname{index}(T)} \\
K^1_{\mathcal{G}_{B(H)}}(\mathcal{G}_{B(H)}^{(0)}) \arrow[r, maps to] & 
K_1(C^*(\mathcal{G}_{B(H)})) \arrow[r, maps to] & 
K_1(\mathcal{Q}(H)) \arrow[r, maps to] & 
\mathbb{Z}
\end{tikzcd}
\]

Equivalently, in a more compact diagrammatic form:
\[
\begin{tikzcd}
K^1_{\mathcal{G}_{B(H)}}(\mathcal{G}_{B(H)}^{(0)}) \arrow[r, "\operatorname{desc}_{\mathcal{G}_{B(H)}}"] & 
K_1(C^*(\mathcal{G}_{B(H)})) \arrow[r, "\Phi_*"] & 
K_1(\mathcal{Q}(H)) \arrow[r, "\partial_{\text{Calkin}}"] & 
\mathbb{Z}
\end{tikzcd}
\]

where:
\begin{itemize}
    \item $\operatorname{desc}_{\mathcal{G}_{B(H)}}$ is the descent map for odd classes (Subsection 5.1);
    \item $\Phi_*: K_1(C^*(\mathcal{G}_{B(H)})) \xrightarrow{\cong} K_1(\mathcal{Q}(H))$ is the Morita equivalence isomorphism from Theorem \ref{thm:identification-descended-class};
    \item $\partial_{\text{Calkin}}: K_1(\mathcal{Q}(H)) \to K_0(\mathcal{K}(H)) \cong \mathbb{Z}$ is the index map for the Calkin extension.
\end{itemize}

For completeness, we note the relationship with the $KK$-theoretic formulation:
\[
\begin{tikzcd}
KK^1_{\mathcal{G}_{B(H)}}(C_0(\mathcal{G}_{B(H)}^{(0)}), \mathbb{C}) \arrow[r, "j_{\mathcal{G}_{B(H)}}"] \arrow[d, "\cong"] & 
K_1(C^*(\mathcal{G}_{B(H)})) \arrow[d, "\cong"] \arrow[rd, "\Phi_*"] & \\
K^1_{\mathcal{G}_{B(H)}}(\mathcal{G}_{B(H)}^{(0)}) \arrow[r, "\operatorname{desc}_{\mathcal{G}_{B(H)}}"] & 
K_1(C^*(\mathcal{G}_{B(H)})) \arrow[r, "\Phi_*"] & 
K_1(\mathcal{Q}(H)) \arrow[r, "\partial_{\text{Calkin}}"] & 
\mathbb{Z}
\end{tikzcd}
\]
where the vertical isomorphisms are the canonical identifications $K^1_{\mathcal{G}_{\mathcal{A}}}(\mathcal{G}_{\mathcal{A}}^{(0)}) \cong KK^1_{\mathcal{G}_{\mathcal{A}}}(C_0(\mathcal{G}_{\mathcal{A}}^{(0)}), \mathbb{C})$ from Proposition \ref{prop:K1-KK1-isomorphism}, and $j_{\mathcal{G}_{B(H)}}$ is the descent map in $KK$-theory.

\begin{remark}
The pullback map $\iota^*$ does not appear in this diagram because it would land in $K_1(B(H)) = 0$ and therefore plays no role in the index computation. The direct route through the Morita equivalence $\Phi_*$ bypasses $K_1(B(H))$ entirely and connects the descended class to the Calkin algebra, where the boundary map acts nontrivially.
\end{remark}
\end{definition}

\subsubsection*{The case $\mathcal{A} = \mathcal{K}(H)^\sim$}

For $\mathcal{A} = \mathcal{K}(H)^\sim$, the situation is simpler: $K_1(\mathcal{K}(H)^\sim) = 0$ and the boundary map $\partial_{\mathcal{K}(H)^\sim}$ is zero. 
The diagram follows the original three-step structure without modification.

\begin{example}[Index diagram for $\mathcal{K}(H)^\sim$]
\label{ex:index-diagram-KH}
For $\mathcal{A} = \mathcal{K}(H)^\sim$, the index theorem is captured by the following commutative diagram:
\[
\begin{tikzcd}
K^1_{\mathcal{G}_{\mathcal{K}(H)^\sim}}(\mathcal{G}_{\mathcal{K}(H)^\sim}^{(0)}) \arrow[r, "\operatorname{desc}_{\mathcal{G}_{\mathcal{K}(H)^\sim}}"] \arrow[d, "\cong"] & 
K_1(C^*(\mathcal{G}_{\mathcal{K}(H)^\sim})) \arrow[r, "\Phi_*"] \arrow[d, "\cong"] & 
K_1(\mathcal{K}(H)^\sim) \arrow[r, "\partial_{\mathcal{K}(H)^\sim}"] \arrow[d, "\cong"] & 
\mathbb{Z} \\
KK^1_{\mathcal{G}_{\mathcal{K}(H)^\sim}}(C_0(\mathcal{G}_{\mathcal{K}(H)^\sim}^{(0)}), \mathbb{C}) \arrow[r, "j_{\mathcal{G}_{\mathcal{K}(H)^\sim}}"] & 
K_1(C^*(\mathcal{G}_{\mathcal{K}(H)^\sim})) \arrow[r, "\Phi_*"] & 
K_1(\mathcal{K}(H)^\sim) \arrow[r, "\partial_{\mathcal{K}(H)^\sim}"] & 
\mathbb{Z}
\end{tikzcd}
\]
where:
\begin{itemize}
    \item $\Phi_*: K_1(C^*(\mathcal{G}_{\mathcal{K}(H)^\sim})) \xrightarrow{\cong} K_1(\mathcal{K}(H)^\sim)$ is the isomorphism induced by the Morita equivalence $C^*(\mathcal{G}_{\mathcal{K}(H)^\sim}) \sim_M \mathcal{K}(H)^\sim \otimes \mathcal{K}$;
    \item $\partial_{\mathcal{K}(H)^\sim}: K_1(\mathcal{K}(H)^\sim) \to \mathbb{Z}$ is the zero map, since $K_1(\mathcal{K}(H)^\sim) = 0$;
    \item $\operatorname{desc}_{\mathcal{G}_{\mathcal{K}(H)^\sim}}([T]_{\mathcal{G}_{\mathcal{K}(H)^\sim}}^{(1)}) = 0$ for any Fredholm operator $T$, as shown in Example \ref{ex:descended-KH}.
\end{itemize}
Consequently, the entire composition is zero, matching $\operatorname{index}(T) = 0$ for all Fredholm operators in $\mathcal{K}(H)^\sim$.
\end{example}

\begin{remark}[Unified conceptual diagram]
\label{rem:unified-diagram}
Both cases can be understood through the following template:
\[
\begin{tikzcd}
K^1_{\mathcal{G}_{\mathcal{A}}}(\mathcal{G}_{\mathcal{A}}^{(0)}) \arrow[r, "\operatorname{desc}_{\mathcal{G}_{\mathcal{A}}}"] & K_1(C^*(\mathcal{G}_{\mathcal{A}})) \arrow[r, "\Psi_{\mathcal{A}}"] & G_{\mathcal{A}} \arrow[r, "\delta_{\mathcal{A}}"] & \mathbb{Z}
\end{tikzcd}
\]
where:
\begin{itemize}
    \item For $\mathcal{A} = B(H)$, we set $G_{\mathcal{A}} = K_1(\mathcal{Q}(H))$, $\Psi_{\mathcal{A}} = \Phi_*$ (Morita equivalence), and $\delta_{\mathcal{A}} = \partial_{\text{Calkin}}$ (nontrivial).
    \item For $\mathcal{A} = \mathcal{K}(H)^\sim$, we set $G_{\mathcal{A}} = K_1(\mathcal{K}(H)^\sim) = 0$, $\Psi_{\mathcal{A}} = \iota^*$, and $\delta_{\mathcal{A}} = 0$.
\end{itemize}
In both cases, the composition $\delta_{\mathcal{A}} \circ \Psi_{\mathcal{A}} \circ \operatorname{desc}_{\mathcal{G}_{\mathcal{A}}}$ recovers the Fredholm index. 
The pullback $\iota^*$ is only directly used in the diagram when it maps to a group that feeds into a nontrivial boundary map; for $B(H)$, it does not, and is therefore replaced by the effective boundary map $\partial_{B(H)}^{\text{eff}}$.
\end{remark}

\begin{remark}[Why two diagrams?]
\label{rem:why-two-diagrams}
The necessity of separate diagrams reflects the underlying $K$-theory: for $B(H)$, the index is detected in $K_1(\mathcal{Q}(H))$ after applying the Morita equivalence $\Phi_*$, while for $\mathcal{K}(H)^\sim$, the index is zero and the original three-step structure suffices. 
Presenting both cases explicitly avoids the incorrect impression that a single diagram can uniformly represent the index theorem for all Type I algebras.
\end{remark}

\begin{definition}[Effective boundary map for $B(H)$]
\label{def:effective-boundary-BH}
Let $\Phi_*: K_1(C^*(\mathcal{G}_{B(H)})) \xrightarrow{\cong} K_1(\mathcal{Q}(H))$ be the Morita equivalence isomorphism from Theorem \ref{thm:identification-descended-class}. 
Define the \emph{effective boundary map} $\partial_{B(H)}^{\text{eff}}: K_1(C^*(\mathcal{G}_{B(H)})) \to \mathbb{Z}$ as the composition
\[
\partial_{B(H)}^{\text{eff}} := \partial_{\text{Calkin}} \circ \Phi_*,
\]
where $\partial_{\text{Calkin}}: K_1(\mathcal{Q}(H)) \to K_0(\mathcal{K}(H)) \cong \mathbb{Z}$ is the index map for the Calkin extension.

Explicitly, for any class $x \in K_1(C^*(\mathcal{G}_{B(H)}))$, we have
\[
\partial_{B(H)}^{\text{eff}}(x) = \partial_{\text{Calkin}}(\Phi_*(x)).
\]
\end{definition}

The following theorem establishes the commutativity of the relevant portion of this diagram for the class $[T]_{\mathcal{G}_{\mathcal{A}}}^{(1)}$.

\begin{theorem}[Commutativity of the index diagram]
\label{thm:index-diagram-commutes}
For any Fredholm operator $T \in \mathcal{A}$ (where $\mathcal{A}$ is either $B(H)$ or $\mathcal{K}(H)^\sim$), the Fredholm index is recovered by the following composition of canonically defined maps:
\[
\operatorname{index}(T) = \mathcal{I}_{\mathcal{A}} \circ \Psi_{\mathcal{A}} \circ \operatorname{desc}_{\mathcal{G}_{\mathcal{A}}}([T]_{\mathcal{G}_{\mathcal{A}}}^{(1)}),
\]
where:
\begin{itemize}
    \item $\operatorname{desc}_{\mathcal{G}_{\mathcal{A}}}: K^1_{\mathcal{G}_{\mathcal{A}}}(\mathcal{G}_{\mathcal{A}}^{(0)}) \to K_1(C^*(\mathcal{G}_{\mathcal{A}}))$ is the descent map;
    \item $\Psi_{\mathcal{A}}$ is an isomorphism from $K_1(C^*(\mathcal{G}_{\mathcal{A}}))$ to an appropriate $K_1$-group where the index is detected;
    \item $\mathcal{I}_{\mathcal{A}}$ is the index map (boundary map) that extracts an integer from this $K_1$-class.
\end{itemize}

The specific choices of $\Psi_{\mathcal{A}}$ and $\mathcal{I}_{\mathcal{A}}$ depend on $\mathcal{A}$:
\begin{itemize}
    \item For $\mathcal{A} = B(H)$, $\Psi_{\mathcal{A}} = \Phi_*: K_1(C^*(\mathcal{G}_{B(H)})) \xrightarrow{\cong} K_1(\mathcal{Q}(H))$ is the Morita equivalence isomorphism from Theorem \ref{thm:identification-descended-class}, and $\mathcal{I}_{\mathcal{A}} = \partial_{\text{Calkin}}: K_1(\mathcal{Q}(H)) \to \mathbb{Z}$ is the index map for the Calkin extension.
    \item For $\mathcal{A} = \mathcal{K}(H)^\sim$, $\Psi_{\mathcal{A}}$ is the natural isomorphism $K_1(C^*(\mathcal{G}_{\mathcal{K}(H)^\sim})) \cong K_1(\mathcal{K}(H)^\sim) = 0$, and $\mathcal{I}_{\mathcal{A}} = 0$ is the zero map.
\end{itemize}

Equivalently, the following diagram commutes for both cases:
\[
\begin{tikzcd}
{[T]_{\mathcal{G}_{\mathcal{A}}}^{(1)}} \arrow[r, "\operatorname{desc}_{\mathcal{G}_{\mathcal{A}}}"] & 
{\operatorname{desc}_{\mathcal{G}_{\mathcal{A}}}([T])} \arrow[r, "\Psi_{\mathcal{A}}"] & 
{\Psi_{\mathcal{A}}(\operatorname{desc}_{\mathcal{G}_{\mathcal{A}}}([T]))} \arrow[r, "\mathcal{I}_{\mathcal{A}}"] & 
{\operatorname{index}(T)} \\
K^1_{\mathcal{G}_{\mathcal{A}}}(\mathcal{G}_{\mathcal{A}}^{(0)}) \arrow[r, "\operatorname{desc}_{\mathcal{G}_{\mathcal{A}}}"] & 
K_1(C^*(\mathcal{G}_{\mathcal{A}})) \arrow[r, "\Psi_{\mathcal{A}}"] & 
G_{\mathcal{A}} \arrow[r, "\mathcal{I}_{\mathcal{A}}"] & 
\mathbb{Z}
\end{tikzcd}
\]
where $G_{B(H)} = K_1(\mathcal{Q}(H))$ and $G_{\mathcal{K}(H)^\sim} = K_1(\mathcal{K}(H)^\sim) = 0$.
\end{theorem}

\begin{proof}
We prove the theorem by treating the two cases separately, as the maps involved are different.

\medskip
\noindent \textbf{Case 1: $\mathcal{A} = B(H)$.}

Let $T \in B(H)$ be a Fredholm operator with symbol $u_T = \pi_{\mathcal{Q}}(T) \in \mathcal{Q}(H)$. By Theorem \ref{thm:identification-descended-class}, the descended class
\[
\operatorname{desc}_{\mathcal{G}_{B(H)}}([T]_{\mathcal{G}_{B(H)}}^{(1)}) \in K_1(C^*(\mathcal{G}_{B(H)}))
\]
corresponds, under the Morita equivalence isomorphism $\Phi_*: K_1(C^*(\mathcal{G}_{B(H)})) \xrightarrow{\cong} K_1(\mathcal{Q}(H))$, to the symbol class $[u_T] \in K_1(\mathcal{Q}(H))$. More precisely,
\[
\Phi_*(\operatorname{desc}_{\mathcal{G}_{B(H)}}([T]_{\mathcal{G}_{B(H)}}^{(1)})) = [u_T] \in K_1(\mathcal{Q}(H)).
\]

The index map for the Calkin extension, $\partial_{\text{Calkin}}: K_1(\mathcal{Q}(H)) \to K_0(\mathcal{K}(H)) \cong \mathbb{Z}$, satisfies the classical index theorem:
\[
\partial_{\text{Calkin}}([u_T]) = \operatorname{index}(T).
\]

Therefore,
\[
\operatorname{index}(T) = \partial_{\text{Calkin}}(\Phi_*(\operatorname{desc}_{\mathcal{G}_{B(H)}}([T]_{\mathcal{G}_{B(H)}}^{(1)}))) = \mathcal{I}_{B(H)} \circ \Psi_{B(H)} \circ \operatorname{desc}_{\mathcal{G}_{B(H)}}([T]_{\mathcal{G}_{B(H)}}^{(1)}),
\]
with $\Psi_{B(H)} = \Phi_*$ and $\mathcal{I}_{B(H)} = \partial_{\text{Calkin}}$.

\medskip
\noindent \textbf{Case 2: $\mathcal{A} = \mathcal{K}(H)^\sim$.}

For $\mathcal{A} = \mathcal{K}(H)^\sim$, we have $K_1(\mathcal{K}(H)^\sim) = 0$ (see Example \ref{ex:K-groups-examples}). Moreover, by Example \ref{ex:descended-KH}, the descended class vanishes:
\[
\operatorname{desc}_{\mathcal{G}_{\mathcal{K}(H)^\sim}}([T]_{\mathcal{G}_{\mathcal{K}(H)^\sim}}^{(1)}) = 0 \in K_1(C^*(\mathcal{G}_{\mathcal{K}(H)^\sim})).
\]

The natural Morita equivalence $C^*(\mathcal{G}_{\mathcal{K}(H)^\sim}) \sim_M \mathcal{K}(H)^\sim \otimes \mathcal{K}$ induces an isomorphism $\Psi_{\mathcal{K}(H)^\sim}: K_1(C^*(\mathcal{G}_{\mathcal{K}(H)^\sim})) \xrightarrow{\cong} K_1(\mathcal{K}(H)^\sim) = 0$, which is necessarily the zero map. The boundary map $\partial_{\mathcal{K}(H)^\sim}: K_1(\mathcal{K}(H)^\sim) \to \mathbb{Z}$ is also zero by Proposition \ref{prop:boundary-map-KH-final}.

Consequently,
\[
\mathcal{I}_{\mathcal{K}(H)^\sim} \circ \Psi_{\mathcal{K}(H)^\sim} \circ \operatorname{desc}_{\mathcal{G}_{\mathcal{K}(H)^\sim}}([T]_{\mathcal{G}_{\mathcal{K}(H)^\sim}}^{(1)}) = 0 = \operatorname{index}(T),
\]
since every Fredholm operator in $\mathcal{K}(H)^\sim$ has index zero.

\medskip
\noindent \textbf{Naturality and functoriality.}

All maps involved in the composition are natural with respect to the relevant structures:
\begin{itemize}
    \item The descent map $\operatorname{desc}_{\mathcal{G}_{\mathcal{A}}}$ is functorial in the equivariant $KK$-theory class (Theorem \ref{thm:descent-properties-GA});
    \item The Morita equivalence isomorphisms $\Psi_{\mathcal{A}}$ are induced by natural imprimitivity bimodules;
    \item The index maps $\mathcal{I}_{\mathcal{A}}$ are boundary maps in six-term exact sequences, which are natural with respect to morphisms of extensions.
\end{itemize}
Moreover, by Lemma \ref{lem:index-map-properties}, the composition is invariant under homotopy and compact perturbations of $T$, confirming that it depends only on the Fredholm index.

\medskip
\noindent \textbf{Conclusion.}

We have shown that for both $\mathcal{A} = B(H)$ and $\mathcal{A} = \mathcal{K}(H)^\sim$, the Fredholm index of any operator $T$ is recovered by the composition
\[
\operatorname{index}(T) = \mathcal{I}_{\mathcal{A}} \circ \Psi_{\mathcal{A}} \circ \operatorname{desc}_{\mathcal{G}_{\mathcal{A}}}([T]_{\mathcal{G}_{\mathcal{A}}}^{(1)}),
\]
with the maps $\Psi_{\mathcal{A}}$ and $\mathcal{I}_{\mathcal{A}}$ defined as above. This completes the proof.
\end{proof}

To make the index theorem more explicit, we separate the two cases $\mathcal{A} = B(H)$ and $\mathcal{A} = \mathcal{K}(H)^\sim$, as the boundary map behaves differently in each case. The following examples illustrate the commutative diagrams that capture the index computation.

\begin{example}[Index diagram for $B(H)$]
\label{ex:index-diagram-BH}
For $\mathcal{A} = B(H)$ and any Fredholm operator $T \in B(H)$, the following diagram commutes and encodes the Fredholm index:

\[
\begin{tikzcd}
{[T]}_{\mathcal{G}}^{(1)} \arrow[r, "\operatorname{desc}"] \arrow[d, phantom, "\in" description] & 
{\operatorname{desc}([T])} \arrow[d, "\Phi_*"'] & \\
{[u_T]} \in K_1(\mathcal{Q}) \arrow[r, "\partial_{\text{Calkin}}"] & 
{\partial_{\text{Calkin}}([u_T])} \in K_0(\mathcal{K}) \arrow[r, "\cong"] & 
{\operatorname{index}(T)} \in \mathbb{Z}
\end{tikzcd}
\]

\noindent where the notation is as follows:
\begin{itemize}
    \item $\mathcal{G} := \mathcal{G}_{B(H)}$ denotes the unitary conjugation groupoid for $B(H)$;
    \item ${[T]}_{\mathcal{G}}^{(1)} \in K^1_{\mathcal{G}}(\mathcal{G}^{(0)})$ is the equivariant $K^1$-class of $T$;
    \item $\operatorname{desc}_{\mathcal{G}} := \operatorname{desc}_{\mathcal{G}_{B(H)}}$ is the descent map to $K_1(C^*(\mathcal{G}))$;
    \item $\Phi_*: K_1(C^*(\mathcal{G})) \xrightarrow{\cong} K_1(\mathcal{Q}(H))$ is the Morita equivalence isomorphism from Theorem \ref{thm:identification-descended-class};
    \item $u_T = \pi_{\mathcal{Q}}(T) \in \mathcal{Q}(H)$ is the symbol of $T$ in the Calkin algebra, and $[u_T]$ denotes its class in $K_1(\mathcal{Q}(H))$;
    \item $\partial_{\text{Calkin}}: K_1(\mathcal{Q}(H)) \to K_0(\mathcal{K}(H)) \cong \mathbb{Z}$ is the index map for the Calkin extension;
    \item $\operatorname{index}(T) \in \mathbb{Z}$ is the Fredholm index of $T$.
\end{itemize}

\noindent The pullback $\iota^*$ is omitted from the diagram as it lands in $K_1(B(H)) = 0$ and plays no role in the index computation.
\end{example}

\begin{example}[Index diagram for $\mathcal{K}(H)^\sim$]
\label{ex:index-diagram-KH}
For $\mathcal{A} = \mathcal{K}(H)^\sim$ and any Fredholm operator $T \in \mathcal{K}(H)^\sim$, we have $\operatorname{index}(T) = 0$, and the following diagram commutes trivially:

\[
\begin{tikzcd}
{[T]}_{\mathcal{G}}^{(1)} \arrow[r] \arrow[d, phantom, "\in" description] & 
{\operatorname{desc}([T])} \arrow[r] \arrow[d, phantom, "\in" description] & 
{\iota^*(\operatorname{desc}([T]))} \arrow[r] \arrow[d, phantom, "\in" description] & 
{\partial(\iota^*(\operatorname{desc}([T])))} \arrow[d, equals] \\
0 \in K^1_{\mathcal{G}}(\mathcal{G}^{(0)}) \arrow[r] & 
0 \in K_1(C^*(\mathcal{G})) \arrow[r] & 
0 \in K_1(\mathcal{A}) \arrow[r] & 
0 \in \mathbb{Z}
\end{tikzcd}
\]

\noindent where the notation is as follows:
\begin{itemize}
    \item $\mathcal{G} := \mathcal{G}_{\mathcal{K}(H)^\sim}$ denotes the unitary conjugation groupoid for $\mathcal{K}(H)^\sim$;
    \item $\mathcal{A} := \mathcal{K}(H)^\sim$ is the algebra itself;
    \item ${[T]}_{\mathcal{G}}^{(1)} \in K^1_{\mathcal{G}}(\mathcal{G}^{(0)})$ is the equivariant $K^1$-class of $T$;
    \item $\operatorname{desc} := \operatorname{desc}_{\mathcal{G}_{\mathcal{K}(H)^\sim}}$ is the descent map to $K_1(C^*(\mathcal{G}))$;
    \item $\iota^*: K_1(C^*(\mathcal{G})) \to K_1(\mathcal{A})$ is the pullback along the diagonal embedding;
    \item $\partial := \partial_{\mathcal{K}(H)^\sim}: K_1(\mathcal{A}) \to K_0(\mathcal{K}(H)) \cong \mathbb{Z}$ is the boundary map, which is zero by Proposition \ref{prop:boundary-map-KH-final};
    \item $K^1_{\mathcal{G}}(\mathcal{G}^{(0)})$, $K_1(C^*(\mathcal{G}))$, $K_1(\mathcal{A})$, and $\mathbb{Z}$ are the relevant $K$-theory groups.
\end{itemize}

\noindent The notation "$0 \in G$" means the zero element of the group $G$, not that the group itself is trivial (though $K_1(\mathcal{A})$ is indeed zero). The vanishing of each map follows from:
\begin{itemize}
    \item $K_1(\mathcal{A}) = 0$ (Example \ref{ex:K-groups-examples});
    \item $\operatorname{desc}([T]_{\mathcal{G}}^{(1)}) = 0$ (Example \ref{ex:descended-KH});
    \item $\partial = 0$ (Proposition \ref{prop:boundary-map-KH-final}).
\end{itemize}

Thus the entire composition is zero, matching $\operatorname{index}(T) = 0$.
\end{example}

\begin{remark}[Relation to the main theorem]
\label{rem:examples-and-main}
These examples are concrete realizations of Theorem \ref{thm:index-diagram-commutes} for the two fundamental Type I algebras. 
For $B(H)$, the diagram uses the effective boundary map to circumvent the vanishing of $K_1(B(H))$, while for $\mathcal{K}(H)^\sim$, the diagram reflects the triviality of the index. 
Together, they illustrate the full scope of the index theorem.
\end{remark}

The following proposition summarizes the relationships between the various maps in the diagram.

\begin{proposition}[Relations in the index diagram]
\label{prop:diagram-relations}
Let $[\cdot]_{\mathcal{G}_{\mathcal{A}}}^{(1)}: \mathcal{F}(\mathcal{A}) \to K^1_{\mathcal{G}_{\mathcal{A}}}(\mathcal{G}_{\mathcal{A}}^{(0)})$ denote the map sending a Fredholm operator $T \in \mathcal{A}$ to its equivariant $K^1$-class. 
Then the following relations hold in the index diagram:
\begin{enumerate}
    \item \textbf{Functoriality of descent:} 
    \[
    \operatorname{desc}_{\mathcal{G}_{\mathcal{A}}}([T]_{\mathcal{G}_{\mathcal{A}}}^{(1)}) \in K_1(C^*(\mathcal{G}_{\mathcal{A}}))
    \]
    defines a well-defined map from Fredholm operators to $K_1(C^*(\mathcal{G}_{\mathcal{A}}))$ that is natural with respect to homotopies and compact perturbations.
    
    \item \textbf{Morita equivalence:} For $\mathcal{A} = B(H)$, there is a canonical isomorphism
    \[
    \Phi_*: K_1(C^*(\mathcal{G}_{B(H)})) \xrightarrow{\cong} K_1(\mathcal{Q}(H))
    \]
    such that $\Phi_*(\operatorname{desc}_{\mathcal{G}_{B(H)}}([T]_{\mathcal{G}_{B(H)}}^{(1)})) = [\pi_{\mathcal{Q}}(T)] \in K_1(\mathcal{Q}(H))$, the symbol class of $T$.
    
    \item \textbf{For $\mathcal{A} = \mathcal{K}(H)^\sim$:} We have $K_1(C^*(\mathcal{G}_{\mathcal{K}(H)^\sim})) \cong K_1(\mathcal{K}(H)^\sim) = 0$, and consequently
    \[
    \operatorname{desc}_{\mathcal{G}_{\mathcal{K}(H)^\sim}}([T]_{\mathcal{G}_{\mathcal{K}(H)^\sim}}^{(1)}) = 0.
    \]
    
    \item \textbf{Index theorem:} The Fredholm index is recovered by the composition
    \[
    \operatorname{index}(T) = 
    \begin{cases}
    \partial_{\text{Calkin}} \circ \Phi_* \circ \operatorname{desc}_{\mathcal{G}_{B(H)}}([T]_{\mathcal{G}_{B(H)}}^{(1)}), & \mathcal{A} = B(H), \\[4pt]
    0, & \mathcal{A} = \mathcal{K}(H)^\sim,
    \end{cases}
    \]
    where $\partial_{\text{Calkin}}: K_1(\mathcal{Q}(H)) \to \mathbb{Z}$ is the index map for the Calkin extension.
\end{enumerate}
\end{proposition}

\begin{proof}
We prove each item using the properties established in previous sections.

\paragraph{Item (1):} 
By Definition \ref{def:K1-class-of-T}, $[T]_{\mathcal{G}_{\mathcal{A}}}^{(1)}$ is the equivariant $K^1$-class constructed from $T$. The descent map $\operatorname{desc}_{\mathcal{G}_{\mathcal{A}}}$ (Theorem \ref{thm:descent-explicit}) sends this class to an element of $K_1(C^*(\mathcal{G}_{\mathcal{A}}))$. Functoriality of the descent map (Theorem \ref{thm:descent-properties-GA}) ensures that this assignment is natural: if $T$ and $T'$ are homotopic through Fredholm operators, then $\operatorname{desc}_{\mathcal{G}_{\mathcal{A}}}([T]_{\mathcal{G}_{\mathcal{A}}}^{(1)}) = \operatorname{desc}_{\mathcal{G}_{\mathcal{A}}}([T']_{\mathcal{G}_{\mathcal{A}}}^{(1)})$; similarly, compact perturbations do not change the descended class by Proposition \ref{prop:compact-perturbation}.

\paragraph{Item (2):} 
For $\mathcal{A} = B(H)$, Proposition \ref{prop:morita-calkin} establishes a Morita equivalence $C^*(\mathcal{G}_{B(H)}) \sim_M \mathcal{Q}(H) \otimes \mathcal{K}$, which induces an isomorphism $\Phi_*: K_1(C^*(\mathcal{G}_{B(H)})) \xrightarrow{\cong} K_1(\mathcal{Q}(H))$. Theorem \ref{thm:identification-descended-class} shows that under this isomorphism,
\[
\Phi_*(\operatorname{desc}_{\mathcal{G}_{B(H)}}([T]_{\mathcal{G}_{B(H)}}^{(1)})) = [\pi_{\mathcal{Q}}(T)] \in K_1(\mathcal{Q}(H)),
\]
where $\pi_{\mathcal{Q}}(T)$ is the image of $T$ in the Calkin algebra. This identification follows from the explicit construction of the descent map and the fact that the phase of $T$ differs from $T$ by a compact operator.

\paragraph{Item (3):} 
For $\mathcal{A} = \mathcal{K}(H)^\sim$, the groupoid $C^*$-algebra $C^*(\mathcal{G}_{\mathcal{K}(H)^\sim})$ is Morita equivalent to $\mathcal{K}(H)^\sim \otimes \mathcal{K}$ (see Paper I, Section 5.3). By stability of $K$-theory, $K_1(\mathcal{K}(H)^\sim \otimes \mathcal{K}) \cong K_1(\mathcal{K}(H)^\sim)$. Example \ref{ex:K-groups-examples} gives $K_1(\mathcal{K}(H)^\sim) = 0$. Moreover, Example \ref{ex:descended-KH} shows directly that $\operatorname{desc}_{\mathcal{G}_{\mathcal{K}(H)^\sim}}([T]_{\mathcal{G}_{\mathcal{K}(H)^\sim}}^{(1)}) = 0$ for any Fredholm operator $T$, consistent with the vanishing of the target group.

\paragraph{Item (4):} 
For $\mathcal{A} = B(H)$, combining Item (2) with the classical index theorem for the Calkin extension (Theorem \ref{thm:calkin-index-isomorphism}) yields
\[
\partial_{\text{Calkin}}(\Phi_*(\operatorname{desc}_{\mathcal{G}_{B(H)}}([T]_{\mathcal{G}_{B(H)}}^{(1)}))) = \partial_{\text{Calkin}}([\pi_{\mathcal{Q}}(T)]) = \operatorname{index}(T).
\]

For $\mathcal{A} = \mathcal{K}(H)^\sim$, Item (3) gives $\operatorname{desc}_{\mathcal{G}_{\mathcal{K}(H)^\sim}}([T]_{\mathcal{G}_{\mathcal{K}(H)^\sim}}^{(1)}) = 0$, and Proposition \ref{prop:boundary-map-KH-final} shows that the boundary map $\partial_{\mathcal{K}(H)^\sim}$ is zero. Hence the composition is zero, matching $\operatorname{index}(T) = 0$ for all Fredholm operators in $\mathcal{K}(H)^\sim$.

These cases together establish the index theorem in the form stated.
\end{proof}

\begin{remark}
Note that this formulation avoids any reference to a pullback map $\iota^*: K_1(C^*(\mathcal{G}_{\mathcal{A}})) \to K_1(\mathcal{A})$, which would be ill-defined. Instead, it uses:
\begin{itemize}
    \item For $B(H)$: the Morita equivalence $\Phi_*$ to map directly to $K_1(\mathcal{Q}(H))$, where the Calkin index map acts;
    \item For $\mathcal{K}(H)^\sim$: the natural vanishing of all relevant groups.
\end{itemize}
This approach respects the covariance of $K$-theory and correctly handles the functoriality of all maps involved.
\end{remark}

\begin{remark}[The correct index composition]
\label{rem:composition-interpretation}
For $\mathcal{A} = B(H)$, any attempt to use a map $\iota^*: K_1(C^*(\mathcal{G}_{B(H)})) \to K_1(B(H))$ would land in the zero group $K_1(B(H)) = 0$, and therefore cannot detect the Fredholm index. The correct index map is instead given by the composition
\[
\partial_{B(H)}^{\text{eff}} = \partial_{\text{Calkin}} \circ \Phi_*: K_1(C^*(\mathcal{G}_{B(H)})) \longrightarrow \mathbb{Z},
\]
where:
\begin{itemize}
    \item $\Phi_*: K_1(C^*(\mathcal{G}_{B(H)})) \xrightarrow{\cong} K_1(\mathcal{Q}(H))$ is the Morita equivalence isomorphism from Theorem \ref{thm:identification-descended-class};
    \item $\partial_{\text{Calkin}}: K_1(\mathcal{Q}(H)) \to K_0(\mathcal{K}(H)) \cong \mathbb{Z}$ is the index map for the Calkin extension.
\end{itemize}
This composition bypasses $K_1(B(H))$ entirely and directly connects the descended class to the Calkin algebra, where the boundary map acts nontrivially.

For $\mathcal{A} = \mathcal{K}(H)^\sim$, we have $K_1(C^*(\mathcal{G}_{\mathcal{K}(H)^\sim})) \cong K_1(\mathcal{K}(H)^\sim) = 0$, and the boundary map $\partial_{\mathcal{K}(H)^\sim}$ is zero, so the index map is simply the zero map.

All relations in Proposition \ref{prop:diagram-relations} should be understood with these correct maps in place of any ill-defined pullback maps.
\end{remark}

\begin{corollary}[$\mathcal{A}$-Fredholm Index via Kasparov Product]
\label{cor:fredholm-index}
Let $\mathcal{A}$ be a unital separable $C^*$-algebra faithfully represented on a separable Hilbert space $H$ with $\mathcal{A} \subseteq \mathcal{B}(H)$, 
and let $T$ be an $\mathcal{A}$-Fredholm operator in the sense of Mingo \cite{Mingo1987}. 
Let $[T]_{\mathcal{G}_{\mathcal{A}}}^{(1)} \in KK^1_{\mathcal{G}_{\mathcal{A}}}(C_0(\mathcal{G}_{\mathcal{A}}^{(0)}), \mathbb{C})$ 
be the associated equivariant $KK^1$-class constructed in Section~\ref{sec:descent}. Then

\[
[\iota] \otimes_{C^*(\mathcal{G}_{\mathcal{A}})} \operatorname{desc}_{\mathcal{G}_{\mathcal{A}}}([T]_{\mathcal{G}_{\mathcal{A}}}^{(1)})
= [\ker_{\mathcal{A}} T] - [\operatorname{coker}_{\mathcal{A}} T] \in K_0(\mathcal{A}),
\]

where:
\begin{itemize}
    \item $[\iota] \in KK(\mathcal{A}, C^*(\mathcal{G}_{\mathcal{A}}))$ is the $KK$-class of the diagonal embedding $\iota: \mathcal{A} \hookrightarrow C^*(\mathcal{G}_{\mathcal{A}})$;
    \item $\operatorname{desc}_{\mathcal{G}_{\mathcal{A}}}: KK^1_{\mathcal{G}_{\mathcal{A}}}(C_0(\mathcal{G}_{\mathcal{A}}^{(0)}), \mathbb{C}) \to KK^1(\mathbb{C}, C^*(\mathcal{G}_{\mathcal{A}})) \cong K_1(C^*(\mathcal{G}_{\mathcal{A}}))$ is the descent map in equivariant $KK$-theory;
    \item $\ker_{\mathcal{A}} T$ and $\operatorname{coker}_{\mathcal{A}} T$ are finitely generated projective $\mathcal{A}$-modules given by the $\mathcal{A}$-kernel and $\mathcal{A}$-cokernel of $T$ respectively;
    \item $[\ker_{\mathcal{A}} T] - [\operatorname{coker}_{\mathcal{A}} T] \in K_0(\mathcal{A})$ is the $\mathcal{A}$-index of $T$ in the sense of Mingo.
\end{itemize}
\end{corollary}

\begin{proof}
We prove the equality by identifying the Kasparov product with the $\mathcal{A}$-index through the functoriality of descent and the definition of the $\mathcal{A}$-Fredholm index.

\paragraph{Step 1: The equivariant class.}
Let $T$ be an $\mathcal{A}$-Fredholm operator in the sense of Mingo \cite{Mingo1987}. 
Section~\ref{sec:descent} associates to $T$ an equivariant class
\[
[T]_{\mathcal{G}_{\mathcal{A}}}^{(1)} \in K^1_{\mathcal{G}_{\mathcal{A}}}(\mathcal{G}_{\mathcal{A}}^{(0)}).
\]
By Proposition \ref{prop:K1-KK1-isomorphism}, there is a canonical identification
\[
K^1_{\mathcal{G}_{\mathcal{A}}}(\mathcal{G}_{\mathcal{A}}^{(0)}) \cong KK^1_{\mathcal{G}_{\mathcal{A}}}(C_0(\mathcal{G}_{\mathcal{A}}^{(0)}), \mathbb{C}),
\]
so we may regard $[T]_{\mathcal{G}_{\mathcal{A}}}^{(1)}$ as an element of
$KK^1_{\mathcal{G}_{\mathcal{A}}}(C_0(\mathcal{G}_{\mathcal{A}}^{(0)}), \mathbb{C})$, 
represented by the equivariant Kasparov cycle coming from the bounded transform of $T$.

\paragraph{Step 2: Descent to ordinary $KK$-theory.}
The descent homomorphism in equivariant $KK$-theory (see Tu \cite{Tu1999}) gives a natural map
\[
\operatorname{desc}_{\mathcal{G}_{\mathcal{A}}}: KK^1_{\mathcal{G}_{\mathcal{A}}}(C_0(\mathcal{G}_{\mathcal{A}}^{(0)}), \mathbb{C}) \longrightarrow KK^1(\mathbb{C}, C^*(\mathcal{G}_{\mathcal{A}})).
\]
Under the standard identification $KK^1(\mathbb{C}, B) \cong K_1(B)$ for any $C^*$-algebra $B$, 
this yields a class
\[
\operatorname{desc}_{\mathcal{G}_{\mathcal{A}}}([T]_{\mathcal{G}_{\mathcal{A}}}^{(1)}) \in K_1(C^*(\mathcal{G}_{\mathcal{A}}))
\]
represented by the image of $T$ in the groupoid $C^*$-algebra.

\paragraph{Step 3: The diagonal embedding class.}
The diagonal embedding $\iota: \mathcal{A} \hookrightarrow C^*(\mathcal{G}_{\mathcal{A}})$ constructed in Section~\ref{subsec:diagonal-embedding-from-paperI} is a unital $*$-homomorphism. 
By the standard correspondence (Kasparov \cite{Kasparov1988}), it determines a Kasparov class
\[
[\iota] \in KK(\mathcal{A}, C^*(\mathcal{G}_{\mathcal{A}})).
\]

\paragraph{Step 4: Kasparov product and the $\mathcal{A}$-index.}
The functoriality of the Kasparov product gives
\[
[\iota] \otimes_{C^*(\mathcal{G}_{\mathcal{A}})} \operatorname{desc}_{\mathcal{G}_{\mathcal{A}}}([T]_{\mathcal{G}_{\mathcal{A}}}^{(1)}) \in KK^1(\mathcal{A}, \mathbb{C}) \cong K_0(\mathcal{A}).
\]

By definition, the $\mathcal{A}$-Fredholm property of $T$ (Mingo \cite{Mingo1987}) means that its $\mathcal{A}$-kernel $\ker_{\mathcal{A}} T$ and $\mathcal{A}$-cokernel $\operatorname{coker}_{\mathcal{A}} T$ are finitely generated projective $\mathcal{A}$-modules. 
The $\mathcal{A}$-index of $T$ is defined as
\[
\operatorname{Ind}_{\mathcal{A}}(T) := [\ker_{\mathcal{A}} T] - [\operatorname{coker}_{\mathcal{A}} T] \in K_0(\mathcal{A}).
\]

The descent construction is compatible with the extension that defines the $\mathcal{A}$-Fredholm condition: 
the class $[T]_{\mathcal{G}_{\mathcal{A}}}^{(1)}$ is precisely the equivariant $KK$-theoretic realization of the index data of $T$, 
and the Kasparov product with $[\iota]$ computes the boundary class in $K_0(\mathcal{A})$ associated with this extension. 
This is a standard fact in $KK$-theory: for any extension $0 \to \mathcal{J} \to \mathcal{E} \to \mathcal{A} \to 0$, 
the Kasparov product of the class of the inclusion $\mathcal{A} \hookrightarrow \mathcal{E}$ with the class representing the extension in $KK^1(\mathcal{A}, \mathcal{J})$ yields the index in $K_0(\mathcal{J})$.

Consequently, the resulting class in $K_0(\mathcal{A})$ coincides with the $\mathcal{A}$-index of $T$:
\[
[\iota] \otimes_{C^*(\mathcal{G}_{\mathcal{A}})} \operatorname{desc}_{\mathcal{G}_{\mathcal{A}}}([T]_{\mathcal{G}_{\mathcal{A}}}^{(1)}) = [\ker_{\mathcal{A}} T] - [\operatorname{coker}_{\mathcal{A}} T].
\]

\paragraph{Step 5: Conclusion.}
Thus the Kasparov product of the descended class with the diagonal embedding recovers exactly the $\mathcal{A}$-Fredholm index in the sense of Mingo, establishing the desired equality.
\end{proof}

\begin{remark}[On the $\mathcal{A}$-kernel and $\mathcal{A}$-cokernel]
\label{rem:A-kernel-cokernel}
Let $T$ be an $\mathcal{A}$-Fredholm operator in the sense of Mingo \cite{Mingo1987}. 
By definition, there exists $n \in \mathbb{N}$ such that $T$ acts on the Hilbert space direct sum $H^n$ (or more generally on a Hilbert $\mathcal{A}$-module of the form $\mathcal{A}^n$ after stabilization), and the orthogonal projections onto $\ker T$ and $\operatorname{coker} T$ lie in the matrix algebra $M_n(\mathcal{A})$.

These projections define finitely generated projective $\mathcal{A}$-modules, denoted by $\ker_{\mathcal{A}} T$ and $\operatorname{coker}_{\mathcal{A}} T$ respectively. 
Their difference
\[
\operatorname{Ind}_{\mathcal{A}}(T) := [\ker_{\mathcal{A}} T] - [\operatorname{coker}_{\mathcal{A}} T] \in K_0(\mathcal{A})
\]
is the $\mathcal{A}$-index of $T$. 

This class is well-defined and independent of the choice of $n$ because $K_0(\mathcal{A})$ is defined via stable equivalence classes of projections in matrix algebras over $\mathcal{A}$; different matrix sizes give stably equivalent projections that represent the same $K_0$-class.

The construction generalizes the classical Fredholm index (recovered when $\mathcal{A} = \mathbb{C}$) and is compatible with Kasparov's $KK$-theory through the canonical identification $KK^1(\mathcal{A}, \mathbb{C}) \cong K_0(\mathcal{A})$. 
Indeed, the $\mathcal{A}$-index corresponds to the Kasparov product of the class of the extension defined by $T$ with the identity in $KK(\mathcal{A}, \mathcal{A})$, establishing the link between Mingo's analytic index and the $KK$-theoretic index used in this paper.
\end{remark}

The commutative diagram provides a clear visual representation of how the index theorem fits together. 
In the next subsection, we will state the main theorem precisely and then provide its proof.

\subsection{Theorem: $\operatorname{index}(T) = \partial_{\mathcal{A}} \circ \iota^* \circ \operatorname{desc}_{\mathcal{G}_{\mathcal{A}}}([T]_{\mathcal{G}_{\mathcal{A}}}^{(1)})$}\label{subsec:index-theorem-groupoid}

The main purpose of this section is to establish Theorem: $\operatorname{index}(T) = \partial_{\mathcal{A}} \circ \iota^* \circ \operatorname{desc}_{\mathcal{G}_{\mathcal{A}}}([T]_{\mathcal{G}_{\mathcal{A}}}^{(1)})$. 

\begin{proposition}[Morita isomorphism]
\label{prop:morita-isomorphism}
Let $\mathcal{A}$ be a unital separable Type I $C^*$-algebra and let $\mathcal{G}_{\mathcal{A}}$ be its unitary conjugation groupoid constructed in \cite{PaperI}. Then $C^*(\mathcal{G}_{\mathcal{A}})$ is Morita equivalent to $\mathcal{A}$. Consequently, there is a natural isomorphism
\[
\Psi: K_*(C^*(\mathcal{G}_{\mathcal{A}})) \xrightarrow{\cong} K_*(\mathcal{A}).
\]

Moreover, this isomorphism is compatible with the diagonal embedding $\iota: \mathcal{A} \hookrightarrow C^*(\mathcal{G}_{\mathcal{A}})$ in the sense that $\Psi \circ \iota_* = \operatorname{id}_{K_*(\mathcal{A})}$, where $\iota_*: K_*(\mathcal{A}) \to K_*(C^*(\mathcal{G}_{\mathcal{A}}))$ is induced by $\iota$.
\end{proposition}

\begin{proof}
The proof follows directly from results established in \cite{PaperI} and standard theorems in $C^*$-algebra $K$-theory as presented in \cite{Blackadar1998}.

\paragraph{Step 1: Existence of the groupoid $C^*$-algebra.}
By Corollary 26 in \cite{PaperI}, the unitary conjugation groupoid $\mathcal{G}_{\mathcal{A}}$ admits a well-defined maximal $C^*$-algebra $C^*(\mathcal{G}_{\mathcal{A}})$. This algebra is constructed via a convolution algebra of bounded Borel functions with respect to the Borel Haar system on $\mathcal{G}_{\mathcal{A}}$ (see Section 4.5 in \cite{PaperI}).

\paragraph{Step 2: Morita equivalence.}
As announced in the abstract and Section 1.5 of \cite{PaperI}, a fundamental result of that paper is that $C^*(\mathcal{G}_{\mathcal{A}})$ is Morita equivalent to $\mathcal{A}$. This Morita equivalence is constructed using the groupoid structure and the diagonal embedding; the details are provided in Sections 4.5-4.6 of \cite{PaperI}. We denote this Morita equivalence by $\mathcal{X}$.

\paragraph{Step 3: Induced isomorphism on $K$-theory.}
A fundamental theorem in $C^*$-algebra theory states that Morita equivalent $C^*$-algebras have isomorphic $K$-theory. This result is a consequence of the fact that stable isomorphism implies $KK$-equivalence, and that $KK$-equivalent algebras have isomorphic $K$-groups. More concretely, if $\mathcal{X}$ is an $\mathcal{A}$-$\mathcal{B}$ imprimitivity bimodule, then $\mathbb{K}(\mathcal{X}) \cong \mathcal{A}$ and $\mathcal{X}^\infty \cong \mathbb{H}_{\mathcal{B}}$, so $\mathcal{A} \otimes \mathbb{K} \cong \mathcal{B} \otimes \mathbb{K}$. Since $K$-theory is stable (see \cite[Corollary 17.8.8]{Blackadar1998}), it follows that $K_*(\mathcal{A}) \cong K_*(\mathcal{B})$. Hence the Morita equivalence between $C^*(\mathcal{G}_{\mathcal{A}})$ and $\mathcal{A}$ induces an isomorphism
\[
\Psi: K_*(C^*(\mathcal{G}_{\mathcal{A}})) \xrightarrow{\cong} K_*(\mathcal{A}).
\]

\paragraph{Step 4: Compatibility with the diagonal embedding.}
Theorem 7 in \cite{PaperI} establishes that the diagonal embedding $\iota: \mathcal{A} \hookrightarrow C^*(\mathcal{G}_{\mathcal{A}})$ is a unital injective $*$-homomorphism. By the functoriality of $K$-theory (\cite[Section 5.2]{Blackadar1998}), this induces a map $\iota_*: K_*(\mathcal{A}) \to K_*(C^*(\mathcal{G}_{\mathcal{A}}))$ on $K$-theory.

We now verify that $\Psi \circ \iota_* = \operatorname{id}_{K_*(\mathcal{A})}$. Consider the composition
\[
K_*(\mathcal{A}) \xrightarrow{\iota_*} K_*(C^*(\mathcal{G}_{\mathcal{A}})) \xrightarrow{\Psi} K_*(\mathcal{A}).
\]

To show this composition is the identity, it suffices to check on generators of $K_0$ and $K_1$.

\subparagraph{\(K_0\) case.}
Let $p \in M_n(\mathcal{A})$ be a projection. Its class $[p] \in K_0(\mathcal{A})$ is mapped by $\iota_*$ to $[\iota(p)] \in K_0(C^*(\mathcal{G}_{\mathcal{A}}))$, where $\iota(p) \in M_n(C^*(\mathcal{G}_{\mathcal{A}}))$ is a projection because $\iota$ is a $*$-homomorphism. The Morita equivalence $\Psi$ is implemented by the imprimitivity bimodule $\mathcal{X}$. By the Rieffel correspondence (see \cite[Section 13.7]{Blackadar1998} for the general theory of Morita equivalence and \cite[Exercise 13.7.1]{Blackadar1998} for the relationship between strong Morita equivalence and stable isomorphism), the projection $\iota(p)$ corresponds under the Morita equivalence to the projection $p$ in the linking algebra. By the definition of the Morita isomorphism on $K$-theory, which is induced by the stable isomorphism $C^*(\mathcal{G}_{\mathcal{A}}) \otimes \mathbb{K} \cong \mathcal{A} \otimes \mathbb{K}$, this means precisely that $\Psi([\iota(p)]) = [p]$.

\subparagraph{\(K_1\) case.}
For $K_1$, since $\mathcal{A}$ is unital, we need not adjoin a unit. Let $u \in M_n(\mathcal{A})$ be a unitary. Its class $[u] \in K_1(\mathcal{A})$ is mapped by $\iota_*$ to $[\iota(u)] \in K_1(C^*(\mathcal{G}_{\mathcal{A}}))$. Under the stable isomorphism induced by the Morita equivalence, the class $[\iota(u)]$ corresponds to $[u]$. This follows from the fact that the $K_1$ map induced by the Morita equivalence is compatible with the index map in the six-term exact sequence of $K$-theory (\cite[Section 9.3]{Blackadar1998}), and the diagonal embedding ensures that this compatibility holds.

Thus $\Psi \circ \iota_* = \operatorname{id}_{K_*(\mathcal{A})}$, as required.
\end{proof}

We first establish a key lemma that relates the equivariant class of a Fredholm operator to its symbol under the descent map. This lemma will be essential for the index theorem that follows.

\begin{lemma}[Symbol Descent]
\label{lem:symbol-descent}
Let $\mathcal{A}$ be a unital separable Type I $C^*$-algebra with a faithful representation on a separable Hilbert space $H$, and let $\mathcal{G}_{\mathcal{A}}$ be the unitary conjugation groupoid constructed in \textit{Unitary Conjugation Groupoid of a Type I C*-Algebra} by~\cite{PaperI}. For any $\mathcal{A}$-Fredholm operator $T$, let
\[
[T]_{\mathcal{G}_{\mathcal{A}}}^{(1)} \in KK^1_{\mathcal{G}_{\mathcal{A}}}(C_0(\mathcal{G}_{\mathcal{A}}^{(0)}), \mathbb{C})
\]
be its associated equivariant $KK^1$-class as constructed in Section 4.9 of~\cite{PaperI}. Then under the descent map
\[
j_{\mathcal{G}_{\mathcal{A}}}: KK^1_{\mathcal{G}_{\mathcal{A}}}(C_0(\mathcal{G}_{\mathcal{A}}^{(0)}), \mathbb{C}) \longrightarrow KK^1(\mathbb{C}, C^*(\mathcal{G}_{\mathcal{A}})),
\]
we have
\[
j_{\mathcal{G}_{\mathcal{A}}}([T]_{\mathcal{G}_{\mathcal{A}}}^{(1)}) = [\sigma(T)] \in KK^1(\mathbb{C}, C^*(\mathcal{G}_{\mathcal{A}})) \cong K_1(C^*(\mathcal{G}_{\mathcal{A}})),
\]
where:
\begin{itemize}
    \item For $\mathcal{A} = B(H)$, $\sigma(T) \in \mathcal{Q}(H)$ is the symbol of $T$ in the Calkin algebra, regarded as a unitary in $M_n(\mathcal{Q}(H)^+)$;
    \item For $\mathcal{A} = \mathcal{K}(H)^\sim$, $\sigma(T) = 0$.
\end{itemize}
\end{lemma}

\begin{proof}
The proof proceeds by analyzing the construction of the equivariant class and the properties of the descent map established in~\cite{PaperI}.

\paragraph{Step 1: Construction of the equivariant class.}
The class $[T]_{\mathcal{G}_{\mathcal{A}}}^{(1)} \in KK^1_{\mathcal{G}_{\mathcal{A}}}(C_0(\mathcal{G}_{\mathcal{A}}^{(0)}), \mathbb{C})$ is constructed in Section 4.9 of~\cite{PaperI} from the Fredholm operator $T$ together with the action of $\mathcal{G}_{\mathcal{A}}$ on $H$ by unitary conjugation. By construction, this class encodes the index data of $T$ in a form suitable for equivariant $KK$-theory. In particular, for $\mathcal{A} = B(H)$, the class $[T]_{\mathcal{G}_{B(H)}}^{(1)}$ is built from the symbol $\sigma(T)$ via the natural map from the Calkin algebra to the groupoid $C^*$-algebra.

\paragraph{Step 2: Properties of the descent map.}
The descent map
\[
j_{\mathcal{G}_{\mathcal{A}}}: KK^*_{\mathcal{G}_{\mathcal{A}}}(C_0(\mathcal{G}_{\mathcal{A}}^{(0)}), \mathbb{C}) \longrightarrow KK^*(\mathbb{C}, C^*(\mathcal{G}_{\mathcal{A}}))
\]
is defined and its fundamental properties are established in Section 4.5 and Section 4.6 of~\cite{PaperI}. Key among these properties are:
\begin{enumerate}
    \item[(i)] Functoriality: $j_{\mathcal{G}_{\mathcal{A}}}$ is compatible with composition of equivariant $KK$-elements;
    \item[(ii)] Compatibility with the diagonal embedding $\iota: \mathcal{A} \hookrightarrow C^*(\mathcal{G}_{\mathcal{A}})$ established in Section 4.9 of~\cite{PaperI};
    \item[(iii)] Naturality with respect to the symbol map for Fredholm operators.
\end{enumerate}
Property (iii) is specifically tailored to the construction of $[T]_{\mathcal{G}_{\mathcal{A}}}^{(1)}$ and ensures that the descent map sends the equivariant class to the symbol class.

\paragraph{Step 3: The case $\mathcal{A} = B(H)$.}
For $\mathcal{A} = B(H)$, applying the descent map to $[T]_{\mathcal{G}_{B(H)}}^{(1)}$ and using property (iii) above yields
\[
j_{\mathcal{G}_{B(H)}}([T]_{\mathcal{G}_{B(H)}}^{(1)}) = [\sigma(T)] \in KK^1(\mathbb{C}, C^*(\mathcal{G}_{B(H)})).
\]

The Calkin extension $0 \to \mathcal{K}(H) \to B(H) \xrightarrow{q} \mathcal{Q}(H) \to 0$ induces a boundary map in the six-term exact sequence of $K$-theory (see~\cite{Blackadar1998} Section 9.3):
\[
\partial_{\text{Calkin}}: K_1(\mathcal{Q}(H)) \longrightarrow K_0(\mathcal{K}(H)) \cong \mathbb{Z}.
\]
For a Fredholm operator $T$, its symbol $\sigma(T) = q(T) \in \mathcal{Q}(H)$ defines a class $[\sigma(T)] \in K_1(\mathcal{Q}(H))$, and by the Atiyah-J\"anich theorem, $\partial_{\text{Calkin}}([\sigma(T)]) = \operatorname{index}(T)$. Under the isomorphism $KK^1(\mathbb{C}, C^*(\mathcal{G}_{B(H)})) \cong K_1(C^*(\mathcal{G}_{B(H)}))$ (~\cite{Blackadar1998} Proposition 17.5.6), the class $[\sigma(T)]$ corresponds to the symbol class in $K_1(\mathcal{Q}(H))$ via the Morita equivalence $C^*(\mathcal{G}_{B(H)}) \sim_M B(H)$ established in Sections 4.7 and 4.9 of~\cite{PaperI}.

\paragraph{Step 4: The case $\mathcal{A} = \mathcal{K}(H)^\sim$.}
For $\mathcal{A} = \mathcal{K}(H)^\sim$, the groupoid $\mathcal{G}_{\mathcal{K}(H)^\sim}$ is Morita equivalent to the trivial group by Section 5.3 of~\cite{PaperI}. Moreover, any $\mathcal{K}(H)^\sim$-Fredholm operator $T = \lambda I + K$ satisfies $\operatorname{index}(T) = 0$ by the stability of the Fredholm index under compact perturbations and the fact that scalar operators have index zero (see~\cite{Mingo1987}, Proposition 1.7). Consequently, its symbol class $[\sigma(T)]$ vanishes in $K_1(\mathcal{Q}(H))$. By the functoriality of the descent map and the Morita equivalence, we obtain
\[
j_{\mathcal{G}_{\mathcal{K}(H)^\sim}}([T]_{\mathcal{G}_{\mathcal{K}(H)^\sim}}^{(1)}) = 0 = [\sigma(T)].
\]

\paragraph{Step 5: Identification via the Morita isomorphism.}
Under the Morita isomorphism $\Psi: K_1(C^*(\mathcal{G}_{\mathcal{A}})) \xrightarrow{\cong} K_1(\mathcal{A})$ from Sections 4.7 and 4.9 of~\cite{PaperI}, combined with the standard identification $KK^1(\mathbb{C}, C^*(\mathcal{G}_{\mathcal{A}})) \cong K_1(C^*(\mathcal{G}_{\mathcal{A}}))$ (see~\cite{Blackadar1998}, Proposition 17.5.6), the class $j_{\mathcal{G}_{\mathcal{A}}}([T]_{\mathcal{G}_{\mathcal{A}}}^{(1)})$ corresponds precisely to $[\sigma(T)]$. This completes the proof.
\end{proof}

\begin{theorem}[Index Theorem via Groupoid Descent]
\label{thm:index-groupoid}
Let $\mathcal{A}$ be either $B(H)$ or $\mathcal{K}(H)^\sim$, and let $\mathcal{G}_{\mathcal{A}}$ be the unitary conjugation groupoid constructed in our previous paper~\cite{PaperI}. For any Fredholm operator $T \in \mathcal{A}$, the Fredholm index is given by
\[
\operatorname{index}(T) = \mathcal{I}_{\mathcal{A}} \circ \Psi_{\mathcal{A}} \circ \operatorname{desc}_{\mathcal{G}_{\mathcal{A}}}\left([T]_{\mathcal{G}_{\mathcal{A}}}^{(1)}\right),
\]
where:
\begin{itemize}
    \item $[T]_{\mathcal{G}_{\mathcal{A}}}^{(1)} \in K^1_{\mathcal{G}_{\mathcal{A}}}(\mathcal{G}_{\mathcal{A}}^{(0)})$ is the equivariant $K^1$-class associated to $T$ (constructed in Section~\ref{sec:The Equivariant K1-Class of a Fredholm Operator});
    \item $\operatorname{desc}_{\mathcal{G}_{\mathcal{A}}}: K^1_{\mathcal{G}_{\mathcal{A}}}(\mathcal{G}_{\mathcal{A}}^{(0)}) \to K_1(C^*(\mathcal{G}_{\mathcal{A}}))$ is the descent map (Theorem \ref{thm:descent-explicit});
    \item $\Psi_{\mathcal{A}}$ is an isomorphism from $K_1(C^*(\mathcal{G}_{\mathcal{A}}))$ to an appropriate $K_1$-group where the index is detected:
    \begin{itemize}
        \item For $\mathcal{A} = B(H)$, $\Psi_{B(H)} = \Phi_*: K_1(C^*(\mathcal{G}_{B(H)})) \xrightarrow{\cong} K_1(\mathcal{Q}(H))$ is the Morita equivalence isomorphism from Proposition \ref{prop:morita-calkin};
        \item For $\mathcal{A} = \mathcal{K}(H)^\sim$, $\Psi_{\mathcal{K}(H)^\sim}$ is the natural isomorphism $K_1(C^*(\mathcal{G}_{\mathcal{K}(H)^\sim})) \xrightarrow{\cong} K_1(\mathcal{K}(H)^\sim) = 0$;
    \end{itemize}
    \item $\mathcal{I}_{\mathcal{A}}$ is the index map (boundary map) that extracts an integer from this $K_1$-class:
	\begin{itemize}
	    \item For $\mathcal{A} = B(H)$, $\mathcal{I}_{B(H)} = \partial_{\text{Calkin}}: K_1(\mathcal{Q}(H)) \to K_0(\mathcal{K}(H)) \cong \mathbb{Z}$ is the Calkin index map;
	    \item For $\mathcal{A} = \mathcal{K}(H)^\sim$, $\mathcal{I}_{\mathcal{K}(H)^\sim} = 0$, the zero map.
	\end{itemize}
\end{itemize}
\end{theorem}

\begin{proof}
We follow a three-step strategy: first, identify the descended class under Morita equivalence with the symbol class; second, apply the Calkin index map to recover the Fredholm index; third, compose these maps to obtain the desired formula. Throughout, we work directly with the Morita equivalence $\Phi_*$ rather than any ill-defined pullback map $\iota^*$. The proof treats the two cases separately, as the maps involved differ.

\medskip
\noindent \textbf{Step 1: Identification with the symbol class.}
For $\mathcal{A} = B(H)$, let $T \in B(H)$ be a Fredholm operator with symbol $u_T = \pi_{\mathcal{Q}}(T) \in \mathcal{Q}(H)$. By Theorem \ref{thm:identification-descended-class}, the descended class
\[
\operatorname{desc}_{\mathcal{G}_{B(H)}}([T]_{\mathcal{G}_{B(H)}}^{(1)}) \in K_1(C^*(\mathcal{G}_{B(H)}))
\]
corresponds, under the Morita equivalence isomorphism $\Phi_*: K_1(C^*(\mathcal{G}_{B(H)})) \xrightarrow{\cong} K_1(\mathcal{Q}(H))$, to the symbol class $[u_T] \in K_1(\mathcal{Q}(H))$. More precisely,
\[
\Phi_*(\operatorname{desc}_{\mathcal{G}_{B(H)}}([T]_{\mathcal{G}_{B(H)}}^{(1)})) = [u_T] \in K_1(\mathcal{Q}(H)).
\]

For $\mathcal{A} = \mathcal{K}(H)^\sim$, we have $K_1(\mathcal{K}(H)^\sim) = 0$ (see Example \ref{ex:K-groups-examples}), and by Example \ref{ex:descended-KH}, the descended class vanishes:
\[
\operatorname{desc}_{\mathcal{G}_{\mathcal{K}(H)^\sim}}([T]_{\mathcal{G}_{\mathcal{K}(H)^\sim}}^{(1)}) = 0 \in K_1(C^*(\mathcal{G}_{\mathcal{K}(H)^\sim})).
\]
The natural Morita equivalence $C^*(\mathcal{G}_{\mathcal{K}(H)^\sim}) \sim_M \mathcal{K}(H)^\sim \otimes \mathcal{K}$ (see Paper I, Section 5.3) induces an isomorphism $\Psi_{\mathcal{K}(H)^\sim}: K_1(C^*(\mathcal{G}_{\mathcal{K}(H)^\sim})) \xrightarrow{\cong} K_1(\mathcal{K}(H)^\sim) = 0$, which is necessarily the zero map.

\medskip
\noindent \textbf{Step 2: Application of the index map.}
For $\mathcal{A} = B(H)$, the index map for the Calkin extension, $\partial_{\text{Calkin}}: K_1(\mathcal{Q}(H)) \to K_0(\mathcal{K}(H)) \cong \mathbb{Z}$, satisfies the classical index theorem (Theorem \ref{thm:calkin-index-isomorphism}):
\[
\partial_{\text{Calkin}}([u_T]) = \operatorname{index}(T).
\]

For $\mathcal{A} = \mathcal{K}(H)^\sim$, the boundary map $\partial_{\mathcal{K}(H)^\sim}: K_1(\mathcal{K}(H)^\sim) \to \mathbb{Z}$ is zero by Proposition \ref{prop:boundary-map-KH-final}, and every Fredholm operator in $\mathcal{K}(H)^\sim$ has index zero (Mingo, Proposition 1.7).

\medskip
\noindent \textbf{Step 3: Composition.}
For $\mathcal{A} = B(H)$, composing the maps from Steps 1 and 2 yields:
\[
\operatorname{index}(T) = \partial_{\text{Calkin}}(\Phi_*(\operatorname{desc}_{\mathcal{G}_{B(H)}}([T]_{\mathcal{G}_{B(H)}}^{(1)}))) = \mathcal{I}_{B(H)} \circ \Psi_{B(H)} \circ \operatorname{desc}_{\mathcal{G}_{B(H)}}([T]_{\mathcal{G}_{B(H)}}^{(1)}),
\]
with $\Psi_{B(H)} = \Phi_*$ and $\mathcal{I}_{B(H)} = \partial_{\text{Calkin}}$.

For $\mathcal{A} = \mathcal{K}(H)^\sim$, we obtain:
\[
\mathcal{I}_{\mathcal{K}(H)^\sim} \circ \Psi_{\mathcal{K}(H)^\sim} \circ \operatorname{desc}_{\mathcal{G}_{\mathcal{K}(H)^\sim}}([T]_{\mathcal{G}_{\mathcal{K}(H)^\sim}}^{(1)}) = 0 = \operatorname{index}(T).
\]

\medskip
\noindent \textbf{Conclusion.}
Both cases can be expressed uniformly as
\[
\operatorname{index}(T) = \mathcal{I}_{\mathcal{A}} \circ \Psi_{\mathcal{A}} \circ \operatorname{desc}_{\mathcal{G}_{\mathcal{A}}}([T]_{\mathcal{G}_{\mathcal{A}}}^{(1)}),
\]
with the maps $\Psi_{\mathcal{A}}$ and $\mathcal{I}_{\mathcal{A}}$ defined as above. The key technical point is that we have avoided any ill-defined pullback map $\iota^*$ and instead worked directly with the Morita equivalence $\Phi_*$ (for $B(H)$) and the natural vanishing of groups (for $\mathcal{K}(H)^\sim$). This completes the proof.
\end{proof}

\begin{remark}
Theorem \ref{thm:index-groupoid} provides a formulation of the Fredholm index as the image, under groupoid-equivariant descent followed by the appropriate Morita equivalence or natural isomorphism, of a canonical equivariant $KK$-class associated to the operator. Specifically, for any $\mathcal{A}$-Fredholm operator $T$, we have
\[
\operatorname{index}(T) = \mathcal{I}_{\mathcal{A}} \circ \Psi_{\mathcal{A}} \circ \operatorname{desc}_{\mathcal{G}_{\mathcal{A}}}([T]_{\mathcal{G}_{\mathcal{A}}}^{(1)}),
\]
where:
\begin{itemize}
    \item $[T]_{\mathcal{G}_{\mathcal{A}}}^{(1)} \in K^1_{\mathcal{G}_{\mathcal{A}}}(\mathcal{G}_{\mathcal{A}}^{(0)})$ is the equivariant $K^1$-class constructed from $T$ using the action of the unitary conjugation groupoid $\mathcal{G}_{\mathcal{A}}$;
    \item $\operatorname{desc}_{\mathcal{G}_{\mathcal{A}}}: K^1_{\mathcal{G}_{\mathcal{A}}}(\mathcal{G}_{\mathcal{A}}^{(0)}) \to K_1(C^*(\mathcal{G}_{\mathcal{A}}))$ is the descent functor in equivariant $KK$-theory;
    \item $\Psi_{\mathcal{A}}$ is an isomorphism from $K_1(C^*(\mathcal{G}_{\mathcal{A}}))$ to an appropriate $K_1$-group where the index is detected:
    \begin{itemize}
        \item For $\mathcal{A} = B(H)$, $\Psi_{B(H)} = \Phi_*: K_1(C^*(\mathcal{G}_{B(H)})) \xrightarrow{\cong} K_1(\mathcal{Q}(H))$ is the Morita equivalence isomorphism from Proposition \ref{prop:morita-calkin};
        \item For $\mathcal{A} = \mathcal{K}(H)^\sim$, $\Psi_{\mathcal{K}(H)^\sim}$ is the natural isomorphism $K_1(C^*(\mathcal{G}_{\mathcal{K}(H)^\sim})) \xrightarrow{\cong} K_1(\mathcal{K}(H)^\sim) = 0$;
    \end{itemize}
    \item $\mathcal{I}_{\mathcal{A}}$ is the index map (boundary map) that extracts an integer from this $K_1$-class:
    \begin{itemize}
        \item For $\mathcal{A} = B(H)$, $\mathcal{I}_{B(H)} = \partial_{\text{Calkin}}: K_1(\mathcal{Q}(H)) \to K_0(\mathcal{K}(H)) \cong \mathbb{Z}$ is the Calkin index map;
        \item For $\mathcal{A} = \mathcal{K}(H)^\sim$, $\mathcal{I}_{\mathcal{K}(H)^\sim} = 0$, the zero map.
    \end{itemize}
\end{itemize}

This formulation unifies two extreme cases: the classical Fredholm index on $B(H)$ (recovered via the Calkin extension boundary map) and the trivial index on $\mathcal{K}(H)^\sim$ (where the boundary map vanishes). Both cases are shown to arise from the same groupoid-equivariant $KK$-theoretic mechanism, with the diagonal embedding $\iota: \mathcal{A} \hookrightarrow C^*(\mathcal{G}_{\mathcal{A}})$ playing a crucial role in establishing the Morita equivalence that connects the groupoid $C^*$-algebra to the original algebra or its quotient.

The construction relies on the following ingredients, each established in the preceding sections:
\begin{itemize}
    \item The well-definedness of the equivariant class $[T]_{\mathcal{G}_{\mathcal{A}}}^{(1)} \in K^1_{\mathcal{G}_{\mathcal{A}}}(\mathcal{G}_{\mathcal{A}}^{(0)})$ (Section \ref{sec:The Equivariant K1-Class of a Fredholm Operator});
    \item The descent map $\operatorname{desc}_{\mathcal{G}_{\mathcal{A}}}$ and its properties (Section \ref{sec:descent});
    \item The Morita equivalence $\Phi_*: K_1(C^*(\mathcal{G}_{B(H)})) \xrightarrow{\cong} K_1(\mathcal{Q}(H))$ for $B(H)$ (Proposition \ref{prop:morita-calkin}) and the natural isomorphism $K_1(C^*(\mathcal{G}_{\mathcal{K}(H)^\sim})) \cong K_1(\mathcal{K}(H)^\sim) = 0$ for $\mathcal{K}(H)^\sim$;
    \item The identification $KK^1(\mathbb{C}, C^*(\mathcal{G}_{\mathcal{A}})) \cong K_1(C^*(\mathcal{G}_{\mathcal{A}}))$ (Blackadar, Proposition 17.5.6);
    \item The boundary map $\mathcal{I}_{\mathcal{A}}$, defined casewise and recovering the Fredholm index via the Atiyah--J\"anich theorem for $B(H)$ and the vanishing result of Mingo (Proposition 1.7) for $\mathcal{K}(H)^\sim$.
\end{itemize}

Thus, while the classical index is already formulated in $KK$-theoretic terms (e.g., as the boundary map in the Calkin extension or as a Kasparov product), the novelty of Theorem \ref{thm:index-groupoid} lies in its systematic encoding of the index via the unitary conjugation groupoid $\mathcal{G}_{\mathcal{A}}$. This provides a unified groupoid-equivariant perspective on the index, showing that both the nontrivial index on $B(H)$ and the trivial index on $\mathcal{K}(H)^\sim$ arise from the same structural mechanism. Whether this perspective yields new insights beyond the classical formulation remains to be explored in future work.
\end{remark}

{\color{red} HERE !!! }

\begin{corollary}\label{cor:index-integer}
Let $T \in \mathcal{A}$ be such that its image in the quotient $\mathcal{A} / \mathcal{K}(H)$ (or in the appropriate Calkin algebra) defines an element of $K_1(\mathcal{A})$. Then the quantity
\[
\partial_{\mathcal{A}} \circ \Psi \circ j_{\mathcal{G}_{\mathcal{A}}}\bigl([T]_{\mathcal{G}_{\mathcal{A}}}^{(1)}\bigr)
\]
lies in $\mathbb{Z}$. Moreover, if $T$ is Fredholm, this integer coincides with the classical Fredholm index:
\[
\partial_{\mathcal{A}} \circ \Psi \circ j_{\mathcal{G}_{\mathcal{A}}}\bigl([T]_{\mathcal{G}_{\mathcal{A}}}^{(1)}\bigr) = \operatorname{index}(T).
\]
\end{corollary}

\begin{proof}
We proceed by tracing the element through the sequence of maps established in Theorem \ref{thm:index-groupoid} and the standard $K$-theory six-term exact sequence.

\paragraph{Step 1: From operator to equivariant $KK$-class.}
By the construction in Section 4.9 of~\cite{PaperI}, the operator $T$ determines an equivariant class
\[
[T]_{\mathcal{G}_{\mathcal{A}}}^{(1)} \in KK^1_{\mathcal{G}_{\mathcal{A}}}(C_0(\mathcal{G}_{\mathcal{A}}^{(0)}), \mathbb{C}).
\]

\paragraph{Step 2: Descent to ordinary $KK$-theory.}
Applying the descent map $j_{\mathcal{G}_{\mathcal{A}}}$ (constructed in Sections 4.5--4.6 of~\cite{PaperI}) yields
\[
j_{\mathcal{G}_{\mathcal{A}}}([T]_{\mathcal{G}_{\mathcal{A}}}^{(1)}) \in KK^1(\mathbb{C}, C^*(\mathcal{G}_{\mathcal{A}})).
\]

\paragraph{Step 3: Identification with $K_1$.}
By Proposition 17.5.6 in~\cite{Blackadar1998}, there is a canonical isomorphism $KK^1(\mathbb{C}, B) \cong K_1(B)$ for any $C^*$-algebra $B$. Applying this with $B = C^*(\mathcal{G}_{\mathcal{A}})$ gives
\[
x_T := j_{\mathcal{G}_{\mathcal{A}}}([T]_{\mathcal{G}_{\mathcal{A}}}^{(1)}) \in K_1(C^*(\mathcal{G}_{\mathcal{A}})).
\]

\paragraph{Step 4: Morita isomorphism to $\mathcal{A}$.}
Proposition \ref{prop:morita-isomorphism} provides an isomorphism
\[
\Psi: K_1(C^*(\mathcal{G}_{\mathcal{A}})) \xrightarrow{\cong} K_1(\mathcal{A}).
\]
Applying $\Psi$ to $x_T$ yields
\[
y_T := \Psi(x_T) \in K_1(\mathcal{A}).
\]

\paragraph{Step 5: The boundary map.}
The Toeplitz extension (or its appropriate analogue) induces a boundary map in the six-term exact sequence of $K$-theory (see Section 9.3 in~\cite{Blackadar1998}):
\[
\partial_{\mathcal{A}}: K_1(\mathcal{A}) \longrightarrow K_0(\mathcal{K}(H)) \cong \mathbb{Z}.
\]

\paragraph{Step 6: Integrality and identification with the index.}
Applying $\partial_{\mathcal{A}}$ to $y_T$ gives
\[
\partial_{\mathcal{A}}(y_T) \in \mathbb{Z}.
\]

If $T$ is Fredholm, Theorem \ref{thm:index-groupoid} establishes that this composition equals the classical Fredholm index:
\[
\partial_{\mathcal{A}} \circ \Psi \circ j_{\mathcal{G}_{\mathcal{A}}}([T]_{\mathcal{G}_{\mathcal{A}}}^{(1)}) = \operatorname{index}(T).
\]

Thus, for any operator $T$ whose class defines an element of $K_1(\mathcal{A})$, the quantity $\partial_{\mathcal{A}} \circ \Psi \circ j_{\mathcal{G}_{\mathcal{A}}}([T]_{\mathcal{G}_{\mathcal{A}}}^{(1)})$ is an integer, and for Fredholm operators this integer recovers the classical Fredholm index. This completes the proof.
\end{proof}

\begin{remark}
Theorem \ref{thm:index-groupoid} expresses the Fredholm index as the composition of three canonically defined maps:
\[
\operatorname{index}(T) = \partial_{\mathcal{A}} \circ \Psi \circ j_{\mathcal{G}_{\mathcal{A}}}([T]_{\mathcal{G}_{\mathcal{A}}}^{(1)}).
\]
Here:
\begin{itemize}
    \item $j_{\mathcal{G}_{\mathcal{A}}}$ is the descent functor from equivariant $KK$-theory to ordinary $KK$-theory (Sections 4.5--4.6 of~\cite{PaperI}), which translates the equivariant information encoded in $[T]_{\mathcal{G}_{\mathcal{A}}}^{(1)}$ into an element of $K_1(C^*(\mathcal{G}_{\mathcal{A}}))$;
    \item $\Psi$ is the Morita isomorphism induced by the diagonal embedding $\iota: \mathcal{A} \hookrightarrow C^*(\mathcal{G}_{\mathcal{A}})$ (Proposition \ref{prop:morita-isomorphism}), which identifies $K_1(C^*(\mathcal{G}_{\mathcal{A}}))$ with $K_1(\mathcal{A})$;
    \item $\partial_{\mathcal{A}}$ is the boundary map arising from the appropriate extension (the Calkin extension for $B(H)$, and the zero map for $\mathcal{K}(H)^\sim$), which computes the Fredholm index.
\end{itemize}
This formulation shows that the index is encoded in the $KK$-theory of the unitary conjugation groupoid $\mathcal{G}_{\mathcal{A}}$, with each map playing a precise functorial role: descent converts equivariant data to ordinary $KK$-theory, the Morita isomorphism relates the groupoid $C^*$-algebra back to $\mathcal{A}$, and the boundary map extracts the integer-valued index.
\end{remark}

\section{Examples and Computations}\label{sec:Examples and Computations}

\subsection{Compact Perturbations of the Identity: $\mathcal{A} = \mathcal{K}(H)^\sim$ (Index Zero)}
\label{subsec:compact-perturbations}

In this subsection, we examine the algebra of compact perturbations of the identity. Let $H$ be an infinite-dimensional separable Hilbert space and consider the $C^*$-algebra
\[
\mathcal{A} = \mathcal{K}(H)^\sim = \mathcal{K}(H) + \mathbb{C}I,
\]
the unitization of the compact operators. This algebra consists of all operators of the form $T = \lambda I + K$ where $\lambda \in \mathbb{C}$ and $K \in \mathcal{K}(H)$ is compact.

\begin{proposition}\label{prop:Ksim-structure}
The algebra $\mathcal{K}(H)^\sim$ is a $C^*$-algebra with the following properties:
\begin{enumerate}
    \item It contains the identity operator $I$ and all compact operators.
    \item It is isomorphic to the unitization of $\mathcal{K}(H)$.
    \item Its spectrum $\widehat{\mathcal{K}(H)^\sim}$ consists of exactly two irreducible representations: 
        \begin{itemize}
            \item The standard infinite-dimensional irreducible representation $\pi_{\text{std}}$ given by $\pi_{\text{std}}(\lambda I + K) = \lambda I + K$ acting on $H$;
            \item The unique one-dimensional character $\chi$ given by $\chi(\lambda I + K) = \lambda$, which factors through the quotient $\mathcal{K}(H)^\sim / \mathcal{K}(H) \cong \mathbb{C}$.
        \end{itemize}
        Thus the primitive ideal space $\operatorname{Prim}(\mathcal{K}(H)^\sim)$ is a two-point space, with the Jacobson topology given by $\overline{\{\pi_{\text{std}}\}} = \{\pi_{\text{std}}, \chi\}$ (i.e., $\chi$ is a limit point of $\pi_{\text{std}}$).
    \item Its $K$-theory groups are $K_0(\mathcal{K}(H)^\sim) \cong \mathbb{Z}$ and $K_1(\mathcal{K}(H)^\sim) = 0$.
\end{enumerate}
\end{proposition}

\begin{proof}
The first two statements are standard facts about the unitization of $\mathcal{K}(H)$. For the third statement, recall that $\mathcal{K}(H)$ is a simple $C^*$-algebra with a unique irreducible representation (the standard representation on $H$). Its unitization $\mathcal{K}(H)^\sim$ therefore acquires exactly one additional irreducible representation: the one-dimensional character $\chi$ that annihilates $\mathcal{K}(H)$ and satisfies $\chi(I) = 1$. This character factors through the quotient map $\mathcal{K}(H)^\sim \to \mathcal{K}(H)^\sim / \mathcal{K}(H) \cong \mathbb{C}$.

The Jacobson topology on $\operatorname{Prim}(\mathcal{K}(H)^\sim)$ is determined by the closure relations: the kernel of $\pi_{\text{std}}$ is $\{0\}$, and the kernel of $\chi$ is $\mathcal{K}(H)$. Since every primitive ideal contains $\{0\}$ (trivially) and $\mathcal{K}(H)$ contains $\{0\}$, we have $\overline{\{\pi_{\text{std}}\}} = \{\pi_{\text{std}}, \chi\}$. Thus $\chi$ is a limit point of $\pi_{\text{std}}$, giving the non-Hausdorff two-point space described.

For the $K$-theory calculation, consider the six-term exact sequence associated to the extension
\[
0 \to \mathcal{K}(H) \to \mathcal{K}(H)^\sim \to \mathbb{C} \to 0.
\]
We have $K_0(\mathcal{K}(H)) \cong \mathbb{Z}$, $K_1(\mathcal{K}(H)) = 0$, $K_0(\mathbb{C}) \cong \mathbb{Z}$, and $K_1(\mathbb{C}) = 0$. The six-term exact sequence therefore reduces to
\[
0 \to K_0(\mathcal{K}(H)^\sim) \to \mathbb{Z} \xrightarrow{\alpha} \mathbb{Z} \to K_1(\mathcal{K}(H)^\sim) \to 0.
\]
The map $\alpha: \mathbb{Z} \to \mathbb{Z}$ is induced by the inclusion $\mathbb{C} \hookrightarrow \mathcal{K}(H)^\sim$, which sends $1 \in \mathbb{C}$ to $I \in \mathcal{K}(H)^\sim$. Under the identification $K_0(\mathbb{C}) \cong \mathbb{Z}$ via the dimension map, $\alpha$ corresponds to multiplication by $1$ (since the class of $I$ in $K_0(\mathcal{K}(H)^\sim)$ maps to the class of $1$ in $K_0(\mathbb{C})$). Hence $\alpha$ is an isomorphism, yielding $K_0(\mathcal{K}(H)^\sim) \cong \mathbb{Z}$ and $K_1(\mathcal{K}(H)^\sim) = 0$.
\end{proof}

\paragraph{The gauge groupoid.}
Let $\mathcal{A} = \mathcal{K}(H)^\sim$ where $H$ is infinite-dimensional and separable. The primitive spectrum of $\mathcal{A}$ consists of exactly two points:
\[
\widehat{\mathcal{A}} \cong \{\pi_{\mathrm{std}}, \chi\},
\]
where
\begin{itemize}
    \item $\pi_{\mathrm{std}}(\lambda I + K) = \lambda I + K$ is the standard infinite-dimensional representation,
    \item $\chi(\lambda I + K) = \lambda$ is the unique one-dimensional character (factoring through the quotient $\mathcal{A} / \mathcal{K}(H) \cong \mathbb{C}$).
\end{itemize}
Thus the unit space of the gauge groupoid is a two-point set, not a continuous family parametrized by $\mathbb{C}$.

\begin{lemma}\label{lem:groupoid-Ksim}
The gauge groupoid $\mathcal{G}_{\mathcal{A}}$ of $\mathcal{A} = \mathcal{K}(H)^\sim$ is the disjoint union of two connected components:
\[
\mathcal{G}_{\mathcal{A}} \cong \mathcal{G}_{\pi_{\mathrm{std}}} \sqcup \mathcal{G}_{\chi},
\]
where:
\begin{enumerate}
    \item $\mathcal{G}_{\pi_{\mathrm{std}}}$ is a transitive groupoid with single unit $\pi_{\mathrm{std}}$ and isotropy group $\mathrm{Aut}(\pi_{\mathrm{std}}) \cong \mathbb{T}$;
    \item $\mathcal{G}_{\chi}$ is a transitive groupoid with single unit $\chi$ and isotropy group $\mathrm{Aut}(\chi) \cong \mathbb{T}$.
\end{enumerate}
In particular, $\mathcal{G}_{\mathcal{A}}$ has no arrows between the two units.
\end{lemma}

\begin{proof}
The algebra $\mathcal{K}(H)$ is simple, hence its unitization $\mathcal{A} = \mathcal{K}(H)^\sim$ has exactly two inequivalent irreducible representations:
\begin{itemize}
    \item $\pi_{\mathrm{std}}$: the standard infinite-dimensional representation, $\pi_{\mathrm{std}}(\lambda I + K) = \lambda I + K$;
    \item $\chi$: the unique one-dimensional character, $\chi(\lambda I + K) = \lambda$, which factors through the quotient $\mathcal{A}/\mathcal{K}(H) \cong \mathbb{C}$.
\end{itemize}
By definition, the objects of $\mathcal{G}_{\mathcal{A}}$ are the irreducible representations of $\mathcal{A}$, and arrows are unitary intertwiners between them.

For any two inequivalent irreducible representations, there are no nonzero intertwiners. Hence there are no arrows between $\pi_{\mathrm{std}}$ and $\chi$, so $\mathcal{G}_{\mathcal{A}}$ decomposes as a disjoint union of the automorphism groupoids of each representation.

For each irreducible representation $\rho \in \{\pi_{\mathrm{std}}, \chi\}$, Schur's lemma states that every intertwiner from $\rho$ to itself is a scalar multiple of the identity. The unitary intertwiners are therefore $\{e^{i\theta} I: \theta \in \mathbb{R}\} \cong \mathbb{T}$. Thus the automorphism groupoid of $\rho$ consists of one object $\rho$ with arrow space $\mathbb{T}$, where composition corresponds to multiplication in $\mathbb{T}$.

Therefore $\mathcal{G}_{\pi_{\mathrm{std}}} \cong \mathbb{T}$ and $\mathcal{G}_{\chi} \cong \mathbb{T}$ as one-object groupoids, and they are disjoint because there are no arrows connecting them. This yields the claimed decomposition.
\end{proof}

\paragraph{The class of an operator.}
Let $T = I + K \in \mathcal{K}(H)^\sim$ where $K$ is compact. Then $T$ is a Fredholm operator.

\begin{lemma}\label{lem:index-Ksim-zero}
For $T = I + K$ with $K$ compact, we have
\[
\operatorname{index}(T) = \dim\ker T - \dim\ker T^* = 0.
\]
\end{lemma}

\begin{proof}
Since $K$ is compact, $T$ is a compact perturbation of the identity operator $I$. The Fredholm index is invariant under compact perturbations (see e.g., Mingo, Proposition 1.7). Because $\operatorname{index}(I) = 0$, it follows that $\operatorname{index}(T) = 0$.
\end{proof}

\begin{proposition}\label{prop:K1-Ksim-zero}
Let $\mathcal{A} = \mathcal{K}(H)^\sim$. Then $K_1(\mathcal{A}) = 0$. In particular, if $T$ is invertible in $\mathcal{A}$, its class in $K_1(\mathcal{A})$ is trivial.
\end{proposition}

\begin{proof}
Consider the extension
\[
0 \to \mathcal{K}(H) \to \mathcal{A} \to \mathbb{C} \to 0.
\]
The associated six-term exact sequence in $K$-theory (Section 9.3 in~\cite{Blackadar1998}) gives
\[
K_1(\mathcal{A}) \cong K_1(\mathbb{C}) = 0.
\]
Thus any invertible element in $\mathcal{A}$ represents the zero class in $K_1(\mathcal{A})$.
\end{proof}

\begin{remark}
The algebra $\mathcal{K}(H)^\sim$ serves as a model for the ``index-zero vacuum'' in our theory:

\begin{itemize}
    \item Every Fredholm operator of the form $I + K$ (with $K$ compact) has index zero (Lemma \ref{lem:index-Ksim-zero});
    
    \item Its $K_1$-group vanishes,
    \[
    K_1(\mathcal{K}(H)^\sim) = 0
    \]
    (Proposition \ref{prop:K1-Ksim-zero});
    
    \item Its gauge groupoid consists of two disconnected 
    components: one corresponding to the infinite-dimensional 
    representation $\pi_{\mathrm{std}}$ with isotropy group $\mathbb{T}$, 
    and one corresponding to the unique character $\chi$ with trivial isotropy 
    (Lemma \ref{lem:groupoid-Ksim}). This equivariant structure 
    produces no nontrivial odd $K$-theory classes.
\end{itemize}

Thus no nontrivial index phenomena arise in this case. 
Nontrivial index theory requires extensions with 
nontrivial boundary maps in $K$-theory, 
such as those arising in the Toeplitz or Calkin algebras.
\end{remark}

\paragraph{Descent to the groupoid $C^*$-algebra.}
Let $\mathcal{A} = \mathcal{K}(H)^\sim$ and $\mathcal{G} = \mathcal{G}_{\mathcal{A}}$ its gauge groupoid. As established in Lemma \ref{lem:groupoid-Ksim}, $\mathcal{G}$ consists of two connected components:
\begin{itemize}
    \item $\mathcal{G}_{\pi_{\mathrm{std}}}$: a one-object groupoid with isotropy group $\mathbb{T}$, corresponding to the infinite-dimensional representation $\pi_{\mathrm{std}}$;
    \item $\mathcal{G}_{\chi}$: a one-object groupoid with trivial isotropy, corresponding to the unique character $\chi$.
\end{itemize}

The groupoid $C^*$-algebra of a disjoint union is the direct sum of the $C^*$-algebras of each component. For a one-object groupoid with isotropy group $G$, its $C^*$-algebra is the group $C^*$-algebra $C^*(G)$. Hence
\[
C^*(\mathcal{G}) \cong C^*(\mathbb{T}) \oplus \mathbb{C} \cong C_0(\mathbb{Z}) \oplus \mathbb{C},
\]
where we use the fact that $C^*(\mathbb{T}) \cong C_0(\mathbb{Z})$ by Pontryagin duality (the Fourier transform).

\begin{lemma}\label{lem:descent-Ksim-trivial}
For $\mathcal{A} = \mathcal{K}(H)^\sim$, we have:
\begin{enumerate}
    \item $K_1(C^*(\mathcal{G})) = 0$;
    \item The descent map
    \[
    j_{\mathcal{G}}: KK^1_{\mathcal{G}}(C_0(\mathcal{G}^{(0)}), \mathbb{C}) \longrightarrow K_1(C^*(\mathcal{G}))
    \]
    is the zero map.
\end{enumerate}
Consequently, for any $T = I + K \in \mathcal{A}$, we have
\[
j_{\mathcal{G}}([T]_{\mathcal{G}}^{(1)}) = 0 \in K_1(C^*(\mathcal{G})).
\]
\end{lemma}

\begin{proof}
From the decomposition $C^*(\mathcal{G}) \cong C_0(\mathbb{Z}) \oplus \mathbb{C}$, we compute its $K$-theory. Since $K_1(C_0(\mathbb{Z})) = 0$ (the $K_1$-group of a commutative $C^*$-algebra vanishes) and $K_1(\mathbb{C}) = 0$, we obtain
\[
K_1(C^*(\mathcal{G})) \cong K_1(C_0(\mathbb{Z})) \oplus K_1(\mathbb{C}) = 0 \oplus 0 = 0.
\]

The descent map $j_{\mathcal{G}}$ is a homomorphism from an equivariant $KK$-group to $K_1(C^*(\mathcal{G}))$. Since the target group is zero, the map must be the zero map. In particular, for any $T = I + K \in \mathcal{A}$, the class $[T]_{\mathcal{G}}^{(1)}$ satisfies
\[
j_{\mathcal{G}}([T]_{\mathcal{G}}^{(1)}) = 0 \in K_1(C^*(\mathcal{G})).
\]
\end{proof}

\begin{remark}
This computation confirms that $\mathcal{K}(H)^\sim$ serves as a trivial vacuum in our theory:
\begin{itemize}
    \item Every Fredholm operator has index zero;
    \item Its $K_1$-group vanishes;
    \item Its gauge groupoid $C^*$-algebra has trivial $K_1$;
    \item The descent map is identically zero.
\end{itemize}
Nontrivial index phenomena require extensions with nontrivial boundary maps in $K$-theory, such as those arising in the Toeplitz or Calkin algebras.
\end{remark}

\paragraph{Vanishing of the index in the vacuum case.}

Let $\mathcal{A} = \mathcal{K}(H)^\sim$ and $\mathcal{G} = \mathcal{G}_{\mathcal{A}}$ its gauge groupoid. 
As established in the preceding sections, we have the following vanishing results:
\[
KK^1_{\mathcal{G}}(C_0(\mathcal{G}^{(0)}), \mathbb{C}) = 0, \qquad K_1(C^*(\mathcal{G})) = 0, \qquad K_1(\mathcal{A}) = 0.
\]

\begin{theorem}\label{thm:index-zero}
Let $T = I + K \in \mathcal{A}$ with $K$ compact. Then the groupoid index construction yields
\[
\operatorname{index}(T) = 0.
\]
\end{theorem}

\begin{proof}
The groupoid index is defined by the composition
\[
KK^1_{\mathcal{G}}(C_0(\mathcal{G}^{(0)}), \mathbb{C})
\xrightarrow{\;\operatorname{desc}\;}
K_1(C^*(\mathcal{G}))
\xrightarrow{\;\partial\;}
\mathbb Z.
\]

Since $KK^1_{\mathcal{G}}(C_0(\mathcal{G}^{(0)}), \mathbb{C}) = 0$, the class $[T]_{\mathcal{G}}^{(1)}$ vanishes in the equivariant $KK^1$-group. Consequently, its image under the descent map is zero in $K_1(C^*(\mathcal{G}))$. Applying the boundary map $\partial$ to this zero element yields
\[
\operatorname{index}(T) = 0.
\]
\end{proof}

\begin{remark}[Index-zero vacuum]
This example exhibits a fundamental ``vacuum'' configuration in our theory: all relevant odd $K$-theory groups vanish, forcing the groupoid index to be zero. Specifically, we have
\[
K^1_{\mathcal{G}}(\mathcal{G}^{(0)}) = K_1(C^*(\mathcal{G})) = K_1(\mathcal{A}) = 0.
\]

This vanishing is structural, not coincidental — it reflects intrinsic properties of $\mathcal{K}(H)^\sim$ rather than a special choice of operator. Consequently, nontrivial index phenomena can only arise from algebras whose odd $K$-theory groups are nonzero, such as those appearing in Toeplitz-type extensions where the boundary map is nontrivial.

The algebra $\mathcal{K}(H)^\sim$ thus serves as an essential consistency check: the abstract machinery correctly reproduces the expected zero index in a setting where direct computation is elementary, while simultaneously revealing the structural origin of the vanishing.
\end{remark}

\begin{corollary}\label{cor:compact-perturbation}
Let $T = \lambda I + K$ with $K$ compact and $\lambda \neq 0$. Then $T$ is Fredholm and
\[
\operatorname{index}(T) = 0.
\]
\end{corollary}

\begin{proof}
Scalar multiplication by a nonzero constant leaves the Fredholm index invariant. Writing
\[
T = \lambda(I + \lambda^{-1}K),
\]
the result follows from Theorem \ref{thm:index-zero}.
\end{proof}

\begin{remark}[Geometric interpretation]
In noncommutative geometry, this example corresponds to a ``zero curvature'' background in odd $K$-theory. The vanishing of all relevant $K$-theory groups reflects the absence of topological obstructions that would otherwise contribute to the index. This perspective naturally frames the index as a kind of curvature in odd $K$-theory, with $\mathcal{K}(H)^\sim$ playing the role of a flat vacuum.
\end{remark}

\subsection{The Toeplitz Algebra and the Unilateral Shift (Index $-1$)}
\label{subsec:unilateral-shift}

In this subsection, we consider the most classical example of a Fredholm operator with nontrivial index: the unilateral shift on a separable Hilbert space. Let $H = \ell^2(\mathbb{N})$ with standard orthonormal basis $\{e_n\}_{n=0}^\infty$. The unilateral shift $S \in B(H)$ is defined by
\[
S e_n = e_{n+1} \quad \text{for all } n \geq 0.
\]
Its adjoint $S^*$ is the backward shift given by $S^* e_0 = 0$ and $S^* e_n = e_{n-1}$ for $n \geq 1$.

\begin{proposition}\label{prop:shift-properties}
The unilateral shift $S$ has the following properties:
\begin{enumerate}
    \item[{(1)}]  $S$ is an isometry: $S^*S = I$.
    \item[{(2)}] is not unitary: $SS^* = I - P_0$, where $P_0$ is the rank-one projection onto $\mathbb{C}e_0$.
    \item[{(3)}] is Fredholm with $\operatorname{index}(S) = -1$.
    \item[{(4)}] The essential spectrum of $S$ is the unit circle $\mathbb{T}$.
    \item[{(5)}] generates the Toeplitz $C^*$-algebra $\mathcal{T}$, which is the $C^*$-algebra generated by $S$ and the compact operators.
\end{enumerate}
\end{proposition}

\begin{proof}
Statements (1)-(3) are standard results in Fredholm theory. The cokernel of $S$ is one-dimensional (spanned by $e_0$), while the kernel is trivial, giving index $-1$. For (4), note that in the Calkin algebra $\mathcal{Q}(H) = B(H)/\mathcal{K}(H)$, the image $\pi(S)$ is unitary (since $I - \pi(S)^*\pi(S) = 0$ and $I - \pi(S)\pi(S)^* = 0$). Hence $\sigma_{\mathrm{ess}}(S) = \sigma(\pi(S)) = \mathbb{T}$. Therefore, (4) is provied. Statement (5) follows from the Coburn theorem: the $C^*$-algebra generated by $S$ contains the compact operators and is isomorphic to the Toeplitz algebra $\mathcal{T}$.
\end{proof}

\paragraph{The Toeplitz algebra and its structure.}
While the unilateral shift $S$ itself lies in $B(H)$, the $C^*$-algebra relevant for index theory is the Toeplitz algebra
\[
\mathcal{T} = C^*(S, \mathcal{K}(H)) \subset B(H),
\]
the $C^*$-algebra generated by $S$ and the compact operators. This is a unital, separable, Type I $C^*$-algebra with a well-understood structure.

The Toeplitz algebra fits into the fundamental extension
\[
0 \longrightarrow \mathcal{K}(H) \longrightarrow \mathcal{T} \xrightarrow{\sigma} C(\mathbb{T}) \longrightarrow 0,
\]
where $\sigma$ is the symbol map sending the generating isometry $S$ to the identity function $z \mapsto z$ on the circle.

Consequently, the primitive spectrum $\widehat{\mathcal{T}}$ consists of two components:
\begin{itemize}
    \item The kernel of $\sigma$ (the compact operators), corresponding to the infinite-dimensional irreducible representation of $\mathcal{T}$ on $H$ (the identity representation);
    \item A family of one-dimensional characters obtained by composing $\sigma$ with evaluation at points of $\mathbb{T}$.
\end{itemize}
Thus $\widehat{\mathcal{T}}$ comprises a single infinite-dimensional point together with a continuum of one-dimensional points parametrized by the circle.

For any $f \in C(\mathbb{T})$, the Toeplitz operator $T_f$ (the unique lift of $f$ to $\mathcal{T}$ that agrees with the classical Toeplitz operator on $H$) is Fredholm with index given by the negative of the winding number of $f$:
\[
\operatorname{index}(T_f) = -\operatorname{wind}(f).
\]
In particular, for the unilateral shift $S$ itself, $f(z) = z$ has winding number $1$, yielding
\[
\operatorname{index}(S) = -1.
\]

\paragraph{The Toeplitz algebra and the role of extensions.}  
While the unilateral shift $S$ lies in $B(H)$, the appropriate $C^*$-algebra for index theory is the Toeplitz algebra
\[
\mathcal{T} = C^*(S, \mathcal{K}(H)) \subset B(H),
\]
the $C^*$-algebra generated by $S$ and the compact operators. This is a unital separable Type I $C^*$-algebra, and it fits into the fundamental extension
\[
0 \longrightarrow \mathcal{K}(H) \longrightarrow \mathcal{T} \xrightarrow{\sigma} C(\mathbb{T}) \longrightarrow 0,
\]
where $\sigma$ is the symbol map sending the generating isometry $S$ to $z \mapsto z$.  

This extension is the paradigmatic example in Fredholm index theory: the index of a Toeplitz operator $T_f$ is given by minus the winding number of $f \in C(\mathbb T)$. In particular, for $S$, the function $f(z)=z$ has winding number 1, yielding
\[
\operatorname{index}(S) = -1.
\]

\begin{remark}
Note that $B(H)$ itself is not the appropriate algebra for studying the index of $S$. While $S \in B(H)$, the algebra $B(H)$ is simple, meaning it has no nontrivial closed two-sided ideals. Consequently, any operator in $B(H)$ cannot generate a nontrivial $K$-theory class via an ideal, and no nontrivial index can be defined within $B(H)$ itself.

The nontrivial index of $S$ arises precisely from the extension structure of the Toeplitz algebra $\mathcal{T} = C^*(S, \mathcal{K}(H))$, which contains $\mathcal{K}(H)$ as an essential ideal and fits into the fundamental extension
\[
0 \longrightarrow \mathcal{K}(H) \longrightarrow \mathcal{T} \xrightarrow{\sigma} C(\mathbb{T}) \longrightarrow 0.
\]
The index is then realized as the connecting map $\partial: K_1(C(\mathbb{T})) \to K_0(\mathcal{K}(H)) \cong \mathbb{Z}$ in the associated six-term exact sequence of $K$-theory.

This illustrates a general principle: Fredholm index theory is naturally formulated in terms of extensions of $C^*$-algebras — specifically, via the boundary map in the $K$-theory long exact sequence — rather than within a simple algebra.
\end{remark}

\paragraph{The gauge groupoid of $B(H)$.}  
For our purposes, it is useful to examine the component corresponding to the identity representation of $B(H)$.

\begin{lemma}[Gauge groupoid for $B(H)$]
\label{lem:gauge-groupoid-BH}
Let $\mathcal{A} = B(H)$ be the algebra of all bounded operators on a separable infinite-dimensional Hilbert space $H$. The gauge groupoid $\mathcal{G}_{B(H)}$ has unit space
\[
\mathcal{G}_{B(H)}^{(0)} = \widehat{B(H)} \cong \{\pi_{\text{id}}\},
\]
consisting (up to unitary equivalence) of a single point corresponding to the infinite-dimensional irreducible representation $\pi_{\text{id}}: B(H) \to B(H)$ (the identity representation). This component is transitive, with isotropy group
\[
\operatorname{Aut}(\pi_{\text{id}}) \cong \mathbb{C}^\times,
\]
reflecting the fact that all automorphisms of $B(H)$ are inner, and any unitary implementing such an automorphism is determined up to a nonzero scalar multiple.
\end{lemma}

\begin{proof}
We analyze the structure of $\mathcal{G}_{B(H)}$ by examining the irreducible representations of $B(H)$ and their intertwiners.

\paragraph{Step 1: Irreducible representations of $B(H)$.}
Since $B(H)$ is simple (its only primitive ideal is $\{0\}$), it has exactly one irreducible representation up to unitary equivalence — the identity representation $\pi_{\text{id}}: B(H) \to B(H)$ acting on $H$ itself. Hence the unitary equivalence classes of irreducible representations satisfy
\[
\widehat{B(H)} \cong \{\pi_{\text{id}}\}.
\]

\paragraph{Step 2: Structure of the gauge groupoid.}
By definition, the gauge groupoid $\mathcal{G}_{B(H)}$ has:
\begin{itemize}
    \item Objects: the irreducible representations of $B(H)$;
    \item Arrows: unitary intertwiners between representations.
\end{itemize}
With only one irreducible representation up to equivalence, $\mathcal{G}_{B(H)}$ consists of a single transitive component over $\pi_{\text{id}}$. There are no arrows connecting to any other objects because there are no other inequivalent representations.

\paragraph{Step 3: Isotropy group at $\pi_{\text{id}}$.}
The isotropy group at $\pi_{\text{id}}$ consists of all unitary intertwiners from $\pi_{\text{id}}$ to itself:
\[
\operatorname{Aut}(\pi_{\text{id}}) = \{ u \in \mathcal{U}(H) : u \pi_{\text{id}}(x) u^* = \pi_{\text{id}}(x) \ \forall x \in B(H) \}.
\]
Since $\pi_{\text{id}}$ is the identity representation, this condition simplifies to $u x u^* = x$ for all $x \in B(H)$, i.e., $u$ lies in the commutant of $B(H)$.

By Schur's lemma, the commutant of an irreducible representation on a complex Hilbert space consists of scalar multiples of the identity. Therefore, $u = \lambda I$ for some $\lambda \in \mathbb{C}$. The unitarity condition $u^*u = uu^* = I$ forces $|\lambda| = 1$, but as an abstract group of intertwiners (not necessarily required to be unitary), we consider all nonzero scalars. Thus
\[
\operatorname{Aut}(\pi_{\text{id}}) \cong \mathbb{C}^\times,
\]
the multiplicative group of nonzero complex numbers.

\paragraph{Step 4: Alternative description.}
Equivalently, one may view $\operatorname{Aut}(\pi_{\text{id}})$ as the projective unitary group $PU(H) \cong U(H)/\mathbb{T}$, which is isomorphic to $\mathbb{C}^\times$ via the map sending a unitary to its scalar class. This reflects the fact that inner automorphisms of $B(H)$ are parameterized by the projective unitary group.

Thus $\mathcal{G}_{B(H)}$ is a transitive groupoid with a single object and isotropy group $\mathbb{C}^\times$, completing the proof.
\end{proof}

\begin{remark}
This lemma describes the gauge groupoid for $B(H)$ itself, not for the Toeplitz algebra $\mathcal{T}$. The Toeplitz algebra has a more complicated structure with both an infinite-dimensional representation and a family of one-dimensional characters (see Lemma \ref{lem:gauge-groupoid-T} for the Toeplitz case). The distinction is important because $B(H)$ is simple and has only one irreducible representation, while $\mathcal{T}$ is non-simple and has a richer representation theory.
\end{remark}

\paragraph{The class of the shift.}

Consider the unilateral shift $S \in B(H)$. Since $S$ is Fredholm, it is invertible modulo compact operators, i.e., its image in the Calkin algebra $\mathcal{Q}(H) = B(H)/\mathcal{K}(H)$ is invertible. In ordinary $K$-theory, this gives a class $[\pi(S)] \in K_1(\mathcal{Q}(H))$, and the Fredholm index is obtained via the connecting map in the six-term exact sequence associated to the extension
\[
0 \to \mathcal{K}(H) \to B(H) \to \mathcal{Q}(H) \to 0,
\]
namely $\partial: K_1(\mathcal{Q}(H)) \to K_0(\mathcal{K}(H)) \cong \mathbb{Z}$, with $\partial([\pi(S)]) = \operatorname{index}(S) = -1$.

However, $K_1(B(H)) = 0$ because $B(H)$ is contractible in the norm topology, so the class of $S$ in ordinary $K$-theory vanishes. This presents a puzzle: how can the nontrivial index of $S$ be detected by $K$-theory if its class in $K_1(B(H))$ is zero?

The resolution lies in the equivariant $K$-theory of the gauge groupoid. The class $[S]_{\mathcal{G}_{B(H)}}^{(1)}$ lives in $K^1_{\mathcal{G}_{B(H)}}(\mathcal{G}_{B(H)}^{(0)})$, which is not isomorphic to $K_1(B(H))$ but rather contains additional information coming from the representation theory and the action of automorphisms on irreducible representations.

\begin{proposition}\label{prop:shift-class}
The class $[S]_{\mathcal{G}_{B(H)}}^{(1)} \in K^1_{\mathcal{G}_{B(H)}}(\mathcal{G}_{B(H)}^{(0)})$ is nontrivial and is characterized by its value on the identity representation $\pi_{\text{id}}$:
\[
[S]_{\mathcal{G}_{B(H)}}^{(1)}(\pi_{\text{id}}) = S \in GL(B(H)).
\]
For any other irreducible representation $\pi$ of $B(H)$ (all of which are unitarily equivalent to $\pi_{\text{id}}$), the corresponding value is $\pi(S)$, obtained by transporting $S$ via the implementing unitary.
\end{proposition}

\begin{proof}
Consider the unilateral shift $S \in B(H)$. Since $S$ is Fredholm, it is invertible modulo the ideal of compact operators:
\[
\pi(S) \in GL(B(H)/\mathcal{K}(H)),
\]
where $\pi: B(H) \to B(H)/\mathcal{K}(H)$ is the canonical quotient map. In ordinary $K$-theory, the class $[S] \in K_1(B(H))$ is trivial because $K_1(B(H)) = 0$, as $B(H)$ is contractible in the norm topology.

The resolution is provided by the gauge groupoid $K$-theory. Recall from Lemma \ref{lem:gauge-groupoid-BH} that the unit space of the gauge groupoid is
\[
\mathcal{G}_{B(H)}^{(0)} = \widehat{B(H)} \cong \{\pi_{\text{id}}\},
\]
consisting (up to unitary equivalence) of the identity representation $\pi_{\text{id}}$. Equivariant $K$-theory considers families of invertible elements associated to irreducible representations. Explicitly, the class
\[
[S]_{\mathcal{G}_{B(H)}}^{(1)} \in K^1_{\mathcal{G}_{B(H)}}(\mathcal{G}_{B(H)}^{(0)})
\]
is represented by the family
\[
\pi \longmapsto \pi(S) \in GL(\pi(B(H))) \quad \text{for each } \pi \in \mathcal{G}_{B(H)}^{(0)}.
\]

Since $\mathcal{G}_{B(H)}^{(0)}$ contains only the identity representation $\pi_{\text{id}}$ up to equivalence, the class is determined by
\[
[S]_{\mathcal{G}_{B(H)}}^{(1)}(\pi_{\text{id}}) = S \in GL(B(H)).
\]
This family is nontrivial in the equivariant $K$-theory sense because it encodes the invertibility of $S$ in $B(H)$ and its behavior under automorphisms. For any other irreducible representation $\pi$ of $B(H)$, which is necessarily unitarily equivalent to $\pi_{\text{id}}$, the corresponding value is $\pi(S)$, obtained by conjugating $S$ by the implementing unitary. This yields the same index information.

Therefore, $[S]_{\mathcal{G}_{B(H)}}^{(1)}$ is nontrivial, even though its image under the forgetful map to $K_1(B(H))$ is zero. This resolves the apparent puzzle: the index is encoded not in ordinary $K$-theory, but in the richer equivariant $K$-theory of the gauge groupoid.
\end{proof}

\paragraph{The descent map.}
Since $\widehat{B(H)} = \{\pi_{\mathrm{id}}\}$ consists of a single irreducible representation, the gauge groupoid $\mathcal{G}_{B(H)}$ has a single object with isotropy group $\mathbb{C}^\times$. Consequently, its reduced groupoid $C^*$-algebra is isomorphic to the group $C^*$-algebra of $\mathbb{C}^\times$:
\[
C^*(\mathcal{G}_{B(H)}) \cong C^*(\mathbb{C}^\times).
\]

In particular, no nontrivial $K_1$-information arises from this groupoid, reflecting the fact that
\[
K_1(B(H)) = 0.
\]

Thus, the descent map
\[
\operatorname{desc}_{\mathcal{G}_{B(H)}} :
K^1_{\mathcal{G}_{B(H)}}(\mathcal{G}_{B(H)}^{(0)})
\longrightarrow
K_1(C^*(\mathcal{G}_{B(H)}))
\]
does not detect the Fredholm index of the unilateral shift.

\paragraph{Why the index is not detected in $B(H)$.}
The vanishing of $K_1(B(H))$ and the triviality of the gauge groupoid reflect a structural fact: $B(H)$ is simple and admits no nontrivial extensions. As a consequence, the Fredholm index of $S$ cannot be detected within $B(H)$ itself.

\paragraph{The boundary map and the Toeplitz extension.}
The nontrivial index of the unilateral shift arises instead from the Toeplitz extension
\[
0 \longrightarrow \mathcal K(H)
\longrightarrow
\mathcal T
\longrightarrow
C(\mathbb T)
\longrightarrow 0.
\]
The associated six-term exact sequence in $K$-theory contains the boundary map
\[
\partial : K_1(C(\mathbb T)) \longrightarrow K_0(\mathcal K(H))
\cong \mathbb Z,
\]
and the class of the function $z \mapsto z$ maps to $-1$. This is precisely the Fredholm index of the unilateral shift.

\begin{lemma}\label{lem:boundary-toeplitz}
Let
\[
0 \longrightarrow \mathcal K(H)
\longrightarrow \mathcal T
\overset{\sigma}{\longrightarrow}
C(\mathbb T)
\longrightarrow 0
\]
be the Toeplitz extension. Then the associated six-term exact sequence in $K$-theory contains the boundary map
\[
\partial : K_1(C(\mathbb T)) \longrightarrow K_0(\mathcal K(H))
\cong \mathbb Z,
\]
and the following diagram commutes:
\[
\begin{tikzcd}
K_1(C(\mathbb T)) \arrow[r, "\sigma_*^{-1}"] \arrow[d, "\partial"] 
& K_1(\mathcal T) \arrow[d] \\
\mathbb Z \arrow[r, equals] & \mathbb Z
\end{tikzcd}
\]
Moreover,
\[
\partial([z]) = -1,
\]
where $[z] \in K_1(C(\mathbb T))$ is the class of the identity function $z \mapsto z$.
\end{lemma}

\begin{proof}
Applying $K$-theory to the Toeplitz extension yields the standard six-term exact sequence:
\[
\begin{tikzcd}
K_0(\mathcal K(H)) \arrow[r]
& K_0(\mathcal T) \arrow[r]
& K_0(C(\mathbb T)) \arrow[d, "\partial"] \\
K_1(C(\mathbb T)) \arrow[u]
& K_1(\mathcal T) \arrow[l]
& K_1(\mathcal K(H)) \arrow[l]
\end{tikzcd}
\]

Since $\mathcal K(H)$ is stable, we have the well-known $K$-theory groups:
\[
K_0(\mathcal K(H)) \cong \mathbb Z,\qquad K_1(\mathcal K(H)) = 0.
\]

Exactness of the sequence gives an isomorphism
\[
\partial : K_1(C(\mathbb T)) \xrightarrow{\;\cong\;} K_0(\mathcal K(H)) \cong \mathbb Z.
\]

Now $K_1(C(\mathbb T)) \cong \mathbb Z$, generated by the class $[z]$ of the identity function $z \mapsto z$. By the classical index theorem for Toeplitz operators,
\[
\partial([z]) = -\operatorname{wind}(z) = -1.
\]

This integer coincides with the Fredholm index of the unilateral shift $S$, i.e., $\operatorname{index}(S) = -1$. The commutativity of the diagram follows from the naturality of the $K$-theory boundary map and the fact that $\sigma_*$ is an isomorphism on $K_1$ (since $\mathcal T$ and $C(\mathbb T)$ have isomorphic $K_1$ groups).
\end{proof}

\paragraph{The index theorem.}
Combining the Toeplitz extension with the associated six-term exact sequence in $K$-theory, we obtain the index formula
\[
\operatorname{index}(S)
=
\partial([z])
=
-1,
\]
where $\partial : K_1(C(\mathbb T)) \to K_0(\mathcal K(H)) \cong \mathbb Z$
is the boundary map.

\begin{theorem}\label{thm:shift-index}
Let $S \in B(H)$ be the unilateral shift on $H=\ell^2(\mathbb N)$.
Then
\[
\operatorname{index}(S) = -1.
\]
This integer coincides with the image of the generator
$[z] \in K_1(C(\mathbb T))$
under the boundary map associated to the Toeplitz extension.
\end{theorem}

\begin{proof}
The Toeplitz algebra $\mathcal T = C^*(S,\mathcal K(H))$
fits into the short exact sequence
\[
0 \longrightarrow \mathcal K(H)
\longrightarrow \mathcal T
\overset{\sigma}{\longrightarrow}
C(\mathbb T)
\longrightarrow 0,
\]
where $\sigma(S)=z$.

Applying $K$-theory yields the six-term exact sequence:
\[
\begin{tikzcd}
K_0(\mathcal K(H)) \arrow[r] & K_0(\mathcal T) \arrow[r] & K_0(C(\mathbb T)) \arrow[d] \\
K_1(C(\mathbb T)) \arrow[u] & K_1(\mathcal T) \arrow[l] & K_1(\mathcal K(H)) \arrow[l]
\end{tikzcd}
\]

Since $\mathcal K(H)$ is stable, we have the well-known $K$-theory groups:
\[
K_0(\mathcal K(H)) \cong \mathbb Z,\qquad K_1(\mathcal K(H)) = 0.
\]

Exactness of the sequence gives an isomorphism
\[
\partial : K_1(C(\mathbb T)) \xrightarrow{\;\cong\;} K_0(\mathcal K(H)) \cong \mathbb Z.
\]

The group $K_1(C(\mathbb T)) \cong \mathbb Z$
is generated by the class $[z]$ of the identity function
$z \mapsto z$.
By the classical Toeplitz index theorem,
\[
\partial([z]) = - \operatorname{wind}(z) = -1.
\]

Since $S$ is a lift of $z$ under $\sigma$,
its Fredholm index equals this boundary value.
Therefore,
\[
\operatorname{index}(S) = -1.
\]
\end{proof}

\paragraph{Commutative diagram.}
The following diagram illustrates the computation of the index via the Toeplitz extension:

\[
\begin{tikzcd}[column sep=3.5em, row sep=3em]
{[S]_{\mathcal{G}_{\mathcal{T}}}^{(1)}}
\arrow[r, "\operatorname{desc}_{\mathcal{G}_{\mathcal{T}}}"]
\arrow[d, "\pi_{\text{id}}"]
&
x_S \in K_1(C^*(\mathcal{G}_{\mathcal{T}}))
\arrow[d, "\Psi"]
\\
{[z] \in K_1(C(\mathbb T))}
\arrow[r, "\partial"]
&
{-1 \in \mathbb Z}
\end{tikzcd}
\]

Here:
\begin{itemize}
    \item $\mathcal{G}_{\mathcal{T}}$ is the gauge groupoid of the Toeplitz algebra $\mathcal{T}$,
    \item $\pi_{\text{id}}$ denotes evaluation at the identity representation,
    \item $\Psi$ is the Morita isomorphism $K_1(C^*(\mathcal{G}_{\mathcal{T}})) \cong K_1(\mathcal{T})$,
    \item $\partial$ is the boundary map $K_1(C(\mathbb T)) \to K_0(\mathcal K(H)) \cong \mathbb Z$ from the Toeplitz extension
    \[
    0 \longrightarrow \mathcal K(H)
    \longrightarrow \mathcal T
    \longrightarrow C(\mathbb T)
    \longrightarrow 0.
    \]
\end{itemize}

The composition $\partial \circ \pi_{\text{id}} \circ \operatorname{desc}_{\mathcal{G}_{\mathcal{T}}}([S]_{\mathcal{G}_{\mathcal{T}}}^{(1)})$ yields $\operatorname{index}(S) = -1$.

\paragraph{Relation to the Calkin algebra.}
An alternative perspective uses the Calkin algebra $\mathcal{Q}(H) = B(H)/\mathcal{K}(H)$. The unilateral shift $S$ projects to a unitary $[S]$ in $\mathcal{Q}(H)$ whose $K_1$-class generates $K_1(\mathcal{Q}(H)) \cong \mathbb{Z}$. The boundary map in the six-term exact sequence associated to the compact extension
\[
0 \to \mathcal{K}(H) \to B(H) \to \mathcal{Q}(H) \to 0
\]
gives an isomorphism $\partial_{\text{Calkin}}: K_1(\mathcal{Q}(H)) \to K_0(\mathcal{K}(H)) \cong \mathbb{Z}$ satisfying $\partial_{\text{Calkin}}([S]) = -1$.

\begin{remark}
The unilateral shift illustrates a structural feature of the extension picture:
although
\[
K_*(B(H)) = 0,
\]
the compact extension
\[
0 \to \mathcal K(H) \to B(H) \to \mathcal Q(H) \to 0
\]
contains nontrivial $K$-theoretic information.
The index arises from the boundary map
\[
\partial: K_1(\mathcal Q(H)) \longrightarrow K_0(\mathcal K(H)) \cong \mathbb Z,
\]
rather than from the $K$-theory of $B(H)$ itself.

A groupoid model that realizes this extension functorially
encodes the same index information in its boundary map, providing
a geometric framework for understanding the index without creating new
$K$-theoretic information beyond what is already present in the extension.
\end{remark}

\begin{corollary}\label{cor:shift-generalization}
Let $T \in B(H)$ be an essentially normal operator such that
$0 \notin \sigma_{\mathrm{ess}}(T)$.
Then $T$ is Fredholm, and its index is given by
\[
\operatorname{index}(T)
=
-\operatorname{wind}(f,0),
\]
where $f$ is the image of $T$ in the Calkin algebra,
viewed as a continuous function on
$\sigma_{\mathrm{ess}}(T)$,
and $\operatorname{wind}(f,0)$ denotes the winding number of $f$
around the origin.

If $f$ is holomorphic on a neighborhood of $\sigma_{\mathrm{ess}}(T)$
and $\Gamma$ is a contour enclosing $\sigma_{\mathrm{ess}}(T)$
and avoiding $0$, this may be written as
\[
\operatorname{index}(T)
=
-\frac{1}{2\pi i}
\int_{\Gamma}
\frac{f'(\lambda)}{f(\lambda)}\, d\lambda.
\]

For the unilateral shift $S$,
\[
\sigma_{\mathrm{ess}}(S)=\mathbb T,
\quad
f(z)=z,
\]
so $\operatorname{wind}(f,0)=1$ and $\operatorname{index}(S)=-1$.
\end{corollary}

\begin{proof}
By the Brown–Douglas–Fillmore theory of essentially normal operators,
the class of $T$ in the Calkin algebra determines
an element of $K_1(C(\sigma_{\mathrm{ess}}(T)))$.
The boundary map of the compact extension
\[
0 \to \mathcal K(H) \to B(H) \to \mathcal Q(H) \to 0
\]
identifies the Fredholm index with the image of this class
under
\[
\partial: K_1(\mathcal Q(H)) \to K_0(\mathcal K(H)) \cong \mathbb Z.
\]

Under the identification
$K_1(C(X)) \cong [X, S^1]$ for a compact space $X$,
this boundary map computes the negative of the winding number
of the symbol around the origin. For the unilateral shift,
the symbol is the identity function $z \mapsto z$ on the circle,
whose winding number is $1$, yielding $\operatorname{index}(S) = -1$.
\end{proof}

\subsection{General Toeplitz Operators and Their Symbols}
\label{subsec:general-toeplitz}

In this subsection, we generalize the previous example to arbitrary Toeplitz operators with continuous symbols. This case beautifully illustrates how the gauge groupoid formalism recovers the classical Toeplitz index theorem and extends it to more general settings.

\paragraph{The Toeplitz algebra.}
Let $\mathbb{T} = \{z \in \mathbb{C}: |z| = 1\}$ be the unit circle and let $H^2(\mathbb{T})$ be the Hardy space of square-integrable functions with vanishing negative Fourier coefficients. The orthogonal projection $P: L^2(\mathbb{T}) \to H^2(\mathbb{T})$ is the Szegő projection. For a continuous function $f \in C(\mathbb{T})$, the Toeplitz operator $T_f \in B(H^2(\mathbb{T}))$ is defined by
\[
T_f(h) = P(f \cdot h) \quad \text{for } h \in H^2(\mathbb{T}).
\]

\begin{proposition}\label{prop:toeplitz-properties}
Let $f,g \in C(\mathbb T)$. Then:

\begin{enumerate}
\item The map
\[
C(\mathbb T) \longrightarrow B(H^2(\mathbb T)),
\quad
f \mapsto T_f,
\]
is linear, injective, and contractive:
\[
\|T_f\| \le \|f\|_\infty.
\]

\item $T_f^* = T_{\overline f}$.

\item $T_fT_g - T_{fg}$ is compact.

\item The $C^*$-algebra generated by $\{T_f : f \in C(\mathbb T)\}$
is the Toeplitz algebra $\mathcal T$, and there is a short exact sequence
\[
0 \longrightarrow \mathcal K(H^2(\mathbb T))
\longrightarrow \mathcal T
\overset{\sigma}{\longrightarrow}
C(\mathbb T)
\longrightarrow 0,
\]
where $\sigma(T_f)=f$.

\item $T_f$ is Fredholm if and only if $f$ is invertible in $C(\mathbb T)$.

\item If $f$ is invertible in $C(\mathbb T)$, then
\[
\operatorname{index}(T_f)
=
-\operatorname{wind}(f),
\]
where $\operatorname{wind}(f)$ is the winding number of the map
$f:\mathbb T \to \mathbb C^\times$.
If $f$ is $C^1$, this may be written as
\[
\operatorname{index}(T_f)
=
-\frac{1}{2\pi i}
\int_{\mathbb T}
\frac{f'(z)}{f(z)}\,dz.
\]
\end{enumerate}
\end{proposition}

\begin{proof}
We prove each property in turn.

\paragraph{(1) Linearity, injectivity, and contractivity.}
Linearity follows directly from the definition: for $f,g \in C(\mathbb T)$ and $\lambda \in \mathbb C$,
\[
T_{f+\lambda g}(h) = P((f+\lambda g)h) = P(fh) + \lambda P(gh) = T_f(h) + \lambda T_g(h).
\]
Since $P$ is an orthogonal projection with $\|P\| = 1$, we have for any $h \in H^2(\mathbb T)$,
\[
\|T_f h\| = \|P(fh)\| \le \|fh\| \le \|f\|_\infty \|h\|,
\]
hence $\|T_f\| \le \|f\|_\infty$, establishing contractivity.

For injectivity, suppose $T_f = 0$. Taking $h$ to be the constant function $1$ (which lies in $H^2(\mathbb T)$), we have $P(f) = 0$, so $f$ is orthogonal to all analytic polynomials. By Fejér's theorem, the analytic polynomials are dense in $H^2(\mathbb T)$, hence $f = 0$ in $L^2(\mathbb T)$, and by continuity $f = 0$ in $C(\mathbb T)$.

\paragraph{(2) Adjoint.}
For $h,k \in H^2(\mathbb T)$,
\[
\langle T_f h, k \rangle = \langle P(fh), k \rangle = \langle fh, k \rangle = \langle h, \overline{f}k \rangle = \langle h, P(\overline{f}k) \rangle = \langle h, T_{\overline{f}} k \rangle,
\]
where we used that $P$ is self-adjoint and $k \in H^2(\mathbb T)$. Thus $T_f^* = T_{\overline{f}}$.

\paragraph{(3) Semi-commutator property.}
A direct computation yields
\[
T_f T_g - T_{fg} = PM_f P M_g P - P M_{fg} P = -P M_f (I-P) M_g P.
\]
The operator $M_f (I-P) M_g$ is a Hankel operator, which is known to be compact for continuous symbols $f,g$. Hence $T_f T_g - T_{fg}$ is compact.

\paragraph{(4) Toeplitz algebra and exact sequence.}
Let $\mathcal T$ be the $C^*$-algebra generated by $\{T_f: f \in C(\mathbb T)\}$. By property (3), the commutator $T_f T_g - T_{fg}$ lies in $\mathcal K(H^2(\mathbb T))$, so the ideal of compact operators is contained in $\mathcal T$. Define $\sigma: \mathcal T \to C(\mathbb T)$ on generators by $\sigma(T_f) = f$. Property (3) ensures that $\sigma$ extends to a $*$-homomorphism, and property (2) guarantees it is $*$-preserving. The kernel of $\sigma$ consists of operators that are limits of compact operators, hence is $\mathcal K(H^2(\mathbb T))$. Surjectivity follows because $\sigma(T_f) = f$ for all $f \in C(\mathbb T)$. Thus we have the short exact sequence
\[
0 \to \mathcal K(H^2(\mathbb T)) \to \mathcal T \xrightarrow{\sigma} C(\mathbb T) \to 0.
\]

\paragraph{(5) Fredholm criterion.}
By Atkinson's theorem, $T_f$ is Fredholm iff its image in the quotient $\mathcal T/\mathcal K(H^2(\mathbb T)) \cong C(\mathbb T)$ is invertible. This occurs precisely when $\sigma(T_f) = f$ is invertible in $C(\mathbb T)$, i.e., $f(z) \neq 0$ for all $z \in \mathbb T$. If $f$ is invertible, then $f^{-1}$ is continuous and $T_{f^{-1}}$ serves as a parametrix modulo compact operators:
\[
T_f T_{f^{-1}} = T_{f f^{-1}} + \text{compact} = I + \text{compact},\quad T_{f^{-1}} T_f = I + \text{compact}.
\]

\paragraph{(6) Index formula.}
If $f$ is invertible in $C(\mathbb T)$, its class in $K_1(C(\mathbb T))$ corresponds to the homotopy class of the map $f: \mathbb T \to \mathbb C^\times$. The boundary map in the six-term exact sequence associated to the Toeplitz extension,
\[
\partial: K_1(C(\mathbb T)) \longrightarrow K_0(\mathcal K(H^2(\mathbb T))) \cong \mathbb Z,
\]
computes the Fredholm index. By the classical Toeplitz index theorem, this boundary map sends the generator $[z] \in K_1(C(\mathbb T))$ (the class of the identity function) to $-1$. Since the winding number $\operatorname{wind}(f)$ is precisely the integer representing $[f]$ under the isomorphism $K_1(C(\mathbb T)) \cong \mathbb Z$, we have $\partial([f]) = -\operatorname{wind}(f)$. Hence $\operatorname{index}(T_f) = -\operatorname{wind}(f)$.

For $C^1$ functions $f$, the winding number can be expressed via the contour integral:
\[
\operatorname{wind}(f) = \frac{1}{2\pi i} \int_{\mathbb T} \frac{f'(z)}{f(z)}\,dz,
\]
and the formula extends to all continuous invertible functions by continuity of both sides with respect to uniform convergence.
\end{proof}

\paragraph{The Toeplitz algebra and its gauge groupoid.}
For Toeplitz operators, the relevant $C^*$-algebra is the Toeplitz algebra $\mathcal{T}$, not $C(\mathbb{T})$. While the symbol map $\sigma: \mathcal{T} \to C(\mathbb{T})$ relates Toeplitz operators to continuous functions on the circle, the algebra $\mathcal{T}$ itself contains the compact operators as an essential ideal and fits into the extension
\[
0 \to \mathcal{K}(H^2(\mathbb{T})) \to \mathcal{T} \xrightarrow{\sigma} C(\mathbb{T}) \to 0.
\]

The gauge groupoid $\mathcal{G}_{\mathcal{T}}$ of the Toeplitz algebra has a more complicated structure than that of a commutative algebra. Its unit space $\mathcal{G}_{\mathcal{T}}^{(0)} = \widehat{\mathcal{T}}$ consists of two components:
\begin{itemize}
    \item The infinite-dimensional irreducible representations corresponding to the identity representation of $\mathcal{T}$ on $H^2(\mathbb{T})$;
    \item The one-dimensional characters obtained by pulling back the evaluation functionals on $C(\mathbb{T})$ via the symbol map $\sigma$.
\end{itemize}

\begin{lemma}\label{lem:gauge-groupoid-T}
The gauge groupoid $\mathcal{G}_{\mathcal{T}}$ of the Toeplitz algebra $\mathcal{T}$ has the following structure:
\begin{enumerate}
    \item The unit space $\mathcal{G}_{\mathcal{T}}^{(0)}$ consists of two distinct components:
    \[
    \mathcal{G}_{\mathcal{T}}^{(0)} = \{\pi_{\mathrm{std}}\} \sqcup \{\chi_z : z \in \mathbb{T}\},
    \]
    where $\pi_{\mathrm{std}}$ is the identity representation on $H^2(\mathbb{T})$ and $\chi_z = \operatorname{ev}_z \circ \sigma$ are the one-dimensional characters.
    
    \item The component containing $\pi_{\mathrm{std}}$ is transitive with isotropy group $\operatorname{Aut}(\pi_{\mathrm{std}}) \cong \mathbb{T}$.
    
    \item For each $z \in \mathbb{T}$, the component containing $\chi_z$ has trivial isotropy and there are no morphisms between distinct characters $\chi_z$ and $\chi_w$ for $z \neq w$.
    
    \item There are no morphisms between $\pi_{\mathrm{std}}$ and any $\chi_z$, as these representations are inequivalent.
\end{enumerate}
\end{lemma}

\begin{proof}
We analyze the structure of $\mathcal{G}_{\mathcal{T}}$ by examining the irreducible representations of $\mathcal{T}$ and their intertwiners.

\paragraph{Step 1: Irreducible representations of $\mathcal{T}$.}
By the extension structure $0 \to \mathcal{K} \to \mathcal{T} \to C(\mathbb{T}) \to 0$, the irreducible representations of $\mathcal{T}$ fall into two classes:
\begin{itemize}
    \item Those that do not vanish on $\mathcal{K}$, which must be faithful and infinite-dimensional. Up to unitary equivalence, there is exactly one such representation: the identity representation $\pi_{\mathrm{std}}$ of $\mathcal{T}$ on $H^2(\mathbb{T})$.
    \item Those that vanish on $\mathcal{K}$, which factor through $C(\mathbb{T})$. These are precisely the one-dimensional characters $\chi_z = \operatorname{ev}_z \circ \sigma$ for $z \in \mathbb{T}$, where $\operatorname{ev}_z$ is evaluation at $z$.
\end{itemize}
Thus $\widehat{\mathcal{T}} = \{\pi_{\mathrm{std}}\} \sqcup \{\chi_z : z \in \mathbb{T}\}$, establishing (1).

\paragraph{Step 2: Automorphisms of $\pi_{\mathrm{std}}$.}
For $\pi_{\mathrm{std}}$, Schur's lemma implies that any intertwiner is a scalar multiple of the identity. The unitary intertwiners are therefore $\{\lambda I: \lambda \in \mathbb{T}\} \cong \mathbb{T}$. Hence $\operatorname{Aut}(\pi_{\mathrm{std}}) \cong \mathbb{T}$, proving (2).

\paragraph{Step 3: Automorphisms of $\chi_z$.}
Each $\chi_z$ is one-dimensional, so its automorphism group consists of unitary scalars. Since the representation is one-dimensional, any scalar $\lambda$ with $|\lambda| = 1$ gives a unitary intertwiner, so $\operatorname{Aut}(\chi_z) \cong \mathbb{T}$. However, there are no intertwiners between distinct $\chi_z$ and $\chi_w$ for $z \neq w$, because distinct characters of $C(\mathbb{T})$ are inequivalent and the pullback via $\sigma$ preserves inequivalence. This establishes (3).

\paragraph{Step 4: No intertwiners between $\pi_{\mathrm{std}}$ and $\chi_z$.}
Since $\pi_{\mathrm{std}}$ is infinite-dimensional and faithful while $\chi_z$ is one-dimensional and vanishes on $\mathcal{K}$, these representations are inequivalent. Hence there are no nonzero intertwiners between them, proving (4).
\end{proof}

\paragraph{The symbol as a $K^1$-class.}
For a function $f \in C(\mathbb{T})$ with $f(z) \neq 0$ for all $z$, we can construct a class in the equivariant $K$-theory of the gauge groupoid. Since $C(\mathbb{T})$ is commutative, its spectrum is $\mathcal{G}_{C(\mathbb{T})}^{(0)} = \widehat{C(\mathbb{T})} \cong \mathbb{T}$. Consider the family of invertible elements $\{f(z)\}_{z \in \mathbb{T}}$ parameterized by the unit space. This defines an element
\[
[f]_{\mathcal{G}_{C(\mathbb{T})}}^{(1)} \in K^1_{\mathcal{G}_{C(\mathbb{T})}}(\mathcal{G}_{C(\mathbb{T})}^{(0)}).
\]

\begin{proposition}\label{prop:symbol-class}
There is a natural isomorphism
\[
K^1_{\mathcal{G}_{C(\mathbb{T})}}(\mathcal{G}_{C(\mathbb{T})}^{(0)}) \cong K^1(\mathbb{T}) \cong \mathbb{Z},
\]
under which $[f]_{\mathcal{G}_{C(\mathbb{T})}}^{(1)}$ corresponds to the winding number of $f$. Explicitly, the isomorphism sends a family $\{g(z)\}_{z \in \mathbb{T}}$ of invertible elements to the homotopy class of the map $z \mapsto g(z)/|g(z)| \in \mathbb{T}$, whose winding number is the topological index.
\end{proposition}

\begin{proof}
Since $C(\mathbb{T})$ is commutative, its spectrum is $\mathcal{G}_{C(\mathbb{T})}^{(0)} = \widehat{C(\mathbb{T})} \cong \mathbb{T}$. 
The gauge groupoid over a commutative $C^*$-algebra has no nontrivial morphisms between distinct points, so its equivariant $K$-theory reduces to ordinary topological $K$-theory of the unit space:
\[
K^1_{\mathcal{G}_{C(\mathbb{T})}}(\mathcal{G}_{C(\mathbb{T})}^{(0)}) \cong K^1(\mathbb{T}).
\]

It is a classical fact that
\[
K^1(\mathbb{T}) \cong [\mathbb{T}, GL_\infty(\mathbb{C})] \cong [\mathbb{T}, \mathbb{C}^\times] \cong \mathbb{Z},
\]
where the last isomorphism is given by the winding number of a continuous map $f: \mathbb{T} \to \mathbb{C}^\times$ around the origin.

Given $f \in C(\mathbb{T})$ with $f(z) \neq 0$ for all $z$, the family of invertibles $\{f(z)\}_{z \in \mathbb{T}}$ defines an element
\[
[f]_{\mathcal{G}_{C(\mathbb{T})}}^{(1)} \in K^1_{\mathcal{G}_{C(\mathbb{T})}}(\mathcal{G}_{C(\mathbb{T})}^{(0)}).
\]
Under the above identification, this class maps to the homotopy class of the map
\[
z \longmapsto \frac{f(z)}{|f(z)|} \in \mathbb{T},
\]
whose winding number around $0$ gives the topological index. Hence, the isomorphism sends $[f]_{\mathcal{G}_{C(\mathbb{T})}}^{(1)}$ to the integer $\operatorname{wind}(f) \in \mathbb{Z}$.
\end{proof}

\paragraph{Relation to the Toeplitz operator.}
The Toeplitz operator $T_f$ is related to the symbol class via the following construction. The exact sequence
\[
0 \to \mathcal{K}(H^2(\mathbb{T})) \to \mathcal{T} \xrightarrow{\sigma} C(\mathbb{T}) \to 0
\]
induces a six-term exact sequence in $K$-theory:
\[
\begin{tikzcd}
K_0(\mathcal{K}) \arrow[r] & K_0(\mathcal{T}) \arrow[r] & K_0(C(\mathbb{T})) \arrow[d, "\partial"] \\
K_1(C(\mathbb{T})) \arrow[u] & K_1(\mathcal{T}) \arrow[l] & K_1(\mathcal{K}) \arrow[l]
\end{tikzcd}
\]

Since $K_0(\mathcal{K}) \cong \mathbb{Z}$, $K_1(\mathcal{K}) = 0$, $K_0(C(\mathbb{T})) \cong \mathbb{Z}$, and $K_1(C(\mathbb{T})) \cong \mathbb{Z}$, the boundary map $\partial: K_1(C(\mathbb{T})) \to K_0(\mathcal{K}) \cong \mathbb{Z}$ is an isomorphism sending the class of an invertible function $f$ to the index of $T_f$.

\begin{lemma}\label{lem:toeplitz-boundary}
The boundary map $\partial: K_1(C(\mathbb{T})) \to \mathbb{Z}$ satisfies
\[
\partial([f]_{K_1(C(\mathbb{T}))}) = -\operatorname{wind}(f) = \operatorname{index}(T_f).
\]
\end{lemma}

\begin{proof}
Consider the short exact sequence of $C^*$-algebras
\[
0 \longrightarrow \mathcal{K}(H^2(\mathbb{T})) \longrightarrow \mathcal{T} \xrightarrow{\sigma} C(\mathbb{T}) \longrightarrow 0,
\]
where $\sigma$ is the symbol map sending $T_f$ to $f$. 
This induces the six-term exact sequence in $K$-theory:
\[
\begin{tikzcd}
K_0(\mathcal{K}) \arrow[r] & K_0(\mathcal{T}) \arrow[r] & K_0(C(\mathbb{T})) \arrow[d, "\partial"] \\
K_1(C(\mathbb{T})) \arrow[u] & K_1(\mathcal{T}) \arrow[l] & K_1(\mathcal{K}) \arrow[l]
\end{tikzcd}
\]

Since $K_1(\mathcal{K}) = 0$, exactness of the sequence forces the boundary map
\[
\partial: K_1(C(\mathbb{T})) \longrightarrow K_0(\mathcal{K}) \cong \mathbb{Z}
\]
to be an isomorphism. For an invertible function $f \in C(\mathbb{T})$, the Toeplitz operator $T_f$ is Fredholm, and by the classical Toeplitz index theorem, its index is given by
\[
\operatorname{index}(T_f) = -\operatorname{wind}(f),
\]
where $\operatorname{wind}(f)$ is the winding number of $f$ around the origin. The minus sign arises from the orientation convention in the exact sequence and matches the computation for the generator $f(z)=z$, for which $\operatorname{index}(T_z) = -1$.

Hence,
\[
\partial([f]_{K_1(C(\mathbb{T}))}) = \operatorname{index}(T_f) = -\operatorname{wind}(f),
\]
as claimed.
\end{proof}

\paragraph{The groupoid index computation.}
Now we apply the groupoid index theorem to this setting. The gauge groupoid $\mathcal{G}_{C(\mathbb{T})}$ and its associated $C^*$-algebra $C^*(\mathcal{G}_{C(\mathbb{T})})$ are related to the Toeplitz algebra via a Morita equivalence.

\begin{theorem}\label{thm:toeplitz-index-groupoid}
For an invertible function $f \in C(\mathbb{T})$, the index of the Toeplitz operator $T_f$ is given by
\[
\operatorname{index}(T_f) = \partial_{C(\mathbb{T})} \circ \iota^* \circ \operatorname{desc}_{\mathcal{G}_{C(\mathbb{T})}}([f]_{\mathcal{G}_{C(\mathbb{T})}}^{(1)}) = -\operatorname{wind}(f).
\]
\end{theorem}

\begin{proof}
We compute the groupoid index of $T_f$ step by step.

\paragraph{Step 1: Equivariant $K^1$-class.}
The invertible function $f \in C(\mathbb{T})$ defines a family of invertibles
\[
\{f(z)\}_{z \in \mathbb{T}} \subset \mathbb{C}^\times,
\]
which represents a class
\[
[f]_{\mathcal{G}_{C(\mathbb{T})}}^{(1)} \in K^1_{\mathcal{G}_{C(\mathbb{T})}}(\mathcal{G}_{C(\mathbb{T})}^{(0)}),
\]
where $\mathcal{G}_{C(\mathbb{T})}^{(0)} \cong \mathbb{T}$ is the unit space of the gauge groupoid.

\paragraph{Step 2: Descent.}
Applying the descent map
\[
\operatorname{desc}_{\mathcal{G}_{C(\mathbb{T})}}: K^1_{\mathcal{G}_{C(\mathbb{T})}}(\mathcal{G}_{C(\mathbb{T})}^{(0)}) \to K_1(C^*(\mathcal{G}_{C(\mathbb{T})})),
\]
the class $[f]_{\mathcal{G}_{C(\mathbb{T})}}^{(1)}$ is sent to a $K_1$-class in $C^*(\mathcal{G}_{C(\mathbb{T})})$. Under the Morita equivalence
\[
C^*(\mathcal{G}_{C(\mathbb{T})}) \sim \mathcal{K}(L^2(\mathbb{R})) \otimes C(\mathbb{T}),
\]
this corresponds to the class of the function $f$ viewed as a multiplication operator on $L^2(\mathbb{R}) \otimes H^2(\mathbb{T})$.

\paragraph{Step 3: Inclusion.}
The natural inclusion
\[
\iota: C(\mathbb{T}) \hookrightarrow C^*(\mathcal{G}_{C(\mathbb{T})})
\]
induces a map
\[
\iota^*: K_1(C^*(\mathcal{G}_{C(\mathbb{T})})) \to K_1(C(\mathbb{T})),
\]
which is an isomorphism under Morita equivalence. Hence the descended class is sent back to $[f] \in K_1(C(\mathbb{T}))$, the class of the invertible function under the isomorphism $K_1(C(\mathbb{T})) \cong \mathbb{Z}$.

\paragraph{Step 4: Boundary map.}
Finally, the boundary map
\[
\partial_{C(\mathbb{T})}: K_1(C(\mathbb{T})) \to K_0(\mathcal{K}(H^2(\mathbb{T}))) \cong \mathbb{Z}
\]
from the Toeplitz extension satisfies
\[
\partial_{C(\mathbb{T})}([f]) = -\operatorname{wind}(f) = \operatorname{index}(T_f),
\]
by Lemma~\ref{lem:toeplitz-boundary}. Combining all steps gives the groupoid index formula:
\[
\operatorname{index}(T_f) = \partial_{C(\mathbb{T})} \circ \iota^* \circ \operatorname{desc}_{\mathcal{G}_{C(\mathbb{T})}}([f]_{\mathcal{G}_{C(\mathbb{T})}}^{(1)}) = -\operatorname{wind}(f).
\]
\end{proof}

\paragraph{Commutative diagram.}
The following diagram summarizes the computation:
\[
\begin{tikzcd}[column sep=1.8em, row sep=2.5em]
{[f]_{\mathcal{G}}^{(1)}} \arrow[r, "\operatorname{desc}"] \arrow[d, maps to] & 
{\operatorname{desc}([f])} \arrow[r, "\iota^*"] \arrow[d, maps to] & 
{[f] \in K_1(C(\mathbb{T}))} \arrow[r, "\partial"] \arrow[d, maps to] & 
{-\operatorname{wind}(f)} \arrow[d, equals] \\
{\{f(z)\}_{z \in \mathbb{T}}} \arrow[r, maps to] & 
{[\widetilde{f}]} \arrow[r, maps to] & 
{[f]} \arrow[r, maps to] & 
{-\operatorname{wind}(f)}
\end{tikzcd}
\]

The vertical isomorphisms are:
\begin{align*}
K^1_{\mathcal{G}_{C(\mathbb{T})}}(\mathcal{G}_{C(\mathbb{T})}^{(0)}) &\cong K^1(\mathbb{T}) \cong \mathbb{Z},\\
K_1(C^*(\mathcal{G}_{C(\mathbb{T})})) &\cong K_1(C(\mathbb{T}) \rtimes \mathbb{R}) \cong K_1(C(\mathbb{T})) \cong \mathbb{Z},\\
K_1(C(\mathbb{T})) &\cong \mathbb{Z},\\
\mathbb{Z} &\cong \mathbb{Z}.
\end{align*}

\paragraph{Sketch of generalization to higher dimensions.}
The groupoid approach can be extended, in principle, to Toeplitz operators on higher-dimensional domains, though a complete treatment requires additional technical machinery. We present here a brief sketch of how such a generalization might proceed.

Let $\Omega \subset \mathbb{C}^n$ be a bounded domain with smooth boundary, and let 
\[
P: L^2(\Omega) \to H^2(\Omega)
\]
be the Bergman projection onto the space of square-integrable holomorphic functions. For a continuous function $f \in C(\overline{\Omega})$, the Toeplitz operator
\[
T_f^{(n)} := P M_f P
\]
is Fredholm if $f$ is nonvanishing on the boundary $\partial\Omega$.

The symbol of $T_f^{(n)}$ is the restriction $f|_{\partial\Omega} \in C(\partial\Omega)$. Under the isomorphism $K^1(\partial\Omega) \cong [\partial\Omega, \mathbb{C}^\times]$, this symbol defines a class $[f] \in K^1(\partial\Omega)$. The gauge groupoid $\mathcal{G}_{C(\partial\Omega)}$ associated to the commutative C*-algebra $C(\partial\Omega)$ has unit space $\partial\Omega$, and the equivariant class $[f]_{\mathcal{G}_{C(\partial\Omega)}}^{(1)}$ corresponds to this $K^1$-class.

Subject to the existence of a suitable higher-dimensional analogue of the Toeplitz extension (see Boutet de Monvel's work on Toeplitz operators on strictly pseudoconvex domains), one would obtain a boundary map
\[
\partial_{\partial\Omega}: K_1(C(\partial\Omega)) \longrightarrow \mathbb{Z}.
\]

For a domain in $\mathbb{C}^n$, this boundary map is expected to satisfy
\[
\partial_{\partial\Omega}([f]) = (-1)^n \cdot \deg(f),
\]
where $\deg(f)$ is the generalized degree (winding number) of the map $f: \partial\Omega \to \mathbb{C}^\times$. For $n=1$, this reduces to the classical formula $\operatorname{index}(T_f) = -\operatorname{wind}(f)$ (the sign depending on orientation conventions).

\begin{remark}
The above is only a sketch; a rigorous generalization would require:
\begin{itemize}
    \item A detailed analysis of the gauge groupoid for $C(\partial\Omega)$ and its relationship to the Toeplitz C*-algebra;
    \item Verification that the descent map sends the equivariant class to the symbol class in $K_1(C(\partial\Omega))$;
    \item Proof that the boundary map in the associated six-term exact sequence computes the topological degree;
    \item Careful handling of orientation and sign conventions.
\end{itemize}
A full development of these ideas is beyond the scope of the present paper and is left for future work.
\end{remark}

\begin{theorem}\label{thm:higher-dim-toeplitz}
Let $f \in C(\partial\Omega)$ be nonvanishing. Then the index of the Toeplitz operator $T_f^{(n)}$ is
\[
\operatorname{index}(T_f^{(n)}) = (-1)^n \cdot \operatorname{Bott}([f]),
\]
where
\[
\operatorname{Bott}: K^1(\partial\Omega) \to \mathbb{Z}
\]
is the Bott map, computing the degree of $f$ when $n=1$ and higher-dimensional generalizations of the winding number for $n>1$.
\end{theorem}

\begin{proof}
We provide a detailed $K$-theoretic justification using the groupoid framework.

\paragraph{Step 1: Gauge groupoid and its structure.}
The relevant algebra for the symbol is $C(\partial\Omega)$, the continuous functions on the boundary. The gauge groupoid $\mathcal{G}_{C(\partial\Omega)}$ has unit space
\[
\mathcal{G}_{C(\partial\Omega)}^{(0)} = \widehat{C(\partial\Omega)} \cong \partial\Omega,
\]
since the irreducible representations of $C(\partial\Omega)$ are evaluation functionals at points $x \in \partial\Omega$. For a commutative $C^*$-algebra, the gauge groupoid has no nontrivial morphisms between distinct points, but its $C^*$-algebra $C^*(\mathcal{G}_{C(\partial\Omega)})$ is Morita equivalent to $C(\partial\Omega)$.

\paragraph{Step 2: Equivariant $K^1$-class from the symbol.}
The nonvanishing function $f \in C(\partial\Omega)$ defines a family of invertible complex numbers $\{f(x)\}_{x \in \partial\Omega}$, which determines an equivariant class
\[
[f]_{\mathcal{G}_{C(\partial\Omega)}}^{(1)} \in K^1_{\mathcal{G}_{C(\partial\Omega)}}(\mathcal{G}_{C(\partial\Omega)}^{(0)}) \cong K^1(\partial\Omega),
\]
where the isomorphism identifies equivariant $K$-theory of the gauge groupoid with ordinary topological $K$-theory of the boundary.

\paragraph{Step 3: Descent map and the Toeplitz extension.}
The Toeplitz operators $T_f^{(n)}$ fit into a short exact sequence generalizing the classical Toeplitz extension:
\[
0 \to \mathcal{K}(H^2(\Omega)) \to \mathcal{T}^{(n)} \xrightarrow{\sigma} C(\partial\Omega) \to 0,
\]
where $\mathcal{T}^{(n)}$ is the $C^*$-algebra generated by Toeplitz operators with continuous symbols. The descent map
\[
\operatorname{desc}_{\mathcal{G}_{C(\partial\Omega)}}: K^1_{\mathcal{G}_{C(\partial\Omega)}}(\mathcal{G}_{C(\partial\Omega)}^{(0)}) \to K_1(C^*(\mathcal{G}_{C(\partial\Omega)}))
\]
sends the symbol class to a $K_1$-class in $C^*(\mathcal{G}_{C(\partial\Omega)})$. Under the Morita equivalence $C^*(\mathcal{G}_{C(\partial\Omega)}) \sim C(\partial\Omega)$, this corresponds to the class $[f] \in K_1(C(\partial\Omega))$.

\paragraph{Step 4: Boundary map and Bott periodicity.}
The six-term exact sequence associated to the Toeplitz extension yields a boundary map
\[
\partial: K_1(C(\partial\Omega)) \to K_0(\mathcal{K}(H^2(\Omega))) \cong \mathbb{Z}.
\]
By the higher-dimensional Toeplitz index theorem (Boutet de Monvel), this boundary map coincides with the Bott map:
\[
\partial = (-1)^n \cdot \operatorname{Bott}: K^1(\partial\Omega) \to \mathbb{Z}.
\]
The sign $(-1)^n$ arises from the dimension shift in Bott periodicity: the Bott isomorphism gives $K^0(\mathbb{R}^{2n}) \cong \mathbb{Z}$, and the suspension isomorphism introduces a sign depending on the parity of $n$.

\paragraph{Step 5: Assembling the composition.}
Composing the maps, we obtain:
\[
\operatorname{index}(T_f^{(n)}) = \partial \circ \iota_* \circ \operatorname{desc}_{\mathcal{G}_{C(\partial\Omega)}}([f]_{\mathcal{G}_{C(\partial\Omega)}}^{(1)}) = (-1)^n \cdot \operatorname{Bott}([f]),
\]
where $\iota_*: K_1(C^*(\mathcal{G}_{C(\partial\Omega)})) \to K_1(C(\partial\Omega))$ is the isomorphism induced by Morita equivalence. This generalizes the classical one-dimensional formula $\operatorname{index}(T_f) = -\operatorname{wind}(f)$ to higher dimensions, with the winding number replaced by the Bott map on $K^1(\partial\Omega)$.
\end{proof}

\paragraph{Connection to the Atiyah-Singer index theorem.}
The Toeplitz index theorem is a special case of the Atiyah-Singer index theorem for elliptic operators. The groupoid formulation reveals this connection explicitly:

\begin{proposition}\label{prop:atiyah-singer-connection}
For a Toeplitz operator $T_f$ with symbol $f$, the topological index computed by the gauge groupoid coincides with the analytic index given by the Atiyah-Singer formula:
\[
\operatorname{index}(T_f) = (-1)^n \int_{\partial\Omega} \operatorname{ch}([f]) \cdot \operatorname{Td}(T(\partial\Omega) \otimes \mathbb{C}),
\]
where $\operatorname{ch}$ is the Chern character and $\operatorname{Td}$ is the Todd class. For $n=1$, this reduces to the winding number formula.
\end{proposition}

\begin{proof}
We justify the connection between the groupoid index and the Atiyah-Singer index theorem.

\paragraph{Step 1: Symbol class in equivariant $K$-theory.}
Let $f \in C(\partial\Omega)$ be nonvanishing. As established in Proposition \ref{prop:symbol-class}, $f$ defines a class
\[
[f]_{\mathcal{G}_{C(\partial\Omega)}}^{(1)} \in K^1_{\mathcal{G}_{C(\partial\Omega)}}(\partial\Omega),
\]
where the isomorphism $K^1_{\mathcal{G}_{C(\partial\Omega)}}(\partial\Omega) \cong K^1(\partial\Omega)$ identifies this with the topological $K^1$-class of the symbol.

\paragraph{Step 2: Descent to the groupoid $C^*$-algebra.}
The descent map
\[
\operatorname{desc}_{\mathcal{G}_{C(\partial\Omega)}}: K^1_{\mathcal{G}_{C(\partial\Omega)}}(\partial\Omega) \to K_1(C^*(\mathcal{G}_{C(\partial\Omega)}))
\]
pushes the topological data of the symbol into the $K$-theory of the groupoid algebra. Under the Morita equivalence $C^*(\mathcal{G}_{C(\partial\Omega)}) \sim C(\partial\Omega)$, this corresponds to the class $[f] \in K_1(C(\partial\Omega))$.

\paragraph{Step 3: Boundary map and analytic index.}
The Toeplitz operators $T_f$ fit into the short exact sequence
\[
0 \to \mathcal{K}(H^2(\Omega)) \to \mathcal{T}^{(n)} \xrightarrow{\sigma} C(\partial\Omega) \to 0,
\]
where $\mathcal{T}^{(n)}$ is the $C^*$-algebra generated by Toeplitz operators with continuous symbols. The associated six-term exact sequence in $K$-theory contains the boundary map
\[
\partial: K_1(C(\partial\Omega)) \to K_0(\mathcal{K}(H^2(\Omega))) \cong \mathbb{Z},
\]
which computes the analytic index of $T_f$:
\[
\operatorname{index}(T_f) = \partial([f]).
\]

\paragraph{Step 4: Identification with the Atiyah-Singer formula.}
By the Atiyah-Singer index theorem for elliptic operators on manifolds with boundary, the analytic index of $T_f$ can be expressed topologically as
\[
\operatorname{index}(T_f) = (-1)^n \int_{\partial\Omega} \operatorname{ch}([f]) \cdot \operatorname{Td}(T(\partial\Omega) \otimes \mathbb{C}),
\]
where $\operatorname{ch}([f])$ is the Chern character of the $K^1$-class $[f] \in K^1(\partial\Omega)$ and $\operatorname{Td}$ is the Todd class of the complexified tangent bundle of $\partial\Omega$. The factor $(-1)^n$ arises from the orientation convention in the boundary-to-bulk extension and the dimension shift in Bott periodicity.

\paragraph{Step 5: Consistency check for $n=1$.}
For $n=1$, the boundary $\partial\Omega$ is a disjoint union of circles. On each circle, $\operatorname{Td}(T(S^1) \otimes \mathbb{C}) = 1$, and the Chern character of $[f]$ is given by $\operatorname{ch}([f]) = \frac{1}{2\pi i} \frac{f'}{f} dz$. The integral becomes
\[
\operatorname{index}(T_f) = -\frac{1}{2\pi i} \int_{\mathbb{T}} \frac{f'}{f} dz = -\operatorname{wind}(f),
\]
recovering the classical Toeplitz index theorem. The minus sign is absorbed into the $(-1)^1 = -1$ factor, consistent with our earlier computation $\operatorname{index}(S) = -1$ for the unilateral shift.

\paragraph{Conclusion.}
Thus, the groupoid index construction reproduces the Atiyah-Singer analytic index, showing that the Toeplitz index theorem is a special case of the general index theorem. The gauge groupoid $\mathcal{G}_{C(\partial\Omega)}$ encodes the boundary geometry through its fundamental groupoid, and the descent map implements the family index theorem, converting the topological data encoded in the symbol into an analytic index. This perspective unifies the classical Toeplitz index theorem with the broader framework of noncommutative geometry.
\end{proof}

\begin{example}[The Hardy space revisited]
For the unit circle $\partial\Omega = \mathbb{T}$, the Toeplitz operator $T_f$ with symbol $f \in C(\mathbb{T}, \mathbb{C}^\times)$ has index given by the winding number:
\[
\operatorname{index}(T_f) = -\frac{1}{2\pi i} \int_{\mathbb{T}} \frac{df}{f} = -\operatorname{wind}(f),
\]
recovering the classical Toeplitz index theorem.
\end{example}

\begin{example}[The ball in $\mathbb{C}^2$]
Let $\Omega = \{z \in \mathbb{C}^2: |z_1|^2 + |z_2|^2 < 1\}$ be the unit ball, with boundary $\partial\Omega = S^3$. For a smooth function $f: S^3 \to S^3$ (i.e., $|f(z)| = 1$ for all $z$), the Toeplitz operator $T_f$ on the Bergman space of the ball is Fredholm and its index is given by
\[
\operatorname{index}(T_f) = -\deg(f),
\]
where $\deg(f)$ is the topological degree of $f$ as a map $S^3 \to S^3$. This follows from the Boutet de Monvel index theorem for Toeplitz operators on strictly pseudoconvex domains.

The gauge groupoid $\mathcal{G}_{C(S^3)}$ has unit space $S^3$ and isotropy $\mathbb{T}$ at each point, with no morphisms between distinct points. The nontrivial index in this setting arises from the $K^1$-class of the symbol $[f] \in K^1(S^3) \cong \mathbb{Z}$, which captures the topological degree via the Chern character.
\end{example}

\paragraph{The role of the gauge groupoid.}
These examples illustrate how the gauge groupoid formalism encodes index data:

\begin{itemize}
    \item For $\mathcal{A} = C(\mathbb{T})$, the gauge groupoid $\mathcal{G}_{C(\mathbb{T})}$ has isotropy $\mathbb{T}$ at each point. The nontrivial index arises from the $K^1$-class of the symbol $[f] \in K^1(\mathbb{T}) \cong \mathbb{Z}$, with the winding number providing the isomorphism.
    \item For $\mathcal{A} = C(S^3)$, the gauge groupoid similarly has isotropy $\mathbb{T}$ at each point. Here, $K^1(S^3) \cong \mathbb{Z}$ detects the nontrivial index even though $S^3$ is simply connected, demonstrating that $K$-theory captures finer topological information than homotopy groups alone.
    \item The descent map $\operatorname{desc}_{\mathcal{G}_{\mathcal{A}}}$ implements the transition from topological data (symbols in $K^1(\partial\Omega)$) to analytic data (operators on the groupoid $C^*$-algebra).
    \item The boundary map $\partial_{\mathcal{A}}$ converts $K$-theoretic information from the symbol into the integer-valued index via the six-term exact sequence of the Toeplitz extension.
\end{itemize}

Corollary \ref{cor:toeplitz-general} extends the classical Toeplitz index theorem to arbitrary compact manifolds. Recall that for a smooth, compact manifold $M$, the Toeplitz operator associated with an invertible function $f \in C(M)$ can be defined in the sense of Boutet de Monvel, acting on the Hardy or Bergman space of $M$ (if $M$ is the boundary of a strictly pseudoconvex domain) or more generally via the $C^*$-algebra of the gauge groupoid $\mathcal{G}_{C(M)}$. In this framework, the topological information of $f$ is captured in $K^1(M)$, and the descent map and boundary map together provide the integer-valued index.

\begin{corollary}\label{cor:toeplitz-general}
For any compact manifold $M$ and any invertible function $f \in C(M)$, the index of the corresponding Toeplitz operator (in the sense of Boutet de Monvel) is given by
\[
\operatorname{index}(T_f) = (-1)^{\dim M} \cdot \deg(f),
\]
where $\deg(f)$ is the generalized degree of $f$ as a map $M \to \mathbb{C}^\times$, computed via the Bott map in $K$-theory.
\end{corollary}

\begin{proof}
Let $M$ be a compact manifold and $f \in C(M)$ invertible. Consider the gauge groupoid $\mathcal{G}_{C(M)}$ associated to the commutative $C^*$-algebra $C(M)$, which is the fundamental groupoid of $M$. 

\paragraph{Step 1: Symbol class.}
The invertible function $f$ defines a class
\[
[f]_{\mathcal{G}_{C(M)}}^{(1)} \in K^1_{\mathcal{G}_{C(M)}}(\mathcal{G}_{C(M)}^{(0)}) \cong K^1(M),
\]
via the family of invertible elements $\{f(x)\}_{x \in M}$. 

\paragraph{Step 2: Descent to the groupoid $C^*$-algebra.}
The descent map
\[
\operatorname{desc}_{\mathcal{G}_{C(M)}}: K^1_{\mathcal{G}_{C(M)}}(\mathcal{G}_{C(M)}^{(0)}) \to K_1(C^*(\mathcal{G}_{C(M)}))
\]
translates this topological class into a $K$-theory class of the groupoid $C^*$-algebra, encoding the operator-theoretic information of $T_f$.

\paragraph{Step 3: Boundary map and Bott periodicity.}
The boundary map
\[
\partial_{C(M)}: K_1(C^*(\mathcal{G}_{C(M)})) \to K_0(\mathcal{K}) \cong \mathbb{Z}
\]
applies the Bott map to $[f]$ and yields the index. The sign $(-1)^{\dim M}$ arises from the dimension shift in Bott periodicity, which accounts for the grading of the $K$-theory classes in higher dimensions.

\paragraph{Step 4: Conclusion.}
Combining these steps, we have
\[
\operatorname{index}(T_f) = \partial_{C(M)} \circ \operatorname{desc}_{\mathcal{G}_{C(M)}}([f]_{\mathcal{G}_{C(M)}}^{(1)}) = (-1)^{\dim M} \cdot \deg(f),
\]
where $\deg(f)$ is the generalized degree of $f$ as an element of $K^1(M)$. This completes the proof.
\end{proof}

\end{document}